\documentclass[11pt]{gsasthesis}
\PassOptionsToPackage{unicode}{hyperref}
\pdfoutput=1
\usepackage{microtype}
\usepackage{textgreek}

\usepackage{etex}
\usepackage[margin={1.2in}]{geometry}

\usepackage{subfiles}
\usepackage{fullpage}

\usepackage{ragged2e}
\usepackage{amsmath}
\usepackage{amssymb}
\usepackage{amsthm}

\usepackage[all]{xy}
\usepackage{tikz-cd}
\usetikzlibrary{matrix, arrows}

\usepackage[shortlabels]{enumitem}

\usepackage{hyperref}
\hypersetup{colorlinks=true,linkcolor={blue!65!black},citecolor={green!50!black},urlcolor={magenta!50!black}}

\usepackage[backend=biber,style=numeric,citestyle=numeric-comp,maxbibnames=100,doi=false,isbn=false,url=false]{biblatex}
\AtEveryBibitem{\iffieldsequal{eid}{pages}{\clearfield{pages}}{}}
\AtEveryBibitem{\clearfield{eid}}

\usepackage{cleveref}

\newtheorem{theorem}{Theorem}[section]
\newtheorem{lemma}[theorem]{Lemma}
\newtheorem{proposition}[theorem]{Proposition}
\newtheorem{corollary}[theorem]{Corollary}

\theoremstyle{definition}
\newtheorem{definition}[theorem]{Definition}
\newtheorem{example}[theorem]{Example}
\newtheorem{nonexample}[theorem]{Non-example}
\newtheorem{remark}[theorem]{Remark}
\newtheorem{convention}[theorem]{Convention}
\newtheorem{notation}[theorem]{Notation}
\newtheorem{warning}[theorem]{Warning}
\newtheorem{question}[theorem]{Question}

\mathchardef\mhyphen="2D 

\setlist[enumerate,1]{label=(\arabic*),font=\normalfont}

\DeclareMathOperator{\Ob}{Ob}
\DeclareMathOperator{\Hom}{Hom}
\DeclareMathOperator{\THom}{\mathbf{Hom}} 
\DeclareMathOperator{\id}{id}
\DeclareMathOperator*{\colim}{colim}
\newcommand{\op}{\mathrm{op}}   
\newcommand{\co}{\mathrm{co}}

\newcommand{\bp}{\mathop{\square}\nolimits}    
\newcommand{\adj}{\rightleftarrows} 
\newcommand{\yo}{\mathbf{y}}        

\newcommand{\iso}{\cong}        
\newcommand{\eqv}{\simeq}       

\newcommand{\Set}{\mathrm{Set}}
\newcommand{\Cat}{\mathrm{Cat}}

\newcommand{\PP}{\mathcal{P}}   
\newcommand{\OO}{\mathcal{O}}   

\newcommand{\llp}{\mathsf{llp}}
\newcommand{\rlp}{\mathsf{rlp}}
\newcommand{\sL}{\mathsf{L}}    
\newcommand{\sR}{\mathsf{R}}    

\newcommand{\All}{\mathsf{All}} 
\newcommand{\Iso}{\mathsf{Iso}} 
\newcommand{\Mono}{\mathsf{Mono}} 
\newcommand{\Epi}{\mathsf{Epi}}   
\newcommand{\Eqv}{\mathsf{Eqv}} 

\newcommand{\sC}{\mathsf{C}}    
\newcommand{\sF}{\mathsf{F}}    
\newcommand{\sW}{\mathsf{W}}    
\newcommand{\sAC}{\mathsf{AC}}  
\newcommand{\sAF}{\mathsf{AF}}  

\newcommand{\cto}{\hookrightarrow} 
\newcommand{\acto}{\stackrel{\sim}{\hookrightarrow}} 
\newcommand{\fto}{\twoheadrightarrow} 
\newcommand{\afto}{\stackrel{\sim}{\twoheadrightarrow}} 
\newcommand{\ancto}{\stackrel{\an}{\hookrightarrow}} 
\newcommand{\anfto}{\stackrel{\an}{\twoheadrightarrow}} 

\newcommand{\htopicl}{\stackrel{l}{\simeq}} 
\newcommand{\htopicr}{\stackrel{r}{\simeq}} 
\newcommand{\htopic}{\simeq}                

\newcommand{\R}{\mathrm{R}}     

\DeclareMathOperator{\Ho}{Ho}

\DeclareMathOperator{\Map}{Map}

\newcommand{\Kan}{\mathrm{Kan}}
\newcommand{\init}{\mathrm{init}} 
\newcommand{\fin}{\mathrm{fin}} 
\newcommand{\tr}{\mathrm{tr}}   
\newcommand{\proj}{\mathrm{proj}} 
\newcommand{\inj}{\mathrm{inj}} 
\newcommand{\Reedy}{\mathrm{Reedy}} 
\newcommand{\cof}{\mathrm{cof}} 
\newcommand{\acof}{\mathrm{acof}} 
\newcommand{\fib}{\mathrm{fib}} 
\newcommand{\cf}{\mathrm{cf}}   
\newcommand{\Ld}{\mathbf{L}}    

\renewcommand{\C}{\mathbf{C}}     

\newcommand{\TCat}{\mathbf{Cat}} 
\newcommand{\LPr}{\mathbf{LPr}}

\newcommand{\Mod}{\mathbf{Mod}}
\newcommand{\PM}{\mathbf{PM}}
\newcommand{\PMA}{\mathbf{PM^A}}
\newcommand{\PMR}{\mathbf{PM^R}}
\newcommand{\CPM}{\mathbf{CPM}}

\newcommand{\QAdj}{\mathrm{QAdj}}
\newcommand{\LQF}{\mathrm{LQF}}
\newcommand{\RQF}{\mathrm{RQF}}
\newcommand{\QBF}{\mathrm{QBF}}

\newcommand{\ltri}{\lhd}        
\newcommand{\rtri}{\rhd}        

\newcommand{\N}{\mathbb{N}}

\newcommand{\an}{\mathbf{A}}    

\newcommand{\II}{\mathbf{I}}
\newcommand{\JJ}{\mathbf{J}}

\newcommand{\Sq}{\mathbf{Sq}}   
\newcommand{\LL}{\mathbf{L}}

\newcommand{\E}{\mathbb{E}}

\DeclareMathOperator{\Spec}{Spec}

\title{A model 2-category of enriched combinatorial premodel categories}
\author{Reid William Barton}
\degreename{Doctor of Philosophy}
\degreefield{Mathematics}
\department{The Department of Mathematics}
\degreemonth{July}
\degreeyear{2019}
\principaladvisor{Professor Michael Hopkins}

\begin{document}


\pagenumbering{roman}

\thesistitlepage
\copyrightpage
\begin{abstract}
  Quillen equivalences induce equivalences of homotopy theories and therefore form a natural choice for the ``weak equivalences'' between model categories.
  In \cite{Ho}, Hovey asked whether the 2-category $\Mod$ of model categories has a ``model 2-category structure'' with these weak equivalences.
  We give an example showing that $\Mod$ does not have pullbacks, so cannot be a model 2-category.

  We can try to repair this lack of limits by generalizing the notion of model category.
  The lack of limits in $\Mod$ is due to the two-out-of-three axiom, so we define a \emph{premodel category} to be a complete and cocomplete category equipped with two nested weak factorization systems.
  Combinatorial premodel categories form a 2-category $\CPM$ with excellent algebraic properties: $\CPM$ has all limits and colimits and is equipped with a tensor product (representing Quillen bifunctors) which is adjoint to an internal Hom.

  The homotopy theory of a model category depends in an essential way on the weak equivalences, so it does not extend directly to premodel categories.
  We build a substitute homotopy theory under an additional axiom on the premodel category, which holds automatically for a premodel category enriched in a monoidal model category $V$.
  The 2-category of combinatorial $V$-premodel categories $V\CPM$ is simply the 2-category of modules over the monoid object $V$, so $V\CPM$ inherits the algebraic structure of $\CPM$.

  We construct a model 2-category structure on $V\CPM$ for $V$ a tractable symmetric monoidal model category, by adapting Szumi\l o's construction of a fibration category of cofibration categories \cite{Sz}.
  For set-theoretic reasons, constructing factorizations for this model 2-category structure requires a technical variant of the small object argument which relies on an analysis of the rank of combinatoriality of a premodel category.
\end{abstract}


\setcounter{tocdepth}{1}
\tableofcontents 

\newpage 

\begin{acknowledgments}
  First, I would like to thank my advisor Mike Hopkins for his support and patience.
  For many, many stimulating discussions, I would like to thank Clark Barwick, Sam Isaacson, and Inna Zakharevich.
  In particular, one day Sam posed the question of whether there exists a tensor product of model categories.
  The answer (``sometimes'') did not satisfy me and this was the starting point of a project that eventually became the present work.
  I would also like to thank Jacob Lurie and Haynes Miller for agreeing to serve on my thesis committee and for taking the time to read this dissertation.
  Finally, many thanks to my family and friends, especially Andrew and Kate and Niki.
\end{acknowledgments}

\pagenumbering{arabic}

\chapter{Introduction}
\label{chap:intro}


\section{Model categories and their homotopy theories}

Ever since their introduction by Quillen \cite{Q}, model categories have been a central part of the language of homotopy theory.
We begin with a brief overview of model categories and the roles they play in homotopy theory today.

A \emph{model category} is a category equipped with a certain kind of additional structure which allows one to carry out the constructions of homotopy theory, such as the formation of mapping cylinders.
In this way, the theory of model categories provides an organizational framework for ``homotopy theories'' in the same way that ordinary categories form an organizational framework used in many other areas of mathematics.
The title of the first chapter of \cite{Q}, ``Axiomatic homotopy theory'', reflects this perspective on model categories.
Examples of homotopy theories include not only ones arising from homotopy theory itself, such as spaces and spectra, but also many of an algebraic nature, the most familiar example being the homotopy theory of chain complexes from homological algebra.

The characteristic feature of a homotopy theory is the existence for each pair of objects (spaces, chain complexes, etc.) $A$ and $B$ of not just a set of maps from $A$ to $B$ but a \emph{space} of such maps.
We may think of points of this space $\Map(A, B)$ as being maps from $A$ to $B$, paths in the space $\Map(A, B)$ as being homotopies between maps, homotopies between paths in $\Map(A, B)$ as being homotopies between homotopies, and so on.
Each connected component of $\Map(A, B)$ corresponds to a homotopy equivalence class of maps from $A$ to $B$, but the space $\Map(A, B)$ also contains higher-order information which encodes, for example, the homotopically distinct homotopies between two given maps from $A$ to $B$.

Given a model category $M$, there are a variety of ways to construct, for any two objects $A$ and $B$ of $M$, a simplicial set $\Map_M(A, B)$ which has the correct homotopy type to represent the space of maps from $A$ to $B$ in the homotopy theory associated to $M$.
Moreover, there are composition maps $\Map_M(B, C) \times \Map_M(A, B) \to \Map_M(A, C)$ which assemble these mapping spaces into a \emph{simplicial category}, that is, a category enriched in simplicial sets.
A simplicial category is the most direct realization of the idea that a homotopy theory can be described in terms of the spaces of maps between its objects.
The simplicial sets $\Map_M(A, B)$ are determined only up to weak homotopy equivalence, so we call a functor between simplicial categories an equivalence (or a Dwyer--Kan equivalence) if it is essentially surjective and induces a weak homotopy equivalence on each mapping space.

Model categories are related to one another by \emph{Quillen adjunctions}, pairs of adjoint functors which respect the model category structures in a particular way.
A Quillen adjunction $F : M \adj N : G$ between two model categories induces an adjunction between their associated simplicial categories.
When this induced adjunction is an equivalence, we say that $F$ and $G$ are \emph{Quillen equivalences} and we think of the model categories they relate as two presentations of the same homotopy theory.



Simplicial categories are themselves the objects of a model category developed in work of Dwyer and Kan \cite{DK1} and Bergner \cite{Ber} whose weak equivalences are the Dwyer--Kan equivalences.
We can think of this model category as describing the \emph{homotopy theory of homotopy theories}.
Other models for homotopy theories include complete Segal spaces \cite{Rezk}, quasicategories \cite{BV,Joy}, and relative categories \cite{RelCat}.
All of these model categories are known to be Quillen equivalent, so all of these model categories are presentations of the homotopy theory of homotopy theories.

Nowadays, ``homotopy theories'' are better known as \emph{$(\infty, 1)$-categories}.
Using the model of quasicategories, Lurie has extended an enormous amount of classical category theory to the $(\infty, 1)$-categorical setting \cite{HTT}.
From this perspective the purpose of a model category is to serve as a presentation of the object of real interest, its associated $(\infty, 1)$-category.
However, model categories are still quite useful for performing calculations.
A popular analogy is that an $(\infty, 1)$-category is like an abstract vector space, while a model category is like a vector space equipped with a choice of basis.

While each model category has an associated $(\infty, 1)$-category, not all $(\infty, 1)$-categories arise from model categories.
Specifically, the $(\infty, 1)$-category associated to any model category always admits all limits and colimits.
Moreover, as mentioned earlier, a Quillen adjunction between model categories induces an adjunction between the associated $(\infty, 1)$-categories.
Thus, the assignment to each model category of its associated $(\infty, 1)$-category lands inside the class of complete and cocomplete homotopy theories, as shown by the vertical arrow in the middle column of the figure below.

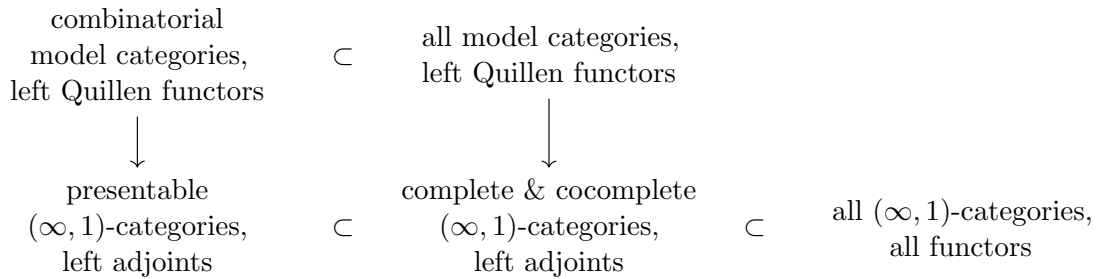
\begin{figure}[b]
\[
  \begin{tikzcd}
    \parbox{10pc}
    {\centering combinatorial\\model categories,\\left Quillen functors}
    \arrow[draw=none]{r}[sloped,auto=false,font=]{\subset} 
    \ar[d] &
    \parbox{10pc}
    {\centering all model categories,\\left Quillen functors} \ar[d] \\
    \parbox{10pc}
    {\centering presentable\\$(\infty, 1)$-categories,\\left adjoints}
    \arrow[draw=none]{r}[sloped,auto=false,font=]{\subset} 
    &
    \parbox{10pc}
    {\centering complete \& cocomplete\\$(\infty, 1)$-categories,\\left adjoints}
    \arrow[draw=none]{r}[sloped,auto=false,font=]{\subset} 
    &
    \parbox{10pc}
    {\centering all $(\infty, 1)$-categories,\\all functors}
  \end{tikzcd}
\]
\caption{The relationship between model categories and their homotopy theories.}
\end{figure}

Almost all model categories of interest are large categories: they have a proper class of objects.
A general large category is a rather unwieldy object.
A more convenient class of model categories is the class of \emph{combinatorial model categories}, which are ones whose structure is determined in a certain sense by a small amount of data.
The corresponding notion for $(\infty, 1)$-categories is that of \emph{presentable $(\infty, 1)$-categories}.
The $(\infty, 1)$-category associated to a combinatorial model category is always presentable, as indicated by the left vertical arrow in the figure.
Moreover, every presentable $(\infty, 1)$-category is the homotopy theory associated to some combinatorial model category.
Hence, we may summarize the relationship between model categories and $(\infty, 1)$-categories by saying that the class of combinatorial model categories provides a model for the class of presentable $(\infty, 1)$-categories, at least in the somewhat weak sense that the left vertical arrow is essentially surjective and takes exactly the Quillen equivalences to equivalences of $(\infty, 1)$-categories.

\section{A model category of combinatorial model categories?}

One might hope that combinatorial model categories actually present the homotopy theory of presentable $(\infty, 1)$-categories in a much stronger sense.
Namely, Hovey posed the following question (paraphrased from \cite[Problem~8.1]{Ho}):
\begin{question}\label{ques:model-model}
  For some reasonable notion of ``model 2-category'', is there a model 2-category of (combinatorial%
  \footnote{Hovey does not include this condition. Its presence or absence will not be important for the purpose of this section.}%
  ) model categories and left Quillen functors whose weak equivalences are the Quillen equivalences?
\end{question}
This hope is reasonable as presentable $(\infty, 1)$-categories and left adjoints between them form a complete and cocomplete $(\infty, 1)$-category \cite[Proposition~5.5.3.13 and Theorem~5.5.3.18]{HTT}.
Furthermore, certain known facts about model categories are suggestive of such a model category structure.
Notably, Dugger \cite{Dug} showed that combinatorial model categories admit ``presentations'': specifically, every combinatorial model category $M$ admits a left Quillen equivalence $F : L_S \Kan^{C^\op}_\proj \to M$ from a left Bousfield localization of the projective model category structure on a category of simplicial presheaves.
Left Quillen functors out of $L_S \Kan^{C^\op}_\proj$ admit a simple description in terms of $C$ and $S$, so the equivalence $F : L_S \Kan^{C^\op}_\proj \to M$ is a plausible candidate for a cofibrant replacement of the model category $M$.


\Cref{ques:model-model} is somewhat imprecise in that it does not specify a notion of ``model 2-category''.
One possible interpretation is to ignore the 2-categorical structure entirely and ask for a model category structure on $\mathrm{Mod}$, the 1-category of model categories and left Quillen functors.
However $\mathrm{Mod}$ fails badly to have limits and colimits.
For example, two parallel left Quillen functors $F : M \to N$ and $F' : M \to N$ might not \emph{strictly} agree on any object, not even the initial object $\emptyset$ of $M$; then no model category can possibly be the equalizer of $F$ and $F'$.
Any positive answer to \cref{ques:model-model} must avoid this difficulty somehow.
Two possibilities suggest themselves.
\begin{itemize}
\item
  We could modify the 1-category $\mathrm{Mod}$ by equipping each model category with a \emph{choice} of colimits of every small diagram, and considering only those functors which preserve the chosen colimits on the nose.
  Then, if $\emptyset$ denotes the chosen colimit of the empty diagram, both $F$ and $F'$ preserve $\emptyset$ strictly and so there at least exists a nonempty category equalizing $F$ and $F'$.
\item
  Alternatively, we could accept that model categories really form a 2-category $\Mod$ and declare that a model 2-category is not required to have \emph{strict} limits but only limits in an appropriate 2-categorical sense, involving diagrams which commute only up to specified isomorphism.
  (We will simply refer to these as ``limits'' or sometimes ``2-limits'', as opposed to ``strict limits''; in the literature they are also known as ``bilimits''.)
\end{itemize}
The former option would require distinguishing categories which are equivalent but not isomorphic, as presumably a cofibrant model category would need to have chosen colimits which are ``freely adjoined'' in a suitable sense.
Philosophically, we can explain our preference for the second option as follows.
One purpose of model categories is to reduce the calculation of homotopy limits and colimits to that of ordinary 1-categorical limits and colimits.
However, for the kinds of categories which appear as the underlying categories of model categories, it is rarely a good idea to compute 1-categorical limits anyways; it is much more sensible to compute 2-categorical limits (e.g., pseudolimits).
Thus, it would be better to work with a notion of ``model 2-category'' which, instead, reduces the calculation of homotopy limits and colimits to that of 2-categorical limits and colimits.

Let us suppose, then, that we have chosen such a notion of model 2-category.
Unfortunately, there is still a more serious problem: the 2-category of combinatorial model categories also lacks limits and colimits even in the 2-categorical sense.
A model category structure is uniquely determined by its cofibrations and acyclic cofibrations, and left Quillen functors preserve cofibrations and acyclic cofibrations.
Thus, the obvious candidate for the limit $\lim_{i \in I} M_i$ of a diagram of model categories is the limit of the underlying categories, equipped with a model category structure in which a morphism is an (acyclic) cofibration if and only if its image in each $M_i$ is an (acyclic) cofibration.
However, there is no reason why the weak equivalences of this candidate model category structure should satisfy the two-out-of-three axiom.
A priori, even if this structure fails to be a model category, it could still have a universal approximation by a model category.
But in fact this does not happen: we give an explicit example of a diagram of model categories which has no limit.

\begin{proposition}\label{prop:model-pullback}
  Suppose $M_1$, $M_2$ and $M_3$ are three model category structures on the same underlying category $M$ such that the identity functors $M_1 \to M_3$ and $M_2 \to M_3$ are left Quillen functors.
  If the pullback of model categories
  \[
  \xymatrix{
    M_0 \ar@{.>}[r] \ar@{.>}[d] & M_1 \ar[d] \\
    M_2 \ar[r] & M_3
  }
  \]
  exists, then (up to equivalence) $M_0$ also has underlying category $M$ and the underlying functors of the left Quillen functors $M_0 \to M_1$ and $M_0 \to M_2$ are the identity of $M$.
\end{proposition}

\begin{proof}
  Let $A$ be the category of sets equipped with the model category structure in which every morphism is an acyclic fibration.
  Then for any model category $N$, giving a left Quillen functor $A \to N$ is the same as giving a left adjoint from the category of sets to $N$, which (up to equivalence) is the same as giving an object of $N$.
  In other words, the Hom-category functor $\THom(A, -)$ sends a model category $N$ to its underlying category $|N|$.

  Now if $M_0 = \lim M_i$ exists, then $\THom(A, M_0) = \lim \THom(A, M_i)$, so the underlying category of $M_0$ fits into a pullback square of categories
  \[
  \xymatrix{
    |M_0| \ar[r] \ar[d] & M \ar^{\id_M}[d] \\
    M \ar_{\id_M}[r] & M
  }
  \]
  hence (up to equivalence) $|M_0| = M$ and the functors in the above square are all the identity of $M$.
\end{proof}

Given an equivalence of categories $M_0 \to M$ and a model category structure on $M_0$, we can uniquely transfer it to $M$ so that the functor $M_0 \to M$ becomes an equivalence of model categories.
Thus, if a pullback square of the form in Proposition~\ref{prop:model-pullback} exists, we may assume that $M_0$ also has underlying category $M$ and the functors $M_0 \to M_1$ and $M_0 \to M_2$ are the identity.
If $M_0'$ is another model category structure on $M$ such
that the identity functors $M_0' \to M_1$ and $M_0' \to M_2$ are left Quillen functors, then it follows from the universal property of the pullback that the identity functor $M_0' \to M_0$ is a left Quillen functor also.
Hence, $M_0$ is the greatest lower bound (or meet) of $M_1$ and $M_2$ in the poset of model structures on $M$, ordered according to inclusion of the cofibrations and the acyclic cofibrations.
This poset always has a maximum object, the model structure in which every morphism is an acyclic cofibration; so we might as well take $M_3$ to be this maximum object, and $M_3$ plays no further role.

\begin{figure}[h]
  \[
  \def\cf#1{\ar[hookrightarrow, two heads, #1]}
  \def\ac#1{\ar[hookrightarrow, "\sim", #1]}
  \def\af#1{\ar[two heads, "\sim", #1]}
  \begin{tikzcd}
    \cdot \af{r} \af{d} & \cdot \af{d} & & & \cdot \cf r \af{d} & \cdot \af{d} \\
    \cdot \ac{r} & \cdot \ar{rd} & & & \cdot \cf{r} \ar{ld} & \cdot \\
    & & \cdot \cf{r, "?"} \af{d} & \cdot \af{d} \\
    & & \cdot \ac{r} \ar{ld} & \cdot \ar{rd} \\
    \cdot \cf{r} \cf{d} & \cdot \cf{d} & & & \cdot \ac{r} \af{d} & \cdot \af{d} \\
    \cdot \ac{r} & \cdot & & & \cdot \ac{r} & \cdot
  \end{tikzcd}
  \]
  \caption{The poset of model category structures on $M$ may not admit
    meets.}
  \label{fig:no-meets}
\end{figure}
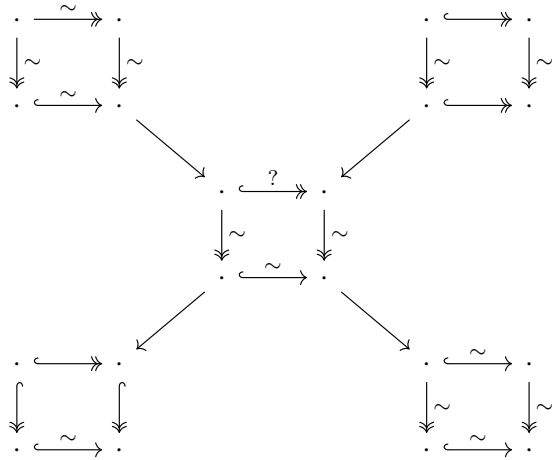

However, the poset of model structures on $M$ does not admit meets in general.
This fails already when $M = \{0 \to 1\} \times \{0 \to 1\}$, as shown in \cref{fig:no-meets}.
One can verify that the structure in each of the four corners of the diagram is indeed a model category, and that each of the four identity functors from one of the top model categories to one of the bottom model categories is a left Quillen functor.
Naming the two bottom model categories $M_1$ and $M_2$, suppose that they admit a meet $M_0$ (shown in the center of the figure).
Then the identity functors from $M_0$ to $M_1$ and $M_2$ are left Quillen functors, and the identity functors from the two top model categories to $M_0$ must also be left Quillen functors.
Since the vertical maps of the bottom right model category structure are acyclic fibrations, they must also be acyclic fibrations in $M_0$.
From the left Quillen functors from the top model categories to $M_0$, we see that the top morphism of $M_0$ must be a cofibration and the bottom morphism an acyclic cofibration.
By two-out-of-three, the top morphism must then be an acyclic cofibration in $M_0$ as well. But its image in the bottom left model category is not an acyclic cofibration, a contradiction.

Thus, we conclude that the diagram formed by $M_1$ and $M_2$ together with the terminal model structure $M_3$ on $M$ and identity functors does not admit a limit, and so the 2-category of combinatorial model categories is not complete.
Hence \emph{\cref{ques:model-model} has a negative answer, even if we only require a model 2-category to have 2-limits}, because $\Mod$ lacks 2-limits.

\begin{remark}
  Suppose we discard the noninvertible natural transformations from $\Mod$, leaving a $(2, 1)$-category, and declare that a model $(2, 1)$-category is only required to have (homotopy) limits.
  Then in \cref{prop:model-pullback}, we could not conclude that $M_0 \to M$ is an equivalence of underlying categories, but only that it is an equivalence of the maximal groupoids contained in $M_0$ and $M$.
  However, we may repeat this argument with $A$ replaced with $A_n$, the category of presheaves of sets on the category $[n] = \{0 \to 1 \to \cdots \to n\}$, to conclude that $M_0 \to M$ induces an equivalence of maximal groupoids of strings of $n$ composable arrows for each $n$.
  As we know from the theory of complete Segal spaces, it follows that $M_0 \to M$ is actually an equivalence of categories.
  Thus, the same argument also shows that the $(2, 1)$-category of combinatorial model categories is not homotopy complete.

  Similarly, if we had chosen to work in the setting of categories with chosen colimits, take $A_n$ to be the free such category on a string of $n$ composable arrows.
  Then we conclude that $M_0$ and $M$ have the same strings of $n$ composable arrows for each $n$, and are therefore isomorphic as categories.
  Hence, we cannot avoid this problem by treating $\mathrm{Mod}$ as a 1-category either.
\end{remark}

\section{Premodel categories}

We have seen that the 2-category of combinatorial model categories does not admit limits, and so cannot be the underlying 2-category of a model 2-category.

However, this is not really cause for concern.
It often occurs that we are primarily interested in only a subclass of the objects of a model category (for example, the fibrant ones).
This subcategory usually will not be closed under limits or colimits, and so the full model category serves as a framework for doing calculations.
In this light, we are led to the question: how might we embed the 2-category of combinatorial model categories in one which is complete and cocomplete?

At this point it is helpful to recall the following concise definition of a model category.
(See for instance \cite{JT}.)

\begin{definition}
  A \emph{model category} is a category $M$ equipped with three classes of maps $\sW$, $\sC$, and $\sF$ such that
  \begin{enumerate}
  \item $(\sC, \sW \cap \sF)$ and $(\sW \cap \sC, \sF)$ are weak factorization systems on $M$, and
  \item $\sW$ satisfies the two-out-of-three condition.
  \end{enumerate}
\end{definition}

In our discussion of limits of model categories, we observed that the obstruction to simply writing down a formula for the limit of a diagram of model categories lies in the weak equivalences.
This suggests the following definition.

\begin{definition}
  A \emph{premodel category} is a category $M$ equipped with four classes of maps $\sC$, $\sAC$, $\sF$, and $\sAF$ such that $(\sC, \sAF)$ and $(\sAC, \sF)$ are weak factorization systems on $M$, and $\sAC \subset \sC$ (equivalently, $\sAF \subset \sF$).
\end{definition}

We call the maps belonging to $\sAC$ \emph{anodyne cofibrations} and those belonging to $\sAF$ \emph{anodyne fibrations}.
Every model category yields a premodel category, by setting $\sAC = \sW \cap \sC$ and $\sAF = \sW \cap \sF$.
In the converse direction, a premodel category can arise in this way from at most model category, the one with $\sW = \sAF \circ \sAC$; but in general this $\sW$ will not satisfy two-out-of-three.

The ``algebraic'' parts of model category theory, such as Quillen functors, Quillen bifunctors, monoidal model categories%
\footnote{For us, a monoidal model category will always have cofibrant unit. We will discuss this point when we introduce monoidal premodel categories.}
and enriched model categories, and the projective, injective, and Reedy model structures on diagram categories, transfer directly to the setting of premodel categories, because these notions do not directly involve the weak equivalences of a model category.
But combinatorial premodel categories turn out to have an even better algebraic theory, like the locally presentable categories they are built on: they form a complete and cocomplete 2-category $\CPM$ with a tensor product and internal Hom.
Monoid objects in $\CPM$ and their modules are precisely the monoidal premodel categories and the premodel categories enriched over them.
The combinatorial \emph{model} categories live inside the combinatorial premodel categories as a full sub-2-category which fails to be closed under most of these operations.
(Even the unit object for the monoidal structure on $\CPM$ is not a model category.)

The trade-off for this rich algebraic structure on the entirety of combinatorial premodel categories is that we have lost the part of the structure needed to define a satisfactory homotopy theory within a single premodel category: the weak equivalences.
For instance, consider the following basic fact about cylinder objects in a model category.

\begin{proposition}
  \label{prop:model-cyl}
  Let $A$ be a cofibrant object of a model category.
  Then there exists a cofibration $A \amalg A \cto C$ such that the two compositions $A \cto A \amalg A \cto C$ are acyclic cofibrations.
\end{proposition}

\begin{proof}
  Factor the fold map $A \amalg A \to A$ into a cofibration $A \amalg A \cto C$ followed by an acyclic fibration $C \afto A$.
  The inclusions $A \to A \amalg A$ are pushouts of the cofibration $\emptyset \cto A$, hence themselves cofibrations.
  Each composition $A \cto A \amalg A \cto C \to A$ is the identity, so by two-out-of-three each $A \cto C$ is a weak equivalence.
\end{proof}

Since each $A \acto C$ is an acyclic cofibration, it has the left lifting property with respect to fibrations.
This kind of fact is used in order to show that the notions of left and right homotopy agree for maps from a cofibrant object to a fibrant one.

The above argument is not available in a premodel category, since anodyne fibrations and anodyne cofibrations are not related via the two-out-of-three property.
In fact, the above proposition has no suitable analogue for a general premodel category.
This lack of cylinder objects in turn means that we cannot construct, for example, homotopy pushouts by the usual method, and so the homotopy theory of a general premodel category ends up looking rather unlike the homotopy theory of a model category.

We call a premodel category \emph{relaxed} if it satisfies certain conditions which roughly amount to the existence of a sufficient supply of cosimplicial and simplicial resolutions.
Then a relaxed premodel category has a homotopy theory which resembles that of a model category.
More precisely, the cofibrant objects of a relaxed premodel category form a \emph{cofibration category}, which by the standard theory of cofibration categories has an associated $(\infty, 1)$-category which is cocomplete.
On the other hand, the fibrant objects of a relaxed premodel category form a fibration category, with an associated homotopy theory which is complete; and these two homotopy theories turn out to be equivalent.
Hence the associated $(\infty, 1)$-category of a relaxed premodel category, like that of a model category, is both complete and cocomplete.
Every model category is relaxed when viewed as a premodel category; this amounts to the existence of resolutions, which is a souped-up version of \cref{prop:model-cyl}.
The $(\infty, 1)$-category associated to a model category can be computed from the cofibration category consisting of the model category's cofibrant objects, so our assignment of an associated $(\infty, 1)$-category to a relaxed premodel category extends the usual one for model categories.

We may then define a left Quillen functor between relaxed premodel categories to be a \emph{Quillen equivalence} if it induces an equivalence of associated homotopy theories.
These functors are the obvious candidate for the class of weak equivalences.
However, if we try to form a model 2-category of relaxed premodel categories, we once more encounter a problem with limits and colimits.
The property of being relaxed amounts to the existence of a sufficient supply of resolutions, but as these resolutions are not part of the structure of a premodel category, there is no apparent reason to think that relaxed premodel categories are closed under the various algebraic operations that general premodel categories enjoy.

To solve this problem, we turn to \emph{enriched} premodel categories.
Recall that, in a simplicial model category, there is a second way to prove Proposition~\ref{prop:model-cyl}: simply take $C = \Delta^1 \otimes A$.
More generally, the tensor and cotensor by $\Delta^n$ yield cosimplicial resolutions of cofibrant objects and simplicial resolutions of fibrant objects.
These constructions also make sense and produce resolutions in $V$-premodel categories for any monoidal model category $V$.
Hence, any $V$-premodel category is automatically relaxed.
Furthermore, the limit or colimit of a diagram of $V$-premodel categories and left Quillen $V$-functors again has the structure of a $V$-premodel category, so combinatorial $V$-premodel categories also form a complete and cocomplete 2-category.

We can now state our main result.

\begin{theorem}
  \label{thm:v-combpremod}
  Let $V$ be a tractable symmetric monoidal model category.
  Then the 2-category of combinatorial $V$-premodel categories has a model 2-category structure in which a left Quillen functor between combinatorial $V$-\emph{model} categories is a weak equivalence if and only if it is a Quillen equivalence.
\end{theorem}

Specializing to $V = \Kan$, the theorem gives a model 2-category structure on combinatorial simplicial premodel categories.

We have already described the underlying 2-category of combinatorial $V$-premodel categories and its weak equivalences at some length.
To complete the outline of the proof of Theorem~\ref{thm:v-combpremod}, we will say something about the cofibrations and fibrations in the model 2-category structure.
We will say just a few words about them here and leave a more detailed overview of their construction and the verification of the model 2-category axioms to \cref{chap:vcpmmodel}.

The fibrations and acyclic fibrations of our model 2-category structure on combinatorial $V$-premodel categories are each characterized by a small number of right lifting properties.
The specific choice of lifting properties is strongly influenced by the fibration category of cofibration categories constructed by Szumi\l o in \cite{Sz}, although some adjustments are necessary to adapt them to the setting of premodel categories.
In particular, a version of \emph{pseudofactorizations} plays an important role in defining the fibrations.

Since the fibrations and acyclic fibrations are each determined by a set of right lifting properties, one might say that the model 2-category is cofibrantly generated, and indeed we will construct factorizations using a variant of the small object argument.
However, there are two problems which arise when trying to apply the small object argument to the 2-category of combinatorial $V$-premodel categories.
First, this 2-category is not locally small.
For example, for a $V$-premodel category $M$, the category $\THom(V, M)$ is equivalent to the full subcategory of $M$ on its cofibrant objects, which is rarely essentially small.
Hence we cannot form a coproduct over all isomorphism classes of squares as in the usual small object argument, because the cardinality of the indexing category would be too large.
Second, almost no object of this 2-category is small.
Even $V$ itself is not a $\mu$-small object for any $\mu$, because a $\mu$-filtered colimit $M = \colim_{i \in I} M_i$ of combinatorial $V$-premodel categories contains arbitrary coproducts of cofibrant objects in the images of the $M_i$, and such a coproduct may not itself belong to the image of any $M_i$.
Here we will mention only that our solution to these problems requires analysis of how the \emph{rank of combinatoriality} behaves under the formation of colimits, limits, tensors and internal Homs of combinatorial premodel categories.


\chapter{Premodel categories}
\label{chap:premodel}

In this chapter we define \emph{premodel categories} and describe how to extend the ``algebraic'' parts of model category theory to premodel categories: Quillen adjunctions, Quillen bifunctors, monoidal premodel categories and their modules, and premodel category structures on diagram categories.
We also introduce the notion of \emph{relaxed} premodel categories.
A relaxed premodel category has a well-behaved homotopy theory, which we will describe in the next chapter.

\section{The 2-category of premodel categories}

\subsection{Weak factorization systems}

The core technical ingredient of model category theory is the notion of a weak factorization system.

\begin{definition}
  Let $C$ be a category and let $i : A \to B$ and $p : X \to Y$ two morphisms of $C$.
  We say that $i$ has the \emph{left lifting property} with respect to $p$,
  or that $p$ has the \emph{right lifting property} with respect to $i$,
  if every (solid) commutative square of the form
  \[
    \begin{tikzcd}
      A \ar[r] \ar[d, "i"'] & X \ar[d, "p"] \\
      B \ar[r] \ar[ru, dotted, "\exists l"] & Y
    \end{tikzcd}
  \]
  admits a lift $l : B \to X$ making both triangles commute, as indicated by the dotted arrow.
\end{definition}

\begin{notation}
  For classes of morphisms $\sL$, $\sR$ of a category $C$, we write $\llp(\sR)$ for the class of morphisms $i$ of $C$ with the left lifting property with respect to every $p \in \sR$, and $\rlp(\sL)$ for the class of morphisms $p$ of $C$ with the right lifting property with respect to every $i \in \sL$.
\end{notation}

\begin{definition}
  Let $C$ be a category.
  A \emph{weak factorization system} on $C$ is a pair $(\sL, \sR)$ of classes of morphisms of $C$ such that:
  \begin{enumerate}
  \item $\sL = \llp(\sR)$ and $\sR = \rlp(\sL)$.
  \item Each morphism of $C$ admits a factorization as a morphism in $\sL$ followed by a morphism in $\sR$.
  \end{enumerate}
  We call $\sL$ the \emph{left class} and $\sR$ the \emph{right class} of the weak factorization system.
\end{definition}

\begin{remark}
  The operations $\llp$ and $\rlp$ are inclusion-reversing,
  that is, if $\sL_1 \subset \sL_2$ then $\rlp(\sL_1) \supset \rlp(\sL_2)$ and similarly for $\llp$.
  It follows that if $(\sL_1, \sR_1)$ and $(\sL_2, \sR_2)$ are two weak factorization systems on the same category $C$,
  then $\sL_1 \subset \sL_2$ if and only if $\sR_1 \supset \sR_2$.
\end{remark}

\begin{convention}\label{conv:wfs-order}
  We regard the collection of all weak factorization systems on $C$ as ordered by inclusion of the \emph{left} class.
  That is, we define the ordering $\le$ by
  \[
    (\sL_1, \sR_1) \le (\sL_2, \sR_2) \iff \sL_1 \subset \sL_2 \iff \sR_1 \supset \sR_2.
  \]
  This relation $\le$ is a partial ordering on weak factorization systems on $C$.
\end{convention}

\begin{remark}
  The axioms for a weak factorization system are self-dual.
  That is, if $(\sL, \sR)$ is a weak factorization system on $C$, then $(\sR^\op, \sL^\op)$ is a weak factorization system on $C^\op$, where $\sR^\op$ (respectively $\sL^\op$) denotes $\sR$ (respectively $\sL$) viewed as a class of morphisms in $C^\op$.
  This relationship lets us transform theorems about the left class of a weak factorization system into dual theorems about the right class and vice versa.
  We will generally write out the full statements for the left class and leave the formulation of the dual statements to the reader.
\end{remark}

\begin{notation}
  Let $C$ be a category (to be inferred from context).
  We will write
  \begin{itemize}
  \item $\Iso$ for the class of isomorphisms of $C$;
  \item $\All$ for the class of all morphisms of $C$.
  \end{itemize}
\end{notation}

\begin{example}
  One easily verifies that $(\Iso, \All)$ and $(\All, \Iso)$ are weak factorization systems on $C$ for any $C$.
  Evidently, these are the minimal and maximal weak factorization systems respectively:
  that is, for any weak factorization system $(\sL, \sR)$, we have $(\Iso, \All) \le (\sL, \sR) \le (\All, \Iso)$.
\end{example}

\begin{example}\label{ex:set-mono-epi}
  As a less trivial example, in the category $\Set$, we will write
  \begin{itemize}
  \item $\Mono$ for the class of monomorphisms, i.e., injective functions;
  \item $\Epi$ for the class of epimorphisms, i.e., surjective functions.
  \end{itemize}
  Then one verifies that
  \begin{itemize}
  \item a function is injective if and only if it has the left lifting property with respect to all surjective functions;
  \item a function is surjective if and only if it has the right lifting property with respect to all injective functions;
  \item any function $f : X \to Y$ can be written as the composition of an injective function and a surjective function,
    for example via the factorization $X \to X \amalg Y \to Y$.
  \end{itemize}
  Therefore $(\Mono, \Epi)$ is a weak factorization system on $\Set$.
  This example will play an important role:
  the category $\Set$ equipped with the weak factorization system $(\Mono, \Epi)$ turns out to be a kind of ``unit object'' (in a sense that will become clear later).
\end{example}

\begin{remark}\label{remark:trivial}
  Of course, the main source of interesting weak factorization systems is from model categories.
  As we will review shortly, each model category gives rise to two weak factorization systems.
  In fact, provided that the ambient category $C$ is complete and cocomplete, every weak factorization system arises from a model category.
  Indeed, a model category on $C$ in which every morphism is a weak equivalence is precisely the same thing as a weak factorization system on $C$.
  We will refer to such model category structures as \emph{trivial}%
  \footnote{
    Some authors use the term ``trivial'' for model categories in which the weak equivalences are just the isomorphisms.
    We prefer to call such model categories \emph{discrete}.}
  because they are the model categories with trivial homotopy category.

  The upshot is that some of the basic properties of model categories are really facts about weak factorization systems and, conversely, questions about weak factorization systems can be ``reduced'' to questions about trivial model categories.
  In particular, proofs of the following facts can be found in any text on model categories.
\end{remark}

\begin{proposition}
  Let $(\sL, \sR)$ be a weak factorization system.
  Then $\sL$ is closed under coproducts, pushouts, transfinite compositions and retracts, and contains all isomorphisms, and dually for $\sR$.
\end{proposition}

\begin{proposition}\label{prop:wfs-of-retract}
  Let $\sL$ and $\sR$ be two classes of morphisms of $C$ such that
  \begin{enumerate}[(i)]
  \item each morphism of $\sL$ has the left lifting property with respect to each morphism of $\sR$;
  \item each morphism of $C$ admits a factorization as a morphism of $\sL$ followed by a morphism of $\sR$;
  \item $\sL$ and $\sR$ are each closed under retracts.
  \end{enumerate}
  Then $(\sL, \sR)$ is a weak factorization system on $C$.
\end{proposition}

\begin{proof}
  The only conditions remaining to be checked are that $\llp(\sR) \subset \sL$ and $\rlp(\sL) \subset \sR$.
  If $f \in \llp(\sR)$, write $f = gh$ with $g \in \sR$, $h \in \sL$.
  By the ``retract argument'' \cite[Proposition~7.2.2]{Hi}, $f$ is a retract of $h$ and so by assumption $f \in \sL$.
  The proof that $\rlp(\sL) \subset \sR$ is dual.
\end{proof}

\begin{definition}
  We say that a class $I$ of morphisms of $M$ \emph{generates} a weak factorization system $(\sL, \sR)$ if $\sR = \rlp(I)$.
\end{definition}

\begin{notation}
  For two categories $C$ and $D$, we write $F : C \adj D : G$ to mean the data of an adjunction between two functors $F : C \to D$ and $G : D \to C$, with $F$ the left and $G$ the right adjoint.
\end{notation}

\begin{proposition}
  Let $C$ and $D$ be two categories equipped with weak factorization systems $(\sL_C, \sR_C)$ and $(\sL_D, \sR_D)$, respectively, and suppose $F : C \adj D : G$ is an adjunction.
  Then the following are equivalent:
  \begin{enumerate}
  \item $F$ sends morphisms of $\sL_C$ to morphisms of $\sL_D$;
  \item $G$ sends morphisms of $\sR_D$ to morphisms of $\sR_C$.
  \end{enumerate}
  If $I$ is a class which generates $(\sL_C, \sR_C)$, then we may add a third equivalent condition:
  \begin{enumerate}
  \item[(3)] $F$ sends morphisms of $I$ to morphisms of $\sL_D$.
  \end{enumerate}
\end{proposition}

\begin{remark}
  In the situation of the preceding Proposition, suppose that $C$ and $D$ are equal as categories but equipped with possibly different weak factorization systems, and let $F : C \adj D : G$ be the identity adjunction.
  Then the equivalent conditions of the Proposition are precisely the definition of the relation $(\sL_C, \sR_C) \le (\sL_D, \sR_D)$.
\end{remark}
  
\begin{convention}\label{conv:adj-direction}
  We will always think of an adjunction $F : C \adj D : G$ as a morphism in the direction of its \emph{left} adjoint $F : C \to D$.
  \Cref{conv:wfs-order} is chosen to be compatible with this convention, the preceding remark, and the usual convention of regarding a partially ordered set as a category in which there is a (unique) morphism $a \to b$ if and only if $a \le b$.
\end{convention}

\subsection{Model categories and premodel categories}

\begin{definition}
  Let $\sW$ be a class of morphisms of a category $C$.
  We say that $\sW$ satisfies the \emph{two-out-of-three condition} if for any morphisms $f : X \to Y$ and $g : Y \to Z$ of $C$, whenever two of $f$, $g$ and $g \circ f$ belong to $\sW$, so does the third.
\end{definition}

\begin{definition}
  Let $M$ be a complete and cocomplete category.
  A \emph{model category structure} on $M$ consists of three classes of morphisms $\sW$, $\sC$, and $\sF$ such that:
  \begin{enumerate}
  \item $(\sC, \sF \cap \sW)$ and $(\sC \cap \sW, \sF)$ are weak factorization systems;
  \item $\sW$ satisfies the two-out-of-three condition.
  \end{enumerate}
  We call $\sC$ the \emph{cofibrations}, $\sF$ the \emph{fibrations} and $\sW$ the \emph{weak equivalences} of the model category structure.
  A \emph{model category} is a complete and cocomplete category equipped with a model category structure.
\end{definition}

\begin{remark}
  We have given an ``optimized'' definition of a model category, apparently due to Joyal and Tierney \cite{JT}.
  Compared to the traditional definition (as presented in \cite[Definition~7.1.3]{Hi}, for example), two conditions may appear to be missing:
  \begin{itemize}
  \item We did not explicitly require $\sW$ to contain the isomorphisms of $M$.
    However, $\sW$ contains $\sC \cap \sW$ which is the left class of a weak factorization system on $M$ and therefore contains all isomorphisms.
  \item We did not explicitly require $\sW$ to be closed under retracts.
    However, this also follows from the given axioms.
    The proof is not trivial; see \cite[Proposition~7.8]{JT}.
  \end{itemize}
\end{remark}

\begin{remark}
  We have chosen not to require the existence of functorial factorizations.
  This choice is not essential, as we will eventually specialize to the combinatorial setting, in which the existence of functorial factorizations is automatic anyways.
  The reader should feel free to assume that all weak factorization systems that appear admit functorial factorizations, which simplifies a few arguments, at the cost of some generality in the earlier parts of the theory.
\end{remark}

\begin{remark}
  The data making up a model category structure is \emph{overdetermined}, in the following sense.
  Suppose the complete and cocomplete category $M$ is equipped with two classes of maps $\sC$ and $\sF$.
  For a model category structure on $M$ with these prescribed classes of cofibrations and fibrations to exist, the following conditions must be satisfied.
  \begin{enumerate}
  \item $\sC$ must be the left class of a weak factorization system $(\sC, \sAF)$ and $\sF$ must be the right class of a weak factorization system $(\sAC, \sF)$.
  \item By the two-out-of-three condition, the class $\sW$ of weak equivalences of $M$ must be precisely the class of maps which can be expressed as a composition of a map of $\sAC$ followed by a map of $\sAF$.
  \end{enumerate}
  Under these conditions, one can show (using the ``retract argument'') that we do have the required equalities $\sAC = \sC \cap \sW$ and $\sAF = \sF \cap \sW$.
  However, there is no reason in general that a class $\sW$ defined in this way should satisfy the two-out-of-three condition.
  For example, $(\sC, \sAF)$ could be the weak factorization system $(\All, \Iso)$, so that $\sW = \sAC$, and of course the left class of a weak factorization system rarely satisfies the two-out-of-three condition.

  This feature of the notion of a model category accounts for the difficulty in constructing new model categories from old ones.
  We will return to this point later in this section; for now, the point is that the two-out-of-three condition is the main obstruction to an algebraically well-behaved theory of model categories.

  This motivates our main definition.
\end{remark}

\begin{definition}
  Let $M$ be a complete and cocomplete category.
  A \emph{premodel category structure} on $M$ is a pair of weak factorization systems $(\sC, \sAF)$ and $(\sAC, \sF)$ on $M$ such that $\sAC \subset \sC$.
  We call
  \begin{itemize}
  \item $\sC$ the \emph{cofibrations} and $\sF$ the \emph{fibrations} of $M$;
  \item $\sAC$ the \emph{anodyne cofibrations} and $\sAF$ the \emph{anodyne fibrations} of $M$.
  \end{itemize}
\end{definition}

\begin{notation}
  As is standard for model categories, we will denote cofibrations and fibrations by arrows $\cto$ and $\fto$ respectively.
  For \emph{anodyne} cofibrations and fibrations, it could be misleading to use arrows $\acto$ and $\afto$ since we do not assume any form of two-out-of-three condition.
  Instead, we denote anodyne cofibrations and fibrations by arrows $\ancto$ and $\anfto$ respectively.
\end{notation}

\begin{example}
  Any model category structure on $M$ determines a premodel category structure on $M$ with the same cofibrations and fibrations.
  In fact, since a model category structure is determined by its cofibrations and fibrations, we may think of a model category structure as a special kind of premodel category structure---one in which the class of morphisms which can be factored as an anodyne cofibration followed by an anodyne fibration satisfies the two-out-of-three condition.
\end{example}

\begin{notation}
  We write $\Kan$ for the category of simplicial sets equipped with its standard (Kan--Quillen) model category structure.
  By the previous example, we may also regard $\Kan$ as a premodel category.
\end{notation}

\begin{notation}
  For a premodel category $M$,
  \begin{itemize}
  \item we call an object $A$ \emph{cofibrant} if the unique map from the initial object to $A$ is a cofibration, and write $M^\cof$ for the full subcategory of $M$ on the cofibrant objects;
  \item we call an object $X$ \emph{fibrant} if the unique map from $X$ to the final object is a fibration, and write $M^\fib$ for the full subcategory of $M$ on the fibrant objects;
  \item we write $M^\cf$ for the full subcategory of $M$ on the objects which are both cofibrant and fibrant.
  \end{itemize}
  This terminology is consistent with the usual terminology for model categories under the above identification of model categories as particular premodel categories.
\end{notation}

Like the axioms of a model category, the axioms of a premodel category are self-dual in a way which interchanges the two factorization systems.

\begin{definition}
  Let $M$ be a premodel category.
  Then $M^\op$ also has the structure of a premodel category.
  A morphism in $M^\op$ is a cofibration (respectively anodyne cofibration, fibration, anodyne fibration) if and only the corresponding morphism in $M$ is a fibration (respectively anodyne fibration, cofibration, anodyne cofibration).
\end{definition}

In order to obtain a theory with good algebraic properties, we will eventually need to restrict to combinatorial premodel categories.
We will say much more about these in \cref{chap:algebra}, but we give the definition now in order to give previews of the algebraic structure of combinatorial premodel categories throughout this chapter.

\begin{definition}
  A premodel category $M$ is \emph{combinatorial} if
  \begin{enumerate}
  \item the underlying category $M$ is locally presentable;
  \item there exist \emph{sets} of morphisms $I$ and $J$ of $M$ such that $\sAF = \rlp(I)$ and $\sF = \rlp(J)$.
  \end{enumerate}
  We call $I$ and $J$ \emph{generating cofibrations} and \emph{generating anodyne cofibrations} of $M$ respectively.
\end{definition}

\begin{remark}
  The notion of combinatorial premodel category is \emph{not} self-dual.
  In fact, the opposite of a locally presentable category is never locally presentable unless the category is a poset \cite[Theorem~1.64]{AR}.
\end{remark}

Suppose $M$ is any locally presentable category.
Then any sets of morphisms $I$ and $J$ such that $\rlp(I) \subset \rlp(J)$ generate a unique combinatorial premodel category structure on $M$; the required factorizations may be constructed by the small object argument.
This trivial ``existence theorem for combinatorial premodel categories'' is part of the reason that combinatorial premodel categories have much better algebraic structure than combinatorial model categories.

\begin{example}\label{ex:premodel-set}
  The central example of a premodel category is the category $\Set$ equipped with the premodel category structure in which
  \begin{itemize}
  \item $(\sC, \sAF)$ is the weak factorization system $(\Mono, \Epi)$ of \cref{ex:set-mono-epi};
  \item $(\sAC, \sF) = (\Iso, \All)$.
  \end{itemize}
  It is combinatorial; we may take $I = \{\emptyset \to *\}$ and $J = \emptyset$.
  It is not a model category; the weak equivalences would have to be $\sAF = \Epi$, but these do not satisfy the two-out-of-three condition.
  We will denote this premodel category simply by $\Set$.
  It turns out to be the ``unit combinatorial premodel category''.
\end{example}


\subsection{Quillen adjunctions}

The definition of a Quillen adjunction between model categories does not directly involve the weak equivalences, only the (acyclic) cofibrations and (acyclic) fibrations.
Therefore, it transfers without difficulty to the setting of premodel categories.

\begin{definition}
  An adjunction $F : M \adj N : G$ between two premodel categories is a \emph{Quillen adjunction} if it satisfies the following two conditions:
  \begin{enumerate}
  \item $F$ preserves cofibrations, or equivalently, $G$ preserves anodyne fibrations.
  \item $F$ preserves anodyne cofibrations, or equivalently, $G$ preserves fibrations.
  \end{enumerate}
  In this situation we also call $F$ a \emph{left Quillen functor} and $G$ a \emph{right Quillen functor}.
\end{definition}

\begin{remark}
  This definition is compatible with the usual definition of a Quillen adjunction between model categories under the identification of model categories as particular premodel categories.
  In particular, any Quillen adjunction between model categories is also an example of a Quillen adjunction between premodel categories.
\end{remark}

\begin{remark}
  As usual, if $I$ and $J$ are sets (or even classes) of maps of $M$ such that $\sAF_M = \rlp(I)$ and $\sF_M = \rlp(J)$, then in order to check that an adjunction $F : M \adj N : G$ is a Quillen adjunction, it suffices to verify that $F$ sends the maps of $I$ to cofibrations of $N$ and the maps of $J$ to anodyne cofibrations of $N$.
\end{remark}

We now turn to the 2-categorical structure of premodel categories.

\begin{definition}
  For premodel categories $M$ and $N$, we define:
  \begin{itemize}
  \item $\QAdj(M, N)$ to be the category whose objects are Quillen adjunctions $F : M \adj N : G$ in which a morphism from $F : M \adj N : G$ to $F' : M \adj N : G'$ is given by a natural transformation $F \to F'$.
  \item $\LQF(M, N)$ to be the full subcategory of the category of functors from $M$ to $N$ on the left Quillen functors.
  \item $\RQF(M, N)$ to be the full subcategory of the category of functors from $M$ to $N$ on the right Quillen functors.
  \end{itemize}
\end{definition}

\begin{remark}
  If $F : M \adj N : G$ and $F' : M \adj N : G'$ are two adjunctions, then each natural transformation $\alpha : F \to F'$ corresponds to a unique natural transformation $\beta : G' \to G$ such that the square below commutes for every choice of objects $X$ of $M$ and $Y$ of $N$.
  \[
    \begin{tikzcd}
      \Hom(F'X, Y) \ar[r, "\sim"] \ar[d, "(\alpha_X)^*"'] & \Hom(X, G'Y) \ar[d, "(\beta_Y)_*"] \\
      \Hom(FX, Y) \ar[r, "\sim"] & \Hom(X, GY)
    \end{tikzcd}
  \]
  Then there are equivalences of categories
  \[
    \begin{tikzcd}
      & \QAdj(M, N) \ar[dl, "\eqv"'] \ar[dr, "\eqv"] & \\
      \LQF(M, N) & & \RQF(N, M)^\op
    \end{tikzcd}
  \]
  which send a Quillen adjunction $F : M \adj N : G$ to $F$ and $G$ respectively.
  By convention, we will primarily work with $\LQF(M, N)$; these equivalences allow us to replace it by $\QAdj(M, N)$ or $\RQF(N, M)^\op$ where convenient.
\end{remark}

\begin{definition}
  The (strict) \emph{2-category of premodel categories} $\PM$ has as objects premodel categories and, for premodel categories $M$ and $N$, the category $\LQF(M, N)$ as the category of morphisms from $M$ to $N$.

  The (strict) \emph{2-category of combinatorial premodel categories} $\CPM$ is the full sub-2-category of $\PM$ containing the premodel categories which are combinatorial.
\end{definition}

\begin{remark}\label{remark:pmr}
  In accordance with \cref{conv:adj-direction}, we always regard an adjunction as a morphism in the direction of its left adjoint.
  We may also define alternative 2-categories of premodel categories
  \begin{itemize}
  \item $\PMA$, with morphism category from $M$ to $N$ given by $\QAdj(M, N)$;
  \item $\PMR$, with morphism category from $M$ to $N$ given by $\RQF(N, M)^\op$.
  \end{itemize}
  By the preceding remark, these 2-categories are related by 2-equivalences
  \[
    \begin{tikzcd}
      & \PMA \ar[dl, "\approx"'] \ar[dr, "\approx"] & \\
      \PM & & \PMR
    \end{tikzcd}
  \]
  which allow us to replace $\PM$ by $\PMA$ or $\PMR$ where convenient.
\end{remark}

\begin{example}
  Let $N$ be a premodel category.
  We will describe left Quillen functors $F : \Set \to N$, where as usual $\Set$ carries the premodel category structure described in \cref{ex:premodel-set}.
  A left Quillen functor $F : \Set \to N$ is in particular a left adjoint, so it is determined up to unique isomorphism by the object $F(*)$ of $N$.
  In order for $F$ to be a left Quillen functor, it must also preserve cofibrations and anodyne cofibrations.
  Recall that $\Set$ has generating cofibrations $I = \{\emptyset \to *\}$ and $J = \emptyset$.
  By the previous remark, it suffices to check that $F$ sends the cofibration $\emptyset \to *$ of $\Set$ to a cofibration $\emptyset = F(\emptyset) \to F(*)$ of $M$.
  In other words, $F(*)$ must be a cofibrant object of $N$.

  We conclude that left Quillen functors $F : \Set \to N$ are the same as cofibrant objects of $N$, or more precisely, the full subcategory of the category of functors from $\Set \to N$ on the left Quillen functors is equivalent to the full category $N^\cof$ of cofibrant objects of $N$.
\end{example}

\begin{remark}
  We will see later that for any combinatorial premodel categories $M$ and $N$, the category of \emph{all} left adjoints from $M$ to $N$ admits a premodel category structure $\CPM(M, N)$ whose cofibrant objects are precisely the left Quillen functors.
  In the case $M = \Set$, the category of left adjoints from $M$ to $N$ can be identified with $N$ itself and the premodel category structure in question is (as one might guess from the preceding example) just the original premodel category structure on $N$.
  This is one manifestation of the role that $\Set$ plays as the unit combinatorial premodel category.
\end{remark}

\begin{remark}
  Let $V$ be a symmetric monoidal category with unit object $1_V$ and $X$ an object of $V$.
  In some contexts, it is appropriate to think of $\Hom_V(1_V, X)$ as the ``underlying set'' of the object $X$.
  For example if $C$ is a $V$-enriched category, then the underlying ordinary category of $C$ is constructed by applying $\Hom_V(1_V, -)$ to each $V$-valued Hom object of $C$.

  If we apply this prescription to $\CPM$ and its unit object $\Set$, we are led to conclude that the ``underlying category'' of a combinatorial premodel category $M$ should be $\THom_\CPM(\Set, M) \eqv M^\cof$, rather than $M$ itself.
  This terminology would obviously be too confusing, and we instead use the phrase ``underlying category'' in its usual sense.
  However, this idea makes sense in some contexts.
  For example, it explains why the objects of the internal Hom $\CPM(M, N)$ of combinatorial premodel categories are not left Quillen functors but actually \emph{all} left adjoints; only the cofibrant objects of $\CPM(M, N)$ are left Quillen functors.
  The homotopy theory of a relaxed premodel category $M$ which we will develop in the next chapter is also defined entirely in terms of $M^\cof$.
\end{remark}

\section{Monoidal premodel categories}

Model category theory has a ``multiplicative structure'' which allows us to define monoidal model categories, enriched model categories, and so on.
The core underlying concept is that of a Quillen bifunctor, which generalizes straightforwardly to the setting of premodel categories.

\subsection{Quillen bifunctors}

We first review some requisite category theory.
A reference for this material is \cite[section~4.1]{Ho}.

\begin{definition}
  Let $C_1$, $C_2$ and $D$ be categories.
  An \emph{adjunction of two variables} from $(C_1, C_2)$ to $D$ consists of
  \begin{enumerate}
  \item a functor $F : C_1 \times C_2 \to D$,
  \item functors $G_1 : C_2^\op \times D \to C_1$ and $G_2 : C_1^\op \times D \to C_2$, and
  \item isomorphisms
    \[
      \Hom(F(X_1, X_2), Y) \simeq \Hom(X_1, G_1(X_2, Y)) \simeq \Hom(X_2, G_2(X_1, Y))
    \]
    natural in $X_1$, $X_2$ and $Y$.
  \end{enumerate}
  In other words,
  \begin{itemize}
  \item $G_1(X_2, -)$ is right adjoint to $F(-, X_2)$ for each $X_2$ in $C_2$, and for each map $X_2 \to X'_2$ the natural transformation $G_1(X'_2, -) \to G_1(X_2, -)$ is the one determined by $F(-, X_2) \to F(-, X'_2)$;
  \item $G_2(X_1, -)$ is right adjoint to $F(X_1, -)$ for each $X_1$ in $C_1$, with the analogous condition on the structure maps of $G_2$.
  \end{itemize}
  We call $F$ the \emph{left part} of the adjunction of two variables $(F, G_1, G_2)$.
  We see that $G_1$ and $G_2$ along with the adjunction isomorphisms are uniquely determined up to unique isomorphism by $F$, so by an abuse of language we will often identify the adjunction of two variables by $F$ alone.
\end{definition}

If $F : C_1 \times C_2 \to D$ is an adjunction of two variables, then each functor $F(X_1, -)$ and $F(-, X_2)$ is a left adjoint.
In particular, $F$ preserves colimits in each variable separately.

\begin{definition}
  Let $C_1$, $C_2$ and $D$ be categories, with $D$ admitting pushouts, and let $F : C_1 \times C_2 \to D$ be a functor.
  For maps $f_1 : A_1 \to B_1$ in $C_1$ and $f_2 : A_2 \to B_2$ in $C_2$, we define the \emph{$F$-pushout product} $f_1 \bp_{F} f_2$ to be the morphism of $D$ indicated by the dotted arrow
  \[
    \begin{tikzcd}
      F(A_1, A_2) \ar[r] \ar[d] & F(A_1, B_2) \ar[d] \ar[rdd] \\
      F(B_1, A_2) \ar[r] \ar[rrd] & \cdot \ar[rd, dotted] \\
      & & F(B_1, B_2)
    \end{tikzcd}
  \]
  where the square is a pushout, so that
  \[
    f_1 \bp_{F} f_2 : F(B_1, A_2) \amalg_{F(A_1, A_2)} F(A_1, B_2) \to F(B_1, B_2).
  \]
  Usually $F$ will be an adjunction of two variables.
  We will omit $F$ from the notation $\bp_F$ when it is clear from context.

  For classes of morphisms $K_1$, $K_2$ of $C_1$, $C_2$ respectively, we also write
  \[
    K_1 \bp_F K_2 = \{\,f_1 \bp_F f_2 \mid f_1 \in K_1, f_2 \in K_2\,\}.
  \]
\end{definition}

\begin{definition}
  Let $M_1$, $M_2$ and $N$ be three premodel categories.
  An adjunction of two variables $(F : M_1 \times M_2 \to N, G_1, G_2)$ is a \emph{Quillen adjunction of two variables} if whenever $f_1 : A_1 \to B_1$ and $f_2 : A_2 \to B_2$ are cofibrations in $M_1$ and $M_2$ respectively, $f_1 \bp_F f_2$ is a cofibration in $N$ which is an anodyne cofibration if either $f_1$ or $f_2$ is.
  A functor $F : M_1 \times M_2 \to N$ is a \emph{Quillen bifunctor} if it is the left part of a Quillen adjunction of two variables.
\end{definition}

If $F : M_1 \times M_2 \to N$ is a Quillen functor, then for any \emph{cofibrant} object $A_1$ of $M_1$ or $A_2$ of $M_2$, the functor $F(A_1, -)$ or $F(-, A_2)$ is a left Quillen functor.

The following basic fact about Quillen bifunctors is proved the same way as in the setting of model categories \cite[Lemma~4.2.2 and Corollary~4.2.5]{Ho}.

\begin{proposition}\label{prop:quillen-bifunctor}
  For an adjunction of two variables $(F : M_1 \times M_2 \to N, G_1, G_2)$, the following conditions are equivalent:
  \begin{enumerate}
  \item $(F, G_1, G_2)$ is a Quillen bifunctor.
  \item For any cofibration $f_2 : A_2 \to B_2$ in $M_2$ and fibration $g : X \to Y$ in $N$, the dotted map to the pullback
    \[
      \begin{tikzcd}
        G_1(B_2, X) \ar[rrd] \ar[rdd] \ar[rd, dotted] \\
        & \cdot \ar[r] \ar[d] & G_1(A_2, X) \ar[d] \\
        & G_1(B_2, Y) \ar[r] & G_1(A_2, Y)
      \end{tikzcd}
    \]
    is a fibration in $M_1$ which is an anodyne fibration if either $f_2$ is an anodyne cofibration or $g$ is an anodyne fibration.
  \item For any cofibration $f_1 : A_1 \to B_1$ in $M_1$ and fibration $g : X \to Y$ in $N$, the dotted map to the pullback
    \[
      \begin{tikzcd}
        G_2(B_1, X) \ar[rrd] \ar[rdd] \ar[rd, dotted] \\
        & \cdot \ar[r] \ar[d] & G_2(A_1, X) \ar[d] \\
        & G_2(B_1, Y) \ar[r] & G_2(A_1, Y)
      \end{tikzcd}
    \]
    is a fibration in $M_2$ which is an anodyne fibration if either $f_1$ is an anodyne cofibration or $g$ is an anodyne fibration.
  \end{enumerate}
  Furthermore, suppose that for $i = 1$ and $2$, $M_i$ has generating cofibrations $I_i$ and generating anodyne cofibrations $J_i$ (which can be classes).
  Then we may add a fourth equivalent condition:
  \begin{enumerate}
  \item[(4)] $I_1 \bp_F I_2 \subset \sC_N$ and $J_1 \bp_F I_2 \cup I_1 \bp_F J_2 \subset \sAC_N$.
  \end{enumerate}
\end{proposition}

\begin{notation}
  We write $\QBF((M_1, M_2), N)$ for the full subcategory of the category of functors $M_1 \times M_2 \to N$ on the Quillen bifunctors.
\end{notation}

\begin{example}\label{ex:qbf-set}
  Let $M$ and $N$ be premodel categories.
  We will determine the Quillen bifunctors $\Set \times M \to N$.
  An adjunction of two variables $F : \Set \times M \to N$ is determined up to unique isomorphism by the functor $F(*, -) : M \to N$, which must be a left adjoint.
  Write $i$ for the unique map $\emptyset \to *$ of $\Set$.
  Then $\Set$ has generating cofibrations $\{i\}$ and generating anodyne cofibrations $\emptyset$, so by the last part of the preceding proposition, $F$ is a Quillen bifunctor if and only if
  \[
    \{i\} \bp_F \sC_M \subset \sC_N \quad\hbox{and}\quad
    \{i\} \bp_F \sAC_M \subset \sAC_N.
  \]
  Now if $f : A \to B$ is any morphism of $M$, then we may compute $i \bp_F f$ by forming the diagram
  \[
    \begin{tikzcd}
      F(\emptyset, A) \ar[r] \ar[d] & F(*, A) \ar[d] \ar[rdd] \\
      F(\emptyset, B) \ar[r] \ar[rrd] & F(*, A) \ar[rd] \\
      & & F(*, B)
    \end{tikzcd}
  \]
  in which the top left square is a pushout because both $F(\emptyset, A)$ and $F(\emptyset, B)$ are initial.
  We see that $i \bp_F f = F(*, f)$.
  Thus, $F$ is a Quillen bifunctor if and only if $F(*, -)$ is a Quillen functor, so we have $\QBF((\Set, M), N) \eqv \LQF(M, N)$.
\end{example}

\begin{remark}
  More generally, for any $n \ge 0$ we may define a notion of $n$-ary \emph{Quillen multifunctor} from an $n$-tuple $(M_1, \dots, M_n)$ of premodel categories to another premodel category $N$.
  For $n = 2$ we recover the notion of a Quillen bifunctor, and for $n = 1$, the notion of a left Quillen functor; and a $0$-ary Quillen multifunctor $() \to N$ is a \emph{cofibrant} object of $N$.
  These Quillen multifunctors assemble into a $\Cat$-valued operad, or 2-multicategory.
  Of course, this construction does not require premodel categories; it works equally well for model categories.
  
  However, a new phenomenon in the setting of premodel categories is that, once we restrict attention to combinatorial premodel categories, Quillen multifunctors become \emph{representable}.
  That is, there is a combinatorial premodel category $M_1 \otimes \dots \otimes M_n$ equipped with a universal Quillen multifunctor $M_1 \times \cdots \times M_n \to M_1 \otimes \cdots \otimes M_n$, inducing an equivalence between the category $\LQF(M_1 \otimes \dots \otimes M_n, N)$ and the category of Quillen multifunctors from $(M_1, \dots, M_n)$ to $N$.
  In the previous example, we effectively computed that $\Set \otimes M \eqv M$, so that $\Set$ is the unit object for the tensor product of combinatorial premodel categories.
  Then $\Set$ must also be the tensor product of zero factors, so we see that a $0$-ary Quillen multifunctor $() \to N$ is the same as a left Quillen functor from $\Set$ to $N$, which we saw earlier amounts to a cofibrant object of $N$, justifying the claim about $0$-ary Quillen multifunctors in the previous paragraph.
  Note that $\Set$ is not a model category, so this is one example of how generalizing from model categories to premodel categories results in an algebraically better behaved theory.
\end{remark}

\subsection{Monoidal and enriched premodel categories}

\begin{definition}\label{def:monoidal-premodel-cat}
  A \emph{monoidal premodel category} is a premodel category $V$ which is also a monoidal category $(V, \otimes, 1_V)$ such that:
  \begin{enumerate}
  \item the tensor product $\otimes : V \times V \to V$ is a Quillen bifunctor;
  \item the unit object $1_V$ of $V$ is cofibrant.
  \end{enumerate}
  A \emph{symmetric} (or \emph{braided}) \emph{monoidal premodel category} is a monoidal premodel category which symmetric (or braided) as a monoidal category.
\end{definition}

\begin{remark}
  Some authors use ``monoidal model category'' as a synonym for ``symmetric monoidal model category''.
  Although we will not have much real need for monoidal premodel categories which are not symmetric monoidal, we prefer to distinguish the two notions in order to clarify what degree of monoidal structure is required for each particular argument.
\end{remark}

\begin{remark}
  There are two different definitions of monoidal model category in the literature.
  One is analogous to the one we give above.
  The other, used for example in \cite{Ho}, replaces the condition on the unit object by the \emph{unit axiom}:
  \begin{enumerate}
  \item[(2$'$)] for any cofibrant replacement $Q \to 1_V$ of the unit object and any cofibrant object $X$, the induced map $Q \otimes X \to 1_V \otimes X \iso X$ is a weak equivalence.
  \end{enumerate}
  (This condition is automatically satisfied if the unit object $1_V$ is cofibrant, since $- \otimes X$ is a left Quillen functor for cofibrant $X$ and therefore preserves the weak equivalence $Q \to 1_V$ between cofibrant objects.)

  Requiring the unit object to be cofibrant fits better into our algebraic framework; \cref{def:monoidal-premodel-cat} makes a monoidal premodel category precisely a pseudomonoid object in the 2-multicategory of premodel categories and Quillen multifunctors.
  Moreover, we do not yet have any notion of weak equivalence in a premodel category, so we cannot even state the alternative unit axiom.
\end{remark}

\begin{example}
  A monoidal model category (with cofibrant unit) is a monoidal premodel category.
  So, for example, $\Kan$ (with its cartesian monoidal structure) is a symmetric monoidal premodel category.
\end{example}

\begin{example}
  We equip the premodel category $\Set$ with the cartesian monoidal structure.
  The product $\times : \Set \times \Set \to \Set$ is an adjunction of two variables because $\Set$ is cartesian closed.
  Moreover, it is a Quillen bifunctor: by \cref{ex:qbf-set}, it suffices to verify that $* \times - : \Set \to \Set$ is a left Quillen functor, and it is the identity functor.
  The unit object $*$ (and indeed every object) of $\Set$ is cofibrant, so $\Set$ is a symmetric monoidal premodel category.
\end{example}

Now for a monoidal premodel category $V$, there is a notion of a \emph{$V$-premodel category}.
There are two essentially equivalent ways to define $V$-premodel categories.
\begin{itemize}
\item
  We may define a $V$-premodel category to be a $V$-enriched category $M$ which is tensored and cotensored over $V$, together with a premodel category structure on the underlying category of $M$ for which the tensor $\otimes : V \times M \to M$ is a Quillen bifunctor.
  In this approach, the compatibility of the tensor $\otimes : V \times M \to M$ with the monoidal structure of $V$ is automatically encoded in the structure of $M$ as a $V$-enriched category.
\item
  Alternatively, we may define a $V$-premodel category to be an ordinary premodel category which is a pseudomodule over $V$ in the 2-multicategory of premodel categories and Quillen multifunctors.
  Such a pseudomodule structure is given by a Quillen bifunctor $\otimes : V \times M \to M$ together with additional coherence data.
  This approach avoids enriched category theory, at the expense of some additional bookkeeping of this coherence data.
\end{itemize}
The second option fits better into our general algebraic approach.
We thus make the following sketch of a definition.
(Compare \cite[Definition~4.2.18]{Ho}, although we work with left modules rather than right ones, and always assume the unit of a monoidal premodel category is cofibrant.)

\begin{definition}
  Let $V$ be a monoidal premodel category.
  A \emph{$V$-premodel category} is a premodel category $M$ equipped with a Quillen bifunctor $\otimes : V \times M \to M$ which is coherently associative and unital with respect to the monoidal structure of $V$.
  In particular, $M$ is equipped with natural isomorphisms $(K \otimes L) \otimes X \iso K \otimes (L \otimes X)$ and $1_V \otimes X \iso X$, which are required to satisfy certain coherence conditions.
\end{definition}

The above definition notwithstanding, we also sometimes use the word ``enriched'' to refer to $V$-premodel categories in general, especially when $V$ is a model category.

\begin{example}
  If $V$ is a monoidal premodel category, then $V$ itself is also a $V$-premodel category, with the action $\otimes : V \times V \to V$ given by the monoidal structure of $V$.
\end{example}

\begin{example}
  Every premodel category is automatically a $\Set$-premodel category in an essentially unique way.
  In fact, we saw that a Quillen bifunctor $\otimes : \Set \times M \to M$ amounts to a left Quillen functor $* \otimes - : M \to M$, and since $*$ is the unit object of $\Set$ this latter functor must be (naturally isomorphic to) the identity.
\end{example}

\begin{example}
  When $V = \Kan$ we call a $V$-premodel category a \emph{simplicial premodel category}.
  Of course, simplicial model categories are examples of simplicial premodel categories.
\end{example}

$V$-premodel categories form a 2-category $V\PM$ whose 1-morphisms are left Quillen $V$-functors and whose 2-morphisms are $V$-natural transformations.
We defer discussion of these concepts to \cref{chap:modules}.
We write $V\CPM$ for the sub-2-category of $V\PM$ consisting of the $V$-premodel categories which are combinatorial.

For our current purposes the only relevant feature of a $V$-premodel category is the following.

\begin{definition}
  Let $V$ be any premodel category, not necessarily monoidal, but with a distinguished cofibrant ``unit'' object $I$.
  We say that a premodel category $M$ \emph{admits a unital action by $V$} if there exists a Quillen bifunctor $\otimes : V \times M \to M$ such that the left Quillen functor $I \otimes - : M \to M$ is naturally isomorphic to the identity.
\end{definition}

If $V$ is a monoidal premodel category, we may take $I$ to be the unit object of the monoidal structure on $V$.
Then any $V$-premodel category $M$ in particular admits a unital action by $V$.
Many of the technical advantages of enriched (e.g., simplicial) model categories over unenriched ones actually only depend on the existence of a unital action and for now we can disregard the rest of the structure of a $V$-premodel category.

\begin{remark}
  Unital actions of simplicial sets and related concepts sometimes arise naturally as well.
  For example, in Morel's homotopy theory of schemes \cite{Mor}, a \emph{$k$-space} is a presheaf of sets on the category of affine smooth schemes over $k$.
  The functor from $\Delta$ to $k$-spaces sending $[n]$ to the presheaf represented by $\Spec k[x_0, \ldots, x_n]/(x_0 + \cdots + x_n - 1)$ extends by colimits to a ``geometric realization'' functor $|{-}|$ from simplicial sets to $k$-spaces.
  Up to isomorphism, the geometric realization of $\Delta^n$ is $\mathbb{A}^n_k$.
  In turn, this defines an action of simplicial sets on $k$-spaces by the formula $K \otimes X = |K| \times X$.
  This action is unital because $|\Delta^0|$ is the terminal space $\Spec k$.
  However, it is not a monoidal action because $|\Delta^1| \times |\Delta^1|$ is the affine plane $\Spec k[x, y]$, while $|\Delta^1 \times \Delta^1| = |\Delta^2 \amalg_{\Delta^1} \Delta^2|$ consists of two planes glued (as $k$-spaces) along a line.
  This means that $k$-spaces do not acquire the structure of a category enriched, tensored and cotensored over simplicial sets.
  Morel calls the actual structure on $k$-spaces a \emph{quasi-simplicial structure}.
  See \cite[paragraph~2.1.3 and section~A.2.1]{Mor}.

\end{remark}

\subsection{Duality and bifunctors}

If $F : M \to N$ is a left Quillen functor with right adjoint $G : N \to M$, then $G^\op : N^\op \to M^\op$ is a left Quillen functor with right adjoint $F^\op : M^\op \to N^\op$.
Quillen bifunctors are also preserved by passage to opposite categories, in a way which we describe next.

Let $(F : M_1 \times M_2 \to N, G_1 : M_2^\op \times N \to M_1, G_2 : M_1^\op \times N \to M_2)$ be an adjunction of two variables.
Then $G_2^\op : M_1 \times N^\op \to M_2^\op$ is the left part of an adjunction of two variables
\begin{align*}
  G_2^\op & : M_1 \times N^\op \to M_2^\op, \\
  G_1 \circ \tau & : N \times M_2^\op \to M_1, \\
  F^\op & : M_1^\op \times M_2^\op \to N^\op
\end{align*}
where $\tau : N \times M_2^\op \to M_2^\op \times N$ swaps the two factors.
Indeed, the natural isomorphisms
\[
  \Hom_N(F(X_1, X_2), Y) \simeq \Hom_{M_1}(X_1, G_1(X_2, Y)) \simeq \Hom_{M_2}(X_2, G_2(X_1, Y))
\]
may be reinterpreted as natural isomorphisms
\[
  \Hom_{M_2^\op}(G_2(X_1, Y), X_2) \simeq \Hom_{M_1}(X_1, (G_1 \circ \tau)(Y, X_2))
  \simeq \Hom_{N^\op}(Y, F(X_1, X_2)).
\]
Moreover, the new adjunction of two variables is a Quillen adjunction of two variables if and only if the original was one.
This can be seen using the second equivalent condition of \cref{prop:quillen-bifunctor}.
The roles of $M_2$ and $N$ are now played by $N^\op$ and $M_2^\op$, so the roles of the cofibration $f_2 : A_2 \to B_2$ in $M_2$ and the fibration $g : X \to Y$ in $N$ are swapped, and thus the new condition on $G_1 \circ \tau$ is the same as the original condition on $G_1$.
We summarize this discussion below.

\begin{proposition}
  A Quillen bifunctor $F : M_1 \times M_2 \to N$ induces a Quillen bifunctor $F' : M_1 \times N^\op \to M_2^\op$ for which $F'(X_1, -) : N^\op \to M_2^\op$ is the opposite of the right adjoint of $F(X_1, -) : M_2 \to N$ for each object $A_1$ of $M_1$.
\end{proposition}

\begin{proof}
  This follows from the above argument, together with the fact that $G_2(X_1, -)$ is the right adjoint of $F(X_1, -)$ for each object $X_1$ of $M_1$.
\end{proof}

\begin{proposition}\label{prop:op-action}
  If the premodel category $M$ admits a unital action of $V$, then so does $M^\op$.
\end{proposition}

\begin{proof}
  The action $\otimes : V \times M \to M$ induces a left Quillen bifunctor $\otimes' : V \times M^\op \to M^\op$ and $I \otimes' -$ is naturally isomorphic to the identity since $I \otimes -$ is.
\end{proof}

\section{Basic constructions on premodel categories}

In this section, we show how to extend various constructions on model categories to the setting of premodel categories.
In fact, the constructions are easier for premodel categories, because we don't have to check anything related to the weak equivalences; all we have to do is construct two weak factorization systems.
The logical order of development would be to describe the constructions for premodel categories, and then prove that when the inputs are model categories, the outputs are as well.
However, the constructions for model categories are already well-known and we would like to reuse them here.
We will employ two strategies for doing so.
\begin{itemize}
\item
  In order to construct a model category, we must in particular produce two weak factorization systems.
  Typically, an inspection of the construction will reveal that it depends directly on only the weak factorization systems of the input model categories, and not their weak equivalences.
  In that case, we can apply the same construction to premodel categories.
\item
  If the construction is difficult, we may not want to rely on an analysis of its dependencies.
  An alternative approach is to use the result for model categories as a black box, as follows.
  A premodel category consists of a pair of weak factorization systems, and each weak factorization system on a complete and cocomplete category corresponds to a trivial model category as described in \cref{remark:trivial}.
  By applying the construction for model categories twice and extracting a weak factorization system from each result, we can ``reduce'' the construction for premodel categories to the (more difficult, but already established) construction for model categories.
\end{itemize}

\subsection{Slice categories}

\begin{notation}
  For a category $C$ and an object $X$ of $C$, we write $C_{X/}$ for the undercategory of objects of $C$ equipped with a map from $X$, and $C_{/X}$ for the overcategory of objects of $C$ equipped with a map to $X$.
\end{notation}

\begin{proposition}
  Let $M$ be a premodel category and $X$ an object of $C$.
  Then $M_{X/}$ and $M_{/X}$ each have a premodel category structure in which a morphism belongs to one of the classes $\sC$, $\sAC$, $\sF$, $\sAF$ if and only if its underlying morphism belongs to the corresponding class of $M$.
\end{proposition}

\begin{proof}
  This can easily be verified directly from the definitions.
  We will also show how to deduce this from the version for model categories \cite[Theorem~7.6.5]{Hi} as an example of the ``black box'' argument.
  We will just treat $M_{X/}$, as the argument for $M_{/X}$ is similar.

  Write $\sC_{X/}$ for the class of morphisms of $M_{X/}$ whose underlying morphism belongs to $\sC$, the class of cofibrations of $M$, and similarly for $\sAC_{X/}$, $\sF_{X/}$, $\sAF_{X/}$.
  We need to show that $(\sC_{X/}, \sAF_{X/})$ is a weak factorization system on $M_{X/}$.
  Now $M$ has a model category structure with cofibrations $\sC$, fibrations $\sAF$ and all maps weak equivalences.
  Applying the result for model categories, $M_{X/}$ has a model category structure with cofibrations $\sC_{X/}$, fibrations $\sAF_{X/}$ and all maps weak equivalences.
  (In particular, $M_{X/}$ is complete and cocomplete.)
  Then $(\sC_{X/}, \sAF_{X/})$ is the weak factorization system given by the cofibrations and acyclic fibrations of this model category structure on $M_{X/}$.
  By a similar argument, $(\sAC_{X/}, \sF_{X/})$ is a weak factorization system on $M_{X/}$, and clearly $\sAC_{X/} \subset \sC_{X/}$.
\end{proof}

\subsection{Products}\label{subsec:premod-products}

\begin{proposition}
  Let $(M_s)_{s \in S}$ be a family%
  \footnote{Unless otherwise specified, a ``family'' means a collection indexed by a \emph{set}, not a proper class.}
  of premodel categories.
  Then the product category $\prod_{s \in S} M_s$ has a premodel category structure in which a morphism belongs to one of the classes $\sC$, $\sAC$, $\sF$, $\sAF$ if and only if each of its components belongs to the corresponding class of $M$.
\end{proposition}

\begin{proof}
  This can be verified directly from the definitions, or by applying the ``black box'' argument to \cite[Proposition~7.1.7]{Hi}.
\end{proof}

\begin{remark}
  For any premodel category $N$, we have evident \emph{isomorphisms} of categories
  \[
    \textstyle
    \LQF(N, \prod_{s \in S} M_s) \iso \prod_{s \in S} \LQF(N, M_s),
  \]
  \[
    \textstyle
    \RQF(N, \prod_{s \in S} M_s) \iso \prod_{s \in S} \RQF(N, M_s).
  \]
  The former isomorphism makes $\prod_{s \in S} M_s$ into a \emph{strict} product object in the 2-category $\PM$.
  Using the second isomorphism, we can write
  \begin{align*}
    \textstyle \LQF(\prod_{s \in S} M_s, N)
    & \textstyle \eqv \RQF(N, \prod_{s \in S} M_s)^\op \\
    & \textstyle \iso \prod_{s \in S} \RQF(N, M_s)^\op \\
    & \textstyle \eqv \prod_{s \in S} \LQF(M_s, N).
  \end{align*}
  This is only an equivalence of categories, so $\prod_{s \in S} M_s$ is also a coproduct object in $\PM$, though not a strict one.
  This asymmetry is a result of our choice to define the morphisms of $\PM$ to be left Quillen functors.
\end{remark}

\begin{remark}
  Of course, the 2-category of model categories also has products and coproducts by the same argument.
  The advantage of premodel categories is that, as we will show in \cref{chap:algebra}, the 2-category of combinatorial premodel categories admits \emph{all} limits and colimits.
\end{remark}

\subsection{Diagram categories}\label{subsec:premod-diagram}

Let $M$ be a premodel category and $K$ a small category.
As for model categories, there are several premodel category structures we could put on the diagram category $M^K$.

\begin{definition}
  The \emph{projective premodel category structure} $M^K_\proj$ on $M^K$ is the one (if it exists) whose fibrations and anodyne fibrations are defined componentwise.
  The \emph{injective premodel category structure} $M^K_\inj$ on $M^K$ is the one (if it exists) whose cofibrations and anodyne cofibrations are defined componentwise.
\end{definition}

\begin{remark}
  The projective and injective premodel category structures may not exist for general $M$ and $K$.
  However, they are uniquely defined if they do exist, because a weak factorization system is determined by either its left or right class.
\end{remark}

\begin{proposition}\label{prop:exists-proj-inj}
  If $M$ is combinatorial, then the premodel category structures $M^K_\proj$ and $M^K_\inj$ exist and are again combinatorial.
\end{proposition}

\begin{proof}
  This follows from the ``black box'' argument applied to \cite[Proposition~A.2.8.2]{HTT}.
  (The corresponding results for model categories are stronger: the weak equivalences are also defined componentwise.
  This implies in particular that the acyclic fibrations or acyclic cofibrations are defined componentwise along with the fibrations or cofibrations.)
\end{proof}

\begin{remark}
  Projective and injective premodel category structures are dual in the sense that $(M^{K^\op}_\proj)^\op$ and $(M^\op)^K_\inj$ are equal whenever either is defined.
  However, the existence of $M^K_\proj$ and the existence of $M^K_\inj$ given by \cref{prop:exists-proj-inj} are not dual because the opposite of a combinatorial premodel category $M$ is never combinatorial unless $M$ is a poset.
  Nevertheless, we can still transfer theorems between projective and injective premodel category structures by duality under the assumption that these structures do exist.
\end{remark}

The projective and injective premodel category structures, whenever they exist, are functorial in $M$: a Quillen adjunction $F : M \adj N : G$ extends componentwise to an adjunction $F^K : M^K \adj N^K : G^K$ which is a Quillen adjunction when we equip both $M^K$ and $N^K$ with the projective or the injective premodel category structure.
We will prove a more refined version of this functoriality use for use in the next section.

\begin{proposition}\label{prop:proj-enriched-functoriality}
  Let $V$, $M$ and $N$ be premodel categories and let $F : V \times M \to N$ be a Quillen adjunction of two variables.
  Then $F$ extends to a Quillen adjunction of two variables $\tilde F : V \times M^K_\proj \to N^K_\proj$ given by the formula $\tilde F(A, X)_k = F(A, X_k)$.
  The same formula also makes $\tilde F$ into a Quillen adjunction of two variables $\tilde F : V \times M^K_\inj \to N^K_\inj$.
\end{proposition}

\begin{proof}
  As an adjunction of two variables, $F$ comes with functors $G_1 : M^\op \times N \to V$ and $G_2 : V^\op \times N \to M$ together with natural isomorphisms
  \[
    \Hom(F(A, X), Y) \iso \Hom(A, G_1(X, Y)) \iso \Hom(X, G_2(A, Y)).
  \]
  Define $\tilde G_2 : V^\op \times N^K \to M^K$ by $\tilde G_2(A, Y)_k = G_2(A, Y_k)$; then there is a natural isomorphism $\Hom(\tilde F(A, X), Y) \iso \Hom(X, \tilde G_2(A, Y))$ given by applying the adjunction between $F$ and $G_2$ componentwise.
  Define $\tilde G_1 : (M^k)^\op \times N^K \to V$ by the formula $\tilde G_1(X, Y) = \int_{k \in K} G_1(X_k, Y_k)$.
  Then we compute
  \begin{align*}
    \textstyle \Hom(A, \tilde G_1(X, Y))
    & \textstyle = \Hom(A, \int_{k \in K} G_1(X_k, Y_k)) \\
    & \textstyle \iso \int_{k \in K} \Hom(A, G_1(X_k, Y_k)) \\
    & \textstyle \iso \int_{k \in K} \Hom(F(A, X_k), Y_k) \\
    & \textstyle = \int_{k \in K} \Hom((\tilde F(A, X))_k, Y_k) \\
    & \textstyle \iso \Hom(\tilde F(A, X), Y).
  \end{align*}
  In the first step we used the representability of ends, and in the last step the fact that natural transformations are computed by an end.
  Therefore $(\tilde F, \tilde G_1, \tilde G_2)$ is an adjunction of two variables.

  Now suppose the projective premodel category structures $M^K_\proj$ and $N^K_\proj$ exist.
  We want to show $\tilde F : V \times M^K_\proj \to N^K_\proj$ is a Quillen bifunctor.
  By \cref{prop:quillen-bifunctor}, it suffices to check that for each cofibration $f : A \to B$ in $V$ and each fibration $g : X \to Y$ in $N^K_\proj$, a certain map built from $\tilde G_2$, $f$ and $g$ is a fibration in $M^K_\proj$ which is anodyne if either $f$ or $g$ is.
  This is immediate from the corresponding statement for $G_2$ because $\tilde G_2$ and the (anodyne) fibrations in $M^K_\proj$ and $N^K_\proj$ are all defined componentwise.

  Similarly, if $M^K_\inj$ and $N^K_\inj$ exist, then we see that $\tilde F$ is a Quillen bifunctor directly from the definitions, because $\tilde F$ and the (anodyne) cofibrations of $M^K_\inj$ and $N^K_\inj$ are all defined componentwise.
\end{proof}

\begin{remark}
  Let $M$ and $N$ premodel categories and $K$ a small category.
  Giving a functor $F : N \to M^K$ is the same as giving a functor from $K$ to the category of functors from $M$ to $N$.
  Using the fact that the diagonal functor $M \to M^K$ and each projection functor $M^K \to M$ are both left and right adjoints, it is not hard to see that $F$ is a left (or right) adjoint if and only if each component $N \to M^K \to M$ is a left (or right) adjoint.

  Now suppose the injective premodel category structure $M^K_\inj$ exists.
  Then for any premodel category $N$, there is an isomorphism of categories
  \[
    \LQF(N, M^K_\inj) \iso \LQF(N, M)^K = \TCat(K, \LQF(N, M))
  \]
  because, by definition of the (anodyne) cofibrations in $M^K_\inj$, a functor from $M$ to $M^K_\inj$ is a left Quillen functor if and only if its composition with each projection $M^K \to M$ is.
  Therefore $M^K_\inj$ is the \emph{cotensor} of $M$ by the category $K$.
  Similarly, if the projective premodel category structure $M^K_\proj$ exists, then
  \[
    \RQF(N, M^K_\proj) \iso \RQF(N, M)^K = \TCat(K, \RQF(N, M))
  \]
  and so
  \begin{align*}
    \LQF(M^K_\proj, N)
    & \eqv \RQF(N, M^K_\proj)^\op \\
    & \iso \TCat(K, \RQF(N, M))^\op \\
    & \iso \TCat(K^\op, \RQF(N, M)^\op) \\
    & \eqv \TCat(K^\op, \LQF(M, N)).
  \end{align*}
  Replacing $K$ by $K^\op$, we conclude that $M^{K^\op}_\proj$ is the \emph{tensor} of $M$ by the category $K$.

  When $M$ is combinatorial, we know that both $M^{K^\op}_\proj$ and $M^K_\inj$ exist and are again combinatorial.
  Therefore, the 2-category $\CPM$ admits tensors and cotensors by (small) categories.
  Of course, the same statement holds for the 2-category of combinatorial model categories as well.
\end{remark}

There is another premodel category structure on $M^K$, the \emph{Reedy premodel category structure}, which exists whenever $K$ is a Reedy category, without any condition on $M$.
We will briefly recall the basic theory of Reedy categories.
The reader may consult \cite[Chapter~15]{Hi} for a fuller treatment.

\begin{definition}[{\cite[Definition~15.1.2]{Hi}}]
  A \emph{Reedy category} is a category $K$ equipped with a function $\deg : \Ob K \to \alpha$ for some ordinal $\alpha$ and two subcategories $K_+$, $K_-$ satisfying the following properties:
  \begin{enumerate}
  \item Each nonidentity morphism of $K_+$ strictly raises degree, and each nonidentity morphism of $K_-$ strictly lowers degree.
  \item Each morphism of $K$ has a unique factorization as a morphism of $K_-$ followed by a morphism of $K_+$.
  \end{enumerate}
\end{definition}

\begin{definition}
  Let $M$ be a complete and cocomplete category and let $X$ be an object of $M^K$.
  For each object $k$ of $K$, we define
  \begin{itemize}
  \item the \emph{latching object} $L_k X = \colim_{k' \stackrel{+}{\to} k} X_{k'}$;
  \item the \emph{matching object} $M_k X = \lim_{k \stackrel{-}{\to} k'} X_{k'}$.
  \end{itemize}
  Here the colimit in the definition of $L_k$ is taken over the category of objects $k'$ equipped with a \emph{nonidentity} morphism $k' \to k$ belonging to $K_+$.
  The construction $L_k X$ is functorial in $X$, and is equipped with a canonical natural transformation $L_k X \to X_k$.
  Similarly, the limit in the definition of $M_k$ is taken over the category of objects $k'$ equipped with a nonidentity morphism $k \to k'$ belonging to $K_-$, and $M_k X$ is functorial in $X$ and equipped with a canonical natural transformation $X \to M_k X$.
\end{definition}

We will write $K_{\le \beta}$ ($K_{< \beta}$) for the full subcategory of $K$ on the objects of degree at most $\beta$ (less than $\beta$).
These categories inherit Reedy category structures from $K$, because the factorizations that appear in the definition of a Reedy category pass through an object of no greater degree than either the domain or codomain of the map to be factored.
Moreover, the latching and matching objects $L_k X$ and $M_k X$ are not affected by replacing $K$ by $K_{\le \beta}$ as long as $\deg k \le \beta$.

The key facts about diagrams indexed on a Reedy category is that objects and morphisms of $M^K$ may be constructed by (possibly transfinite) induction over $\alpha$.
Specifically, suppose we have constructed a diagram $X : K_{< \beta} \to M$.
Then giving an extension of $X$ to a diagram $K_{\le \beta} \to M$ is equivalent to giving, for each object $k$ with $\deg k = \beta$, a factorization of the map $L_k X \to M_k X$ as $L_k X \to X_k \to M_k X$.
Similarly, given diagrams $X$, $Y : K \to M$ and a natural transformation $t$ between their restrictions to $K_{< \beta}$, giving an extension of $t$ to their restrictions to $K_{\le \beta}$ is equivalent to giving, for each object $k$ with $\deg k = \beta$, a dotted arrow making the following diagram commute:
\[
  \begin{tikzcd}
    L_k X \ar[r] \ar[d, "L_k t"'] & X_k \ar[r] \ar[d, dotted] & M_k X \ar[d, "M_k t"] \\
    L_k Y \ar[r] & Y_k \ar[r] & M_k Y
  \end{tikzcd}
\]
Note that the vertical maps on the two sides of the diagram are already determined by $t$, because $L_k X$ and $M_k X$ only depend on $X_{k'}$ for objects $k'$ with $\deg k' < \deg k = \beta$.

For use in a future technical argument, we prove a slight generalization of the expected result on the existence of Reedy premodel category structures.

\begin{proposition}\label{prop:reedy-wfs}
  Let $(\sL_k, \sR_k)_{k \in K}$ be a family of weak factorization systems on $M$ indexed by the objects of the Reedy category $K$.
  Let $\sL$ be the class of morphisms $X \to Y$ such that for each $k$, the map $X \amalg_{L_k X} L_k Y \to Y_k$ belongs to $\sL_k$, and let $\sR$ be the class of morphisms such that for each $k$, the map $X_k \to Y_k \times_{M_k Y} M_k X$ belongs to $\sR_k$.
  Then $(\sL, \sR)$ is a weak factorization system on $M^K$.
\end{proposition}

\begin{proof}
  By \cref{prop:wfs-of-retract}, to show that $(\sL, \sR)$ is a weak factorization system, it suffices to show that maps of $\sL$ have the left lifting property with respect to maps of $\sR$, that any map can be factored as a map of $\sL$ followed by a map of $\sR$, and that $\sL$ and $\sR$ are closed under retracts.
  Lifts and factorizations may be constructed inductively by the usual argument, and the closure of $\sL$ and $\sR$ under retracts is clear.
\end{proof}


\begin{definition}
  Let $M$ be a premodel category and $K$ a Reedy category.
  The \emph{Reedy premodel category structure} $M^K_\Reedy$ on $M^K$ is the one in which a morphism $X \to Y$ is
  \begin{itemize}
  \item a cofibration or anodyne cofibration if $X \amalg_{L_k X} L_k Y \to Y_k$ is one for every $k \in K$;
  \item a fibration or anodyne fibration if $X_k \to Y_k \times_{M_k Y} M_k X$ is one for every $k \in K$.
  \end{itemize}
\end{definition}

\begin{proposition}
  This structure $M^K_\Reedy$ is indeed a premodel category structure on $M^K$.
\end{proposition}

\begin{proof}
  Applying \cref{prop:reedy-wfs} with $(\sL_k, \sR_k) = (\sC_M, \sAF_M)$ for every $k$, we see that the cofibrations and anodyne fibrations of $M^K$ do indeed form a weak factorization system, and similarly for the anodyne cofibrations and fibrations.
\end{proof}

For particular classes of categories $K$, the Reedy premodel category structure on $M^K$ agrees with either the projective or injective one.

\begin{definition}
  A \emph{direct category} is a Reedy category $K$ in which $K_-$ contains only isomorphisms, so that $K_+$ is all of $K$.
  Equivalently, $K$ is a category which admits a function $\deg : \Ob K \to \alpha$ for which every nonidentity morphism strictly increases degree.
  An \emph{inverse category} is the opposite of a direct category.
\end{definition}

\begin{proposition}\label{prop:reedy-direct-proj}
  Let $K$ be a direct category.
  Then the (anodyne) fibrations in $M^K_\Reedy$ are the componentwise (anodyne) fibrations, so $M^K_\Reedy = M^K_\proj$.
  Dually, if $K$ is an inverse category, then $M^K_\Reedy = M^K_\inj$.
\end{proposition}

\begin{proof}
  When $K$ is a direct category, the indexing category for $M_k X$ is empty for every $k$, so the morphism $X_k \to Y_k \times_{M_k Y} M_k X$ used to define the (anodyne) fibrations is just $X_k \to Y_k$.
  The argument for an inverse category is dual.
\end{proof}

In particular, when $K$ is a direct (respectively, inverse) category, $M^K_\proj$ (respectively, $M^K_\inj$) exists even when $M$ is not assumed to be combinatorial.

\begin{example}
  Let $K = [1] = \{0 \to 1\}$, with degree function $\deg : \Ob K \to \N$ given by $\deg 0 = 0$ and $\deg 1 = 1$.
  Then $K$ is a direct category, so for any premodel category $M$, there is a Reedy (or projective) premodel category structure $M^{[1]}_\Reedy = M^{[1]}_\proj$.
  The (anodyne) fibrations are defined componentwise, of course.
  We can give an explicit description of the (anodyne) cofibrations from the definition of the Reedy premodel category structure.
  By inspection, $L_0 X = \emptyset$ while $L_1 X = X_0$.
  Therefore, a morphism $X \to Y$ is a cofibration or anodyne cofibration if and only if each of the maps $X_0 \to Y_0$ and $X_1 \amalg_{X_0} Y_0 \to Y_1$ is one.
\end{example}

\begin{example}\label{ex:reedy-arrow}
  Let us consider the same category $K$, but in the more general setting of \cref{prop:reedy-wfs}.
  Let $(\sL_0, \sR_0)$ and $(\sL_1, \sR_1)$ be two weak factorization systems on $M$.
  Then in the weak factorization system $(\sL, \sR)$ produced by \cref{prop:reedy-wfs},
  \begin{itemize}
  \item a map $X \to Y$ belongs to $\sL$ if $X_0 \to Y_0$ belongs to $\sL_0$ and $X_1 \amalg_{X_0} Y_0 \to Y_1$ belongs to $\sL_1$;
  \item a map $X \to Y$ belongs to $\sR$ if $X_0 \to Y_0$ belongs to $\sR_0$ and $X_1 \to Y_1$ belongs to $\sR_1$.
  \end{itemize}
  We will use this weak factorization system in the construction of the internal Hom of combinatorial premodel categories in \cref{chap:algebra}.
\end{example}

To describe the generating (anodyne) cofibrations of a Reedy premodel category structure, or the more general structures produced by \cref{prop:reedy-wfs}, it is more convenient to work with diagrams indexed on $K^\op$.
The category $K^\op$ is again a Reedy category, with the same degree function as $K$ and the roles of $K_+$ and $K_-$ reversed.
Suppose $(\sL_k, \sR_k)$ is a weak factorization system on $M$ for each $k \in \Ob K$ which is generated by a class of morphisms $I_k$ (so that $\sR_k = \rlp(I_k)$).
We want to describe generators for the weak factorization system on $M^{K^\op}_\Reedy$ constructed in \cref{prop:reedy-wfs}.

We write $\yo : K \to \Set^{K^\op}$ for the Yoneda embedding.
For each object $k$ of $K$, we also define $\partial (\yo k) = \colim_{k' \stackrel{+}{\to} k} \yo k'$.
Here the colimit is taken over all nonidentity morphisms of $K_+$ for the original Reedy structure on $K$.
The object $\partial (\yo k)$ is equipped with a canonical map $d_k : \partial (\yo k) \to \yo k$ which is induced by the maps $\yo k' \to \yo k$ for each $k' \stackrel{+}{\to} k$.

There is an adjunction of two variables $\otimes : \Set^{K^\op}_\Reedy \times M \to M^{K^\op}_\Reedy$ given by the formula $(S \otimes A)_k = A^{\amalg S_k}$.
We write $\bp$ for the corresponding pushout product.

\begin{lemma}\label{lemma:reedy-lifting}
  For any morphism $i : A \to B$ of $M$, a morphism $X \to Y$ of $M^{K^\op}$ has the right lifting property with respect to $d_k \bp i$ if and only if the induced map $X_k \to Y_k \times_{M_k Y} M_k X$ has the right lifting property with respect to $i$.
\end{lemma}

Here $M_k$ denotes the matching functor for $M^{K^\op}$, so it is computed using the morphisms of $(K^\op)^-$, which correspond to the morphisms of $K^+$.

\begin{proof}
  With some notational changes this is the statement of \cite[Corollary~15.6.21]{Hi}.
\end{proof}

\begin{proposition}\label{prop:reedy-gen}
  The weak factorization system $(\sL, \sR)$ on $M^{K^\op}$ produced by \cref{prop:reedy-wfs} is generated by the class $I = \{\,d_k \bp i \mid k \in \Ob K, i \in I_k\,\}$.
  In particular, it is generated by a set if each $(\sL_k, \sR_k)$ is, and so $M^{K^\op}_\Reedy$ is combinatorial if $M$ is.
\end{proposition}

\begin{proof}
  This follows immediately from the definition of $\sR$ and \cref{lemma:reedy-lifting}.
\end{proof}

In particular, taking $M = \Set$, the premodel category $\Set^{K^\op}_\Reedy$ has generating cofibrations $I = \{\,d_k \mid k \in \Ob K\,\}$ and $J = \emptyset$.

\begin{example}
  Take $K = \Delta$ with its usual Reedy structure.
  Then $d_{[n]}$ is the usual boundary inclusion $d_{[n]} : \partial \Delta^n \to \Delta^n$.
  Hence the generating cofibrations of $\Set^{\Delta^\op}_\Reedy$ are the standard generating cofibrations of $\Kan$, so the cofibrations of $\Set^{\Delta^\op}_\Reedy$ are the monomorphisms.
  Meanwhile $\Set^{\Delta^\op}_\Reedy$ has no generating anodyne cofibrations, so its anodyne cofibrations are just the isomorphisms.
\end{example}

\subsection{The algebra of diagram categories}

We have already mentioned the functoriality of $M^K_\proj$ and $M^K_\inj$ in $M$.
These constructions are also functorial in $K$.
In fact, this functoriality is a formal consequence of the identification of $M^K_\proj$ and $M^K_\inj$ as the tensor of $M$ by $K^\op$ and the cotensor of $M$ by $K$ respectively.

Let $F : K \to L$ be a functor between small categories.
Then $F$ induces a left Quillen functor between the projective premodel categories
\[
  F_! : M^K_\proj = K^\op \otimes M \xrightarrow{F^\op \otimes M} L^\op \otimes M = M^L_\proj.
\]
For any premodel category $N$, this left Quillen functor induces the functor
\[
  \LQF(M^L_\proj, N) \eqv \TCat(L^\op, \LQF(M, N)) \to \TCat(K^\op, \LQF(M, N)) \eqv \LQF(M^K_\proj, N)
\]
given by precomposition by $F^\op : K^\op \to L^\op$.
Passing to right adjoints produces
\[
  \RQF(N, M^L_\proj) \iso \TCat(L, \RQF(N, M)) \to \TCat(K, \RQF(N, M)) \iso \RQF(N, M^K_\proj)
\]
where the middle functor is given by precomposition by $F$.
Now, taking $N = M^L_\proj$, the identity functor of $M^L_\proj$ is sent by this composition to the restriction functor $F^* : M^L_\proj \to M^K_\proj$.
Hence the original left Quillen functor $F_! : M^K_\proj \to M^L_\proj$ induced by $F$ is the left adjoint of restriction, namely left Kan extension.
If $t : F \to F'$ is a natural transformation between the functors $F : K \to L$ and $F' : K \to L$, then $t$ induces an evident natural transformation $t^* : F^* \to (F')^*$ and hence a corresponding natural transformation $t_! : (F')_! \to F_!$.

The injective premodel category structure $M^K_\inj$ is the cotensor of $M$ by $K$, so it has the opposite variance in $K$.
The functor $F$ induces a left Quillen functor
\[
  F^* : M^L_\inj = M^L \xrightarrow{M^F} M^K = M^K_\inj
\]
which induces the composition
\[
  \LQF(N, M^L_\inj) \iso \TCat(L, \LQF(N, M)) \to \TCat(K, \LQF(N, M)) \iso \LQF(N, M^K_\inj)
\]
for any premodel category $N$, where the middle functor is given by precomposition by $F$.
In particular taking $N = M^L_\inj$ we conclude that $F^* : M^L_\inj \to M^K_\inj$ is simply the restriction.
If $t : F \to F'$ is a natural transformation then $t$ induces a natural transformation $t^* : F^* \to (F')^*$.

The identifications $M^K_\proj = K^\op \otimes M$ and $M^K_\inj = M^K$ also imply the identities
\[
  (M^L_\proj)^K_\proj = M^{K \times L}_\proj, \quad
  (M^L_\inj)^K_\inj = M^{K \times L}_\inj, \quad
  M^{\{*\}}_\proj = M^{\{*\}}_\inj = M
\]
where $\{*\}$ denotes the terminal category.
Here we have identified the underlying categories $(M^L)^K$ and $M^{K \times L}$ as well as $M^{\{*\}}$ and $M$ in the obvious way, so that the above identities become equalities of premodel category structures.
These equalities are also immediate from the definitions of the projective and injective premodel category structures.
When $K$ and $L$ are Reedy categories, $K \times L$ also has a natural Reedy category structure and the corresponding formula $(M^L_\Reedy)^K_\Reedy = M^{K \times L}_\Reedy$ also holds by \cite[Theorem~15.5.2]{Hi}.
(Hirschhorn only considers Reedy categories whose degree function takes values in the natural numbers, but the argument can be extended to ordinal-valued degree functions by equipping $K \times L$ with the degree function $\deg {(k, l)} = \deg k \oplus \deg l$ where $\oplus$ denotes the natural sum of ordinals.)

For premodel categories $M$ and $N$, the equivalences $\LQF(M^{K^\op}_\proj, N) \eqv \TCat(K, \LQF(M, N))$ and $\LQF(M, N^K_\inj) \eqv \TCat(K, \LQF(M, N))$ compose to form a ``tensor--cotensor'' or ``projective--injective'' adjunction
\[
  \LQF(M^{K^\op}_\proj, N) \eqv \LQF(M, N^K_\inj).
\]
Even though the Reedy premodel structure $M^K_\Reedy$ is not involved in an adjunction with $\TCat(K, -)$, there is also a ``Reedy--Reedy'' adjunction
\[
  \LQF(M^{K^\op}_\Reedy, N) \eqv \LQF(M, N^K_\Reedy)
\]
which we establish next.
Both of these adjunctions are restrictions of the equivalence between left adjoints $F : M^{K^\op} \to N$ and left adjoints $F' : M \to N^K$ given by the formula
\[
  (F'A)_k = F(\yo k \otimes A)
\]
where $\otimes : \Set^{K^\op} \times M \to M^{K^\op}$ is the adjunction of two variables mentioned earlier.

It will be convenient to establish the Reedy--Reedy adjunction in the following form.
\begin{proposition}
  Let $M_1$, $M_2$ and $N$ be premodel categories and $K$ a Reedy category.
  Then there is an equivalence
  \[
    \QBF(((M_1)^{K^\op}_\Reedy, M_2), N) \eqv \QBF((M_1, M_2), N^K_\Reedy)
  \]
  given by sending $F : (M_1)^{K^\op} \times M_2 \to N$ to $F' : M_1 \times M_2 \to N^K$ defined by the formula $(F'(A_1, A_2))_k = F(\yo k \otimes A_1, A_2)$.
\end{proposition}

\begin{proof}
  The point is that the generating cofibrations of $(M_1)^{K^\op}_\Reedy$ and the property of being a cofibration of $N^K_\Reedy$ both involve latching-type constructions.
  By \cref{prop:quillen-bifunctor,prop:reedy-gen}, an adjunction of two variables $F : (M_1)^{K^\op}_\Reedy \times M_2 \to N$ is a Quillen bifunctor if, for every cofibrations $f_1$ of $M_1$ and $f_2$ of $M_2$, $(d_k \bp f_1) \bp_F f_2$ is a cofibration of $N$ which is anodyne if either $f_1$ or $f_2$ is.
  A calculation shows that $(d_k \bp f_1) \bp_F f_2$ can be identified with the induced map $X \amalg_{L_k X} L_k Y \to Y$ for $X \to Y$ the morphism $f_1 \bp_{F'} f_2$ of $N^K$.
  The result then follows from the definition of the (anodyne) cofibrations of $N^K_\Reedy$.
\end{proof}

\begin{proposition}\label{prop:reedy-reedy}
  Let $M$ and $N$ be premodel categories and $K$ a Reedy category.
  Then there are canonical equivalences
  \[
    \LQF(M^{K^\op}_\Reedy, N) \eqv \QBF((\Set^{K^\op}_\Reedy, M), N)
    \eqv \LQF(M, N^K_\Reedy).
  \]
\end{proposition}

\begin{proof}
  This follows from the preceding proposition together with the equivalences
  \[
    \LQF(M^{K^\op}_\Reedy, N) \eqv \QBF((\Set, M^{K^\op}_\Reedy), N)
  \]
  and
  \[
    \LQF(M, N^K_\Reedy) \eqv \QBF((\Set, M), N^K_\Reedy). \qedhere
  \]
\end{proof}

In particular, we can take $M = \Set$.
Then a left Quillen functor from $\Set^{K^\op}_\Reedy$ to $N$ is the same as a Reedy cofibrant diagram in $N$.

\begin{remark}
  Once we define the tensor product and internal Hom of combinatorial premodel categories, we may regard \cref{prop:reedy-reedy} as computing
  \[
    \Set^{K^\op}_\Reedy \otimes M \eqv M^{K^\op}_\Reedy, \quad
    \CPM(\Set^{K^\op}_\Reedy, N) \eqv N^K_\Reedy.
  \]
  In particular, there is an adjunction
  \[
    \LQF(\Set^{K^\op}_\Reedy \otimes M, N) \eqv \LQF(M, \Set^K_\Reedy \otimes N).
  \]
  It follows that $\Set^{K^\op}_\Reedy$ and $\Set^K_\Reedy$ are \emph{dual} objects of $\CPM$.
  This duality is a consequence of the close relationship between the fibrations of $M^{K^\op}_\Reedy$ and the cofibrations of $M^K_\Reedy$.
\end{remark}

\section{Relaxed premodel categories}\label{sec:relaxed-premodel}

We will see in \cref{chap:algebra} that the 2-category $\CPM$ of combinatorial premodel categories has excellent algebraic properties.
In particular, it admits all limits and colimits.
In order to equip $\CPM$ with a model 2-category structure, our next task should be to define what it means for a left Quillen functor between premodel categories to be a weak equivalence.
However, here we encounter a serious issue.

Recall that each \emph{model} category has an associated homotopy category and also an associated $(\infty, 1)$-category (given for example by its simplicial localization).
A left Quillen functor is a \emph{Quillen equivalence} when it induces an equivalence of homotopy categories or, equivalently, an equivalence of associated $(\infty, 1)$-categories.
This is the natural class of functors to take as the weak equivalences between model categories.
However, all of the concepts we have just mentioned appear to depend on an essential way on the weak equivalences of a model category.
(Indeed, they can all be defined in terms of the weak equivalences alone, and do not depend on the choices of cofibrations or fibrations.)

Our challenge, then, is to find an alternative way to associate an $(\infty, 1)$-category to each \emph{premodel} category in a manner which extends the associated $(\infty, 1)$-categories of model categories.
In general, we cannot expect this to work out as neatly as for model categories.
For example, one basic feature of the associated $(\infty, 1)$-category of a model category $M$ is that each object of the $(\infty, 1)$-category is ``modeled'' by (that is, equivalent to) some object of $M$.
Let's consider our basic example of a premodel category, $\Set$.
Suppose $N$ is a model category, so that we know how to define its associated $(\infty, 1)$-category $\mathbf{N}$.
A left Quillen functor from $\Set$ to $N$ ought to induce a left adjoint between the associated $(\infty, 1)$-categories.
We saw that giving a left Quillen functor from $\Set$ to any premodel category $N$ is the same as giving a cofibrant object of $N$, so it is reasonable to assume that giving a left adjoint from the $(\infty, 1)$-category associated to the premodel category $\Set$ to $\mathbf{N}$ is the same as giving an object of $\mathbf{N}$.
This is possible only if the $(\infty, 1)$-category associated to $\Set$ is the $(\infty, 1)$-category of \emph{spaces}.
In particular, the left Quillen functor $\Set \to \Kan$ corresponding to the object $\Delta^0$ ought to be a Quillen equivalence.
It is clearly unreasonable to expect a simple way to extract the $(\infty, 1)$-category of spaces from the premodel category $\Set$ or to recognize $\Set \to \Kan$ as a Quillen equivalence.

The problem is somehow that $\Set$ lacks cylinder objects, or to put it another way, not every morphism is equivalent to a cofibration.
In effect, we could say that as a premodel category, $\Set$ does not have all homotopy pushouts and therefore we cannot build up objects which represent general spaces.

There are two ways we could attempt to surmount these difficulties and associate an $(\infty, 1)$-category to a premodel category $M$.
\begin{itemize}
\item We could try to ``close up'' $M$ under homotopy colimits, either in the world of premodel categories (say by replacing $M$ by the category of simplicial objects of $M$), or via some universal construction in the world of $(\infty, 1)$-categories.
\item Alternatively, we could restrict our attention to a class of well-behaved premodel categories (excluding examples like $\Set$) over which we have sufficient control to imitate the usual constructions of model category theory.
\end{itemize}

Here, we adopt the second approach.
In this section we will introduce the notion of a \emph{relaxed} premodel category.
In \cref{chap:homotopy}, we will show how to associate an $(\infty, 1)$-category to each relaxed premodel category.
Each model category will be relaxed when viewed as a premodel category, and in this case the associated $(\infty, 1)$-category will be the usual one.
Furthermore, the associated $(\infty, 1)$-categories of relaxed premodel categories have the same general features as those of model categories:
for example, they are complete and cocomplete $(\infty, 1)$-categories.
Hence we may define a \emph{Quillen equivalence} between relaxed premodel categories to be a left Quillen functor which induces an equivalence of associated $(\infty, 1)$-categories.
This defines a class of functors extending the Quillen equivalences between model categories and with the expected properties of a class of weak equivalences.

The definition of a relaxed premodel category that we will give below may appear a bit artificial:
it is essentially just the condition needed to make our arguments work.
It may seem that we are now no better off than we started, as there is no apparent reason that relaxed combinatorial premodel categories should be closed under limits or colimits.
The key idea is that the property of being relaxed can be guaranteed by the presence of additional \emph{algebraic} structure.
Specifically, if $V$ is a monoidal model category, then a $V$-premodel category will automatically be relaxed.
On the other hand, when $V$ is combinatorial, the 2-category $V\CPM$ of combinatorial $V$-premodel categories still admits all limits and colimits, thanks to the formal algebraic properties of $\CPM$.
Since $V$-premodel categories are relaxed, $V\CPM$ is also equipped with a class of weak equivalences which extends the class of Quillen equivalences between combinatorial $V$-model categories.
Thus, once we have established all these properties of $V\CPM$, we will have a 2-category which is at least a candidate for being equipped with a model category structure.

\subsection{Resolutions in model categories}

We begin by recasting the classical theory of resolutions in model categories (following \cite[Chapter~16]{Hi}) in a form that will be suitable for generalization to premodel categories.
Throughout this section, let $M$ be a model category.

\begin{definition}\label{def:model-cosimplicial-resolution}
  A \emph{cosimplicial object} $A^\bullet$ of $M$ is a functor $\Delta \to M$.
  A cosimplicial object $A^\bullet$ is \emph{Reedy cofibrant} if it is cofibrant when viewed as an object of $M^\Delta$ equipped with the Reedy model category structure.
  A \emph{cosimplicial resolution} in $M$ is a Reedy cofibrant cosimplicial object $A^\bullet$ for which the structure map $A^n \to A^0$ is a weak equivalence for each $n$.
\end{definition}

\begin{remark}
  Our definition of a cosimplicial resolution agrees with \cite[Definition~16.1.26]{Hi} by the two-out-of-three condition.
  A possibly more common notion is that of a \emph{cosimplicial resolution of an object} \cite[Definition~16.1.2]{Hi}; we won't have a use for this notion, so we have chosen a definition of a cosimplicial resolution which avoids it.
\end{remark}

Any model category has a sufficient supply of cosimplicial resolutions.

\begin{proposition}\label{prop:enough-resolutions}
  Let $M$ be a model category.
  Then for every cofibrant object $A$ of $M$, there exists a cosimplicial resolution $A^\bullet$ such that $A^0$ is (isomorphic to) $A$.
\end{proposition}

\begin{proof}
  See \cite[Proposition~16.6.8]{Hi}.
  There it is assumed that a model category admits a functorial factorization, but without this hypothesis the same proof still constructs a single cosimplicial resolution with a specified cofibrant object in degree $0$.
\end{proof}

Now we adopt a slightly different point of view on cosimplicial objects.
The category of simplicial sets $\Set^{\Delta^\op}$ is the free cocompletion of the category $\Delta$, so the category of cosimplicial objects of $M$ is equivalent to the category of left adjoints from $\Set^{\Delta^\op}$ to $M$.

\begin{notation}[{\cite[Definition~16.3.1]{Hi}}]
  For $A^\bullet$ a cosimplicial object of $M$, we write $A^\bullet \otimes - : \Set^{\Delta^\op} \to M$ for the corresponding left adjoint, so that for example $A^\bullet \otimes \Delta^n = A^n$ (up to canonical isomorphism).
\end{notation}

We can then reformulate the notion of cosimplicial resolution as follows.

\begin{proposition}
  Let $A^\bullet$ be a cosimplicial object of $M$.
  \begin{enumerate}
  \item The object $A^\bullet$ is Reedy cofibrant if and only if $A^\bullet \otimes - : \Set^{\Delta^\op} \to M$ takes the map $\partial \Delta^n \to \Delta^n$ to a cofibration for each $n \ge 0$.
  \item The object $A^\bullet$ is a cosimplicial resolution if and only if, in addition, $A^\bullet \otimes - : \Set^{\Delta^\op} \to M$ takes the map $\Lambda^n_k \to \Delta^n$ to an acyclic cofibration for each $n \ge 1$ and each $0 \le k \le n$.
  \end{enumerate}
\end{proposition}

\begin{proof}
  By \cite[Proposition~16.3.8]{Hi}, the map $A^\bullet \otimes \partial\Delta^n \to A^\bullet \otimes \Delta^n$ is naturally isomorphic to the latching map $L_n A^\bullet \to A^n$.
  This proves part (1).
  For part (2), the forward direction is \cite[Proposition~16.4.11]{Hi}.
  Conversely, suppose that $A^\bullet$ is a cosimplicial object such that $A^\bullet \otimes -$ takes the maps $\partial \Delta^n \to \Delta^n$ for $n \ge 0$ to cofibrations and the maps $\Lambda^n_k \to \Delta^n$ for $n \ge 1$ and $0 \le k \le n$ to acyclic cofibrations.
  Then $A^\bullet \otimes -$ is a left Quillen functor from the standard (Kan--Quillen) model category structure on simplicial sets to $M$.
  The map $A^n \to A^0$ is naturally isomorphic to $A^\bullet \otimes \Delta^n \to A^\bullet \otimes \Delta^0$, and the latter map is a weak equivalence because the map $\Delta^n \to \Delta^0$ is a weak equivalence between cofibrant objects.
\end{proof}

As we already noted in the course of the preceding proof, we therefore have the following alternative description of cofibrant resolutions.

\begin{corollary}\label{cor:model-resolution-quillen}
  Under the correspondence between cosimplicial objects of $M$ and left adjoints from $\Set^{\Delta^\op}$ to $M$, a cosimplicial resolution $A^\bullet$ corresponds precisely to a left Quillen functor from $\Kan$ to $M$, namely the functor $A^\bullet \otimes -$.
  The inverse correspondence takes a left Quillen functor $F : \Kan \to M$ to the cosimplicial object $F(\Delta^\bullet)$ of $M$.
\end{corollary}

We may then reformulate the existence of cofibrant resolutions as follows.

\begin{proposition}\label{prop:enough-quillen}
  Let $M$ be a model category.
  Then every cofibrant object $A$ of $M$ is (isomorphic to) $F(\Delta^0)$ for some left Quillen functor $F : \Kan \to M$.
\end{proposition}

\begin{remark}
  We may thus interpret the existence of cosimplicial resolutions as asserting that the pair $(\Kan, \Delta^0)$ is a ``weakly universal model category equipped with a cofibrant object''.
  It may be of interest to find other pairs $(N, A)$ with the same property or to find necessary conditions on such a pair.
\end{remark}

So far we have discussed cosimplicial resolutions.
By duality, there are corresponding descriptions of simplicial objects, Reedy fibrant simplicial objects, and simplicial resolutions in terms of left Quillen functors $\Kan \to M^\op$, or equivalently, right Quillen functors $\Kan^\op \to M$ or left Quillen functors $M \to \Kan^\op$.
We summarize them below and leave the details to the reader.

\begin{notation}[{\cite[Definition~16.3.1]{Hi}}]
  For $X_\bullet : \Delta^\op \to M$ a simplicial object of $M$, we write $(X_\bullet)^- : (\Set^{\Delta^\op})^\op \to M$ for the corresponding right adjoint, so that for example $(X_\bullet)^{\Delta^n} = X_n$ (up to canonical isomorphism).
\end{notation}

\begin{proposition}
  Under the correspondence between simplicial objects of $M$ and right adjoints from $(\Set^{\Delta^\op})^\op$ to $M$, simplicial resolutions correspond precisely to right Quillen functors from $\Kan^\op$ to $M$.
\end{proposition}

\begin{proposition}
  Let $M$ be a model category and $X$ a fibrant object of $M$.
  Then there exists a left Quillen functor $F : M \to \Kan^\op$ whose right adjoint $G : \Kan^\op \to M$ sends $(\Delta^0)^\op$ to (an object isomorphic to) $X$.
\end{proposition}

\subsection{Resolutions in premodel categories}

\Cref{def:model-cosimplicial-resolution} mentions the weak equivalences of $M$, so it cannot be adapted directly to the setting of premodel categories.
However, \cref{cor:model-resolution-quillen} gives an alternate description of cosimplicial resolutions in terms of left Quillen functors into $M$, which are defined in terms of the cofibrations and acyclic cofibrations of $M$ only.
Therefore, we may make the following definition.

\begin{definition}
  Let $M$ be a premodel category.
  A cosimplicial object $A^\bullet$ of $M$ is a \emph{cosimplicial resolution} if $A^\bullet \otimes -$, viewed as a functor from $\Kan$ to $M$, is a left Quillen functor.
\end{definition}

In practice we tend to dispense with the cosimplicial object $A^\bullet$ and instead work directly with left Quillen functors from $\Kan$ to $M$.

The analogue of \cref{prop:enough-quillen}, however, does \emph{not} hold for a general premodel category.

\begin{example}
  Let $F : \Kan \to \Set$ be a left Quillen functor, where as usual $\Set$ carries the premodel category structure in which the cofibrations are the monomorphisms and the anodyne cofibrations are the isomorphisms.
  Write $A = F(\Delta^0)$ and $B = F(\Delta^1)$.
  Then $F$ sends the cofibration $\Delta^0 \amalg \Delta^0 = \partial\Delta^0 \to \Delta^1$ of $\Kan$ to a monomorphism $A \amalg A \to B$.
  Precomposing with either inclusion $\Delta^0 \to \Delta^0 \amalg \Delta^0$ gives an acyclic cofibration $\Delta^0 \to \Delta^1$, so the two compositions $A \to A \amalg A \to B$ are anodyne cofibrations in $\Set$, that is, isomorphisms.
  Then $A \amalg A \to B$ must be isomorphic to the fold map $A \amalg A \to A$.
  But this is only a monomorphism when $A$ is empty, so for any nonempty set $A$, there is no left Quillen functor $F : \Kan \to \Set$ with $F(\Delta^0) = A$.
\end{example}

We therefore introduce the following terminology.

\begin{definition}
  Let $M$ be a premodel category.
  We say $M$ \emph{has enough cosimplicial resolutions} if for every cofibrant object $A$ of $M$, there exists a left Quillen functor $F : \Kan \to M$ with $F(\Delta^0) = A$.
  We say $M$ \emph{has enough simplicial resolutions} if $M^\op$ has enough cosimplicial resolutions.
\end{definition}

While the definition above (and the term ``cosimplicial resolution'' itself) makes reference to the model category of simplicial sets, we can also give an equivalent ``unbiased'' characterization.

\begin{proposition}
  Let $M$ be a premodel category.
  The following conditions are equivalent:
  \begin{enumerate}
  \item $M$ has enough cosimplicial resolutions.
  \item For every cofibrant object $A$ of $M$, there exists a model category $N$, a cofibrant object $B$ of $N$ and a left Quillen functor $F : N \to M$ with $FB = A$.
  \end{enumerate}
\end{proposition}

\begin{proof}
  Condition (1) implies condition (2) by taking $N = \Kan$ and $B = \Delta^0$ for any $A$.
  Suppose $M$ satisfies condition (2) and let $A$ be a cofibrant object of $M$.
  Then there is a model category $N$, a cofibrant object $B$ of $N$ and a left Quillen functor $F$ with $FB = A$.
  By the existence of cosimplicial resolutions in model categories, there is a left Quillen functor $F' : \Kan \to N$ with $F'(\Delta^0) = B$.
  Then $FF' : \Kan \to M$ is a left Quillen functor with $FF'(\Delta^0) = A$.
\end{proof}

\begin{example}\label{ex:model-cosimplicial}
  Any model category has enough cosimplicial resolutions and, dually, enough simplicial resolutions.
\end{example}

Of course, if this were the only way to produce premodel categories with enough cosimplicial resolutions, there would be no point in generalizing to premodel categories in the first place.
The key point is that we can also produce cosimplicial resolutions by working in an \emph{enriched} setting.
More generally, we have the following result.

\begin{proposition}\label{prop:action-resolutions}
  If $V$ has enough cosimplicial resolutions and $M$ admits a unital action by $V$, then $M$ has enough cosimplicial resolutions and enough simplicial resolutions.
\end{proposition}

\begin{proof}
  Let $F : V \times M \to M$ with $F(I, -) \iso \id_M$ be a unital action by $V$ on $M$.
  By definition, we may choose a left Quillen functor $F' : \Kan \to V$ with $F'(\Delta^0) = I$.
  Then for any cofibrant object $A$ of $M$, $F(F'(-), A) : \Kan \to M$ is a left Quillen functor and $F(F'(\Delta^0), A) = F(I, A) = A$.
  So $M$ has enough cosimplicial resolutions.
  
  By \cref{prop:op-action}, $M^\op$ also admits a unital action by $V$ and so $M$ has enough simplicial resolutions as well.
\end{proof}

\begin{proposition}\label{prop:enriched-enough-cosimplicial}
  Let $V$ be a symmetric monoidal model category and $M$ a $V$-premodel category.
  Then $M$ has enough cosimplicial resolutions and enough simplicial resolutions.
\end{proposition}

\begin{proof}
  $M$ admits a unital action by $V$, and $V$ has enough cosimplicial resolutions by \cref{ex:model-cosimplicial}.
\end{proof}

\subsection{Relaxed premodel categories}

We can now introduce the condition on a premodel category that we will use to ensure that it has a well-behaved homotopy theory.

Recall that for any premodel category $M$, the arrow category of $M$ has a projective (or Reedy) premodel category structure $M^{[1]}_\proj$ in which a morphism $A \to B$ is a cofibration (anodyne cofibration) if and only if both $A_0 \to B_0$ and the induced map $A_1 \amalg_{A_0} B_0 \to B_1$ are cofibrations (anodyne cofibrations).

\begin{definition}
  We say that a premodel category $M$ is \emph{left relaxed} if the premodel category $M^{[1]}_\proj$ has enough cosimplicial resolutions.
  We say that $M$ is \emph{right relaxed} if $M^\op$ is left relaxed, and that $M$ is \emph{relaxed} if $M$ is both left relaxed and right relaxed.
\end{definition}

\begin{remark}
  If $M$ is left relaxed, then $M$ also has enough cosimplicial resolutions, as for any cofibrant object of $A$, the object $\emptyset \to A$ is cofibrant in $M^{[1]}_\proj$ and its image under the left Quillen functor from $M^{[1]}_\proj$ to $M$ taking a morphism to its codomain equals $A$.
\end{remark}

\begin{remark}\label{remark:concrete-relaxed}
  Concretely, $M$ is left relaxed if it satisfies the following ``relative axiom'' on cosimplicial resolutions.
  An object $f : A \to B$ of $M^{[1]}_\proj$ is cofibrant when $A$ is cofibrant and $f$ is a cofibration.
  For each such $f$, we require the existence of a diagram between cosimplicial objects of $M$
  \[
    \begin{tikzcd}
      A \ar[r, equal] \ar[d, "f"'] &
      A^0 \ar[r, shift left] \ar[r, shift right] \ar[d] &
      A^1 \ar[l] \ar[r, shift left=2] \ar[r] \ar[r, shift right=2] \ar[d] &
      \cdots \ar[l, shift left=1] \ar[l, shift right=1] \\
      B \ar[r, equal] &
      B^0 \ar[r, shift left] \ar[r, shift right] &
      B^1 \ar[l] \ar[r, shift left=2] \ar[r] \ar[r, shift right=2] &
      \cdots \ar[l, shift left=1] \ar[l, shift right=1]
    \end{tikzcd}
  \]
  satisfying the following properties, where for brevity we write $A[K]$ for $A^\bullet \otimes K$ and $B[K]$ for $B^\bullet \otimes K$.
  (Note that $A[K]$ and $B[K]$ depend not just on $A$, $B$ and $K$ but also on the original map $f : A \to B$ and the choice of cosimplicial resolution.)
  \begin{enumerate}
  \item For each $n \ge 0$, the maps $A[\partial \Delta^n] \to A[\Delta^n] = A^n$ and $A[\Delta^n] \amalg_{A[\partial \Delta^n]} B[\partial \Delta^n] \to B[\Delta^n] = B^n$ are cofibrations in $M$.
  \item For each $n \ge 1$ and $0 \le k \le n$, the maps $A[\Lambda^n_k] \to A[\Delta^n] = A^n$ and $A[\Delta^n] \amalg_{A[\Lambda^n_k]} B[\Lambda^n_k] \to B[\Delta^n] = B^n$ are anodyne cofibrations in $M$.
  \end{enumerate}
  In brief, these are the properties we would expect of the functors $A[K] = A \otimes K$ and $B[K] = B \otimes K$ for a cofibration $A \to B$ between cofibrant objects of a simplicial model category $M$, except that $A[K]$ is not globally defined as a functor in $A$.

  In particular, taking $n = 1$, the two maps $A^0 \to A^1$ are each required to be anodyne cofibrations, and the induced map $A^0 \amalg A^0 \to A^1$ is required to be a cofibration.
  Thus a left relaxed premodel category has a sufficient supply of ``cylinder objects'' on cofibrant objects, and moreover, when $f : A \to B$ is a cofibration, we may choose cylinder objects on $A$ and $B$ which are compatible in the usual sense.
  (Compare the ``relative cylinder axiom'' of a Baues $I$-category \cite[Definition~I.3.1]{Bau}.)
\end{remark}

\begin{proposition}
  A model category is relaxed.
\end{proposition}

\begin{proof}
  If $M$ is a model category, then $M^{[1]}_\proj$ is also a model category and so has enough cosimplicial resolutions.
  Thus $M$ is left relaxed.
  Now $M^\op$ is also a model category, and so $M$ is right relaxed as well.
\end{proof}

\begin{proposition}
  If $M$ admits a unital action by a premodel category $V$ with enough cosimplicial resolutions, then $M$ is relaxed.
\end{proposition}

\begin{proof}
  Let $F : V \times M \to M$ with $F(I, -) = \id_M$ be a unital action by $V$.
  By \cref{prop:proj-enriched-functoriality}, the extension ${\tilde F}(K, (A_0 \to A_1)) = (F(K, A_0) \to F(K, A_1))$ is again a Quillen adjunction of two variables, and ${\tilde F}(I, -)$ is naturally isomorphic to the identity on $M^{[1]}_\proj$.
  Thus \cref{prop:action-resolutions} implies that $M^{[1]}_\proj$ has enough cosimplicial resolutions and so $M$ is left relaxed.

  By \cref{prop:op-action}, $M^\op$ also admits a unital action by a model category, so $M^\op$ is also left relaxed and therefore $M$ is relaxed.
\end{proof}

\begin{corollary}
  Let $V$ be a model category.
  Then any $V$-premodel category is relaxed.
\end{corollary}

\begin{remark}
  In the next chapter we will develop the homotopy theory of a relaxed premodel category $M$.
  Our eventual intention is to apply this theory to define the weak equivalences of $V\CPM$.
  In that context, we could assume that $M$ comes equipped with a $V$-premodel category structure, and this assumption would actually simplify some of our arguments.
  However, we prefer to work under only the relaxedness assumption.
  One reason is that being relaxed is a \emph{property}, not additional structure.
  This implies in particular that the homotopical notions we will develop are intrinsic to $M$ itself, and do not depend on the $V$-module structure.
  For example, they are unchanged under the ``base change'' operation induced by a strong monoidal Quillen adjunction $V' \adj V$.
  Assuming only that $M$ is relaxed also lets us work in slightly greater generality, as it also allows $M$ to be an ordinary (unenriched) model category.
\end{remark}



\chapter{The homotopy theory of a relaxed premodel category}
\label{chap:homotopy}

Our next task is to understand the ``homotopy theory'' of a premodel category.
The homotopy-theoretic aspects of model category theory (unlike the algebraic aspects, which we were able to extend to premodel categories in the previous chapter) depend in an essential way on the weak equivalences and the two-out-of-three property, and so we cannot reuse the same methods for premodel categories.
As we explained in \cref{sec:relaxed-premodel}, we cannot expect a \emph{general} premodel category to have a homotopy theory which resembles that of a model category and so in this chapter we work only with \emph{relaxed} premodel categories.

Specifically, a major goal of this chapter is to define the \emph{homotopy category} of a relaxed premodel category and show that it is functorial with respect to left Quillen functors.
We then call a left Quillen functor a \emph{Quillen equivalence} when it induces an equivalence of homotopy categories.
We will also need to understand the homotopy category and its relation to the original relaxed premodel category in enough detail in order to obtain a useful criterion for detecting Quillen equivalences.
(When we define the model 2-category structure on combinatorial $V$-premodel categories in \cref{chap:vcpmmodel}, we will define the fibrations and acyclic fibrations by lifting conditions and we will need to be able to check that a fibration is an acyclic fibration if and only if it is a Quillen equivalence.)

A relaxed premodel category has an adequate supply of well-behaved cylinder and path objects which can be used to emulate the traditional left and right homotopy relations on maps from a cofibrant object to a fibrant object.
This is enough to define the \emph{classical homotopy category} of a relaxed premodel category, as the category of objects which are both cofibrant and fibrant with morphisms given by homotopy classes of maps.
This construction has the advantages that it obviously agrees with the corresponding construction for model categories and it is \emph{self-dual}---the classical homotopy category of $M^\op$ is (essentially by definition) the opposite of the classical homotopy of $M$.
However, it is not obviously functorial because left Quillen functors generally do not preserve fibrant objects.

In model category theory, the homotopy category is made functorial in left Quillen functors by using cofibrant replacement to construct left derived functors.
This is sensible because any two cofibrant replacements for the same object are weakly equivalent, by the two-out-of-three condition, and left Quillen functors preserve weak equivalences between cofibrant objects.
Since a relaxed premodel category does not come equipped with a class of weak equivalences, we cannot apply this strategy directly.
In fact, having a notion of weak equivalence (unsurprisingly) turns out to be generally useful for homotopy theory, and so we will define one in terms of the premodel category structure.

An observation that simplifies matters is that we only really need to invent a notion of weak equivalence between cofibrant objects.
After all, every object of a model category is weakly equivalent to a cofibrant one and left Quillen functors preserve cofibrant objects, so for the purpose of defining a homotopy category which is functorial in left Quillen functors, we might as well restrict our attention to cofibrant objects.
Now, weak equivalences between cofibrant objects are also preserved by left Quillen functors so, by extension, we call a morphism between cofibrant objects of a relaxed premodel category a \emph{left weak equivalence} if every left Quillen functor to a model category sends it to a weak equivalence.
This class of morphisms automatically inherits many good properties of the weak equivalences between cofibrant objects of a model category, including (tautologically) being preserved by left Quillen functors.
Accordingly, we define the \emph{homotopy category} $\Ho^L M$ of a relaxed premodel category $M$ as the localization of $M^\cof$ at the class of left weak equivalences.
This construction also obviously agrees (up to equivalence) with the homotopy category of a model category, and is manifestly functorial in left Quillen functors.
It is not obviously self-dual, but we will show that it is equivalent to the classical homotopy category of $M$ and therefore self-dual up to equivalence.
In fact the left homotopy equivalences agree with the dual ``right weak equivalences'' when both are defined, that is, for maps between cofibrant and fibrant objects.

We will show that the cofibrations and left weak equivalences of $M$ make $M^\cof$ into a \emph{cofibration category}, a structure which resembles the ``left half'' of a model category.
This allows us to apply a criterion due to Cisinski for detecting the functors which induce an equivalence of homotopy categories, which are our Quillen equivalences.
It also implies, for example, that the simplicial localization of $M^\cof$ at the left weak equivalences forms a cocomplete $(\infty, 1)$-category.
This simplicial localization is also Dwyer--Kan equivalent\footnote{
  At least when $M$ has functorial factorizations; the author expects that this is also true in general.
}
to the simplicial localization of $M^\cf$ at the homotopy equivalences, and therefore complete as well by duality.
We will also make use of the homotopy theory of a cofibration category to show that a cofibration $f : A \to B$ is a left weak equivalence if and only if there exists an anodyne cofibration $g : B \to C$ such that $gf : A \to C$ is also an anodyne cofibration.
This remarkably useful fact relates the homotopy category of $M$, which is defined in terms of the left weak equivalences, back to the anodyne cofibrations, which can be detected directly using lifting properties.
This is the key tool for establishing the expected relationship between the fibrations, acyclic fibrations, and weak equivalences (Quillen equivalences) of the model 2-category structure on combinatorial $V$-premodel categories.


\section{Homotopy of maps}\label{sec:homotopy}

In this section we introduce anodyne cylinder and path objects and the left and right homotopy relations in a relaxed premodel category.

In a model category, a cylinder object on an object $A$ is usually defined to be a factorization of the fold map $A \amalg A \to A$ as a cofibration followed by a weak equivalence.
This definition cannot be directly adapted to premodel categories, since a premodel category lacks a notion of weak equivalence.
However, when $A$ is cofibrant, there is an alternate definition which does make sense in a premodel category.

\begin{notation}
  We write $\iota_0 : A \to A \amalg B$, $\iota_1 : B \to A \amalg B$ for the inclusions into a coproduct, and $\pi_0 : X \times Y \to X$, $\pi_1 : X \times Y \to Y$ for the projections from a product.
\end{notation}

\begin{definition}\label{def:cylinder-path}
  Let $M$ be a premodel category and $A$ a cofibrant object of $M$.
  An \emph{anodyne cylinder object} for $A$ is a factorization $A \amalg A \xrightarrow{i} C \xrightarrow{q} A$ of the fold map on $A$ such that
  \begin{enumerate}
  \item $i : A \amalg A \to C$ is a cofibration;
  \item $i \iota_0 : A \to C$ and $i \iota_1 : A \to C$ are anodyne cofibrations.
  \end{enumerate}
  Dually, if $X$ is a fibrant object, an \emph{anodyne path object} for $X$ is a factorization $X \xrightarrow{j} Z \xrightarrow{p} X \times X$ of the diagonal map on $X$ such that
  \begin{enumerate}
  \item $p : Z \to X \times X$ is a fibration;
  \item $\pi_0 p : Z \to X$ and $\pi_1 p : Z \to X$ are anodyne fibrations.
  \end{enumerate}
\end{definition}

The definition of an anodyne cylinder object on $A$ requires that $A$ be cofibrant.
This means that we will only be able to define the left homotopy relation between maps with cofibrant domain, but as the left homotopy relation is really only well-behaved in this situation anyways, this is not a serious disadvantage.
Dual comments apply to anodyne path objects.

\begin{notation}
  If $A \amalg A \xrightarrow{i} C \xrightarrow{q} A$ is an anodyne cylinder object, we will use the shorthand notation $i_0 = i \iota_0$, $i_1 = i \iota_1$.
  Similarly, if $X \xrightarrow{j} Z \xrightarrow{p} X \times X$ is an anodyne path object, we will write $p_0 = \pi_0 p$, $p_1 = \pi_1 p$.
\end{notation}

\begin{remark}
  We saw in \cref{remark:concrete-relaxed} that a left relaxed premodel category in particular admits an anodyne cylinder object on each cofibrant object $A$, and a right relaxed premodel category admits an anodyne path object on each fibrant object $X$.
\end{remark}

Henceforth, let $M$ be a relaxed premodel category.

\begin{remark}
  When $M$ is a model category, the notion of an anodyne cylinder object for $A$ agrees with the usual notion of a cylinder object for $A$ as long as $A$ is cofibrant.
  Specifically, we claim that for any factorization $A \amalg A \xrightarrow{i} C \xrightarrow{p} A$ of the fold map with $i$ a cofibration, $p$ is a weak equivalence if and only if both $i_0$ and $i_1$ are acyclic cofibrations.
  Indeed, the inclusions $\iota_0$, $\iota_1 : A \to A \amalg A$ are cofibrations (since $A$ is cofibrant) and hence so are $i_0$ and $i_1$.
  The maps $i_0$ and $i_1$ are each one-sided inverses of the map $p$, so by the two-out-of-three axiom, all three of these maps are weak equivalences whenever any one is.

  Since the notion of a anodyne cylinder object agrees (when it is defined) with the classical notion of a cylinder object in a model category, one might wonder why we bother including the word ``anodyne''.
  The reason is that we will later give $M^\cof$ the structure of a cofibration category, and a cofibration category has a genuinely different notion of a cylinder object: $i_0$ and $i_1$ will only be cofibrations that are left weak equivalences, and not necessarily anodyne cofibrations.
  The prefix ``anodyne'' indicates the requirement that $i_0$ and $i_1$ actually be anodyne cofibrations.
\end{remark}

\begin{definition}
  Let $A$ be a cofibrant object of $M$, $X$ any object of $M$, and $f_0$, $f_1 : A \to X$ two maps.
  A \emph{left homotopy} from $f_0$ to $f_1$ consists of an anodyne cylinder object $A \amalg A \xrightarrow{i} C \xrightarrow{q} A$ together with a map $H : C \to X$ such that $H i_0 = f_0$ and $H i_1 = f_1$.
  We say that $f_0$ and $f_1$ are \emph{left homotopic}, and write $f_0 \htopicl f_1$, if there exists a left homotopy from $f_0$ to $f_1$.

  Dually, suppose instead that $X$ is fibrant while $A$ is an arbitrary object.
  A \emph{right homotopy} from $f_0$ to $f_1$ consists of an anodyne path object $X \xrightarrow{j} Z \xrightarrow{p} X \times X$ together with a map $H : A \to Z$ such that $p_0 H = f_0$ and $p_1 H = f_1$.
  We say that $f_0$ and $f_1$ are \emph{right homotopic}, and write $f_0 \htopicr f_1$, if there exists a right homotopy from $f_0$ to $f_1$.
\end{definition}

The aim of the remainder of this section is to prove that these left and right homotopy relations behave the same way as in a model category.

\begin{proposition}\label{prop:homotopy}
  Let $A$ and $X$ be objects of a relaxed premodel category $M$, and consider the left and right homotopy relations on maps from $A$ to $X$.
  \begin{enumerate}
  \item When $A$ is cofibrant, left homotopy is an equivalence relation.
  \item When $X$ is fibrant, right homotopy is an equivalence relation.
  \item When $A$ is cofibrant and $X$ is fibrant, left and right homotopy agree.
    Moreover, in this situation, left (right) homotopy can be tested on any choice of anodyne cylinder (path) objects.
    That is, if $f_0$, $f_1 : A \to X$ are left homotopic, then for \emph{any} anodyne cylinder object $A \amalg A \xrightarrow{i} C \xrightarrow{q} A$, there exists $H : C \to X$ with $H i_0 = f_0$ and $H i_1 = f_1$; and dually for right homotopy.
  \end{enumerate}
\end{proposition}

In the last case, we say that maps $f_0$, $f_1 : A \to X$ are \emph{homotopic}, and write $f_0 \htopic f_1$, when $f_0$ is left (or equivalently right) homotopic to $f_1$, and call $\simeq$ the \emph{homotopy relation}.

\begin{proof}
  We can reuse the proof for model categories \cite[section~7.4]{Hi}, because \cref{def:cylinder-path} isolates the properties of cylinder and path objects which are actually required.
  We will sketch the argument and fill in the key points, leaving a few details to the reader.

  Suppose $A$ is cofibrant.
  Then $A$ admits some anodyne cylinder object $A \amalg A \xrightarrow{i} C \xrightarrow{q} A$, because $M$ is left relaxed.
  For any map $f : A \to X$, the map $H = f q : C \to X$ defines a left homotopy from $f$ to $f$, so $f \htopicl f$.
  Given an anodyne cylinder object $A \amalg A \xrightarrow{i} C \xrightarrow{q} A$, we may define a new anodyne cylinder object $A \amalg A \xrightarrow{i'} C \xrightarrow{q} A$ where $i'$ is obtained from $i$ by precomposing with the map $A \amalg A \to A \amalg A$ which switches the copies of $A$.
  It is easy to check that this is again an anodyne cylinder object, and that if $H : C \to X$ defines a homotopy from $f_0 : A \to X$ to $f_1 : A \to X$ with respect to the original anodyne cylinder object, then $H$ also defines a homotopy from $f_1$ to $f_0$ with respect to this new anodyne cylinder object.
  Thus, $\htopicl$ is a symmetric relation.
  For transitivity, we must explain how to glue left homotopies.
  Suppose we are given two anodyne cylinder objects $A \amalg A \xrightarrow{i} C \xrightarrow{q} A$ and $A \amalg A \xrightarrow{i'} C' \xrightarrow{q'} A$ on $A$.
  Form the pushout of $C$ and $C'$ and the induced map $q''$ as shown in the diagram below.
  \[
    \begin{tikzcd}
      & & A \ar[d, "i_0"] \\
      & A \ar[r, "i_1"'] \ar[d, "i'_0"] & C \ar[d, "j_0"] \ar[rdd, bend left, "q"] \\
      A \ar[r, "i'_1"'] & C' \ar[r, "j_1"'] \ar[rrd, bend right, "q'"'] & C'' \ar[rd, "q''"] \\
      & & & A
    \end{tikzcd}
  \]
  Let $i'' : A \amalg A \to C''$ be the map induced by $j_0 i_0$ and $j_1 i_1$.
  We claim that $A \amalg A \xrightarrow{i''} C' \xrightarrow{q''} A$ is another anodyne cylinder object on $A$.
  The proof is the same as the proof of \cite[Lemma~7.4.2]{Hi}, with the steps related to the two-out-of-three axiom removed.
  For example, we need to show that $i''_0 = j_0 i_0$ is an anodyne cofibration; this is true because $i_0$ is an anodyne cofibration and $j_0$ is a pushout of the anodyne cofibration $i'_0$.
  Since we can glue anodyne cylinder objects in this fashion, we can also glue left homotopies, and so $\htopicl$ is transitive.
  This proves part (1), and part (2) is dual.

  Now suppose $A$ is cofibrant and $X$ is fibrant.
  Following \cite[Proposition~7.4.7]{Hi}, we will prove that if the maps $f_0$, $f_1 : A \to X$ are left homotopic and $X \xrightarrow{j} Z \xrightarrow{p} X \times X$ is any anodyne path object for $X$, then there exists a right homotopy $H$ from $f_0$ to $f_1$ defined on the given path object, that is, a map $H : A \to Z$ with $p_0 H = f_0$ and $p_1 H = f_1$.
  This statement, together with its dual and the existence of anodyne cylinder (path) objects on anodyne cofibrant (fibrant) objects, implies part (3).
  Again, we can use the same proof as in \cite{Hi}.
  Let $A \amalg A \xrightarrow{i} C \xrightarrow{q} A$ be an anodyne cylinder object for $A$ and $G : C \to X$ a left homotopy from $f_0$ to $f_1$, and form the diagram below.
  \[
    \begin{tikzcd}
      & A \ar[r, "jf_0"] \ar[d, "i_0"'] & Z \ar[d, "p"] \\
      A \ar[r, "i_1"'] & C \ar[r, "{(f_0q, G)}"'] \ar[ru, dotted, "h"] & X \times X
    \end{tikzcd}
  \]
  As $i_0$ is an anodyne cofibration and $p$ is a fibration, the square admits a lift $h$ as shown by the dotted arrow.
  Then setting $H = hi_1$ yields the desired right homotopy.
\end{proof}

\section{The classical homotopy category}

Still following the case of model categories, we can use the results of the previous section to define the ``classical homotopy category'' of a relaxed premodel category.

\begin{proposition}\label{prop:homotopy-congr}
  Let $M$ be a relaxed premodel category and $f_0$, $f_1 : A \to X$ two parallel maps of $M$.
  \begin{enumerate}
  \item If $A$ is cofibrant and $f_0 \htopicl f_1$, then for any map $g : X \to X'$, we have $gf_0 \htopicl gf_1$.
  \item If $X$ is fibrant and $f_0 \htopicr f_1$, then for any map $h : A' \to A$, we have $f_0h \htopicr f_1h$.
  \end{enumerate}
\end{proposition}

\begin{proof}
  Suppose $A$ is cofibrant.
  If $A \amalg A \xrightarrow{i} C \xrightarrow{q} A$ is an anodyne cylinder object for $A$ and $H : C \to X$ is a left homotopy from $f_0$ to $f_1$, then $gH : C \to X'$ is also a left homotopy from $gf_0$ to $gf_1$.
  This proves part (1), and part (2) is dual.
\end{proof}

Recall that $M^\cf$ denotes the full subcategory of $M$ on the objects which are both cofibrant and fibrant.
By \cref{prop:homotopy}, when we restrict to $M^\cf$, the left and right homotopy relations agree.
\Cref{prop:homotopy-congr} then implies that the homotopy relation is a congruence on $M^\cf$.
Therefore, we may form a new category by quotienting each $\Hom_M(X, Y)$ by the homotopy relation.

\begin{definition}
  Let $M$ be a relaxed premodel category.
  The \emph{classical homotopy category} $\pi M^\cf$ of $M$ is the category with the same objects as $M^\cf$ and with
  \[
    \Hom_{\pi M^\cf}(X, Y) = \Hom_M(X, Y) / {\htopic}.
  \]
  The category $\pi M^\cf$ is equipped with a canonical functor $M^\cf \to \pi M^\cf$.
\end{definition}

\begin{definition}
  Let $X$ and $Y$ be cofibrant and fibrant objects of $M$.
  A map $f : X \to Y$ is a \emph{homotopy equivalence} if its image in $\pi M^\cf$ is an isomorphism.
\end{definition}

Concretely, $f : X \to Y$ is a homotopy equivalence if there exists a map $g : Y \to X$ such that $fg \htopic \id_Y$ and $\id_X \htopic gf$.

\begin{proposition}
  In $M^\cf$, the homotopy equivalences satisfy the two-out-of-three condition.
\end{proposition}

\begin{proof}
  This follows automatically from the fact that the class of homotopy equivalences is the preimage under the functor $M^\cf \to \pi M^\cf$ of the class of isomorphisms of $\pi M^\cf$.
\end{proof}

\begin{remark}\label{remark:fib-cylinder}
  Let $A$ be a cofibrant and fibrant object of $M$, and let $A \amalg A \xrightarrow{i} C \xrightarrow{q} A$ be an anodyne cylinder object for $A$.
  Then $C$ is cofibrant, but not necessarily fibrant.
  However, we can factor $q : C \to A$ as an anodyne cofibration $C \to C'$ followed by a fibration $q' : C' \to A$.
  Then $A \amalg A \to C \to C' \xrightarrow{q'} A$ is another anodyne cylinder object for $A$, and $C'$ is fibrant.
  By \cref{prop:homotopy}, for fibrant $X$, we can detect the homotopy relation on maps from $A$ to $X$ using this new anodyne cylinder object.

  It follows that within $M^\cf$, the homotopy relation is not changed by requiring anodyne cylinder (or path) objects to be both cofibrant and fibrant.
  In particular, the classical homotopy category $\pi M^\cf$ depends only on the restriction of the premodel category structure of $M$ to $M^\cf$.
\end{remark}

\begin{remark}
  We will show later that an anodyne cofibration or anodyne fibration between cofibrant and fibrant objects is a homotopy equivalence.
  From this and \cref{remark:fib-cylinder}, it is not hard to see that a functor from $M^\cf$ to a category $C$ which inverts homotopy equivalences also sends homotopic maps to equal maps of $C$.
  Therefore such a functor factors through $\pi M^\cf$.
  Conversely homotopy equivalences in $M^\cf$ become isomorphisms in $\pi M^\cf$ by definition, and so $\pi M^\cf$ is also the localization of $M^\cf$ at the class of homotopy equivalences.
\end{remark}

\begin{remark}
  Duality interchanges anodyne cylinder objects and anodyne path objects, and left and right homotopy.
  If $M$ is a relaxed premodel category, then so is $M^\op$, and we can identify $(M^\cf)^\op$ with $(M^\op)^\cf$.
  Because left and right homotopy agree in $M^\cf$, we can also identify $(\pi M^\cf)^\op$ with $\pi (M^\op)^\cf$.
  In particular, a morphism $f : X \to Y$ in $M^\cf$ is a homotopy equivalence if and only if the corresponding morphism of $(M^\op)^\cf$ is one.
  We summarize the situation by saying that the classical homotopy category construction is \emph{self-dual}.
\end{remark}

\section{Cofibration categories}\label{sec:cof-cats}

The classical homotopy category $\pi M^\cf$ of a relaxed premodel category $M$ is easy to define, and we know that at least when $M$ is a model category, $\pi M^\cf$ correctly captures the intended homotopy theory of $M$.
However, $M^\cf$ has the drawback that it is not functorial with respect to Quillen adjunctions: if $F : M \to N$ is a left Quillen functor, then $F$ preserves cofibrant objects but usually will not preserve fibrant objects, so $F$ does not directly induce a functor $\pi M^\cf \to \pi N^\cf$.

We want to associate a homotopy category (or, better yet, an $(\infty, 1)$-category) to each relaxed premodel category in a way which is functorial with respect to left Quillen functors, and which extends (at least up to equivalence) the usual definition of the homotopy category (or associated $(\infty, 1)$-category) of a model category.
This will allow us to define Quillen equivalences between relaxed premodel categories as those functors which induce equivalences on the level of homotopy categories.
In order to do so, we will make use of the theory of \emph{cofibration categories}.
In this section, we review the parts of this theory that we'll need.

\begin{definition}
  A \emph{cofibration category} is a category $C$ equipped with two subcategories, the \emph{cofibrations} and the \emph{weak equivalences}.
  These are required to satisfy the following conditions.
  (We call a morphism which is both a cofibration and a weak equivalence a \emph{trivial cofibration}.)
  \begin{enumerate}
  \item[C1] Isomorphisms are both cofibrations and weak equivalences.
  \item[C2] The weak equivalences satisfy the \emph{two-out-of-six condition}:
    Suppose $f : W \to X$, $g : X \to Y$ and $h : Y \to Z$ are three composable morphisms.
    If $gf$ and $hg$ are weak equivalences, then so are all of $f$, $g$, and $h$.
  \item[C3] The category $C$ has an initial object $\emptyset$, and for every object $X$, the unique map $\emptyset \to X$ is a cofibration.
  \item[C4] Pushouts of cofibrations exist and are again cofibrations.
    Pushouts of trivial cofibrations are trivial cofibrations.
  \item[C5] Any morphism of $C$ can be factored as a cofibration followed by a weak equivalence.
  \item[C6] Say that an object $X$ of $C$ is \emph{cc-fibrant} if each trivial cofibration $u : X \to Y$ admits a retract $r : Y \to X$ (so $ru = \id_X$).
    Then any object $A$ of $C$ admits a \emph{cc-fibrant approximation}, that is, a weak equivalence to a fibrant object.
  \item[C7] The colimit of a countable sequence $X_0 \xrightarrow{i_0} X_1 \xrightarrow{i_1} X_2 \xrightarrow{i_2} \cdots$ of cofibrations exists, and the induced map from $X_0$ to the colimit is a cofibration which is trivial if each $i_n$ is.
  \item[C8] $C$ admits arbitrary (small) coproducts, and coproducts of arbitrary families of (trivial) cofibrations are again (trivial) cofibrations.
  \end{enumerate}
  For cofibration categories $C$ and $D$, a \emph{cofibration functor} is a functor $F : C \to D$ which preserves cofibrations, trivial cofibrations, pushouts by cofibrations, colimits of countable direct sequences of cofibrations, and arbitrary coproducts.
\end{definition}

\begin{remark}
  We will be working in categories equipped with distinct classes of anodyne cofibrations and trivial cofibrations.
  We choose the term ``trivial cofibrations'' instead of ``acyclic cofibrations'' in order to reduce the risk of confusion between these classes.

  We will also have two notions of fibrant objects: the one that appears in axiom C6 as well as the usual notion in a premodel category.
  We've renamed the condition in axiom C6 to ``cc-fibrant'' in order to distinguish these.
  However, it will turn out that a cofibrant object is cc-fibrant if and only if it is fibrant in the usual sense, so we'll drop the distinction after proving that.
\end{remark}

\begin{remark}
  There are several definitions of cofibration category in the literature and, in keeping with tradition, we have chosen a definition which is different from any of the ones in our sources (\cite{Bau}, \cite{CisCD}, \cite{RB}, \cite{Sz}).
  Our definition is based on that of a ``(homotopy) cocomplete cofibration category'' of \cite{Sz}, combined with the ``axiom on fibrant models'' of \cite{Bau}.
  Using such a strong definition makes sense in our situation because we are interested in using cofibration categories as a tool to understand the homotopy theory of a relaxed premodel category, and the cofibration categories which arise from relaxed premodel categories are particularly well-behaved ones.
  We will take a moment to describe the relationship between our definition and the other ones found in the literature.

  First, \cite{Bau}, \cite{CisCD} and \cite{RB} only require the weak equivalences to satisfy the two-out-of-three axiom, not the two-out-of-six axiom.
  We will refer to this axiom as C2$'$.
  These sources also do not require all objects of a cofibration category to be cofibrant.
  Some of the other axioms are then only required to hold when certain objects involved are cofibrant; we won't describe these differences in detail.

  The definition of cofibration category in \cite{Bau} does not even require the existence of an initial object, let alone that all objects are cofibrant.
  It does include our axioms C1, C2$'$, C4, C5, and C6.
  \cite{Bau} also requires that the pushout of a weak equivalence by a cofibration is again a weak equivalence (``left properness'').
  However, as noted in \cite[Lemma~1.4]{Bau}, this condition is redundant when all objects are cofibrant.

  The definition of an ``ABC cofibration category'' in \cite{RB} includes our axioms C1, C2$'$, C4, C5, C7, and C8, with the caveat that many of the axioms only apply when the objects involved are cofibrant.
  The definition of a ``cat\'egorie d\'erivable \`a droite'' in \cite{CisCD} is the same except that it omits axioms C7 and C8.

  The definition of a cofibration category in \cite{Sz} includes our axioms C1, C2, C3, C4, and C5, and a cofibration category is called ``(homotopy) cocomplete'' if it also satisfies axioms C7 and C8.

  The main point is that we can apply theorems from any of our sources, although we must take care with variations in definitions of concepts such as ``homotopic'' or ``fibrant''.
\end{remark}

\begin{remark}
  We have chosen the term ``cofibration functor'' for the type of functor between cofibration categories that we consider here in preference to the term ``cocontinuous functor'' used in \cite{Sz}, because we want to reserve the latter term for its ordinary meaning of functors which preserve all colimits.
\end{remark}

\begin{definition}
  Let $C$ be a cofibration category.
  The \emph{homotopy category} $\Ho C$ of $C$ is the localization $C[\sW^{-1}]$ of $C$ at its class of weak equivalences $\sW$.
  The \emph{associated $(\infty, 1)$-category} of $C$ is the simplicial localization of $C$ at the class $\sW$.
  (For definiteness, we take the simplicial localization to be defined by the hammock localization $L^H C = L^H(C, \sW)$ of \cite{DK2}.)
\end{definition}

The category $\Ho C$ and the simplicial category $L^H C$ have the same objects as $C$.
They are related by $\Ho C(X, Y) = \pi_0 L^H C(X, Y)$ for any objects $X$, $Y$ of $C$ \cite[Proposition~3.1]{DK2}.
A priori, the sets $\Ho C(X, Y)$ could be large, but due to axiom C6, they are actually small, as we will see below.
The simplicial sets $L^H C(X, Y)$ are large but weakly homotopy equivalent to small simplicial sets.

Because we require every object of a cofibration category to be cofibrant, Ken Brown's factorization lemma implies that a cofibration functor $F : C \to D$ preserves all weak equivalences.
Therefore $F$ induces a functor $\Ho F : \Ho C = C[\sW_C^{-1}] \to D[\sW_D^{-1}] = \Ho D$ on homotopy categories and also a simplicial functor $L^H F : L^H C \to L^H D$ on simplicial localizations.
These constructions $\Ho$ and $L^H$ preserve compositions and identities strictly.

\begin{example}\label{ex:cc-of-model-cat}
  Let $M$ be a model category.
  Then the full subcategory $M^\cof$ of cofibrant objects of $M$ has a cofibration category structure in which the cofibrations and weak equivalences are those of $M$.
  The axioms C1--C8 are all well-known properties of cofibrations and weak equivalences in a model category.
  The homotopy category $\Ho M^\cof$ is equivalent to $\Ho M = M[\sW_M^{-1}]$ via the functor $\Ho M^\cof \to \Ho M$ induced by the inclusion $M^\cof \to M$ and, in fact, this inclusion also induces a Dwyer--Kan equivalence of simplicial categories $L^H M^\cof \to L^H M$ by \cite[Proposition~5.2]{DK3}.
  Thus, we may identify (up to equivalence) the homotopy category and the associated $(\infty, 1)$-category of $M$ with those of the cofibration category $M^\cof$.

  If $F : M \to N$ is a left Quillen functor, then $F$ preserves cofibrant objects and the restriction $F^\cof : M^\cof \to N^\cof$ is a cofibration functor.
  Thus $F$ induces a functor $\Ho F^\cof : \Ho M^\cof \to \Ho N^\cof$ on homotopy categories, which we may identify with the left derived functor $\Ld F$, as well as a simplicial functor $L^H F^\cof : L^H M^\cof \to L^H N^\cof$ between the associated $(\infty, 1)$-categories.
\end{example}

\begin{remark}
  There is a technical subtlety hidden in the preceding argument.
  The category $M^\cof$ is a full subcategory of the cocomplete category $M$.
  It follows that if a diagram $D : K \to M^\cof$ has a colimit in $M$ which belongs to $M^\cof$, then the colimit in $M$ is also the colimit in $M^\cof$.
  However, it is possible that the category $M^\cof$ could have additional colimits which are not colimits in $M$.
  Now suppose $F : M \to N$ is a left Quillen functor; then $F^\cof : M^\cof \to N^\cof$ need not preserve these additional colimits.

  In order to show that $F^\cof$ preserves, say, arbitrary coproducts, we must use the fact that arbitrary coproducts \emph{in $M$} of cofibrant objects are again cofibrant and, therefore, any coproduct in $M^\cof$ is actually a coproduct in $M$ and therefore sent by $F$ to a coproduct in $N$, which is in turn a coproduct in $N^\cof$.
  Similar remarks apply to the preservation of pushouts of cofibrations and countable transfinite compositions of cofibrations.
  We will gloss over this point in future arguments of this type.
\end{remark}

Our aim over the next several sections will be to equip the subcategory $M^\cof$ of each relaxed premodel category $M$ with the structure of a cofibration category in a way which extends the above construction for model categories, and which is again functorial with respect to left Quillen functors.
This cofibration structure will have the same cofibrations as $M$, and a new class of a weak equivalences which we will call the \emph{left weak equivalences}.


In the remainder of this section, we describe the consequences of such a construction.
First and foremost, we obtain a homotopy category $\Ho^L M$, defined as $\Ho M^\cof$, which is (strictly!) \emph{functorial} with respect to left Quillen functors, addressing the major drawback of the classical homotopy category $\pi M^\cf$.
We will see that, in fact, there is an equivalence of categories $\pi M^\cf \to \Ho^L M$, so that $\pi M^\cf$ can be made (pseudo)functorial with respect to left Quillen functors as well.
However, for our primary objective of constructing a model 2-category structure on combinatorial $V$-premodel categories, it will be more convenient to work directly with $\Ho^L M$.

This functoriality allows us to define a \emph{Quillen equivalence} $F : M \to N$ between relaxed premodel categories to be a left Quillen functor which induces an equivalence of homotopy categories.
By construction, the class of Quillen equivalences satisfies the two-out-of-three condition.
When $M$ and $N$ are model categories, this definition is equivalent to the classical one \cite[Proposition~1.3.13]{Ho}.

The following criterion due to Cisinski \cite[Th\'eor\`eme~3.19]{CisCD} allows us to understand the Quillen equivalences between relaxed premodel categories in more concrete terms.

\begin{theorem}
  Let $F : C \to D$ be a cofibration functor between cofibration categories.
  Then $F$ induces an equivalence $\Ho F : \Ho C \to \Ho D$ if and only if it satisfies both of the following conditions.
  \begin{enumerate}
  \item[(AP1)]
    For each morphism $f$ of $C$, $Ff$ is a weak equivalence if and only if $f$ is.
  \item[(AP2)]
    Let $A$ be an object of $C$ and $g : FA \to B$ a morphism of $D$.
    Then there exists a morphism $f : A \to A'$ in $C$ and weak equivalences $B \to B'$ and $FA' \to B'$ making the diagram below commute.
    \[
      \begin{tikzcd}
        FA \ar[r, "Ff"] \ar[d, "g"'] & FA' \ar[d, "\sim"] \\
        B \ar[r, "\sim"'] & B'
      \end{tikzcd}
    \]
  \end{enumerate}
\end{theorem}

Another theorem of Cisinski \cite[Th\'eor\`eme~3.25]{CisInv} gives a third characterization of the Quillen equivalences.

\begin{theorem}\label{thm:eqv-iff-dk}
  Let $F : C \to D$ be a cofibration functor between cofibration categories.
  Then $F$ induces an equivalence $\Ho F : \Ho C \to \Ho D$ if and only if it it induces a Dwyer--Kan equivalence $L^H F : L^H C \to L^H D$.
\end{theorem}

We will not make use of this result other than to note that it offers further justification for considering the Quillen equivalences to be the weak equivalences between relaxed premodel categories.
(Compare the discussion after \cite[Definition~1.7]{Sz}.)

A cofibration category has its own notion of homotopy between maps from a (cofibrant) object to a cc-fibrant object and a ``classical homotopy category'' which is related to its homotopy category $C[\sW^{-1}]$ in the following way.

\begin{theorem}
  Let $C$ be a cofibration category.
  \begin{enumerate}
  \item The homotopy relation is a congruence on the full subcategory $C^\fib$ of cc-fibrant objects of $C$.
  \item A map of $C^\fib$ is a weak equivalence if and only if it is a homotopy equivalence.
  \item
    Write $\pi C^\fib$ for the quotient of $C^\fib$ by the homotopy relation.
    Then the induced functor $\Ho C^\fib \to \pi C^\fib$ (which exists by (2)) is an isomorphism of categories.
  \item The induced functor $\Ho C^\fib \to \Ho C$ is an equivalence of categories.
  \end{enumerate}
\end{theorem}

\begin{proof}
  See \cite[section~II.3]{Bau}.
\end{proof}

In particular, $\Ho C$ is a locally small category because it is equivalent to $\pi C^\fib$.

We will show that if $M$ is a relaxed premodel category then an object $X$ of $M^\cof$ is cc-fibrant if and only if it is fibrant in $M$, and also that the homotopy relation of $M^\cof$ as a cofibration category agrees with the one we defined in \cref{sec:homotopy}.
We may therefore identify $\pi C^\fib$ with the classical homotopy category of $M$, and by part (b), also identify $\Ho C^\fib$ with the localization of $M^\cf$ at the class of homotopy equivalences.

Now, $M^\cf$ and its class of homotopy equivalences are manifestly \emph{self-dual} constructions.
Thus, suppose we replace $M$ by $M^\op$ and repeat the construction of the associated cofibration category.
This will yield a \emph{fibration category} structure on $M^\fib$, satisfying axioms dual to C1--C8, with a class of \emph{right weak equivalences}.
Its homotopy category $\Ho^R M = \Ho M^\fib$ and simplicial localization $L^H M^\fib$ will satisfy results dual to the ones described above.
In particular a morphism in $(M^\fib)^\cof$ is a right weak equivalence if and only if it is a homotopy equivalence (in $M^\op$), and the induced functor $\Ho {(M^\fib)^\cof} \to M^\fib$ is an equivalence.
But $(M^\fib)^\cof$ equals $(M^\cof)^\fib = M^\cf$, and the homotopy equivalences of $M$ and $M^\op$ are the same.
Thus, we obtain a diagram
\[
  \begin{tikzcd}
    & \pi M^\cf \\
    \Ho^L M & \Ho M^\cf \ar[l, "\eqv"'] \ar[u, "\iso"'] \ar[r, "\eqv"] & \Ho^R M
  \end{tikzcd}
\]
and, in particular, we obtain an equivalence $\Ho^L M \eqv \Ho^R M$.

Now suppose that $M$ also admits functorial factorizations, as is the case for instance when $M$ is combinatorial.
Then using \cite[Proposition~3.5]{DK2}, one can show that the inclusion $M^\cf \to M^\cof$ induces not only an equivalence $\Ho M^\cf \to \Ho M^\cof$ but also a Dwyer--Kan equivalence $L^H M^\cf \to L^H M^\cof$.
(We expect that the argument of \cite[(8.1)]{DK3} can be adapted to show that this holds even in the absence of functorial factorizations.)
We therefore conclude that $L^H M^\cof$ and $L^H M^\fib$ are Dwyer--Kan equivalent.

A model category has an associated $(\infty, 1)$-category which is both cocomplete and complete.
The same holds for a relaxed premodel category.
We won't need to make use of this fact, so we just briefly mention it here.

\begin{theorem}\label{thm:cofibration-cocomplete}
  A cofibration category admits all homotopy colimits.
  In particular, if $M$ is a relaxed premodel category, then $L^H M^\cof$ is cocomplete as an $(\infty, 1)$-category.
\end{theorem}

\begin{proof}
  See \cite[Theorem~9.6.3]{RB}.
\end{proof}

Since $L^H M^\cof$ is Dwyer--Kan equivalent to $L^H M^\fib$, the $(\infty, 1)$-category $L^H M^\cof$ is also complete.




\section{Left weak equivalences}

Our current aim is to construct a cofibration category structure on $M^\cof$ as described in the previous section.
More specifically, we will prove the following result.

\begin{theorem}\label{thm:cof}
  There is a unique way to assign to each relaxed premodel category $M$ a class of morphisms of $M^\cof$, the \emph{left weak equivalences}, with the following properties.
  \begin{enumerate}
  \item[(a)] For each relaxed premodel category $M$, the cofibrations of $M$ and the left weak equivalences make $M^\cof$ into a cofibration category.
  \item[(b)] Each left Quillen functor $F : M \to N$ between relaxed premodel categories restricts to a cofibration functor $F^\cof : M^\cof \to N^\cof$.
  \item[(c)] For each relaxed premodel category $M$, the anodyne cofibrations of $M^\cof$ are left weak equivalences.
  \item[(d)] For each \emph{model} category $M$, the left weak equivalences are precisely the weak equivalences of $M$ between cofibrant objects.
  \end{enumerate}
\end{theorem}

In this section we will give several equivalent definitions of the class of left weak equivalences and verify most parts of the above theorem.

\begin{proposition}\label{prop:left-weq-def}
  Let $M$ be a relaxed premodel category and let $f : A \to B$ be a morphism between cofibrant objects of $M$.
  Then the following conditions are equivalent.
  \begin{enumerate}
  \item Every left Quillen functor from $M$ to a model category $N$ takes $f$ to a weak equivalence.
  \item Every left Quillen functor from $M$ to $\Kan^\op$ takes $f$ to a weak equivalence.
  \item For every simplicial resolution $X_\bullet$ in $M$, the induced map of simplicial sets $f^* : \Hom(B, X_\bullet) \to \Hom(A, X_\bullet)$ is a weak equivalence.
  \end{enumerate}
\end{proposition}

\begin{proof}
  Condition (3) is simply a reformulation of condition (2), since a simplicial resolution $X_\bullet$ in $M$ is the same as a left Quillen functor from $\Kan$ to $M^\op$.
  The associated left Quillen functor from $M$ to $\Kan^\op$ sends an object $A$ to the simplicial set $\Hom(A, X_\bullet)$.

  (1) implies (2) is obvious.
  To show (2) implies (1), it suffices to consider the case when $M$ is a model category.
  Using the equivalence of (2) and (3), the assumption is that $f^*$ induces a weak equivalence on all mapping spaces to fibrant objects and this is one characterization of the weak equivalences of a model category.
  (We do not need to derive the mapping spaces $\Hom(A, X_\bullet)$ and $\Hom(B, X_\bullet)$ because $A$ and $B$ are already cofibrant and $X$ is a simplicial resolution.)
\end{proof}

\begin{definition}\label{def:left-weq}
  A map between cofibrant objects of a relaxed premodel category is a \emph{left weak equivalence} if it satisfies the above equivalent conditions.
\end{definition}

\begin{remark}
  \Cref{prop:left-weq-def} holds even when $M$ is not relaxed.
  However, we have stated it only for relaxed premodel categories because we have no strong reason to believe that the class of morphisms it describes is the correct notion of weak equivalence except for a premodel category that is relaxed, or at least one that has enough simplicial resolutions.
\end{remark}

We call a morphism a \emph{trivial cofibration} if it is both a cofibration and a left weak equivalence.

With this choice of the left weak equivalences, we can immediately verify parts (b), (c) and (d) of \cref{thm:cof}.

\begin{proposition}\label{prop:cc-functorial}
  Any left Quillen functor $F : M \to N$ between relaxed premodel categories restricts to a cofibration functor between the associated cofibration categories.
\end{proposition}

\begin{proof}
  Once we verify that $F$ preserves trivial cofibrations, we can use the same argument that we used for model categories in \cref{ex:cc-of-model-cat}.
  In fact, we'll show directly that $F$ preserves all left weak equivalences.
  This follows immediately from characterization (1) of the left weak equivalences, because the composition of left Quillen functors is a left Quillen functor.
\end{proof}

\begin{proposition}\label{prop:acof-is-weq}
  An anodyne cofibration between cofibrant objects is a left weak equivalence.
\end{proposition}

\begin{proof}
  Follows from characterization (1), because left Quillen functors preserve anodyne cofibrations and an anodyne cofibration in a model category is a weak equivalence.
\end{proof}

\begin{warning}
  The converse implication is not true in general: not every trivial cofibration is an anodyne cofibration.
  This is a major point of departure from the theory of model categories.
  In a model category, the anodyne cofibrations and trivial cofibrations are in fact the same; this is a consequence of the following fact.
\end{warning}

\begin{proposition}\label{prop:model-left-weq-iff-weq}
  If $M$ is a model category, then the left weak equivalences of $M$ are exactly the weak equivalences between cofibrant objects.
\end{proposition}

\begin{proof}
  Let $f : A \to B$ be a map of $M$ between cofibrant objects.
  If $f$ is a left weak equivalence, then the identity functor $M \to M$ is a left Quillen functor so $f$ is also a weak equivalence.
  Conversely, if $f$ is a weak equivalence, then the image of $f$ under any left Quillen functor to another model category is also a weak equivalence.
\end{proof}

\begin{remark}\label{remark:left-weq-motivation}
  The fact that left Quillen functors between model categories preserve weak equivalences between cofibrant objects gives a kind of formal justification for characterization (1) of the left weak equivalences: it is the largest class of maps preserved by left Quillen functors with the correct value for model categories.

  However, we may also arrive at characterization (3) by following a common strategy for defining the weak equivalences in a model category, namely, defining $f : A \to B$ to be a weak equivalence when $f$ induces some kind of equivalence on mapping spaces or homotopy classes of maps into each fibrant object $X$.
  This strategy relies on a preexisting notion of mapping space or homotopy, which may only be well-behaved when the objects $A$ and $B$ are cofibrant.
  For example, the weak equivalences in a left Bousfield localization are the maps which induce weak equivalences on the (derived) spaces of maps to each local object.
  In this case, the original model category structure provides a notion of mapping space which is well-behaved for all objects.
  In Cisinski model structures \cite{CisMC} and their generalization by Olschok to locally presentable categories \cite{Ol}, the weak equivalences are the maps which define isomorphisms on homotopy classes of maps to each fibrant object, where homotopy is defined in terms of a cylinder functor which is part of the input data.
  This homotopy relation is only well-behaved for cofibrant domains, which is why these constructions require all objects to be cofibrant in order to produce a model category structure.

  In our case, we obtain a notion of mapping spaces from simplicial resolutions of fibrant objects.
  These mapping spaces are also only well-behaved for cofibrant domains.
  But this is not a problem because our objective is only to define a notion of weak equivalence between cofibrant objects, not on all of $M$.
\end{remark}

Most of the axioms of a cofibration category are also formal consequences of the definition of the left weak equivalences.

\begin{proposition}\label{prop:cof-weq}
  The left weak equivalences of $M^\cof$ form a subcategory which includes the isomorphisms and satisfies the two-out-of-six condition.
\end{proposition}

\begin{proof}
  Clearly isomorphisms are left weak equivalences.
  For the two-out-of-six condition, let $f : W \to X$, $g : X \to Y$, $h : Y \to Z$ be three composable morphisms between cofibrant objects and suppose that $gf$ and $hg$ are left weak equivalences.
  Then for any model category $N$ and any left Quillen functor $F : M \to N$, the morphisms $F(gf) = (Fg)(Ff)$ and $F(hg) = (Fh)(Fg)$ are weak equivalences.
  The weak equivalences of a model category satisfy the two-out-of-six condition, so all of $Ff$, $Fg$, $Fh$ are weak equivalences.
  Since $F$ was arbitrary, we conclude that all of $f$, $g$, $h$ are left weak equivalences in $M$.
  Similarly, left weak equivalences are closed under composition because the weak equivalences of a model category are.
\end{proof}

\begin{proposition}\label{prop:cof-tcof}
  Trivial cofibrations are closed under pushout, countable transfinite compositions, and arbitrary coproducts.
\end{proposition}

\begin{proof}
  Let $f : A \to B$ be a trivial cofibration and $g : A \to A'$ an arbitrary map in $M^\cof$, and form the pushout $f'$ of $f$ along $g$.
  We need to show that $f'$ is again a trivial cofibration, so let $F : M \to N$ be a left Quillen functor from $M$ to a model category $N$.
  Since $f$ is a trivial cofibration, its image $Ff$ is a cofibration and also a weak equivalence, hence an acyclic cofibration of $N$.
  Now $F$ preserves pushouts, so $Ff'$ is a pushout of $Ff$, hence also an acyclic cofibration.
  Therefore $f'$ is again a left weak equivalence, and it is a cofibration because the cofibrations of $M$ are closed under pushouts.

  The other parts follow by similar arguments, since a left Quillen functor also preserves transfinite compositions and arbitrary coproducts.
\end{proof}

By a standard argument \cite{Bro}, the factorization axiom C5 follows from the existence of cylinder objects.

\begin{proposition}\label{prop:factor-cof-weq}
  Any morphism $f : A \to B$ between cofibrant objects can be factored as a cofibration followed by a left weak equivalence.
\end{proposition}

\begin{proof}
  Choose an anodyne cylinder object $A \amalg A \xrightarrow{i} C \xrightarrow{p} A$ for $A$ and form the diagram
  \[
    \begin{tikzcd}
      A \ar[r, hookrightarrow, "\iota_1"] \ar[d, "f"'] &
      A \amalg A \ar[r, hookrightarrow, "i"] \ar[d, "\id_A \amalg f"'] &
      C \ar[r, "p"] \ar[d] &
      A \ar[d, "f"] \\
      B \ar[r, hookrightarrow, "\iota_1"'] &
      A \amalg B \ar[r, hookrightarrow, "i'"'] &
      C \amalg_A B \ar[r, "p'"'] &
      B
    \end{tikzcd}
  \]
  in which the left two squares (and therefore also the third square) are pushouts.
  The bottom composition is the identity, and $i'\iota_1$ is a pushout of the anodyne cofibration $i\iota_1 = i_1$, hence an anodyne cofibration and in particular a left weak equivalence.
  Therefore $p'$ is a left weak equivalence by the two-out-of-three property.
  Then $A \xrightarrow{i'\iota_0} C \amalg_A B \xrightarrow{p'} B$ is a factorization of $f$ of the required form.
\end{proof}

The only difficult part of proving that the left weak equivalences make $M^\cof$ into a cofibration category is identifying the cc-fibrant objects with the objects which are fibrant in the usual sense, so that axiom C6 can be verified using the factorization in $M$.
We summarize the results of this section below.

\begin{proposition}\label{prop:almost-cof-cat}
  Let $M$ be a relaxed premodel category and equip $M^\cof$ with the class of left weak equivalences of \cref{def:left-weq}.
  Assume that any cofibrant and fibrant object of $M$ is also cc-fibrant in $M^\cof$.
  Then $M^\cof$ is a cofibration category.
\end{proposition}

\begin{proof}
  We verified axiom C5 in \cref{prop:factor-cof-weq}.
  For C6, suppose $A$ is an object of $M^\cof$ and use the factorization system of $M$ to factor the map $A \to *$ as an anodyne cofibration $A \to X$ followed by a fibration $X \to *$.
  Then $X$ is cofibrant and fibrant in $M$, so by assumption $X$ is also cc-fibrant, and by \cref{prop:acof-is-weq} the anodyne cofibration $A \to X$ is a left weak equivalence.

  We split the remaining axioms into groups depending on whether they pertain to the weak equivalences (C1, C2), cofibrations (C1, C3, C4, C7, C8), or trivial cofibrations (C4, C7, C8).
  The axioms for the weak equivalences were checked in \cref{prop:cof-weq}.
  The axioms for the cofibrations follow from basic properties of weak factorization systems.
  Finally, the axioms for the trivial cofibrations were checked in \cref{prop:cof-tcof}.
\end{proof}

We will verify the assumption of \cref{prop:almost-cof-cat} in the next section.
This will complete the proof of the existence part of \cref{thm:cof}; we will prove the uniqueness part in \cref{sec:lwe-revisited}.

\section{Trivial cofibrations and lifting}

A morphism $f : A \to B$ in a model category is an acyclic cofibration if and only if it has the left lifting property with respect to all fibrations.
If we know a priori that $f$ is a cofibration, then it suffices to check the lifting property for fibrations between fibrant objects.

\begin{proposition}
  Let $M$ be a model category and $f : A \to B$ a cofibration in $M$.
  Then $f$ is an acyclic cofibration if and only if $f$ has the left lifting property with respect to every fibration between fibrant objects.
\end{proposition}

\begin{proof}
  The ``if'' direction is clear, so suppose $f : A \to B$ is a cofibration with the left lifting property with respect to fibrations between fibrant objects.
  Factor $B \to *$ into an acyclic cofibration $B \to \hat B$ followed by a fibration $\hat B \to *$, and then factor the composition $A \xrightarrow{f} B \to \hat B$ into an acyclic cofibration $A \to \hat A$ followed by a fibration $\hat A \xrightarrow{p} \hat B$, producing the square below.
  \[
    \begin{tikzcd}
      A \ar[r, hookrightarrow, "\sim"] \ar[d, hookrightarrow, "f"'] & \hat A \ar[d, two heads, "p"] \\
      B \ar[r, hookrightarrow, "\sim"] \ar[ru, dotted, "\exists l"] & \hat B \ar[r, two heads] & *
    \end{tikzcd}
  \]
  The map $p$ is a fibration between fibrant objects, so by assumption $f$ has the left lifting property with respect to $p$ and so the square admits a lift $l$ as shown by the dotted arrow.
  Now applying the two-out-of-six condition to $f$, $l$ and $p$, we conclude that all three maps are weak equivalences, and in particular $f$ is an acyclic cofibration.
\end{proof}

In a relaxed premodel category, the same condition instead characterizes the \emph{trivial} cofibrations.

\begin{proposition}\label{prop:tcof-lift}
  Let $M$ be a relaxed premodel category and $f : A \to B$ a cofibration between cofibrant objects of $M$.
  Then $f$ is a trivial cofibration if and only if $f$ has the left lifting property with respect to every fibration between fibrant objects.
\end{proposition}

\begin{proof}
  The same argument proves the ``only if'' direction, because anodyne cofibrations are left weak equivalences and the left weak equivalences satisfy the two-out-of-six property.

  Conversely, suppose that $f : A \to B$ is a trivial cofibration between cofibrant objects and $p : X \to Y$ is a fibration between fibrant objects.
  We must show that any square
  \[
    \begin{tikzcd}
      A \ar[d, "f"'] \ar[r] & X \ar[d, "p"] \\
      B \ar[r] & Y
    \end{tikzcd}
  \]
  admits a lift.
  Choose a simplicial resolution of $p$.
  That is, $p$ being a fibration between fibrant objects of $M$, it represents a cofibration between cofibrant objects of $M^\op$ and therefore a cofibrant object of $(M^\op)^{[1]}_\proj$.
  Since $M$ is right relaxed, $M^\op$ is left relaxed and so $(M^\op)^{[1]}_\proj$ has enough cosimplicial resolutions.
  Therefore, we may choose a left Quillen functor $F : \Kan \to (M^\op)^{[1]}_\proj$ sending $\Delta^0$ to the object corresponding to $p$.

  Now since $[1]$ is a direct category, $(M^\op)^{[1]}_\proj$ is also the Reedy premodel category structure $(M^\op)^{[1]}_\Reedy$, where $[1]$ has the Reedy structure in which the nonidentity morphism increases degree.
  Dually, as $[1]$ is also an inverse category, $\Kan^{[1]}_\inj$ is the Reedy model category $\Kan^{[1]}_\Reedy$ where this time the nonidentity morphism of $[1]$ decreases degree.
  Thus the Reedy--Reedy adjunction takes the form
  \[
    \LQF(\Kan^{[1]}_\inj, M^\op) \eqv \LQF(\Kan, (M^\op)^{[1]}_\proj).
  \]
  Let $F' : \Kan^{[1]}_\inj \to M^\op$ be the preimage of $F$ under this equivalence.
  Then the right adjoint of $F'$ is a right Quillen functor $G' : M^\op \to \Kan^{[1]}_\inj$.
  Unpack $G'$ into components $G'_0 : M^\op \to \Kan$ and $G'_1 : M^\op \to \Kan$ and a natural transformation $g : G'_0 \to G'_1$.
  Since $G'$ is a right Quillen functor, it sends the trivial cofibration $f : A \to B$ to an acyclic fibration in $\Kan^{[1]}_\inj$.
  In particular, the induced map $s : G'_0 B \to G'_0 A \times_{G'_1 A} G'_1 B$ is an acyclic Kan fibration.
  
  Now unfolding the adjunction relationships and using the fact that $F$ sends $\Delta^0$ to the object $p$, we compute that for any object $A$ of $M$, the map on $0$-simplices of the component $g_A$ is given by $p_* : \Hom(A, X) \to \Hom(A, Y)$.
  Hence the map on $0$-simplices of $s$ is the map $\Hom(B, X) \to \Hom(A, X) \times_{\Hom(A, Y)} \Hom(B, Y)$.
  Since $s$ is an acyclic Kan fibration, this map is a surjection and therefore any lifting problem as above has a solution.
\end{proof}

We can now complete the proof that the left weak equivalences make $M^\cof$ into a cofibration category.

\begin{proposition}\label{prop:cc-fibrant-iff-fibrant}
  Let $M$ be a relaxed premodel category and $X$ a cofibrant object of $M$.
  Then $X$ is cc-fibrant in $M^\cof$ if and only if $X$ is fibrant in $M$.
\end{proposition}

\begin{proof}
  Suppose that $X$ is fibrant in $M$.
  Then the map $X \to *$ is a fibration between fibrant objects, so by \cref{prop:tcof-lift} it has the right lifting property with respect to any trivial cofibration $f : X \to Y$, and so $X$ is cc-fibrant.

  Conversely, suppose that $X$ is cc-fibrant.
  Consider a lifting problem
  \[
    \begin{tikzcd}
      A \ar[r] \ar[d, "f"'] & X \ar[d] \\
      B \ar[r] & *
    \end{tikzcd}
  \]
  with $f : A \to B$ an arbitrary anodyne cofibration.
  Form the pushout square as shown on the left below.
  \[
    \begin{tikzcd}
      A \ar[r] \ar[d, "f"'] & X \ar[r, equal] \ar[d, "f'"'] & X \ar[d] \\
      B \ar[r] & Y \ar[r] \ar[ru, dotted] & *
    \end{tikzcd}
  \]
  The map $f' : X \to Y$ is an anodyne cofibration (because it is a pushout of $f$), hence a trivial cofibration.
  Moreover, its domain $X$ is cofibrant.
  Since $X$ is cc-fibrant, we can fill in the dotted arrow, and then the composition $B \to Y \to X$ provides the desired lift of the original square.
\end{proof}

\begin{proposition}\label{prop:exists-cof-cat}
  The cofibrations of $M$ and the left weak equivalences make $M^\cof$ into a cofibration category.
\end{proposition}

\begin{proof}
  Follows from \cref{prop:almost-cof-cat} and \cref{prop:cc-fibrant-iff-fibrant}.
\end{proof}

Now that we have constructed the cofibration category structure on $M^\cof$, we officially make some definitions which were previewed in \cref{sec:cof-cats}.

\begin{definition}
  The \emph{homotopy category} $\Ho^L M$ of a relaxed premodel category $M$ is the localization of $M^\cof$ at the class of weak equivalences.
  If $F : M \to N$ is a left Quillen functor between relaxed premodel categories, then $F$ preserves cofibrant objects and left weak equivalences so it induces a functor $\Ho^L F : \Ho^L M \to \Ho^L N$.
  We call $F$ a \emph{Quillen equivalence} if $\Ho^L F$ is an equivalence of categories.
\end{definition}

We can also show that the various notions of homotopy in a cofibration category considered in the literature agree with the notion of left homotopy that we defined in \cref{sec:homotopy}.

\begin{proposition}
  Let $A$ and $X$ be objects of $M^\cof$ with $X$ fibrant and $f_0$, $f_1 : A \to X$ two maps.
  Then the following conditions are equivalent.
  \begin{enumerate}
  \item $f_0$ and $f_1$ are homotopic in the sense of \cref{sec:homotopy}.
  \item $f_0$ and $f_1$ are homotopic in the sense of \cite[section~I.1]{Bau}.
    That is, there exists a factorization of the fold map $A \amalg A \to A$ into a cofibration $(i_0, i_1) : A \amalg A \to C$ followed by a left weak equivalence $C \to A$, together with a map $H : C \to X$ with $Hi_0 = f_0$ and $Hi_1 = f_1$.
  \item $f_0$ and $f_1$ are left homotopic in the sense of \cite[section~6.3]{RB}.
    That is, there exists a weak equivalence $w : X \to X'$ such that $wf_0$ and $wf_1$ are homotopic in the sense of (2).
  \end{enumerate}
\end{proposition}

\begin{proof}
  By \cref{prop:homotopy}, homotopy of $f_0$ and $f_1$ may be checked on any choice of anodyne cylinder object for $A$.
  By \cite[Proposition~II.2.2]{Bau}, condition (2) satisfies the analogous property.
  Hence, to show the equivalence of (1) and (2), it suffices to exhibit some factorization $A \amalg A \to C \to A$ which is both an anodyne cylinder object and a cylinder object in the sense of (2).
  In fact, any anodyne cylinder object will do, because if $A \amalg A \xrightarrow{i} C \xrightarrow{p} A$ is a cylinder object for $A$, then $p$ is a one-sided inverse to the anodyne cofibration $i_0 : A \to C$, and so $p$ is a left weak equivalence by the two-out-of-three condition.

  Clearly (2) implies (3) since we can take $w = \id_X$.
  Conversely, suppose (3) holds, so there exists a commutative diagram of the form below.
  \[
    \begin{tikzcd}
      A \amalg A \ar[r, "{(f_0, f_1)}"] \ar[d, hookrightarrow, "i"'] & X \ar[d, "w", "\sim"'] \\
      C \ar[r, "H"'] & X'
    \end{tikzcd}
  \]
  By factoring $X' \to *$ into an anodyne cofibration $j : X' \to X''$ followed by a fibration $X'' \to *$ and replacing $w$ by $jw$ and $H$ by $jH$, we may assume that $X'$ is fibrant.
  Then $w$ induces a map $w_* : \Hom(A, X)/{\htopic} \to \Hom(A, X')/{\htopic}$.
  By \cite[Proposition~II.2.11]{Bau}, $w_*$ is an isomorphism.
  Then since $wf_0$ and $wf_1$ are homotopic, $f_0$ and $f_1$ must be as well.
\end{proof}

\begin{remark}
  Condition (3) has the advantage of being well-behaved even when $X$ is not fibrant, which is crucial if one does not assume the existence of fibrant approximations.
\end{remark}

Now that we have verified that our notions of fibrant objects, homotopy and homotopy equivalence agree with those of \cite{Bau}, we can freely apply the homotopy theory of a cofibration category developed therein.

\begin{definition}[{\cite[section~I.1]{Bau}}]
  Let $u : A \to B$ be a cofibration in $M^\cof$.
  A \emph{relative cylinder} on $u$ is a factorization of $B \amalg_A B \to B$ into a cofibration followed by a weak equivalence.

  Let $X$ be a fibrant object of $M^\cof$.
  Two morphisms $f_0$, $f_1 : B \to X$ are \emph{homotopic rel $u$} if there exists a relative cylinder $B \amalg_A B \xrightarrow{i} C \xrightarrow{\sim} B$ and a map $H : C \to X$ with $Hi_0 = f_0$ and $Hi_1 = f_1$.
\end{definition}

\begin{remark}\label{remark:homotopic-rel-id}
  Here $i_0$ and $i_1$ denote the compositions of $i : B \amalg_A B \to C$ with the two inclusions $B \to B \amalg_A B$.
  These two inclusions agree on $A$, so $i_0u = i_1u$ and therefore if $f_0$ and $f_1$ are homotopic rel $u$, we must have $f_0u = f_1u$.

  In particular, suppose $u : B \to B$ is the identity map.
  Then two maps $f_0$, $f_1 : B \to X$ are homotopic rel $u$ if and only if they are equal.

  When $u : \emptyset \to B$ is the unique morphism from the initial object, homotopy rel $u$ is the same as the ordinary homotopy relation.
\end{remark}

\begin{proposition}[{\cite[Proposition~II.2.2]{Bau}}]\label{prop:homotopy-rel-any}
  Homotopy rel $u$ is an equivalence relation on maps $B \to X$, and can be detected on any choice of relative cylinder on $u$.
\end{proposition}

\begin{proposition}[{\cite[Proposition~II.2.12]{Bau}}]\label{prop:dold}
  In the commutative diagram
  \[
    \begin{tikzcd}[column sep=tiny]
      & A \ar[ld, "u"'] \ar[rd, "v"] \\
      U \ar[rr, "g"'] & & V
    \end{tikzcd}
  \]
  suppose $u$ and $v$ are cofibrations, $U$ and $V$ are fibrant and $g$ is a left weak equivalence.
  Then $g$ is a homotopy equivalence under $A$.
  In other words, there exists $f : V \to U$ with $fv = u$, $fg$ homotopic to $\id_U$ rel $u$ and $gf$ homotopic to $\id_V$ rel $v$.
\end{proposition}

\begin{proposition}\label{prop:cf-weq}
  Let $X$ and $Y$ be objects of $M^\cf$ and $f : X \to Y$ a map.
  Then the following are equivalent.
  \begin{enumerate}
  \item $f$ is a left weak equivalence (as a morphism in $M^\cof$).
  \item $f$ is a homotopy equivalence.
  \item $f$ is a right weak equivalence (as a morphism in $M^\fib$).
  \end{enumerate}
  In particular, an anodyne fibration in $M^\cof$ \emph{between fibrant objects} is also a left weak equivalence (since it is a right weak equivalence in $M^\fib$).
\end{proposition}

\begin{proof}
  The equivalence between (1) and (2) is \cref{prop:dold} with $A = \emptyset$.
  The equivalence between (2) and (3) is dual.
\end{proof}

\begin{proposition}
  The localization $\Ho M^\cf$ at the class of homotopy equivalences is isomorphic to the classical homotopy category $\pi M^\cf$.
\end{proposition}

\begin{proof}
  See \cite[section~II.3]{Bau}.
\end{proof}

\section{Multiplicative structure}

Suppose $F : M_1 \times M_2 \to N$ is a Quillen bifunctor.
Then one of the conditions on $F$ is that if $f_1 : A_1 \to B_1$ is a cofibration and $f_2 : A_2 \to B_2$ is an anodyne cofibration, then $f_1 \bp_F f_2$ is an anodyne cofibration.
In this section we show that an analogous property also holds with ``anodyne'' replaced by ``trivial''.

\begin{proposition}\label{prop:bifunctor-trivial}
  Let $(F : M_1 \times M_2 \to N, G_1, G_2)$ be a Quillen adjunction of two variables and assume that $M_2$ and $N$ are relaxed.
  Let $f_1 : A_1 \to B_1$ be a cofibration in $M_1^\cof$ and $f_2 : A_2 \to B_2$ a trivial cofibration in $M_2^\cof$.
  Then $f_1 \bp_F f_2$ is a trivial cofibration in $N^\cof$.
\end{proposition}

\begin{proof}
  By \cref{prop:tcof-lift}, it suffices to verify that $f_1 \bp_F f_2$ has the left lifting property with respect to every fibration $p : X \to Y$ between fibrant objects of $N$.
  For each such $p$, by the usual adjunction argument, this is equivalent to $f_2$ having the left lifting property with respect to the induced morphism $G_2(B_1, X) \to G_2(A_1, X) \times_{G_2(A_1, Y)} G_2(B_1, X)$.
  This latter morphism is a fibration between fibrant objects because $f_1$ is a cofibration between cofibrant objects and $p : X \to Y$ is a fibration between fibrant objects.
  Therefore, we conclude that the required lifting property holds by using \cref{prop:tcof-lift} again.

  (The author learned this efficient argument from the proof of \cite[Theorem~3.2]{He}.)
\end{proof}

From this it follows for example that the homotopy category of a monoidal relaxed premodel category inherits a monoidal structure.
We will not need to make use of results of this nature, so we do not pursue this further here.

\section{Strong deformation retracts}


In \cref{prop:cf-weq} we showed that in $M^\cf$, the classes of left weak equivalences, right weak equivalences, and homotopy equivalences agree.
In this section we'll show that in $M^\cf$ the classes of anodyne cofibrations and trivial cofibrations also agree, as do the classes of anodyne fibrations and trivial fibrations.
Informally, we could say that $M^\cf$ resembles the full subcategory on the cofibrant and fibrant objects of a \emph{model} category.

Our strategy is as follows.
By the general homotopy theory of cofibration categories, any trivial cofibration $u$ between cofibrant and fibrant objects is the inclusion of a \emph{strong deformation retract}.
The relative homotopy in a strong deformation retraction can be detected on any choice of relative cylinder.
Using the condition that $M$ is left relaxed, we can choose a relative cylinder in which certain structural maps are anodyne cofibrations.
We will then be able to express $u$ as a retract of one of these structural maps.
Because anodyne cofibrations are closed under retracts, it follows that $u$ itself is also an anodyne cofibration.

Only \cref{prop:tcof-to-acof} of this section will be used subsequently, and everything else introduced in this section is needed only in its proof.
We first show how homotopy rel $u$ can be detected using only the cofibrations and anodyne cofibrations of $M$.

\begin{lemma}\label{lemma:homotopy-rel-anodyne}
  Let $u : A \to B$ be a cofibration in $M^\cof$ and $X$ a fibrant object of $M^\cof$.
  Suppose $f_0$, $f_1 : B \to X$ are two maps which are homotopic rel $u$.
  Then we can construct the following:
  \begin{enumerate}
  \item a diagram
    \[
      \begin{tikzcd}
        A \amalg A \ar[r, hookrightarrow, "i"] \ar[d, hookrightarrow, "u \amalg u"'] &
        C \ar[r, "p"] \ar[d, hookrightarrow, "v"] &
        A \ar[d, hookrightarrow, "u"] \\
        B \amalg B \ar[r, hookrightarrow, "j"] &
        D \ar[r, "q"] &
        B
      \end{tikzcd}
    \]
    making $v : C \to D$ an anodyne cylinder object on the cofibrant object $u : A \to B$ of the category $M^{[1]}_\proj$;
  \item a map $H : D \to X$ with $Hj_0 = f_0$, $Hj_1 = f_1$, and $Hv = f_0up = f_1up$.
  \end{enumerate}
\end{lemma}

\begin{remark}
  The condition that $v : C \to D$ is an anodyne cylinder object on $u : A \to B$ encodes a total of six conditions that certain maps are (anodyne) cofibrations, most of which are not encoded in the above diagram.
  For example, the ``pushout corner map'' in the left square is a cofibration, and the pushout corner maps of the similar squares obtained by replacing $i$ and $j$ by $i_0$ and $j_0$ or $i_1$ and $j_1$ are \emph{anodyne} cofibrations and not just trivial cofibrations.
\end{remark}

\begin{proof}
  Since $M$ is left relaxed, $M^{[1]}_\proj$ has enough cosimplicial resolutions, so we may choose an anodyne cylinder object $v : C \to D$ on $u : A \to B$.
  This gives a diagram as shown in (1).
  By forming pushout squares repeatedly, we may expand this to the solid part of the diagram below, in which each of the three squares is a pushout.
  \[
    \begin{tikzcd}
      A \amalg A \ar[r, hookrightarrow, "i"] \ar[d, hookrightarrow, "u \amalg u"'] &
      C \ar[r, "p", "\sim"'] \ar[d, hookrightarrow] &
      A \ar[d, hookrightarrow] \\
      B \amalg B \ar[r, hookrightarrow] \ar[rd, hookrightarrow, "j"'] &
      B \amalg_A C \amalg_A B \ar[r, "\sim"'] \ar[d, hookrightarrow] &
      B \amalg_A B \ar[d, hookrightarrow, "j'"] \\
      &
      D \ar[r, "p'", "\sim"'] \ar[rd, "q"', "\sim"] &
      E \ar[d, "q'"] \ar[r, dotted, "H'"] &
      X \\
      & & B
    \end{tikzcd}
  \]
  The composition $B \amalg_A B \to E \to B$ sends each copy of $B$ to $B$ via the identity, via inspection of the triangle with diagonal $qj$.
  The map $B \amalg_A B \to E$ is a cofibration because it is a pushout of the map $B \amalg_A C \amalg_A B \to D$, which is in turn a cofibration because $v$ is an anodyne cylinder object on $u$ in $M^{[1]}_\proj$.
  The maps $p$ and $q$ are left weak equivalences by two-out-of-three, because $A \amalg A \xrightarrow{i} C \xrightarrow{p} A$ and $B \amalg B \xrightarrow{j} D \xrightarrow{q} B$ are in particular anodyne cylinder objects on $A$ and $B$ respectively.
  The map $p'$ is a pushout of $p$, hence a left weak equivalence, because a cofibration category (with all objects cofibrant) is left proper.
  Therefore the induced map $q'$ is also a left weak equivalence.
  So, we conclude that $B \amalg_A B \xrightarrow{j'} E \xrightarrow{q'} B$ is a relative cylinder on $u : A \to B$.

  Now suppose $f_0$, $f_1 : B \to X$ are two maps which are homotopic rel $u$.
  By \cref{prop:homotopy-rel-any}, the homotopy can be detected by this relative cylinder, that is, there exists a map $H' : E \to X$ with $H'j'_0 = f_0$ and $H'j'_1 = f_1$.
  Define $H : D \to X$ by $H = H'p'$.
  Then $Hj_0 = H'p'j_0 = H'j'_0 = f_0$ and similarly $Hj_1 = f_1$.
  Moreover, $Hv = H'p'v = H'j'_0up = f_0up$ and also $Hv = f_1up$ because the map $A \to B \amalg_A B$ in the above diagram equals the composition of $u : A \to B$ with either inclusion $B \to B \amalg_A B$.
\end{proof}

We therefore introduce the following version of a strong deformation retraction.

\begin{definition}
  Let $u : A \to B$ be a cofibration in $M^\cof$.
  An \emph{anodyne strong deformation retraction} of $u$ consists of the following data.
  \begin{enumerate}
  \item A map $r : B \to A$ with $ru = \id_A$.
  \item An anodyne cylinder object $v : C \to D$ on $u$, forming a diagram as in \cref{lemma:homotopy-rel-anodyne}.
    \[
      \begin{tikzcd}
        A \amalg A \ar[r, hookrightarrow, "i"] \ar[d, hookrightarrow, "u \amalg u"'] &
        C \ar[r, "p"] \ar[d, hookrightarrow, "v"] &
        A \ar[d, hookrightarrow, "u"] \\
        B \amalg B \ar[r, hookrightarrow, "j"] &
        D \ar[r, "q"] &
        B
      \end{tikzcd}
    \]
  \item A map $H : D \to B$ with $Hj_0 = ur$, $Hj_1 = \id_B$, and $Hv = up$.
  \end{enumerate}
\end{definition}

\begin{lemma}\label{lemma:tcof-to-asdr}
  Let $u : A \to B$ be a trivial cofibration between fibrant objects of $M^\cof$.
  Then $u$ admits an anodyne strong deformation retraction.
\end{lemma}

\begin{proof}
  Apply \cref{prop:dold} to the diagram
  \[
    \begin{tikzcd}[column sep=tiny]
      & A \ar[ld, "\id_A"'] \ar[rd, "u"] \\
      A \ar[rr, "u"'] & & B
    \end{tikzcd}
  \]
  to obtain a map $r : B \to A$ with $ru = \id_A$ and $ur$ homotopic to $\id_B$ rel $u$.
  Now applying \cref{lemma:homotopy-rel-anodyne} to the latter condition produces an anodyne cylinder object $v : C \to D$ and a map $H : D \to B$ with the required properties.
\end{proof}

\begin{lemma}\label{lemma:asdr-to-acof}
  Let $u : A \to B$ be a cofibration in $M^\cof$ which admits an anodyne strong deformation retraction.
  Then $u$ is an anodyne cofibration.
\end{lemma}

\begin{proof}
  Let $r : B \to A$, $v : C \to D$, $H : D \to B$ be an anodyne strong deformation retraction of $u$.
  Form the pushout of $u$ along $i_0 : A \to C$ as shown below.
  \[
    \begin{tikzcd}
      A \ar[d, "u"'] \ar[r, "i_0"] & C \ar[d, "u'"] \ar[rdd, bend left=10, "v"] \\
      B \ar[r, "i'_0"] \ar[rrd, bend right=10, "j_0"'] & B \amalg_A C \ar[rd, "w"] \\
      & & D
    \end{tikzcd}
  \]
  Because $v$ is an anodyne cylinder object on $u$, the induced map $w : B \amalg_A C \to D$ is an anodyne cofibration.
  Now, form the diagram below.
  \[
    \begin{tikzcd}
      A \ar[r, "u' i_1"] \ar[d, "u"'] & B \amalg_A C \ar[r, "{(r, p)}"] \ar[d, "w"] & A \ar[d, "u"] \\
      B \ar[r, "j_1"] & D \ar[r, "H"] & B
    \end{tikzcd}
  \]
  The map $(r, p) : B \amalg_A C \to A$ is well-defined because $ru = \id_A = pi_0$.
  We claim that this diagram exhibits $u$ as a retract of $w$.
  \begin{enumerate}
  \item The left square commutes because $j_1u = vi_1 = wu'i_1$.
  \item The right square commutes because $Hwi'_0 = Hj_0 = ur$ and $Hwu' = Hv = up$.
  \item The top composition equals $pi_0 = \id_A$.
  \item The bottom composition is $Hj_1 = \id_B$.
  \end{enumerate}
  Therefore $u$ is a retract of the anodyne cofibration $w$, and is itself an anodyne cofibration.
\end{proof}

\begin{proposition}\label{prop:tcof-to-acof}
  Let $u : A \to B$ be a trivial cofibration between fibrant objects of $M^\cof$.
  Then $u$ is an anodyne cofibration.
  Dually, if $u$ is a trivial fibration, then it is an anodyne fibration.
\end{proposition}

\begin{proof}
  Follows from \cref{lemma:tcof-to-asdr} and \cref{lemma:asdr-to-acof}.
\end{proof}

\section{Left weak equivalences revisited}\label{sec:lwe-revisited}

Our original definition of the left weak equivalences of $M$ involved quantification over all left Quillen functors $M \to \Kan^\op$ or, equivalently, over all simplicial resolutions in $M$.
In this section, we will give an ``elementary'' criterion (\cref{prop:left-weq-elt}) to detect left weak equivalences in terms of the anodyne cofibrations of $M$.

\begin{lemma}\label{lemma:cf-cylinder}
  Let $f : A \to B$ be a map between fibrant objects of $M^\cof$.
  Then $f$ admits a factorization as a cofibration followed by an anodyne fibration which has a section which is an anodyne cofibration.
\end{lemma}

\begin{proof}
  Factor $(f, \id_B) : A \amalg B \to B$ into a cofibration $j : A \amalg B \to C$ followed by an anodyne fibration $p : C \to B$.
  Then $C$ is again cofibrant and fibrant, and we can form the diagram below.
  \[
    \begin{tikzcd}
      & B \ar[d, hookrightarrow] \\
      A \ar[r, hookrightarrow] & A \amalg B \ar[r, hookrightarrow, "j"] & C \ar[r, two heads, "p", "\an"'] & B
    \end{tikzcd}
  \]
  The sequence $A \to C \to B$ will be our chosen factorization of $f$.
  The composition $s : B \to A \amalg B \to C$ is a section of $C \to B$.
  Because $B$ and $C$ are both cofibrant and fibrant, $p$ is a left weak equivalence, hence $s$ is also one by the two-out-of-three property, and therefore $s$ is an anodyne cofibration by \cref{prop:tcof-to-acof}.
\end{proof}

\begin{remark}
  If $f : A \to B$ was a left weak equivalence, then the first map of such a factorization will in fact be an \emph{anodyne} cofibration, by the two-out-of-three property and \cref{prop:tcof-to-acof}.
\end{remark}

\begin{proposition}\label{prop:left-weq-elt}
  Let $f : A \to B$ be a morphism of $M^\cof$.
  Then $f$ is a left weak equivalence if and only if there is a diagram
  \[
    \begin{tikzcd}
      A \ar[r, hookrightarrow, "\an"] \ar[d, "f"'] & A' \ar[d, "p"] \\
      B \ar[r, hookrightarrow, "\an"] & B'
    \end{tikzcd}
  \]
  in which the horizontal maps are anodyne cofibrations and the map $p$ has a section which is an anodyne cofibration.
\end{proposition}

\begin{proof}
  Suppose there is a diagram of the above form.
  Then $p$ is a left weak equivalence by the two-out-of-three property, because it has a section which is an anodyne cofibration.
  By the two-out-of-three property again, $f$ is also a left weak equivalence.

  Conversely, suppose $f$ is a left weak equivalence.
  Factor $B \to *$ into an anodyne cofibration $B \to \hat B$ followed by a fibration $\hat B \to *$, and then factor the composition $A \to B \to \hat B$ into an anodyne cofibration $A \to \hat A$ followed by a fibration $\hat f : \hat A \to \hat B$.
  \[
    \begin{tikzcd}
      A \ar[r, hookrightarrow, "\an"] \ar[d, "f"'] & \hat A \ar[d, two heads, "\hat f"] \\
      B \ar[r, hookrightarrow, "\an"] & \hat B \ar[r, two heads] & *
    \end{tikzcd}
  \]
  By the two-out-of-three property, $\hat f$ is a left weak equivalence, and the objects $\hat A$ and $\hat B$ are cofibrant and fibrant.
  Therefore, by \cref{lemma:cf-cylinder} and the following remark, we may factor $\hat f$ as an anodyne cofibration followed by a map $p$ which has a section which is an anodyne cofibration, as shown below.
  \[
    \begin{tikzcd}
      A \ar[r, hookrightarrow, "\an"] \ar[d, "f"'] &
      \hat A \ar[r, hookrightarrow, "\an"] & C \ar[d, "p"] \\
      B \ar[rr, hookrightarrow, "\an"] & & \hat B
    \end{tikzcd}
  \]
  Taking $A' = C$, $B' = \hat B$ we arrive at a diagram of the above form.
\end{proof}

We can now complete the proof of \cref{thm:cof}.

\begin{proof}[Proof of \cref{thm:cof}]
  We proved that the left weak equivalences (as defined in \cref{def:left-weq}) make $M^\cof$ into a cofibration category (\cref{prop:exists-cof-cat}) which is functorial with respect to left Quillen functors (\cref{prop:cc-functorial}), that the left weak equivalences contain the anodyne cofibrations (\cref{prop:acof-is-weq}), and that they agree with the weak equivalences in a model category (\cref{prop:model-left-weq-iff-weq}).

  Conversely, suppose we have some other assignment $\sW'_M$ of weak equivalences in $M^\cof$ to each relaxed premodel category $M$ which also satisfies these conditions.
  Then as observed in \cref{remark:left-weq-motivation}, for any morphism $f : A \to B$ belonging to $\sW'_M$ and any left Quillen functor $F : M \to N$ to a model category $N$, $Ff$ must be a weak equivalence of $N$, and so $\sW'_M$ is contained in the left weak equivalences of $M$ for each $M$.
  Conversely, if $f : A \to B$ is a left weak equivalence in a relaxed premodel category $M$, apply \cref{prop:left-weq-elt} to $f$.
  The anodyne cofibrations in the resulting diagram belong to $\sW'_M$ by hypothesis, so by the two-out-of-three condition, $f$ belongs to $\sW'_M$ as well.
\end{proof}

\begin{proposition}\label{prop:ho-acof}
  The localization of $M^\cof$ at the class of anodyne cofibrations equals $\Ho^L M$ (the localization of $M^\cof$ at the class of left weak equivalences).
\end{proposition}

\begin{proof}
  Anodyne cofibrations are left weak equivalences, so it suffices to verify that each left weak equivalence is inverted in the localization of $M^\cof$ at the anodyne cofibrations, and this follows from \cref{prop:left-weq-elt}.
\end{proof}

If $f$ is known a priori to be a cofibration, we can give a simpler criterion to determine when it is also a left weak equivalence.

\begin{proposition}\label{prop:tcof-iff-comp-acof}
  Let $f : A \to B$ be a cofibration in $M^\cof$.
  Then $f$ is a left weak equivalence if and only if there exists an anodyne cofibration $g : B \to C$ such that $gf$ is also an anodyne cofibration.
\end{proposition}

\begin{proof}
  If such a $g$ exists, then $f$ is a left weak equivalence by the two-out-of-three property.
  Conversely, suppose $f$ is a left weak equivalence.
  Factor $A \to *$ into an anodyne cofibration $A \to \hat A$ followed by a fibration $\hat A \to *$.
  Then form the pushout $B \amalg_A \hat A$, and factor the map $B \amalg_A \hat A \to *$ into an anodyne cofibration $B \amalg_A \hat A \to \hat B$ followed by a fibration $\hat B \to *$.
  \[
    \begin{tikzcd}
      A \ar[r, hookrightarrow, "\an"] \ar[d, hookrightarrow, "f"', "\sim"] & \hat A \ar[d, hookrightarrow] \ar[rd, hookrightarrow] \\
      B \ar[r, hookrightarrow, "\an"'] & B \amalg_A \hat A \ar[r, hookrightarrow, "\an"'] & \hat B
    \end{tikzcd}
  \]
  All the maps displayed are cofibrations, and $\hat A$ and $\hat B$ are fibrant.
  The map $\hat A \to \hat B$ is a left weak equivalence by the two-out-of-three property, and therefore an anodyne cofibration by \cref{prop:tcof-to-acof}.
  We may therefore take $g$ to be the composition $B \to B \amalg_A \hat A \to \hat B$.
\end{proof}

\begin{remark}
  The above statement also appears as \cite[Corollaire~A.2.25]{Mor} where it is used to show that the trivial cofibrations are compatible with the quasi-simplicial structure (our \cref{prop:bifunctor-trivial}).
  As mentioned in the introduction to this chapter, we will use it to relate the left weak equivalences of $M^\cof$ (and thus the homotopy theory of $M$) to the anodyne cofibrations of $M$, which can be detected using lifting conditions.
\end{remark}




\chapter{The algebra of combinatorial premodel categories}
\label{chap:algebra}

Our original motivation for introducing premodel categories was to rectify the lack of limits and colimits in the 2-category of model categories.
In this chapter we show that, if we restrict our attention to \emph{combinatorial} premodel categories, this goal is achieved: the 2-category $\CPM$ admits all limits and colimits.

We also show that any two combinatorial premodel categories $M_1$ and $M_2$ admit a \emph{tensor product} $M_1 \otimes M_2$, for which a left Quillen functor from $M_1 \otimes M_2$ to $N$ is the same as a Quillen bifunctor from $M_1 \times M_2$ to $N$.
Furthermore, this tensor product is left adjoint to an \emph{internal Hom} $\CPM(-, -)$, so that there are equivalences
\[
  \THom_\CPM(M_1 \otimes M_2, N) \eqv \QBF((M_1, M_2), N) \eqv \THom_\CPM(M_1, \CPM(M_2, N)).
\]
The underlying category of $\CPM(M, N)$ is the category $\LPr(M, N)$ of \emph{all} left adjoints from $M$ to $N$, and its cofibrant objects are the left Quillen functors from $M$ to $N$.

Later, we will define a 2-category $V\CPM$ of combinatorial $V$-premodel categories for any monoidal combinatorial premodel category $V$.
Our eventual goal is to (under suitable hypotheses on $V$) equip $V\CPM$ with a model 2-category structure whose weak equivalences are the Quillen equivalences.
In order to produce the required factorizations, we would like to use some version of the small object argument.
One obvious difficulty is that, already in $\CPM$ itself, the Hom categories $\THom_\CPM(M, N) = \LQF(M, N)$ are usually not essentially small.
This means that we cannot adjoin solutions to all possible lifting problems as one usually does at each stage of the small object argument.
We will eventually argue that (under favorable conditions) for sufficiently large $\lambda$, it suffices to adjoin solutions to all lifting problems involving functors that preserve $\lambda$-compact objects.
To lay the groundwork for this argument, in this chapter we begin the investigation of the \emph{rank of combinatoriality} of a combinatorial premodel category.
More precisely, we introduce a filtration of the 2-category $\CPM$ by the sub-2-categories $\CPM_\lambda$ for each regular cardinal $\lambda$, consisting of the $\lambda$-combinatorial premodel categories and the left Quillen functors which preserve $\lambda$-compact objects, and we study the extent to which these sub-2-categories are closed under the aforementioned algebraic structure on $\CPM$.

In order to obtain good control over the internal Hom $\CPM(M, N)$ it turns out that we need to bound not only the ranks of combinatoriality of $M$ and $N$ but also the ``size'' of $M$.
To this end, we introduce the notion of a \emph{$\mu$-small $\lambda$-combinatorial} premodel category.
In this chapter, we show that if $M$ is $\lambda$-small $\lambda$-combinatorial and $N$ is $\lambda$-combinatorial, then $\CPM(M, N)$ is $\lambda$-combinatorial and, moreover, its $\lambda$-compact objects are precisely those left adjoints which preserve $\lambda$-compact objects.
In the next chapter we will show that if $M$ is $\mu$-small $\lambda$-combinatorial then $\THom_{\CPM_\lambda}(M, -) : \CPM_\lambda \to \TCat$ preserves $\mu$-directed colimits; this allows us to carry out the small object argument in $\CPM_\lambda$.

\section{The subcategories $\CPM_\lambda$}

Recall that a premodel category is combinatorial if its underlying category is locally presentable and it admits \emph{sets} $I$ and $J$ of generating cofibrations and anodyne cofibrations respectively (that is, $\sAF = \rlp(I)$ and $\sF = \rlp(J)$).

\begin{definition}
  Let $\lambda$ be a regular cardinal.
  A premodel category is \emph{$\lambda$-combinatorial} if its underlying category is locally $\lambda$-presentable and it admits generating cofibrations and anodyne cofibrations that are morphisms between $\lambda$-compact objects.
\end{definition}

\begin{example}
  The model category $\Kan$ is $\lambda$-combinatorial for all $\lambda$, because its underlying category $\Set^{\Delta^\op}$ is a presheaf category and its standard generating cofibrations and acyclic cofibrations are morphisms between simplicial sets built out of finitely many simplices.
\end{example}

A $\lambda$-combinatorial premodel category is $\lambda'$-combinatorial for all $\lambda' \ge \lambda$.
Any combinatorial premodel category is $\lambda$-combinatorial for some (and hence all sufficiently large) $\lambda$.
We call the smallest $\lambda$ for which a combinatorial premodel category $M$ is $\lambda$-combinatorial the \emph{rank of combinatoriality} of $M$.

\begin{proposition}\label{prop:mk-lambda-comb}
  Let $M$ be a locally $\lambda$-presentable category and let $I$ and $J$ be sets of morphisms between $\lambda$-compact objects of $M$ such that $\rlp(I) \subset \rlp(J)$.
  Then there exists a unique $\lambda$-combinatorial premodel category structure on $M$ with generating cofibrations $I$ and generating anodyne cofibrations $J$.
\end{proposition}

\begin{proof}
  We must take $\sAF = \rlp(I)$, $\sC = \llp(\rlp(I))$, $\sF = \rlp(J)$, $\sAC = \llp(\rlp(J))$.
  The required factorizations are provided by the small object argument.
  The weak factorization systems $(\sC, \sAF)$ and $(\sAC, \sF)$ then make $M$ into a $\lambda$-combinatorial premodel category.
\end{proof}

The simplicity of this statement is part of the reason that combinatorial premodel categories inherit so much algebraic structure of the underlying locally presentable categories.
Left Quillen functors out of a premodel category have a simple description in terms of its generating cofibrations and anodyne cofibrations.
Thus, in order to perform a ``left adjoint'' type construction (colimits, projective premodel structures and tensor products) in $\CPM$, it will suffice to perform the corresponding construction in $\LPr$, write down appropriate generating (anodyne) cofibrations and use the preceding result to produce a combinatorial premodel category with the correct universal property.
For ``right adjoint'' type constructions (limits, injective premodel structures and internal Homs) we will instead have a direct characterization of the desired (anodyne) cofibrations, and then need to show that they are generated by sets in order to obtain a combinatorial premodel category.
Typically we will have no direct characterization of these generating sets (other than ``all cofibrations between $\lambda$-compact objects'' for some sufficiently large $\lambda$) and thus no direct description of the (anodyne) fibrations.

The following basic fact is a simple (and not particularly representative) example of this phenomenon.

\begin{proposition}\label{prop:init-fin-comb}
  Let $M$ be a locally $\lambda$-presentable category.
  Then $M_\init$ and $M_\fin$ are $\lambda$-combinatorial premodel categories.
\end{proposition}

Recall that $M_\init$ (respectively, $M_\fin$) denotes the premodel category structure on $M$ in which every morphism is an anodyne fibration (respectively, an anodyne cofibration) and so only the isomorphisms are cofibrations (respectively, fibrations).
For any premodel category $N$, a left Quillen functor from $M_\init$ to $N$ is the same as a left adjoint from $M$ to $N$, while a left Quillen functor from $N$ to $M_\fin$ is the same as a left adjoint from $N$ to $M$.

\begin{proof}
  The premodel category $M_\init$ has generating (anodyne) cofibrations $I = J = \emptyset$, so it is $\lambda$-combinatorial.
  For $M_\fin$, we must find generating (anodyne) cofibrations which are morphisms between $\lambda$-compact objects.
  We claim that we can take both $I$ and $J$ to be the set of maps of the forms $\emptyset \to A$ and $A \amalg A \to A$ as $A$ runs over a set of representatives of all $\lambda$-compact objects of $M$.
  Indeed, suppose $f : X \to Y$ is a morphism with the right lifting property with respect to all of these maps.
  For any $\lambda$-compact object $A$, the right lifting property
  \[
    \begin{tikzcd}
      \emptyset \ar[r] \ar[d] & X \ar[d, "f"] \\
      A \ar[r] \ar[ru, dotted] & Y
    \end{tikzcd}
  \]
  means that $f_* : \Hom(A, X) \to \Hom(A, Y)$ is surjective, while the right lifting property
  \[
    \begin{tikzcd}
      A \amalg A \ar[r] \ar[d] & X \ar[d, "f"] \\
      A \ar[r] \ar[ru, dotted] & Y
    \end{tikzcd}
  \]
  means that $f_* : \Hom(A, X) \to \Hom(A, Y)$ is injective.
  Thus $f_* : \Hom(A, X) \to \Hom(A, Y)$ is an isomorphism for every $\lambda$-compact object $A$.
  Since the $\lambda$-compact objects of a $\lambda$-presentable category are dense, the functors $\Hom(A, -)$ are in particular jointly conservative and so $f$ is an isomorphism.
\end{proof}

\begin{remark}
  The condition that a premodel category is $\lambda$-combinatorial for some particular $\lambda$ should not be thought of as a ``size'' condition, but as a bound on its ``complexity''.
  For example, the premodel category $\Set^{A^\op}_\init$ is $\aleph_0$-combinatorial for a (small) category $A$ of any cardinality.
  We will discuss a notion of the ``size'' of a combinatorial premodel category in \cref{sec:algebra-size}.
\end{remark}

It turns out to be useful to consider a restricted class of morphisms between $\lambda$-combinatorial premodel categories.
The following result is well-known.

\begin{proposition}\label{prop:strongly-accessible-iff}
  Let $F : C \adj D : G$ be an adjunction with $C$ locally $\lambda$-presentable and $D$ cocomplete.
  Then the following conditions are equivalent:
  \begin{enumerate}
  \item $F$ preserves $\lambda$-compact objects.
  \item $G$ preserves $\lambda$-filtered colimits.
  \end{enumerate}
\end{proposition}


\begin{proof}
  Let $A$ be a $\lambda$-compact object of $C$ and $Y = (Y_i)_{i \in I}$ be a $\lambda$-filtered diagram in $D$; then in the commutative diagram
  \[
    \begin{tikzcd}
      \colim_{i \in I} \Hom_D(FA, Y_i) \ar[r, "\sim"] \ar[dd, "\varphi"'] &
      \colim_{i \in I} \Hom_C(A, GY_i) \ar[d, "\sim"] \\
      & \Hom_C(A, \colim_{i \in I} GY_i) \ar[d, "\psi"] \\
      \Hom_D(FA, \colim_{i \in I} Y_i) \ar[r, "\sim"'] &
      \Hom_C(A, G(\colim_{i \in I} Y_i))
    \end{tikzcd}
  \]
  $\varphi$ is an isomorphism if and only if $\psi$ is.
  Now $F$ preserves $\lambda$-compact objects if and only if $\varphi$ is an isomorphism for every choice of $A$ and $Y$.
  Since the $\lambda$-compact objects of $C$ form a dense subcategory, for fixed $Y$, the map $\psi : \Hom_C(A, \colim_{i \in I} GY_i) \to \Hom_C(A, G(\colim_{i \in I} Y_i))$ is an isomorphism for all choices of $A$ if and only if $\colim_{i \in I} GY_i \to G(\colim_{i \in I} Y_i)$ is an isomorphism; this holds for all choices of $Y$ exactly when $G$ preserves $\lambda$-filtered colimits.
\end{proof}

\begin{definition}
  A left adjoint $F : C \to D$ between locally $\lambda$-presentable categories is called \emph{strongly $\lambda$-accessible} if it preserves $\lambda$-compact objects.
\end{definition}

Although this terminology will sometimes be useful, we will also quite often use the phrase ``preserves $\lambda$-compact objects'' directly.
Condition (2) of \cref{prop:strongly-accessible-iff} makes it clear that a strongly $\lambda$-accessible functor is also strongly $\lambda'$-accessible for any $\lambda' \ge \lambda$.

\begin{proposition}\label{prop:exists-strongly-accessible}
  Every left adjoint $F : C \to D$ between locally presentable categories is strongly $\lambda$-accessible for some (and hence all sufficiently large) $\lambda$.
\end{proposition}

\begin{proof}
  This is proved in the more general setting of accessible functors between accessible categories as \cite[Theorem~2.19]{AR}.
\end{proof}

\begin{definition}
  We write $\LPr_\lambda$ for the sub-2-category of $\LPr$ whose objects are the locally $\lambda$-presentable categories and with $\THom_{\LPr_\lambda}(C, D)$ the full subcategory of $\THom_\LPr(C, D)$ on the strongly $\lambda$-accessible functors.
\end{definition}

\begin{definition}
  $\CPM_\lambda$ is the sub-2-category of $\CPM$ whose objects are the $\lambda$-combinatorial premodel categories and with $\THom_{\CPM_\lambda}(M, N)$ the full subcategory of $\THom_\CPM(M, N)$ on the strongly $\lambda$-accessible functors.
\end{definition}

The sub-2-categories $\LPr_\lambda$ form a filtration $\LPr_{\aleph_0} \subset \LPr_{\aleph_1} \subset \cdots$ of $\LPr$.
Every locally presentable category is locally $\lambda$-presentable for some $\lambda$, so by \cref{prop:exists-strongly-accessible}, every small diagram in $\LPr$ factors through $\LPr_\lambda$ for sufficiently large $\lambda$.
Similarly, the sub-2-categories $\CPM_\lambda$ form a filtration $\CPM_{\aleph_0} \subset \CPM_{\aleph_1} \subset \cdots$ of $\CPM$ and every small diagram in $\CPM$ factors through $\CPM_\lambda$ for sufficiently large $\lambda$.

\section{Coproducts and products}

Coproducts and products are of course a special case of colimits and limits, which will be treated in generality later on.
We treat them specifically in this section for a few reasons.
As simple examples of colimits and limits, they provide a gentle introduction to some of the phenomena which will recur throughout this chapter.
On the other hand, there are also some features unique to coproducts and products which are not shared by general colimits and limits.
Finally, coproducts and products will be one of the ingredients in the subsequent construction of general colimits and limits, so we must treat them separately anyways.

One distinctive feature of coproducts and products in $\CPM$ is that they agree.
Recall from \cref{subsec:premod-products} that, for any family $(M_s)_{s \in S}$ of premodel categories, the product category $M = \prod_{s \in S} M_s$ has a premodel category structure in which all of the classes are defined componentwise.
Each projection $\pi_s : M \to M_s$ is then both a left and right Quillen functor.
The left adjoint of $\pi_s : M \to M_s$ is the functor $\iota_s : M_s \to M$ given by $(\iota_s A)_s = A$ and $(\iota_s A)_{s'} = \emptyset$ for $s' \ne s$, while the right adjoint of $\pi_s$ is the functor $\tau_s : M_s \to M$ given by $(\tau_s A)_s = A$ and $(\tau_s A)_{s'} = *$ for $s' \ne s$.
As we explained in \cref{subsec:premod-products}, the functors $\pi_s : M \to M_s$ make $M$ the (strict) product of the $M_s$ in the 2-category $\PM$ of all premodel categories and left Quillen functors.
The functors $\pi_s : M \to M_s$ also make $M$ the (strict) product of the $M_s$ in the 2-category of premodel categories and \emph{right} Quillen functors, and therefore their left adjoints $\iota_s : M_s \to M$ make $M$ the (non-strict) \emph{coproduct} of the $M_s$ in $\PM$.

\begin{proposition}\label{prop:prod-lambda-comb}
  Let $(M_s)_{s \in S}$ be a family of $\lambda$-combinatorial premodel categories.
  Then $\prod_{s \in S} M_s$ is also $\lambda$-combinatorial.
\end{proposition}

\begin{proof}
  The product of locally $\lambda$-presentable categories is locally $\lambda$-presentable.
  For each $s$, let $I_s$ ($J_s$) be a set of generating (anodyne) cofibrations for $M_s$ made up of morphisms between $\lambda$-compact objects.
  Then one easily verifies that
  \[
    I = \{\,\iota_s f \mid s \in S, f \in I_s\,\}, \quad
    J = \{\,\iota_s f \mid s \in S, f \in J_s\,\}
  \]
  are generating cofibrations and generating anodyne cofibrations for $\prod_{s \in S} M_s$.
  The functors $\iota_s$ are strongly $\lambda$-accessible by \cref{prop:strongly-accessible-iff}, because their right adjoints $\pi_s$ are again left adjoints and therefore preserve \emph{all} colimits.
  Hence $I$ and $J$ are again sets of morphisms between $\lambda$-compact objects.
\end{proof}

In particular, the product of a family $(M_s)_{s \in S}$ of combinatorial premodel categories is again combinatorial.
Since $\CPM$ is a full sub-2-category of $\PM$, we conclude that $\prod_{s \in S} M_s$ is also both a coproduct and product of the $M_s$ in $\CPM$.
We summarize this discussion below.

\begin{proposition}\label{prop:cpm-prod}
  $\CPM$ has arbitrary coproducts and products.
  Both the coproduct and the product of a family $(M_s)_{s \in S}$ are computed as the product $M = \prod_{s \in S} M_s$ of premodel categories.
  The projection functors $\pi_s : M \to M_s$ make $M$ into the product of the $M_s$, while their left adjoints $\iota_s : M_s \to M$ make $M$ into the coproduct of the $M_s$.
\end{proposition}

\begin{proof}
  This is just a restatement of the preceding discussion.
\end{proof}

The special case of the empty family merits its own notation.

\begin{notation}
  We write $0$ for the terminal category equipped with its unique premodel category structure.
\end{notation}

Of course $0$ is actually a model category, and is $\lambda$-combinatorial for every $\lambda$.

\begin{proposition}
  $0$ is a zero object in $\CPM$.
  That is, $0$ is both initial and final.
\end{proposition}

\begin{proof}
  This follows from \cref{prop:cpm-prod} by taking $S = \emptyset$.
  (It is also obvious by inspection, since there is a unique left Quillen functor $F : N \to 0$ for any $N$, while a left Quillen functor $F : 0 \to N$ must send the unique object of $0$ to an initial object of $N$ and any two such functors are uniquely isomorphic.)
\end{proof}

\begin{remark}
  Since coproducts and products of (combinatorial) premodel categories agree, we could call $\prod_{s \in S} M_s$ the \emph{direct sum} of the $M_s$ and denote it by $\bigoplus_{s \in S} M_s$.
  This is analogous to the direct sum of abelian groups (or better, of commutative monoids), with the operations of forming colimits playing the role corresponding to addition.
  The situation is even better for premodel categories because \emph{arbitrary} coproducts and products agree, not just finite ones; this is because we are allowed to form colimits indexed on arbitrary (small) categories.

  As we will see next, the situation is a bit more delicate in $\CPM_\lambda$.
  Thus, in order to avoid confusion, rather than using the notation $\bigoplus_{s \in S} M_s$, we will tend to write $\coprod_{s \in S} M_s$ when we want to emphasize the role of the product premodel category as a coproduct, and $\prod_{s \in S} M_s$ when we think of it as a product or want to do calculations involving its underlying category.
\end{remark}

We now turn to the relationship between coproducts and products and the sub-2-categories $\CPM_\lambda$.
Specifically, we seek to understand under what conditions a coproduct or product in $\CPM$ is also one in $\CPM_\lambda$.
Because $\CPM_\lambda$ imposes restrictions on both objects (they must be $\lambda$-combinatorial) and 1-morphisms (they must be strongly $\lambda$-accessible), but not 2-morphisms, there are in general three reasons why a colimit or limit in $\CPM$ of a diagram in $\CPM_\lambda$ might not also be a colimit or limit in $\CPM_\lambda$.
First, the colimit or limit object itself might not belong to $\CPM_\lambda$.
Second, even if it does belong to $\CPM_\lambda$, the colimit or limit morphisms connecting it to the original diagram might not be 1-morphisms of $\CPM_\lambda$.
Finally, even if the entire colimit or limit diagram belongs to $\CPM_\lambda$, it might fail to have the correct universal property in $\CPM_\lambda$.

Specializing to the case of coproducts and products, let $(M_s)_{s \in S}$ be a family of $\lambda$-combinatorial premodel categories.
By \cref{prop:prod-lambda-comb}, the coproduct and product object $M = \prod_{s \in S} M_s$ is again $\lambda$-combinatorial.
Moreover, the projections $\pi_s : M \to M_s$ and $\iota_s : M_s \to M$ are each strongly $\lambda$-accessible, hence morphisms of $\CPM_\lambda$.
For $\iota_s$, we verified this in the proof of \cref{prop:prod-lambda-comb}.
For $\pi_s$, the formula for its right adjoint $\tau_s$ shows that $\tau_s$ preserves filtered (indeed, connected) colimits and so $\pi_s$ is strongly $\lambda$-accessible by \cref{prop:strongly-accessible-iff}.
It remains to check whether the equivalences
\[
  \THom_\CPM(M, N) \eqv \prod_{s \in S} \THom_\CPM(M_s, N), \quad
  \THom_\CPM(N, M) \eqv \prod_{s \in S} \THom_\CPM(N, M_s)
\]
restrict to equivalences
\[
  \THom_{\CPM_\lambda}(M, N) \eqv \prod_{s \in S} \THom_{\CPM_\lambda}(M_s, N), \quad
  \THom_{\CPM_\lambda}(N, M) \eqv \prod_{s \in S} \THom_{\CPM_\lambda}(N, M_s)
\]
for $N$ an object of $\CPM_\lambda$.

We first consider the equivalence $E : \THom_\CPM(M, N) \eqv \prod_{s \in S} \THom_\CPM(M_s, N)$ which makes $M$ into the coproduct of the $M_s$.
This equivalence is the bottom functor in the commutative square
\[
  \begin{tikzcd}
    \THom_{\CPM_\lambda}(M, N) \ar[r] \ar[d] & \prod_{s \in S} \THom_{\CPM_\lambda}(M_s, N) \ar[d] \\
    \THom_\CPM(M, N) \ar[r, "E", "\eqv"'] & \prod_{s \in S} \THom_\CPM(M_s, N)
  \end{tikzcd}
\]
whose vertical functors are inclusions of full (and replete) subcategories.
The top functor is therefore fully faithful, so it is an equivalence if and only if the inverse image under $E$ of any $S$-tuple of morphisms in $\prod_{s \in S} \THom_{\CPM_\lambda}(M_s, N)$ belongs to $\THom_{\CPM_\lambda}(M, N)$.

The forward direction of the equivalence $E$ sends a left Quillen functor $F : M \to N$ to the family $(F \circ \iota_s)_{s \in S}$.
Now any object $A$ of $M = \prod_{s \in S} M_s$ can be expressed in the form $A = \coprod_{s \in S} \iota_s A_s$; any left Quillen functor $F : M \to N$ preserves this coproduct, so that $FA = \coprod_{s \in S} (F \circ \iota_s) A_s$.
Therefore, the inverse of $E$ sends a family $(F_s : M_s \to N)_{s \in S}$ of left Quillen functors to the left Quillen functor $F : M \to N$ defined by the formula $FA = \coprod_{s \in S} F_s A_s$.

\begin{lemma}\label{lemma:compact-in-prod}
  The $\lambda$-compact objects of the product $M = \prod_{s \in S} M_s$ of locally $\lambda$-presentable categories are the objects of the form $\coprod_{s \in S'} \iota_s A_s$ for $S'$ a $\lambda$-small subset of $S$ and $A_s$ a $\lambda$-compact object of $M_s$ for each $s \in S'$.
\end{lemma}

\begin{proof}
  The functors $\iota_s : M_s \to M$ preserve $\lambda$-compact objects and so do $\lambda$-small coproducts, so every object of this form is $\lambda$-compact.
  Conversely, suppose $A$ is a $\lambda$-compact object of $M$.
  As noted earlier, $A$ can be written as $\coprod_{s \in S} \iota_s A_s$ where the objects $A_s = \pi_s A$ are the components of $A$.
  The coproduct $A = \coprod_{s \in S} \iota_s A_s$ is the $\lambda$-filtered colimit of the objects $\coprod_{s \in S'} \iota_s A_s$ as $S'$ ranges over all $\lambda$-small subsets of $S$.
  Since $A$ is $\lambda$-compact, $A$ must be a retract of $\coprod_{s \in S'} \iota_s A_s$ for some particular $\lambda$-small subset $S' \subset S$.
  Then for any $s \notin S'$, the component $A_s$ is a retract of the initial object of $M_s$, hence initial; so $A$ is actually isomorphic to $\coprod_{s \in S'} \iota_s A_s$.
  Finally, the strongly $\lambda$-accessible left adjoint $\pi_s : M_s \to M$ takes $A$ to $A_s$, so $A_s$ is a $\lambda$-compact object of $M_s$ for each $s \in S'$.
\end{proof}

Suppose now that $(F_s : M_s \to N)_{s \in S}$ is a family of morphisms of $\CPM_\lambda$, that is, of strongly $\lambda$-accessible left Quillen functors.
We claim that the corresponding $F : M \to N$ is always strongly $\lambda$-accessible.
Indeed, by the lemma, if $A$ is any $\lambda$-compact object of $M$, we can express $A$ as $\coprod_{s \in S'} \iota_s A_s$ for a $\lambda$-small subset $S' \subset S$ and $\lambda$-compact objects $A_s$ of $M_s$.
Then $FA = \coprod_{s \in S'} F_s A_s$ is $\lambda$-compact because each $F_s$ preserves $\lambda$-compact objects and so does the coproduct over $S'$.
Hence, we conclude that $\THom_{\CPM_\lambda}(M, N) \eqv \prod_{s \in S} \THom_{\CPM_\lambda}(M_s, N)$ is indeed an equivalence.
We summarize this result as follows.

\begin{proposition}
  $\CPM_\lambda$ is closed under all coproducts in $\CPM$.
  That is, $\CPM_\lambda$ admits all coproducts and the inclusion of $\CPM_\lambda$ in $\CPM$ preserves all coproducts.
\end{proposition}

\begin{proof}
  As we have just verified, given any family $(M_s)_{s \in S}$ of objects of $\CPM_\lambda$, the coproduct of the family in $\CPM$ is also a coproduct in $\CPM_\lambda$.
\end{proof}

We now turn to products.
The object $M = \prod_{s \in S} M_s$ is also the product of the $M_s$ in $\CPM_\lambda$ if and only if the top functor in the diagram
\[
  \begin{tikzcd}
    \THom_{\CPM_\lambda}(N, M) \ar[r] \ar[d] & \prod_{s \in S} \THom_{\CPM_\lambda}(N, M_s) \ar[d] \\
    \THom_\CPM(N, M) \ar[r, "E", "\eqv"'] & \prod_{s \in S} \THom_\CPM(N, M_s)
  \end{tikzcd}
\]
is an equivalence.
In this case the bottom equivalence $E$ is the isomorphism sending $F : N \to M$ to the family $(\pi_s \circ F : N \to M_s)_{s \in S}$ and so its inverse sends a family $(F_s : N \to M_s)_{s \in S}$ to the functor $F : N \to M$ given by $(FB)_s = F_s B$.
Suppose that each $F_s : N \to M_s$ is strongly $\lambda$-accessible.
Then for any $\lambda$-compact object $B$ of $N$, each component $(FB)_s = F_s B$ of $FB$ is $\lambda$-compact.
However, this does not imply that $FB$ is $\lambda$-compact.
Indeed, as we showed in \cref{lemma:compact-in-prod}, only objects for which all but a $\lambda$-small subset of their components are initial objects are $\lambda$-compact in $M$.
Thus $M$ (together with the morphisms $\pi_s : M \to M_s$) is usually \emph{not} a product of the $M_s$ in $\CPM_\lambda$ when $S$ has cardinality $\lambda$ or greater.

However, for $S$ a $\lambda$-small set, this difficulty does not arise and we obtain the following result.

\begin{proposition}
  $\CPM_\lambda$ is closed in $\CPM$ under $\lambda$-small products.
\end{proposition}

\begin{proof}
  By the preceding discussion, it remains to verify that when $S$ has cardinality less than $\lambda$, an object $A$ of $M = \prod_{s \in S} M_s$ is $\lambda$-compact if each of its components is.
  This follows from the formula $A = \coprod_{s \in S} \iota_s A_s$.
\end{proof}

This is the first example of a general phenomenon: $\CPM_\lambda$ is closed in $\CPM$ under all constructions of ``left adjoint'' type, but closed under constructions of ``right adjoint'' type under an additional cardinality assumption of some kind.

\begin{remark}
  Most likely $\CPM_\lambda$ does admit products of families of cardinality $\lambda$ or greater.
  However, as we have seen, these products cannot be preserved by the inclusion of $\CPM_\lambda$ in $\CPM$.
  We will have no need for such products and so we do not pursue this matter here.
\end{remark}

\section{Left- and right-induced premodel structures}

A standard technique for constructing model categories (which dates back to Quillen's original work \cite{Q}) is by transferring an existing model category structure across an adjunction.
Suppose $F : M \adj N : G$ is an adjunction between complete and cocomplete categories and $M$ is already equipped with a model category structure.
We may then attempt to make $N$ into a model category and $F$ into a left Quillen functor in a universal way, as follows.
Define a morphism of $N$ to be a weak equivalence or a fibration if its image under $G$ is one in $M$.
Under certain conditions, this defines a model category structure on $N$.
Then, by definition, $G$ preserves fibrations and acyclic fibrations, so $F : M \to N$ is a left Quillen functor.
Moreover, $F$ has the following universal property: if $H : N \to N'$ is any left adjoint to a model category $N'$, then $H$ is a left Quillen functor if and only if $H \circ F$ is one.
Indeed, if $H \circ F$ is a left Quillen functor, then the right adjoint of $H$ sends (acyclic) fibrations to morphisms of $N$ which are sent by $G$ to (acyclic) fibrations, which are by definition the (acyclic) fibrations of the new model category structure on $N$.
In particular, considering those left Quillen functors $H : N \to N'$ whose underlying functor is the identity, we see that $N$ is equipped with the minimal (in the sense of having the smallest classes cofibrations and acyclic cofibrations) model structure such that $F : M \to N$ is a left Quillen functor.

The conditions needed for the existence of this transferred model category structure are nontrivial.
Let us assume that $M$ is combinatorial and $N$ is locally presentable, so that there is no difficulty in constructing factorizations.
There is still an additional consistency condition required to construct the transferred model category structure: any morphism with the left lifting property with respect to all maps sent by $G$ to fibrations has to be an acyclic cofibration in the new model category structure on $N$, so in particular it must itself be sent by $G$ to a weak equivalence.
In general this condition may not be easy to verify.

In contrast, in the setting of combinatorial premodel categories there is no additional condition required.

\begin{proposition}\label{prop:right-induced}
  Let $F : M \adj N : G$ be an adjunction between a combinatorial premodel category $M$ and a locally presentable category $N$.
  Then there exists a (unique) combinatorial premodel category structure on $N$ in which a morphism is a fibration or an anodyne fibration if and only if its image under $G$ is one.
  Moreover, if $M$ is $\lambda$-combinatorial, $N$ is locally $\lambda$-presentable and $F : M \to N$ is strongly $\lambda$-accessible, then this premodel category structure on $N$ is also $\lambda$-combinatorial.
\end{proposition}

\begin{proof}
  Under the original conditions, we may always choose $\lambda$ large enough so that the additional conditions related to $\lambda$ are also satisfied.
  So assume $M$ is $\lambda$-combinatorial, $N$ is locally $\lambda$-presentable and $F$ preserves $\lambda$-compact objects.
  Choose generating cofibrations $I$ and generating anodyne cofibrations $J$ for $M$ which are sets of morphisms between $\lambda$-compact objects.
  For any morphism $f : X \to Y$ of $N$, $Gf$ has the right lifting property with respect to $I$ (respectively, $J$) if and only if $f$ has the right lifting property with respect to $FI$ (respectively, $FJ$).
  Therefore $FI$ and $FJ$ form generating cofibrations and generating anodyne cofibrations for the required premodel category structure on $N$, and this structure is $\lambda$-combinatorial because $N$ is locally $\lambda$-presentable and $F$ preserves $\lambda$-compact objects.
\end{proof}

\begin{proposition}\label{prop:right-induced-universal}
  In the setting of the previous proposition, $F : M \to N$ is a left Quillen functor and for any left adjoint $H : N \to N'$ to a premodel category $N'$, $H$ is a left Quillen functor if and only if $H \circ F$ is.
\end{proposition}

\begin{proof}
  By definition $G$ preserves fibrations and anodyne fibrations, so $F$ is a left Quillen functor.
  If $H : N \to N'$ is any left adjoint to a premodel category and $H \circ F$ is a left Quillen functor, then the right adjoint of $H$ sends (anodyne) fibrations to morphisms sent by $G$ to (anodyne) fibrations of $M$, which by definition are the (anodyne) fibrations of $N$; hence $H$ is also a left Quillen functor.
\end{proof}

We call a premodel category structure constructed in the above fashion a \emph{right-induced premodel structure}, since its (anodyne) fibrations are induced by the right adjoint $G$.
For us, the relevant feature of a right-induced premodel structure is the universal property described in the previous proposition.
For constructions of ``left adjoint'' type, we need to produce a combinatorial premodel category $N$ for which $\THom_\CPM(N, -)$ is given by some particular prescription.
Our general strategy is to perform the corresponding construction in $\LPr$ to obtain the underlying locally presentable category of $N$, and then to find a premodel category structure on $N$ so that $\THom_\CPM(N, -)$ is the correct full subcategory of $\THom_\LPr(N, -)$.
Often we can describe this subcategory as consisting of those functors whose composition with some functor $F$ from an already-constructed combinatorial premodel category $M$ to $N$ is a left Quillen functor, and in this case the right-induced premodel structure on $N$ then has the required universal property.

There is a dual theory of left-induced model structures \cite{LI}.
In this case, one begins with an adjunction $F : M \adj N : G$ from a locally presentable category $M$ to a model category $N$, and wishes to equip $M$ with a model category structure in which a morphism is a weak equivalence or a cofibration if and only if its image under $F$ is one in $N$.
In place of the easy argument used in \cref{prop:right-induced} this theory relies on a more difficult result of Makkai and Rosick\'y, the statement of which we have slightly augmented below.

\begin{theorem}\label{thm:left-induced-wfs}
  Let $F : M \to N$ be a left adjoint between locally presentable categories and suppose that $N$ is equipped with a cofibrantly generated weak factorization system.
  Then there is a ``left-induced'' cofibrantly generated weak factorization system on $M$ whose left class consists of the morphisms sent by $F$ to the left class of the given weak factorization system on $N$.

  Moreover, suppose that $\lambda$ is an uncountable regular cardinal, $M$ and $N$ are locally $\lambda$-presentable categories, $F$ is strongly $\lambda$-accessible, and the weak factorization system on $N$ is generated by a set of morphisms between $\lambda$-compact objects.
  Then the left class of the left-induced weak factorization system on $M$ is also generated by a set of morphisms between $\lambda$-compact objects.
\end{theorem}

\begin{proof}
  As explained in \cite[Remark~3.8]{MR}, we may obtain the first statement by applying \cite[Theorem~3.2]{MR} to the pseudopullback square
  \[
    \begin{tikzcd}
      M \ar["F", r] \ar["\id"', d] & N \ar["\id", d] \\
      M_\tr \ar["F"', r] & N_\tr
    \end{tikzcd}
  \]
  in which $M_\tr$ and $N_\tr$ are the corresponding categories equipped with the weak factorization system in which every morphism belongs to the left class.
  These weak factorization systems are generated by morphisms between $\lambda$-compact objects when $M$ and $N$ are locally $\lambda$-presentable.
  (We proved this in \cref{prop:init-fin-comb}.)
  For the second statement, we examine the outline of the proof of \cite[Theorem~3.2]{MR}.
  Our conditions on $\lambda$ are precisely the conditions introduced on the regular cardinal $\kappa$ in the second paragraph of this proof.
  The proof proceeds to show that if $\mathcal{X}$ is the collection of all morphisms between $\lambda$-compact objects of $M$ whose images under $\id$ and $F$ belong to the left classes of the weak factorization systems on $M_\tr$ and $N$, then the left class of the induced weak factorization system on $M$ is generated by $\mathcal{X}$.
  Hence, in particular, $\mathcal{X}$ forms a set of morphisms between $\lambda$-compact objects which generates the left-induced weak factorization system on $M$.
\end{proof}

Using this result, we obtain dual versions of \cref{prop:right-induced,prop:right-induced-universal}.

\begin{proposition}\label{prop:left-induced}
  Let $F : M \adj N : G$ be an adjunction between a locally presentable category $M$ and a combinatorial premodel category $N$.
  Then there exists a (unique) combinatorial premodel category structure on $M$ in which a morphism is a cofibration or an anodyne cofibration if and only if its image under $F$ is one.
  Moreover, if $\lambda$ is an uncountable regular cardinal for which $M$ is locally $\lambda$-presentable, $N$ is $\lambda$-combinatorial and $F : M \to N$ is strongly $\lambda$-accessible, then this premodel category structure on $M$ is also $\lambda$-combinatorial.
\end{proposition}

\begin{proof}
  Again, it suffices to prove the second statement, which follows immediately by applying \cref{thm:left-induced-wfs} to the two weak factorization systems which make up the premodel category structure of $N$.
\end{proof}

\begin{proposition}\label{prop:left-induced-universal}
  In the setting of the previous proposition, $F : M \to N$ is a left Quillen functor and for any left adjoint $H : M' \to M$ from a premodel category $M'$, $H$ is a left Quillen functor if and only if $F \circ H$ is.
\end{proposition}

\begin{proof}
  Dual to the proof of \cref{prop:right-induced-universal}.
\end{proof}

We call a premodel category structure constructed in the above fashion a \emph{left-induced premodel structure}.
We will use left-induced premodel structures to perform constructions of ``right adjoint'' type, in which we seek a combinatorial premodel category $M$ for which $\THom_\CPM(-, M)$ is given by some particular prescription, in a manner dual to the one described earlier.

\begin{remark}\label{remark:multi-induced-structures}
  In the setting of right-induced premodel structures, suppose that, rather than a single left adjoint $F : M \to N$ from a combinatorial premodel category $M$ to the locally presentable category $N$, we are given a (small) family $F_s : M_s \to N$ of such left adjoints.
  This family induces a left adjoint $F$ from the coproduct $M = \coprod_{s \in S} M_s$ to $N$.
  Applying \cref{prop:right-induced} to $F : M \to N$, we obtain a right-induced combinatorial premodel category structure on $N$, in which the (anodyne) fibrationts are the morphisms which are sent to (anodyne) fibrations by the right adjoint of every $F_s$.
  Combining \cref{prop:right-induced-universal} with the universal property of the coproduct, we see that a left adjoint $H : N \to N'$ is a left Quillen functor if and only if each composition $M_s \xrightarrow{\iota_s} M \xrightarrow{F} N \xrightarrow{H} N'$ is one.
  In particular, this right-induced premodel structure is the one with the smallest classes of (anodyne) cofibrations which makes each $F_s$ a left Quillen functor.
  Dual comments apply to the case of left-induced premodel category structures, using the product $M = \prod_{s \in S} M_s$.

  For example, the projective and injective premodel category structures on $M^D$ can be produced in this way, using the family of functors $M^D \to M$ given by evaluation at each object of $D$.
  The hard work of constructing generating (anodyne) cofibrations for the injective premodel category structure is contained in the proof of \cref{thm:left-induced-wfs}.
\end{remark}

\section{Conical colimits and limits}

Our next goal is to describe the construction of colimits and limits in $\CPM$ and to determine under what conditions $\CPM_\lambda$ is closed under colimits and limits.

The most general colimit notion in a 2-category is that of a weighted colimit.
Two special kinds of weighted colimits are conical colimits and tensors (by small categories).
A 2-category which has all conical colimits and all tensors by small categories has all colimits.
We have no particular use for general weighted colimits, so we will treat conical colimits and tensors independently.
Dual comments apply to weighted limits, which can be built out of conical limits and cotensors.
We will treat conical colimits and limits in this section and tensors and cotensors in the next section.

We first briefly review the notions of conical limits and colimits.
First, suppose that $X : K \to \C$ is a diagram (i.e., a pseudofunctor) from a (small, ordinary) category $K$ to a 2-category $\C$.
We write $X_k$ for the value of $X$ on an object $k$ of $K$ and $X_f : X_k \to X_l$ for the value of $X$ on a morphism $f : k \to l$.
A \emph{cone} on $X$ with vertex $Y \in \Ob \C$ is a pseudonatural transformation from the constant diagram $K \to \C$ with value $Y$ to $X$.
We can also describe this data as an extension of $X$ to a diagram $X^+ : K^\ltri \to \C$, where $K^\ltri$ denotes the category $K$ with a new initial object $\bot$ adjoined, which is normalized to send the identity of $\bot$ to the identity functor of $Y = X^+(\bot)$.
Concretely, a cone on $X$ with vertex $Y$ consists of \emph{cone 1-morphisms} $F_k : Y \to X_k$ for each $k$, and invertible 2-morphisms between $Y \xrightarrow{F_k} X_k \xrightarrow{X_f} X_l$ and $Y \xrightarrow{f_l} X_l$ for each morphism $f : k \to l$ of $K$, satisfying coherence conditions involving the structural 2-morphisms of $X$.
As a shorthand, we will represent the data of a cone on $X$ with vertex $Y$ by the notation $Y \to (X_k)_{k \in K}$.
Dually, a \emph{cocone} on $X$ amounts to an extension of $X$ to a (normalized) diagram $X^+ : K^\rtri \to \C$, where $K^\rtri$ denotes $K$ with a new terminal object $\top$ adjoined; we represent a cocone on $X$ with vertex $Y$ by the notation $(X_k)_{k \in K} \to Y$.

When $\C = \TCat$, each diagram $X : K \to \TCat$ has a \emph{pseudolimit}, a specific cone $\lim_K X \to (X_k)_{k \in K}$ defined explicitly as follows.
The objects of $\lim_K X$ are the cones on $X$ with vertex the terminal category $*$.
Concretely, such a cone $A$ consists of an object $A_k$ of each category $X_k$ together with specified coherent isomorphisms $\alpha_f : X_f A_k \iso A_l$ for each morphism $f : k \to l$ of $K$.
A morphism $\varphi : A \to B$ of $\lim_K X$ is a family of morphisms $\varphi_k : A_k \to B_k$ which are compatible with the coherence isomorphisms of $A$ and $B$.
For each object $k$ of $K$ there is a functor $\pi_k : \lim_K X \to X_k$ which sends $A$ to $A_k$, and the coherence isomorphism $\alpha_f$ of an object $A$ provides the value on $A$ of an invertible natural transformation between $\lim_K X \to X_k \to X_l$ and $\lim_K X \to X_l$.
This data assembles to form the pseudolimit cone $\lim_K X \to (X_k)_{k \in K}$.
By construction, this pseudolimit cone is the universal cone on $X$: giving a cone on $X$ with vertex $Y$ is the same thing as giving a functor $Y \to \lim_K X$.

A cone $Y \to (X_k)_{k \in K}$ in $\TCat$ is a \emph{limit cone} if the induced functor $Y \to \lim_K X$ is an equivalence of categories.
In this case, for any category $Z$, the induced functor $\TCat(Z, Y) \to \lim_{k \in K} \TCat(Z, X_k)$ is also an equivalence of categories.
More generally, suppose $X : K \to \C$ is a diagram valued in an arbitrary 2-category $\C$.
A cone $Y \to (X_k)_{k \in K}$ is a \emph{limit cone} if the functor $\THom_\C(Z, -)$ sends it to a limit cone for every object $Z$ of $\C$, that is, if the induced functor $\THom_\C(Z, Y) \to \lim_{k \in K} \THom_\C(Z, X_k)$ is an equivalence for every $Z$.
Dually, a cocone $(X_k)_{k \in K} \to Y$ is a \emph{colimit cocone} if the functor $\THom_\C(-, Z)$ sends it to a limit cone for every object $Z$ of $\C$, that is, if $\THom_\C(Y, Z) \to \lim_{k \in K} \THom_\C(X_k, Z)$ is an equivalence for every $Z$.
We will often summarize these situations informally by calling $Y$ the limit or colimit of the diagram $X$, and writing $\THom_\C(Z, \lim_K X) = \lim_{k \in K} \THom_\C(Z, X_k)$ or $\THom_\C(\colim_K X, Z) = \lim_{k \in K} \THom_\C(X_k, Z)$.
The limit or colimit $Y$ of a diagram is described by a universal property which describes the category of morphisms into or out of $Y$ up to \emph{equivalence}; and therefore the limit or colimit is determined only up to equivalence in $\C$.

A stronger condition on a cone $Y \to (X_k)_{k \in K}$ is that the induced functor $\THom_\C(Z, Y) \to \lim_{k \in K} \THom_\C(Z, X_k)$ is an \emph{isomorphism} of categories for every object $Z$ of $\C$.
When $\C = \TCat$, this condition is equivalent to the condition that the induced functor $Y \to \lim_K X$ is an isomorphism, and so in general such a cone is called a \emph{pseudolimit cone}.
(Note that, since we use unadorned terms such as ``limit'' for the weakest, equivalence-invariant notions, the property of being a pseudolimit is \emph{stronger} than the property of being a limit.)
Dually, a cocone $(X_k)_{k \in K} \to Y$ is a \emph{pseudocolimit} if the induced functor $\THom_\C(Y, Z) \to \lim_{k \in K} \THom_\C(X_k, Z)$ is an isomorphism for every object $Z$ of $\C$.

\begin{example}
  Let $\C$ be a 2-category and let $K$ be a discrete category with set of objects $S$.
  Then an $S$-indexed family of objects $(X_s)_{s \in S}$ determines a (strict) functor $X : K \to \C$.
  A cone on $X$ with vertex $Y$ is simply a family of morphisms $(F_s : Y \to X_s)_{s \in S}$.
  (A cone also contains invertible 2-morphisms between $Y \xrightarrow{F_s} X_s \xrightarrow{X_{\id_s}} X_s$ (where $X_{\id_s} = \id_{X_s}$) and $Y \xrightarrow{F_s} X_s$, but the coherence conditions require these 2-morphisms to be identities.)

  When $\C = \TCat$, the pseudolimit $\lim_K X$ is (isomorphic to) the product $\prod_{s \in S} X_s$.
  In general, then, a cone $(F_s : Y \to X_s)_{s \in S}$ on $X : K \to \C$ is a limit cone if the induced functor $\THom_\C(Z, Y) \to \prod_{s \in S} \THom_\C(Z, X_s)$ is an equivalence for every object $Z$ of $\C$.
  We have already used this as the definition of products in a 2-category in \cref{subsec:premod-products}.
  Dually, a cocone $(F_s : X_s \to Y)_{s \in S}$ is a colimit cocone if the induced functor $\THom_\C(Y, Z) \to \prod_{s \in S} \THom_\C(X_s, Z)$ is an equivalence for every object $Z$ of $\C$.
\end{example}

\begin{remark}
  As we have seen, the product $M = \prod_{s \in S} M_s$ of premodel categories (together with the projections $\pi_s : M \to M_s$) makes $M$ a ``pseudoproduct'' of the $M_s$ in $\PM$, but in general $\PM$ does not have ``pseudocoproducts'', only coproducts.
  This asymmetry is a result of our choice to use left adjoints as the morphisms of $\PM$ and discard the right adjoints.
  If instead we took the morphisms of $\PM$ to be Quillen adjunctions, $\PM$ would not have either pseudoproducts or pseudocoproducts.

  We are primarily interested in the equivalence-invariant notions of limits and colimits; for us pseudolimits, when they exist, amount to a convenient way to construct limits.
\end{remark}

We now turn to the problem of constructing colimits and limits in $\CPM$.
As with all the algebraic constructions of this chapter, the plan is to perform the corresponding construction in $\LPr$ and then equip the resulting locally presentable category with a combinatorial premodel category structure which gives it the correct universal property in $\CPM$.
Hence, we first need to understand colimits and limits in $\LPr$.
These were constructed by Bird \cite{Bird}.

\begin{proposition}
  $\LPr$ admits all colimits, and for any regular cardinal $\lambda$, the sub-2-category $\LPr_\lambda$ is closed under all colimits.
\end{proposition}

\begin{proof}
  In the notation of \cite{Bird}, $\mathit{Ladj}$ is our $\LPr$ and $\lambda\mhyphen\mathit{Ladj}$ our $\LPr_\lambda$.
  $\mathit{Loc}$ denotes the 2-category of locally presentable categories, \emph{right} adjoints and natural transformations, and $\lambda\mhyphen\mathit{Loc}$ its sub-2-category of locally $\lambda$-presentable categories and right adjoints which preserve $\lambda$-filtered colimits.
  The 2-categories $\mathit{Loc}$ and $\lambda\mhyphen\mathit{Loc}$ are biequivalent to $\mathit{Ladj}^{\co\op}$ and $\lambda\mhyphen\mathit{Ladj}^{\co\op}$ respectively, by \cref{prop:strongly-accessible-iff}.

  By \cite[Theorem~2.17 and Theorem~2.18]{Bird}, $\lambda\mhyphen\mathit{Loc}$ and $\mathit{Loc}$ both admit all ``limits of retract type'' and the inclusions $\lambda\mhyphen\mathit{Loc} \to \TCat$ and $\mathit{Loc} \to \TCat$ preserve them.
  Pseudolimits are of retract type \cite[Proposition~4.2]{FlexLim} and so in particular $\lambda\mhyphen\mathit{Loc}$ and $\mathit{Loc}$ are complete and the inclusion $\lambda\mhyphen\mathit{Loc} \to \mathit{Loc}$ preserves preserves limits.
  The claim follows because limits in $\lambda\mhyphen\mathit{Loc}$ and $\mathit{Loc}$ are colimits in $\mathit{Ladj}$ and $\lambda\mhyphen\mathit{Ladj}$.
\end{proof}

Later in this chapter we will give a second proof of this fact which gives an alternative description of colimits in $\LPr$.

\begin{proposition}\label{prop:lpr-limits}
  $\LPr$ is closed in $\TCat$ under all limits.
  For an uncountable regular cardinal $\lambda$, the sub-2-category $\LPr_\lambda$ is closed in $\LPr$ under $\lambda$-small limits.
\end{proposition}

\begin{proof}
  With notation as above, by \cite[Proposition~3.14 and Theorem~3.15]{Bird} $\lambda\mhyphen\mathit{Ladj}$ admits $\lambda$-small limits of retract type and the inclusion $\lambda\mhyphen\mathit{Ladj} \to \TCat$ preserves them, and $\mathit{Ladj}$ admits all limits of retract type and the inclusion $\mathit{Ladj} \to \TCat$ preserves them.
\end{proof}

Armed with these facts, we can now construct colimits and limits in $\CPM$.

\begin{proposition}\label{prop:cpm-colimits}
  $\CPM$ admits colimits and $\CPM_\lambda$ is closed under all colimits.
  These colimits are preserved by the forgetful functors $\CPM \to \LPr$ and $\CPM_\lambda \to \LPr_\lambda$.
\end{proposition}

\begin{proof}
  Let $M : K \to \CPM_\lambda$ be a diagram of $\lambda$-combinatorial premodel categories and strongly $\lambda$-accessible left Quillen functors.
  Construct the colimit $Y$ of the diagram $K \xrightarrow{M} \CPM_\lambda \to \LPr_\lambda$, which is also the colimit in $\LPr$.
  Then for any locally presentable category $N$, giving a cocone of left adjoints on $M$ with vertex $N$ is the same as giving a left adjoint from $Y$ to $N$.
  That is, there is an induced equivalence
  \[
    \THom_\LPr(Y, N) \eqv \lim_{k \in K} \THom_\LPr(M_k, N).
  \]

  We now carry out the plan outlined in \cref{remark:multi-induced-structures}.
  Write $\coprod M$ for the premodel category $\coprod_{k \in \Ob K} M_k$.
  By \cref{prop:prod-lambda-comb}, $\coprod M$ is $\lambda$-combinatorial.
  The object $Y$ of $\LPr_\lambda$ is equipped with a cone morphism from each $M_k$, which jointly induce a functor $F : \coprod M \to Y$ in $\LPr_\lambda$ (because the underlying category of $\coprod M$ is also the coproduct of the $M_k$ in $\LPr_\lambda$).
  Apply \cref{prop:right-induced} to the functor $F$ to obtain a $\lambda$-combinatorial right-induced premodel category structure on $Y$.
  In this premodel category structure, a map of $Y$ is a fibration or anodyne fibration if and only if its image under the right adjoint of each colimit cone morphism $M_k \to Y$ is one.
  We must show that $Y$ is the colimit of the diagram $M$ in $\CPM_\lambda$ and also in $\CPM$.
  We know already that the underlying category of $Y$ is the colimit of the diagram $M$ in $\LPr_\lambda$ and also in $\LPr$.
  The Hom categories of $\CPM_\lambda$ are full subcategories of the Hom categories of $\LPr_\lambda$.
  Thus for any object $N$ of $\LPr_\lambda$, there is a diagram
  \[
    \begin{tikzcd}
      \THom_{\CPM_\lambda}(Y, N) \ar[r] \ar[d] &
      \lim\limits_{k \in K} \THom_{\CPM_\lambda}(M_k, N) \ar[d] \\
      \THom_{\LPr_\lambda}(Y, N) \ar[r, "\eqv"'] &
      \lim\limits_{k \in K} \THom_{\LPr_\lambda}(M_k, N)
    \end{tikzcd}
  \]
  in which the bottom morphism is an equivalence and the vertical morphisms are fully faithful.
  Hence, it remains only to show that a left adjoint $Y \to N$ is a left Quillen functor if each composition $M_k \to Y \to N$ is one.
  But by \cref{prop:right-induced-universal}, a left adjoint $Y \to N$ is a left Quillen functor if and only if the composition $\coprod M \to Y \to N$ is one, and by the universal property of $\coprod M$, the latter holds if and only if each $M_k \to Y \to N$ is one.
  The same argument shows that $Y$ is also the colimit in $\CPM$.
\end{proof}

\begin{remark}\label{remark:colimit-gen}
  By examining the proofs of \cref{prop:prod-lambda-comb,prop:right-induced-universal}, we see that $\colim M$ has generating (anodyne) cofibrations given by the images of the generating (anodyne) cofibrations of the $M_k$ under the colimit cone morphisms $M_k \to \colim M$.
\end{remark}

\begin{proposition}
  $\CPM$ admits limits and if $\lambda$ is an uncountable regular cardinal, then $\CPM_\lambda$ is closed under $\lambda$-small limits in $\CPM$.
  These limits are preserved by the forgetful functors $\CPM \to \LPr$ and (for $\lambda$-small limits with $\lambda$ uncountable) $\CPM_\lambda \to \LPr_\lambda$.
\end{proposition}

\begin{proof}
  We imitate the proof of \cref{prop:cpm-colimits}, using \cref{prop:left-induced} in place of \cref{prop:right-induced} and \cref{prop:left-induced-universal} in place of \cref{prop:right-induced-universal}.
  Suppose $M : K \to \CPM_\lambda$ is a $\lambda$-small diagram and let $Y$ be the limit of the diagram $K \xrightarrow{M} \CPM_\lambda \to \LPr_\lambda$, which is also the limit in $\LPr$.
  Write $\prod M$ for the product $\prod_{k \in \Ob K} M_k$ and $F : Y \to \prod M$ for the left adjoint induced by the limit cone morphisms $Y \to M_k$.
  Since $K$ has fewer than $\lambda$ objects, $\prod M$ is also the product of the $M_k$ in $\LPr_\lambda$, so $F$ is strongly $\lambda$-accessible.
  Then the conditions of \cref{prop:right-induced} are satisfied, so the left-induced structure on $Y$ exists and is $\lambda$-combinatorial.
  The rest of the argument is the same as the proof of \cref{prop:cpm-colimits}.
  (Specifically, this argument shows that $\CPM$ admits $\lambda$-small limits for any uncountable regular cardinal $\lambda$, and therefore admits all small limits.)
\end{proof}

\begin{remark}
  The (anodyne) cofibrations in $\lim M$ are the morphisms which are sent by every limit cone morphism $\lim M \to M_k$ to a cofibration (or anodyne cofibration).
  In general, one cannot expect to give an explicit description of the (anodyne) fibrations of $\lim M$.
\end{remark}

\section{Tensors and cotensors by categories}

In this section, we finish the construction of weighted colimits and limits by treating the cases of tensors and cotensors by small categories.
In fact, we showed already in \cref{subsec:premod-diagram} that when $M$ is combinatorial, the projective premodel category $M^{K^\op}_\proj$ and the injective premodel category $M^K_\inj$ exist and are again combinatorial, and are equipped with equivalences
\[
  \THom_\CPM(M^{K^\op}_\proj, N) \eqv \TCat(K, \THom_\CPM(M, N))
  \eqv \THom_\CPM(N, M^{K}_\inj),
\]
making them the tensor and cotensor of $M$ by $K$ in $\CPM$ respectively.
Our only remaining task in this section is to determine when these objects are also tensors and cotensors in $\CPM_\lambda$.

\begin{proposition}
  Let $M$ be a $\lambda$-combinatorial premodel category and $K$ a small category.
  Then $M^{K^\op}_\proj$ is $\lambda$-combinatorial and it is the tensor of $M$ by $K$ in $\CPM_\lambda$.
  In other words, $\CPM_\lambda$ is closed in $\CPM$ under all tensors.
\end{proposition}

\begin{proof}
  The category $M^{K^\op}$ is locally $\lambda$-presentable.
  We already know it is the tensor of $M$ by $K$ in $\LPr$, so for any locally $\lambda$-presentable category $N$, there is an equivalence
  \begin{equation*}
    \THom_\LPr(M^{K^\op}, N) \eqv \TCat(K, \THom_\LPr(M, N)).
    \tag{$*$}
  \end{equation*}
  After passing from left to right adjoints, this equivalence is given simply by sending a right adjoint $N \to M^{K^\op}$ to the corresponding $K^\op$-indexed diagram of right adjoints $N \to M$.
  By \cref{prop:strongly-accessible-iff}, a right adjoint is the right adjoint of a strongly $\lambda$-accessible functor if and only if it preserves $\lambda$-filtered colimits.
  Since colimits in $M^{K^\op}$ are computed componentwise it follows that the equivalence $(*)$ restricts to an equivalence
  \[
    \THom_{\LPr_\lambda}(M^{K^\op}, N) \eqv \TCat(K, \THom_{\LPr_\lambda}(M, N))
  \]
  making $M^{K^\op}$ also the tensor of $M$ by $K$ in $\LPr_\lambda$.

  For each object $k$ of $K$, the evaluation $\pi_k : M^{K^\op} \to M$ has a left adjoint $\iota_k : M \to M^{K^\op}$, which is strongly $\lambda$-accessible because $\pi_k$ also has a right adjoint.
  By the definition of the projective premodel category structure, a left adjoint $M^{K^\op}_\proj \to N$ is a left Quillen functor if and only if its composition with each $\iota_k$ is one.
  Therefore $M^{K^\op}_\proj$ is the right-induced premodel category structure produced by \cref{prop:right-induced} from the functor $F : \coprod_{k \in \Ob K} M \to M^{K^\op}$ induced by each $\iota_k$.
  Now $\coprod_{k \in \Ob K} M$ is $\lambda$-combinatorial and $F$ is strongly $\lambda$-accessible, so $M^{K^\op}_\proj$ is also $\lambda$-combinatorial.
  It has the correct universal property to be the tensor of $M$ by $K$ in $\CPM_\lambda$ because the same is true in $\CPM$ and in $\LPr_\lambda$.
\end{proof}

The corresponding statement for injective premodel structures requires a cardinality assumption on $K$ so that we can control the $\lambda$-compact objects of $M^K$.

\begin{proposition}\label{prop:compact-diagram}
  Let $M$ be a locally $\lambda$-presentable category and let $K$ be a $\lambda$-small category.
  Then an object $X$ of $M^K$ is $\lambda$-compact if and only if each component of $X$ is $\lambda$-compact in $M$.
  In particular, the projection functors $\pi_k : M^K \to M$ are strongly $\lambda$-accessible.
\end{proposition}

\begin{proof}
  Each evaluation functor $\pi_k : M^K \to M$ has a right adjoint $\tau_k : M \to M^K$ given by $(\tau_k X)_l = X^{\Hom_K(k, l)}$.
  This functor $\tau_k$ preserves $\lambda$-filtered colimits because $K$ is $\lambda$-small and therefore so is the set $\Hom_K(k, l)$.
  Thus, $\pi_k$ is strongly $\lambda$-accessible by \cref{prop:strongly-accessible-iff}.
  Conversely, suppose every component of the object $X$ of $M^K$ is $\lambda$-compact.
  For an object $Y$ of $M^K$, $\Hom_{M^K}(X, Y)$ can be computed as the end or equalizer
  \[
    \Hom_{M^K}(X, Y) = \int_{k \in K} \Hom_M(X_k, Y_k)
    = \lim \big[ \prod_{k \in \Ob K} \Hom_M(X_k, Y_k) \rightrightarrows
    \prod_{k \to l} \Hom_M(X_k, Y_l) \big]
  \]
  and this formula commutes with $\lambda$-filtered colimits in $Y$ because the objects $X_k$ are $\lambda$-compact and the limit is over a $\lambda$-small diagram.
\end{proof}

\begin{proposition}
  Let $\lambda$ be an uncountable regular cardinal, $M$ a $\lambda$-combinatorial premodel category and $K$ a $\lambda$-small category.
  Then $M^K_\inj$ is $\lambda$-combinatorial and it is the cotensor of $M$ by $K$ in $\CPM_\lambda$.
  In other words, $\CPM_\lambda$ is closed in $\CPM$ under cotensors by $\lambda$-small categories.
\end{proposition}

\begin{proof}
  The category $M^K$ is locally $\lambda$-presentable and it is the cotensor of $M$ by $K$ in $\LPr$, so for any locally $\lambda$-presentable category $N$ there is an equivalence
  \begin{equation*}
    \THom_\LPr(N, M^K) \eqv \TCat(K, \THom_\LPr(N, M))
  \end{equation*}
  sending a left adjoint $N \to M^K$ to the corresponding $K$-indexed diagram of left adjoints $N \to M$.

  To show that $M^K$ is also the cotensor of $M$ by $K$ in $\LPr_\lambda$, we must show that a left adjoint $N \to M^K$ is strongly $\lambda$-accessible if and only if each component $N \to M^K \to M$ is; this follows from \cref{prop:compact-diagram}.

  It remains to show that $M^K_\inj$ is $\lambda$-combinatorial.
  This follows by applying \cref{prop:left-induced} to the functor $M^K \to \prod_{k \in \Ob K} M$, which is a strongly $\lambda$-accessible functor to a $\lambda$-combinatorial premodel category because (as $K$ is $\lambda$-small) $\prod_{k \in \Ob K} M$ is also the product of $\Ob K$ copies of $M$ in $\LPr_\lambda$ and in $\CPM_\lambda$.
\end{proof}

\section{Orthogonality classes}

Before continuing on to the tensor product and internal Hom of combinatorial premodel categories, we need some more background on locally presentable categories.

\begin{definition}
  Let $C$ be a category and $e : A \to B$ a morphism of $C$.
  A morphism $f : X \to Y$ is \emph{orthogonal} to $e$ if every square
  \[
    \begin{tikzcd}
      A \ar[r] \ar[d, "e"'] & X \ar[d, "f"] \\
      B \ar[r] \ar[ru, dotted, "\exists!"] & Y
    \end{tikzcd}
  \]
  admits a \emph{unique} lift as shown by the dotted arrow.
  An object $X$ is orthogonal to $e$ if the morphism $X \to *$ is orthogonal to $e$.

  If $E$ is a class of morphisms of $C$, we say that a morphism $f : X \to Y$ or an object $X$ is orthogonal to $E$ if $f$ or $X$ is orthogonal to every member of $e$.
\end{definition}

\begin{lemma}\label{lemma:orth-all}
  Suppose that every object $X$ of $C$ is orthogonal to $e : A \to B$.
  Then $e$ is an isomorphism.
\end{lemma}

\begin{proof}
  Since $A$ is orthogonal to $e$, there exists $l : B \to A$ with $le = \id_A$.
  Then $ele = e$ so we can form the diagram below.
  \[
    \begin{tikzcd}
      A \ar[r, "e"] \ar[d, "e"'] & B \ar[d] \\
      B \ar[r] \ar[ru, "el"'] & *
    \end{tikzcd}
  \]
  Since $B$ is also orthogonal to $e$, $el$ must equal $\id_B$ and so $l$ is an inverse for $e$.
\end{proof}

Suppose now that $C$ is a locally presentable category.

\begin{definition}
  For a morphism $e : A \to B$, let $\nabla_e : B \amalg_A B \to B$ be the morphism induced by the identity on each copy of $B$.
  For a morphism $f : X \to Y$, let $\Delta_f : X \to X \times_Y X$ be the morphism induced by the identity on each copy of $X$.
\end{definition}

\begin{proposition}\label{prop:orthogonal-iff}
  Let $e : A \to B$ and $f : X \to Y$ be morphisms of $C$.
  Then the following are equivalent.
  \begin{enumerate}
  \item $f$ is orthogonal to $e$.
  \item $f$ has the right lifting property with respect to both $e$ and $\nabla_e$.
  \item $e$ has the left lifting property with respect to both $f$ and $\Delta_f$.
  \item $f$ is orthogonal to both $e$ and $\nabla_e$.
  \item Both $f$ and $\Delta_f$ are orthogonal to $e$.
  \end{enumerate}
\end{proposition}
  
\begin{proof}
  The existence of lifts in the square
  \[
    \begin{tikzcd}
      B \amalg_A B \ar[r] \ar[d, "\nabla_e"'] & X \ar[d, "f"] \\
      B \ar[r] \ar[ru, dotted] & Y
    \end{tikzcd}
    \quad\mbox{or}\quad
    \begin{tikzcd}
      A \ar[r] \ar[d, "e"'] & X \ar[d, "\Delta_f"] \\
      B \ar[r] \ar[ru, dotted] & X \times_Y X
    \end{tikzcd}
  \]
  is equivalent to the uniqueness of lifts in the square
  \[
    \begin{tikzcd}
      A \ar[r] \ar[d, "e"'] & X \ar[d, "f"] \\
      B \ar[r] \ar[ru, dotted] & Y
    \end{tikzcd}
  \]
  and therefore (1), (2), and (3) are equivalent.
  To show that (4) is equivalent to (1), it suffices to show that if $f$ is orthogonal to $e$ then $f$ is orthogonal to $\nabla_e$.
  We already know that the square
  \[
    \begin{tikzcd}
      B \amalg_A B \ar[r] \ar[d, "\nabla_e"'] & X \ar[d, "f"] \\
      B \ar[r] \ar[ru, dotted] & Y
    \end{tikzcd}
  \]
  admits lifts, by (2).
  These lifts are also unique because $B \amalg_A B \to B$ is an epimorphism.
  The equivalence of (1) and (5) is dual.
\end{proof}

\begin{definition}
  If $E$ is a class of morphisms of $C$, we define $E^+ = E \cup \{\,\nabla_e \mid e \in E\,\}$.
\end{definition}

If $E$ is a set of morphisms between $\lambda$-compact objects of $C$, then so is $E^+$.
A morphism $f : X \to Y$ is orthogonal to $E$ if and only if it has the right lifting property with respect to $E^+$.
Hence we can use the tools of weak factorization systems, such as the small object argument and its variants, in order to study orthogonality.
If $E$ is any \emph{set} of morphisms of $C$, then $E^+$ generates a weak factorization system $(\llp(\rlp(E^+)), \rlp(E^+))$ whose right class is the class of morphisms orthogonal to $E$.
This turns out to be an \emph{orthogonal factorization system}: every morphism of the right class is orthogonal to every morphism of the left class, and hence the factorizations are unique up to isomorphism.

\begin{proposition}\label{prop:ortho-fs}
  If $e \in \llp(\rlp(E^+))$ and $f \in \rlp(E^+)$, then $f$ is orthogonal to $e$.
\end{proposition}

\begin{proof}
  By definition $e$ has the left lifting property with respect to $f$, so by the equivalence of (a) and (c) in \cref{prop:orthogonal-iff} it suffices to show that $e$ also has the left lifting property with respect to $\Delta_f$.
  But $f \in \rlp(E^+)$ means $f$ is orthogonal to $E$, and so by the equivalence of (a) and (e) in \cref{prop:orthogonal-iff}, $\Delta_f$ is also orthogonal to $E$ and so $\Delta_f \in \rlp(E^+)$.
\end{proof}

\begin{proposition}
  The factorizations in this weak factorization system are unique (up to isomorphism).
\end{proposition}

\begin{proof}
  Suppose $X \to Y \to Z$ and $X \to Y' \to Z$ are two factorizations of the same morphism.
  Then $Y' \to Z$ is orthogonal to $X \to Y$ and $Y \to Z$ is orthogonal to $X \to Y'$, so there are lifts in both directions between $Y$ and $Y'$ as shown below.
  \[
    \begin{tikzcd}
      X \ar[r] \ar[d] & Y' \ar[d] \ar[dl, shift left, "g"] \\
      Y \ar[r] \ar[ru, shift left, "f"] & Z
    \end{tikzcd}
  \]
  Also $Y \to Z$ is orthogonal to $X \to Y$, so as the composition $gf : Y \to Y$ is a lift in the square below, it must equal $\id_Y$.
  \[
    \begin{tikzcd}
      X \ar[r] \ar[d] & Y \ar[d] \\
      Y \ar[r] \ar[ru, shift left, "gf"] \ar[ru, shift right, "\id_Y"'] & Z
    \end{tikzcd}
  \]
  Similarly $fg = \id_{Y'}$ and so the factorizations $X \to Y \to Z$ and $X \to Y' \to Z$ are isomorphic.
\end{proof}

\begin{proposition}
  The full subcategory $E^\perp$ of $C$ consisting of the objects orthogonal to $E$ is reflective.
  Writing $R : C \to C$ for the associated monad, the unit $X \to RX$ is an $E^+$-cell complex.
\end{proposition}

\begin{proof}
  Let us instead write $X \to RX \to *$ for the functorial factorization obtained by applying the small object argument to $E^+$.
  Then $X \to RX$ is an $E^+$-cell complex and $RX \to *$ has the right lifting property with respect to $E^+$, so $RX$ belongs to $E^\perp$.
  We need to show that $\Hom(RX, Y) = \Hom(X, Y)$ for any $Y$ belonging to $E^\perp$.
  This is equivalent to the existence of unique lifts in the square
  \[
    \begin{tikzcd}
      X \ar[r] \ar[d] & Y \ar[d] \\
      RX \ar[r] \ar[ru, dotted] & *
    \end{tikzcd}
  \]
  and this follows from \cref{prop:ortho-fs}.
\end{proof}

We write $L : C \to E^\perp$ for the factorization of $R : C \to C$ through the inclusion of $E^\perp$ in $C$.
Since $\Hom(LX, Y) = \Hom(X, Y)$ for any object $Y$ of $E^\perp$, the functor $L$ takes each morphism $X \to RX$ to an isomorphism of $E^\perp$.

By \cite[Theorem~1.39]{AR}, $E^\perp$ is a locally presentable category.
More specifically, if $C$ is locally $\lambda$-presentable and $E$ consists of morphisms between $\lambda$-compact objects, then $E^\perp$ again consists of morphisms between $\lambda$-compact objects and $L : C \to E^\perp$ is strongly $\lambda$-accessible (this amounts to the fact that $R : C \to C$, as constructed by the small object argument, preserves $\lambda$-filtered colimits).

\begin{proposition}\label{prop:out-of-perp}
  Let $D$ be a cocomplete category.
  Then precomposition with $L : C \to E^\perp$ defines an equivalence of categories between left adjoints from $E^\perp$ to $D$ and the full subcategory of left adjoints from $C$ to $D$ which send the morphisms of $E$ to isomorphisms.
  Moreover, if $C$ is locally $\lambda$-presentable and $E$ consists of morphisms between $\lambda$-compact objects of $C$, then the same also holds if we restrict to strongly $\lambda$-accessible left adjoints on both sides.
\end{proposition}

\begin{proof}
  Passing to right adjoints, we must show that composition with the inclusion $E^\perp \to C$ induces an equivalence of categories between right adjoints from $D$ to $E^\perp$ and right adjoints from $D$ to $C$ whose left adjoint sends each morphism of $C$ to an isomorphism.
  In fact, we claim that for an adjunction $F : C \adj D : G$, $G$ factors through $E^\perp$ if and only if $F$ sends the morphisms of $E$ to isomorphisms.
  Indeed, the image $GY$ of an object $Y$ belongs to $E^\perp$ if and only if $GY \to * = G(*)$ is orthogonal to each $e$ in $E$; in turn this holds if and only if $Y \to *$ is orthogonal to $Fe$ for each $e$ in $E$.
  By \cref{lemma:orth-all} this is equivalent to $F$ sending each morphism of $E$ to an isomorphism.
  Furthermore, when $G$ does factor through $E^\perp \subset C$, the induced functor $G' : D \to E^\perp$ is again a right adjoint, with left adjoint $F' : E^\perp \to C \xrightarrow{F} D$, since $\Hom_D(F'A, Y) = \Hom_C(A, GY) = \Hom_{E^\perp}(A, G'Y)$.

  Now suppose $C$ is locally $\lambda$-presentable and $E$ consists of morphisms between $\lambda$-compact objects.
  Then $L : C \to E^\perp$ is strongly $\lambda$-accessible and so $E^\perp$ is closed in $C$ under $\lambda$-filtered colimits.
  It follows that a right adjoint $G : D \to E^\perp$ preserves $\lambda$-filtered colimits if and only if its composition with the inclusion of $E^\perp$ in $C$ does so.
\end{proof}

In other words, $L : C \to E^\perp$ is the universal (in a 2-categorical sense) left adjoint which inverts the morphisms of $E$.

We are mainly interested in the case that $C$ is a presheaf category.
Let $A$ be a small category and recall that we write $\PP(A)$ for the category of presheaves of sets on $A$.

\begin{notation}
  If $E$ is a set of morphisms of $\PP(A)$, we write $\OO(A, E)$ for the full subcategory of $\PP(A)$ on the objects orthogonal to $E$.
  We write $L : \PP(A) \to \OO(A, E)$ for the left adjoint to the inclusion of $\OO(A, E)$ in $\PP(A)$, and $R : \PP(A) \to \PP(A)$ for the associated monad.
\end{notation}

We call $\OO(A, E)$ an \emph{orthogonality class}, and a \emph{$\lambda$-orthogonality class} when $E$ consists of morphisms between $\lambda$-compact objects; in this case $\OO(A, E)$ is a locally $\lambda$-presentable category.
By \cite[Theorem~1.46]{AR}, every locally $\lambda$-presentable category $C$ can be expressed in this form.
In fact, we can take $A$ to be the full subcategory on (a set of representatives for) the $\lambda$-compact objects of $C$, and then choose $E$ such that the inclusion of $A$ in $C$ induces $L : \PP(A) \to \OO(A, E) \eqv C$.

The utility of this description of locally presentable categories is that the form $\OO(A, E)$ is really a ``presentation'' of a locally presentable category, in that it is easy to describe the left adjoints out of a category of this form.
Recall that $\PP(A)$ is the free cocompletion of $A$, so that giving a left adjoint out of $\PP(A)$ is the same as giving an ordinary functor out of $A$.
More precisely, for any cocomplete category $N$, composition with the Yoneda embedding $\yo : A \to \PP(A)$ induces an equivalence between the category of left adjoints from $\PP(A)$ to $N$ and the category of all functors from $A$ to $N$.
Then, by \cref{prop:out-of-perp}, left adjoints from $\OO(A, E)$ to $N$ are the left adjoints from $\PP(A)$ to $N$ which send the morphisms of $E$ to isomorphisms.

This description of left adjoints out of $\OO(A, E)$ makes it a useful form for studying constructions such as colimits.

\begin{proposition}\label{prop:o-form}
  Any diagram $C : K \to \LPr_\lambda$ is equivalent to one of the form $C_k = \OO(A_k, E_k)$, where $A : K \to \TCat$ is a diagram of small categories and $E_k$ is a set of morphisms between $\lambda$-compact objects of $\PP(A_k)$, and the functors $\OO(A_k, E_k) \to \OO(A_l, E_l)$ are the ones induced by $\PP(A_k) \to \PP(A_l)$ given by the left Kan extension of $A_k \to A_l$.

  \[
    \begin{tikzcd}
      A_k \ar[r] \ar[d] & A_l \ar[d] \\
      \PP(A_k) \ar[r] \ar[d] & \PP(A_l) \ar[d] \\
      \OO(A_k, E_k) \ar[r] \ar[d, "\eqv"'] & \OO(A_l, E_l) \ar[d, "\eqv"] \\
      C_k \ar[r] & C_l
    \end{tikzcd}
  \]
\end{proposition}

\begin{proof}
  Take $A_k$ to be the full subcategory of $C_k$ on a set of representatives for the $\lambda$-compact objects of $C_k$.
  Since the functors $C_k \to C_l$ preserve $\lambda$-compact objects, they induce functors $A_k \to A_l$.
  Then choose $E_k$ to be a set of morphisms between $\lambda$-compact objects of $\PP(A_k)$ such that the induced functors $\PP(A_k) \to C_k$ factor through an equivalence $\PP(A_k) \to \OO(A_k, E_k) \xrightarrow{\eqv} C_k$.
\end{proof}

\begin{proposition}\label{prop:colimit-of-o-form}
  Let $(\OO(A_k, E_k))_{k \in K}$ be a diagram in $\LPr_\lambda$ of the form described in \cref{prop:o-form}.
  Let $A$ be the pseudocolimit of the diagram $(A_k)_{k \in K}$ in $\TCat$, and let $E$ be the union of the images of the sets $E_k$ under the left Kan extensions $\PP(A_k) \to \PP(A)$.
  Then $\OO(A, E)$ is the colimit of the diagram $(\OO(A_k, E_k))_{k \in K}$ in $\LPr_\lambda$, and also in $\LPr$.
\end{proposition}

\begin{proof}
  First, note that $E$ is a set of morphisms between $\lambda$-compact objects of $\PP(A)$ because the left Kan extensions $\PP(A_k) \to \PP(A)$ preserve $\lambda$-compact objects.
  Thus $\OO(A, E)$ is a locally $\lambda$-presentable category.
  The functor $\PP(A_k) \to \PP(A)$ sends the morphisms of $E_k$ to morphisms of $E$ by definition, and therefore the composition $\PP(A_k) \to \PP(A) \to \OO(A, E)$ inverts the morphisms of $E_k$ and so factors essentially uniquely through $\OO(A_k, E_k)$ as a strongly $\lambda$-accessible left adjoint.
  These functors $\OO(A_k, E_k) \to \OO(A, E)$ define a cocone on the original diagram.
  Let $D$ be any object of $\LPr_\lambda$.
  Then there is a diagram
  \[
    \begin{tikzcd}
      \THom_{\LPr_\lambda}(\OO(A, E), D) \ar[r] \ar[d] &
      \lim\limits_{k \in K} \THom_{\LPr_\lambda}(\OO(A_k, E_k), D) \ar[d] \\
      \THom_{\LPr_\lambda}(\PP(A), D) \ar[r] \ar[d, "\eqv"'] &
      \lim\limits_{k \in K} \THom_{\LPr_\lambda}(\PP(A_k), D) \ar[d, "\eqv"] \\
      \TCat(A, D^\lambda) \ar[r, "\eqv"'] &
      \lim\limits_{k \in K} \TCat(A_k, D^\lambda)
    \end{tikzcd}
  \]
  in which the bottom functor is an equivalence because $A$ is the (pseudo)colimit of the diagram $(A_k)_{k \in K}$ in $\TCat$.
  Thus, the middle horizontal functor is an equivalence.
  It sends a left adjoint $F : \PP(A) \to D$ to the cone whose $k$ component is the composition $\PP(A_k) \to \PP(A) \xrightarrow{F} D$.
  Hence, each component belongs to the subcategory $\THom_{\LPr_\lambda}(\OO(A_k, E_k), D)$ of left adjoints which send the morphisms of $E_k$ to isomorphisms of $D$ if and only if the original functor sends the morphisms of $E$ to isomorphisms of $D$, by the definition of $E$.
  So, by \cref{prop:out-of-perp}, the top horizontal functor is also an equivalence.
  Hence $\OO(A, E)$ is the colimit of the diagram $\OO(A_k, E_k)$ in $\LPr_\lambda$.
  The same argument with $\LPr_\lambda$ replaced by $\LPr$ and $D^\lambda$ replaced by all of $D$ shows that $\OO(A, E)$ is also the colimit in $\LPr$.
\end{proof}

Together \cref{prop:o-form,prop:colimit-of-o-form} imply that $\LPr$ has all colimits and $\LPr_\lambda$ is closed in $\LPr$ under all colimits.
We will use \cref{prop:colimit-of-o-form} in the next chapter, when studying directed colimits in $\LPr_\lambda$.

\section{Presentations of combinatorial premodel categories}

The morphisms $E$ that define an orthogonality class $\OO(A, E)$ play a role analogous to that of the generating cofibrations $I$ and anodyne cofibrations $J$ of a premodel category; both impose conditions on a morphism (left adjoint or left Quillen functor) out of the object in question.
Accordingly, it makes sense to combine the data of $A$, $E$, $I$, and $J$ as we describe next.

\begin{definition}
  A \emph{presentation} of a combinatorial premodel category consists of
  \begin{enumerate}
  \item a small category $A$,
  \item a set $E$ of morphisms of $\PP(A)$, and
  \item sets $I$ and $J$ of morphisms of $\OO(A, E)$ such that $\rlp(I) \subset \rlp(J)$.
  \end{enumerate}
  The combinatorial premodel category ``presented'' by this presentation has underlying category $\OO(A, E)$, generating cofibrations $I$ and generating anodyne cofibrations $J$.

  We call $(A, E, I, J)$ a \emph{$\lambda$-presentation} if $E$ consists of morphisms between $\lambda$-compact objects of $\PP(A)$ and $I$ and $J$ consist of morphisms between $\lambda$-compact objects of $\OO(A, E)$.
\end{definition}

Let $M$ be a premodel category with $\lambda$-presentation $(A, E, I, J)$ and let $N$ be a locally $\lambda$-presentable category.
Then we can describe the category $\THom_{\CPM_\lambda}(M, N)$ as follows.
Any functor $A \to N^\lambda$ induces an essentially unique morphism $\PP(A) \to N$ of $\LPr_\lambda$.
This functor factors essentially uniquely through $\OO(A, E)$ (as a morphism of $\LPr_\lambda$) if and only if it sends the morphisms of $E$ to isomorphisms of $N$.
In turn, this left adjoint is a left Quillen functor and hence belongs to $\CPM_\lambda$ if and only if it sends the morphisms of $I$ to cofibrations and the morphisms of $J$ to anodyne cofibrations.
Thus, we obtain a sequence
\begin{align*}
  \THom_{\CPM_\lambda}(M, N)
  & \to \THom_{\LPr_\lambda}(\OO(A, E), N) \\
  & \to \THom_{\LPr_\lambda}(\PP(A), N) \\
  & \eqv \TCat(A, N^\lambda)
\end{align*}
in which the first two functors are fully faithful.
The same applies for $\THom_\CPM(M, N)$, without the restriction that the original functor out of $A$ takes values in the $\lambda$-compact objects of $N$.

A premodel category is $\lambda$-combinatorial if and only if it admits a $\lambda$-presentation.
More generally, any diagram in $\CPM_\lambda$ admits a ``diagram of $\lambda$-presentations'' in which the diagram of underlying locally $\lambda$-presentable categories has the form described in \cref{prop:o-form}.
By \cref{remark:colimit-gen,prop:colimit-of-o-form}, the colimit in $\CPM_\lambda$ or $\CPM$ of a diagram presented in this way is computed as follows:
\begin{enumerate}
\item First, form the colimit $\OO(A, E)$ of the underlying locally presentable categories $\OO(A_k, E_k)$ by taking $A$ to be $\colim_{k \in K} A_k$ and $E$ the union of the images of the $E_k$ in $\PP(A)$.
\item Second, take the union of the images of the generating (anodyne) cofibrations of each original combinatorial premodel category and use these as the generating (anodyne) cofibrations for the colimit $\OO(A, E)$.
\end{enumerate}

\section{Tensor products and the internal Hom}

For locally presentable categories $M$ and $N$, the category of all left adjoints from $M$ to $N$ is again locally presentable.
We denote it by $\LPr(M, N)$.
By the adjoint functor theorem a left adjoint from $M$ to $N$ is the same as a colimit-preserving functor, and because colimits commute with colimits it follows that colimits in $\LPr(M, N)$ are computed componentwise.

Now let $M_1$, $M_2$ and $N$ be locally presentable categories.
A functor from $M_1$ to $\LPr(M_2, N)$ is a left adjoint if and only if it preserves colimits, by the adjoint functor theorem.
Because colimits in $\LPr(M_2, N)$ are computed componentwise, such a functor is the same as a functor from $M_1 \times M_2$ to $N$ which preserves colimits in each variable separately, that is, an adjunction of two variables (using the adjoint functor theorem again).

Suppose that $M_1 = \OO(A_1, E_1)$ and $M_2 = \OO(A_2, E_2)$.
There is an ``external tensor product'' $\boxtimes : \PP(A_1) \times \PP(A_2) \to \PP(A_1 \times A_2)$ which is the adjunction of two variables sending $(\yo a_1, \yo a_2)$ to $\yo (a_1, a_2)$.
Set
\[
  E = \{\,e_1 \boxtimes \yo a_2 \mid e_1 \in E_1, a_2 \in \Ob A_2\,\} \cup
  \{\,\yo a_1 \boxtimes e_2 \mid a_1 \in \Ob A_1, e_2 \in E_2\,\}.
\]
The functor $- \boxtimes \yo a_2 : \PP(A_1) \to \PP(A_1 \times A_2)$ is just the left Kan extension of $(-, a_2) : A_1 \to A_1 \times A_2$ and similarly for $\yo a_1 \boxtimes -$, so we could define $E$ without reference to $\boxtimes$.
In particular $E$ is a set of morphisms between $\lambda$-compact objects if both $E_1$ and $E_2$ are.

\begin{proposition}\label{prop:tensor-o}
  There is a canonical equivalence
  \[
    \LPr(\OO(A_1, E_1), \LPr(\OO(A_2, E_2), N)) \eqv
    \LPr(\OO(A_1 \times A_2, E), N).
  \]
  Hence $\LPr(\OO(A_1 \times A_2, E), N)$ is equivalent to the category of adjunctions of two variables from $\OO(A_1, E_1) \times \OO(A_2, E_2)$ to $N$.
\end{proposition}

\begin{proof}
  First, there is a canonical equivalence
  \[
    \LPr(\PP(A_1), \LPr(\PP(A_2), N)) \eqv (N^{A_2})^{A_1}
    \eqv N^{A_1 \times A_2}
    \eqv \LPr(\PP(A_1 \times A_2), N).
  \]
  The full subcategory $\LPr(\OO(A_2, E_2), N)$ of $\LPr(\PP(A_2), N)$ is closed under colimits because colimits preserve isomorphisms.
  Then one can compute that a left adjoint from $\PP(A_1)$ to $\LPr(\PP(A_2), N)$ lands in $\LPr(\OO(A_2, E_2), N)$ and sends $E_1$ to isomorphisms if and only if the corresponding left adjoint from $\PP(A_1 \times A_2)$ to $N$ sends $E$ to isomorphisms.
\end{proof}

Thus, for any locally presentable categories $M_1$ and $M_2$, there exists a tensor product $M_1 \otimes M_2$ equipped with a universal adjunction of two variables $\otimes : M_1 \times M_2 \to M_1 \otimes M_2$.
Namely, we write $M_1 = \OO(A_1, E_1)$ and $M_2 = \OO(A_2, E_2)$ and then set $M_1 \otimes M_2 = \OO(A_1 \times A_2, E)$ with $E$ defined as above.
Then for any $N$, $\LPr(M_1 \otimes M_2, N)$ is equivalent to the category of adjunctions of two variables from $M_1 \times M_2$ to $N$.
The tensor product $M_1 \otimes M_2$ is determined up to equivalence by this universal property, and is therefore independent of the choices of presentations of $M_1$ and $M_2$ and automatically functorial.
By the above computation, $M_1 \otimes M_2$ is locally $\lambda$-presentable if both $M_1$ and $M_2$ are.

Suppose $M_1 = \OO(A_1, E_1)$, $M_2 = \OO(A_2, E_2)$, $M_3 = \OO(A_3, E_3)$ are three locally presentable categories presented as orthogonality classes.
The formula we have given for the tensor product makes it clear that both $(M_1 \otimes M_2) \otimes M_3$ and $M_1 \otimes (M_2 \otimes M_3)$ are canonically equivalent to $M_1 \otimes M_2 \otimes M_3 = \OO(A_1 \times A_2 \times A_3, E)$ for
\begin{align*}
  E & = \{\,e_1 \boxtimes \yo a_2 \boxtimes \yo a_3 \mid e_1 \in E_1, a_2 \in \Ob A_2, a_3 \in \Ob A_3\,\} \\
    & \qquad \cup \{\,\yo a_1 \boxtimes e_2 \boxtimes \yo a_3 \mid a_1 \in \Ob A_1, e_2 \in E_2, a_3 \in \Ob A_3\,\} \\
    & \qquad \cup \{\,\yo a_1 \boxtimes \yo a_2 \boxtimes e_3 \mid a_1 \in \Ob A_1, a_2 \in \Ob A_2, e_3 \in E_3\,\}
\end{align*}
and an argument similar to that of \cref{prop:tensor-o} shows that $M_1 \otimes M_2 \otimes M_3$ is the universal recipient of an adjunction of \emph{three} variables out of $M_1 \times M_2 \times M_3$.
In particular, $\otimes$ is associative up to equivalence.
The object $\Set = \OO(\{*\}, \emptyset)$ is the unit object for the tensor product.

For any locally presentable category $M$ and small category $A$ the tensor $A \otimes M$ can be identified with $\PP(A) \otimes M$, since for any $N$,
\begin{align*}
  \THom_\LPr(A \otimes M, N)
  & \eqv \TCat(A, \LPr(M, N)) \\
  & \eqv \THom_\LPr(\PP(A), \LPr(M, N)) \\
  & \eqv \THom_\LPr(\PP(A) \otimes M, N).
\end{align*}
Thus tensors by categories are a special case of tensor products in $\LPr$.

\begin{definition}
  An adjunction of two variables $F : M_1 \times M_2 \to N$ is \emph{strongly $\lambda$-accessible} if $F(X_1, X_2)$ is $\lambda$-compact whenever both $X_1$ and $X_2$ are $\lambda$-compact.
\end{definition}

\begin{proposition}\label{prop:strongly-accessible-bifunctor}
  Suppose $M_1$, $M_2$ and $N$ are locally $\lambda$-presentable.
  Then the subcategory $\LPr_\lambda(M_1 \otimes M_2, N)$ of $\LPr(M_1 \otimes M_2, N)$ corresponds to the subcategory of adjunctions of two variables from $M_1 \times M_2$ to $N$ which are strongly $\lambda$-accessible.
  In particular, the universal adjunction of two variables $\otimes : M_1 \times M_2 \to M_1 \otimes M_2$ is itself strongly $\lambda$-accessible.
\end{proposition}

\begin{proof}
  We can choose presentations $M_i = \OO(A_i, E_i)$ where $A_i$ is the full subcategory consisting of the $\lambda$-compact objects of $M_i$.
  Then $M_1 \otimes M_2 = \OO(A_1 \times A_2, E)$.
  For an adjunction of two variables $F : M_1 \times M_2 \to N$, let $F' : M_1 \otimes M_2 \to N$ be the corresponding left adjoint.
  Then the composition
  \begin{equation*}
    A_1 \times A_2 \to \PP(A_1 \times A_2) \to \OO(A_1 \times A_2, E) \eqv M_1 \otimes M_2 \xrightarrow{F'} N
    \tag{$*$}
  \end{equation*}
  agrees up to isomorphism with the restriction of $F : M_1 \times M_2 \to N$ to the subcategory $A_1 \times A_2$.
  The result follows because $F'$ is strongly $\lambda$-accessible if and only if the composition $(*)$ sends $A_1 \times A_2$ into $N^\lambda$.
\end{proof}

We now extend these facts to combinatorial premodel categories.

\begin{definition}
  Let $M_1$ and $M_2$ be combinatorial premodel categories.
  The \emph{tensor product} $M_1 \otimes M_2$ is a combinatorial premodel category equipped with a Quillen bifunctor $\otimes : M_1 \times M_2 \to M_1 \otimes M_2$ such that for any combinatorial premodel category $N$, composition with $\otimes$ induces an equivalence between $\THom_\CPM(M_1 \otimes M_2, N)$ and the category of Quillen bifunctors $M_1 \times M_2 \to N$.
\end{definition}

\begin{proposition}\label{prop:cpm-tensor}
  For any combinatorial premodel categories $M_1$ and $M_2$, the tensor product $M_1 \otimes M_2$ exists, and is $\lambda$-combinatorial if both $M_1$ and $M_2$ are.
  This tensor product is preserved by the forgetful functor $\CPM \to \LPr$.
\end{proposition}

\begin{proof}
  Choose generating cofibrations $I_i$ and anodyne cofibrations $J_i$ of $M_i$ for $i = 1$ and $2$, which are sets of morphisms between $\lambda$-compact objects if both $M_i$ are $\lambda$-combinatorial.
  Define the tensor product $M_1 \otimes M_2$ to have underlying category the tensor product of the locally presentable categories $M_1$ and $M_2$, and generating (anodyne) cofibrations
  \[
    I = I_1 \bp I_2, \quad J = J_1 \bp I_2 \cup I_1 \bp J_2
  \]
  where $\bp$ is taken with respect to $\otimes : M_1 \times M_2 \to M_1 \otimes M_2$.
  The correctness of this construction follows from the universal property of the tensor product of locally presentable categories together with \cref{prop:quillen-bifunctor}.
  When $M_1$ and $M_2$ are both $\lambda$-combinatorial, the category $M_1 \otimes M_2$ is locally $\lambda$-presentable, $\otimes : M_1 \times M_2 \to M_1 \otimes M_2$ is strongly $\lambda$-accessible, and $I_1$, $J_1$, $I_2$, $J_2$ are sets of morphisms between $\lambda$-compact objects; therefore $M_1 \otimes M_2$ is $\lambda$-combinatorial.
\end{proof}

Suppose $M_1$, $M_2$, $M_3$ are three combinatorial premodel categories and choose generating cofibrations $I_i$ and anodyne cofibrations $J_i$ of $M_i$ for $i = 1$, $2$, $3$.
Under the identifications of the underlying categories
\[
  (M_1 \otimes M_2) \otimes M_3 \eqv M_1 \otimes M_2 \otimes M_3 \eqv
  M_1 \otimes (M_2 \otimes M_3),
\]
both premodel category structures $(M_1 \otimes M_2) \otimes M_3$ and $M_1 \otimes (M_2 \otimes M_3)$ are generated by
\[
  I = I_1 \bp I_2 \bp I_3, \quad
  J = J_1 \bp I_2 \bp I_3 \cup
  I_1 \bp J_2 \bp I_3 \cup
  I_1 \bp I_2 \bp J_3.
\]
Hence the tensor product of combinatorial premodel categories is associative up to equivalence.
Its unit object is $\Set$ with its usual premodel category structure generated by $I = \{\emptyset \to *\}$ and $J = \emptyset$, because $\emptyset \to *$ is the unit for $\bp$ up to the identification $\Set \otimes M \eqv M$.

We now turn to the internal Hom $\CPM(M, N)$, which will be a combinatorial premodel category with underlying category $\LPr(M, N)$, the category of all left adjoints from $M$ to $N$.
(Note that the category $\THom_\CPM(M, N) = \LQF(M, N)$ could not be the underlying category of a premodel category, since $\LQF(M, N)$ generally does not have limits or colimits.)
The cofibrations and anodyne cofibrations are defined explicitly as described below.

\begin{definition}\label{def:cpm-cof}
  Let $M$ and $N$ be premodel categories and let $t : F \to F'$ a natural transformation between two left adjoints $F : M \to N$ and $F' : M \to N$.
  Then we call $t$
  \begin{itemize}
  \item a \emph{cofibration} if for every cofibration $f : A \to B$ of $M$, the induced map $F'A \amalg_{FA} FB \to F'B$ is a cofibration which is anodyne if $f$ is;
  \item an \emph{anodyne cofibration} if for every cofibration $f : A \to B$ of $M$, the induced map $F'A \amalg_{FA} FB \to F'B$ is an anodyne cofibration.
  \end{itemize}
\end{definition}

We will show below that these classes of morphisms do indeed define a combinatorial premodel category structure $\CPM(M, N)$ on the category $\LPr(M, N)$.
More or less by definition, an object $F : M \to N$ of $\CPM(M, N)$ is cofibrant if and only if it is a left Quillen functor.
In particular,
\begin{align*}
  \THom_\CPM(\Set, \CPM(M, N))
  & \eqv \CPM(M, N)^\cof \\
  & = \THom_\CPM(M, N) \\
  & \eqv \THom_\CPM(\Set \otimes M, N).
\end{align*}
More generally, we will show that $\CPM(M, N)$ has the expected adjunction relationship with the tensor product of $\CPM$.

For now, since we have not yet shown that $\CPM(M, N)$ is a combinatorial premodel category, we will continue to denote its underlying category by $\LPr(M, N)$.

The property of being a left Quillen functor can be checked on generating (anodyne) cofibrations.
We show that, more generally, the classes of (anodyne) cofibrations in $\LPr(M, N)$ have the same property.

\begin{lemma}\label{lemma:cpm-cof-gen}
  Suppose $I$ and $J$ are classes of generating cofibrations and anodyne cofibrations for $M$.
  Then the property of being a cofibration or an anodyne cofibration in $\LPr(M, N)$ as defined in \cref{def:cpm-cof} can be checked on just the cofibrations $f : A \to B$ belonging to $I$ and the anodyne cofibrations belonging to $J$.
\end{lemma}

\begin{proof}
  Let $t : F \to F'$ be a morphism of $\LPr(M, N)$.
  We will show that the following are equivalent:
  \begin{enumerate}
  \item For every cofibration $f : A \to B$ of $M$, the induced map $F'A \amalg_{FA} FB \to F'B$ is a cofibration.
  \item For every cofibration $f : A \to B$ belonging to $I$, the induced map $F'A \amalg_{FA} FB \to F'B$ is a cofibration.
  \end{enumerate}
  Similar arguments will apply to the second condition for $t$ to be a cofibration, and for the condition for $t$ to be an anodyne cofibration.

  Clearly it suffices to show that (2) implies (1).
  Write $G : N \to M$ and $G' : N \to M$ for the right adjoints of $F$ and $F'$ respectively.
  The natural transformation $t : F \to F'$ induces a corresponding natural transformation $u : G' \to G$.
  Consider a lifting problem of the form below, where $p : X \to Y$ is an anodyne fibration of $N$.
  \[
    \begin{tikzcd}
      F'A \amalg_{FA} FB \ar[r] \ar[d] & X \ar[d, "p"] \\
      F'B \ar[r] \ar[ru, dotted] & Y
    \end{tikzcd}
  \]
  An adjointness argument shows that finding lifts in such squares is equivalent to finding lifts in the squares
  \[
    \begin{tikzcd}
      A \ar[r] \ar[d, "f"'] & G' X \ar[d] \\
      B \ar[r] & G' Y \times_{G Y} G X
    \end{tikzcd}
  \]
  If all such squares admits lifts for all $f \in I$ and all anodyne fibrations $p$, then the same is true without the condition on $f$ because $I$ generates the cofibrations of $M$.
  Therefore if $F'A \amalg_{FA} FB \to F'B$ is a cofibration for all $f \in I$, then it is a cofibration for all cofibrations $f$ of $M$.
\end{proof}

\begin{proposition}\label{prop:cpm-internal-hom}
  Suppose $M$ and $N$ are combinatorial premodel categories.
  Then the above cofibrations and anodyne cofibrations define a combinatorial premodel category structure on $\LPr(M, N)$.
\end{proposition}

\begin{proof}
  The strategy is to construct a left adjoint $E : \LPr(M, N) \to Z$ for some combinatorial premodel category $Z$ such that the resulting left-induced premodel category structure on $\LPr(M, N)$ has the classes of cofibrations and anodyne cofibrations defined by \cref{def:cpm-cof}.
  More specifically, we will prove the following result.
  Let $\lambda$ be an uncountable regular cardinal such that:
  \begin{enumerate}
  \item The category $\LPr(M, N)$ is locally $\lambda$-presentable.
  \item For each $\lambda$-compact object $A$ of $M$, the functor from $\LPr(M, N)$ to $N$ given by evaluation at $A$ is strongly $\lambda$-accessible.
  \item $M$ admits generating cofibrations $I$ and anodyne cofibrations $J$ which are $\lambda$-small sets of morphisms between $\lambda$-compact objects.
  \item $N$ is $\lambda$-combinatorial.
  \end{enumerate}
  Then $\CPM(M, N)$ is a $\lambda$-combinatorial premodel category.

  Clearly conditions (1), (3) and (4) are satisfied for all sufficiently large $\lambda$.
  We will study condition (2) in the next section, after which it will also be clear that it holds for sufficiently large $\lambda$.
  This can also be checked directly using the fact that $\LPr(M, N)$ has only a set of isomorphism classes of $\lambda$-compact objects for any $\lambda$.

  Suppose that $\lambda$ satisfies all of the above conditions.
  We set $Z = \prod_{i \in I} N^{[1]}_\cof \times \prod_{j \in J} N^{[1]}_\acof$ where
  \begin{itemize}
  \item $N^{[1]}_\cof$ has the premodel category structure in which a morphism $X \to Y$ is a cofibration or an anodyne cofibration if and only if $Y_0 \amalg_{X_0} X_1 \to Y_1$ is one in $N$;
  \item $N^{[1]}_\acof$ has the premodel category structure in which a morphism $X \to Y$ is both a cofibration and an anodyne cofibration if and only if $Y_0 \amalg_{X_0} X_1 \to Y_1$ is an anodyne cofibration in $N$.
  \end{itemize}
  By \cref{ex:reedy-arrow,prop:reedy-gen} these classes do determine $\lambda$-combinatorial premodel category structures $N^{[1]}_\cof$ and $N^{[1]}_\acof$ on the arrow category $N^{[1]}$.

  For each member $i : A \to B$ of $I$, we define a functor $E_i : \LPr(M, N) \to N^{[1]}_\cof$ which sends a left adjoint $F : M \to N$ to the morphism $Fi : FA \to FB$.
  In the same way, we define $E_j : \LPr(M, N) \to N^{[1]}_\acof$ by $E_j(F) = Fj$.
  (The only difference between the $E_i$ and the $E_j$ is the chosen premodel category structure on the target category.)
  These functors are colimit-preserving and strongly $\lambda$-accessible by condition (2), since $I$ and $J$ consist of morphisms between $\lambda$-compact objects of $M$.
  They are the components of a strongly $\lambda$-accessible left adjoint $E : \LPr(M, N) \to Z$ because $I$ and $J$ each have fewer than $\lambda$ elements and so $Z = \prod_{i \in I} N^{[1]}_\cof \times \prod_{j \in J} N^{[1]}_\acof$ is also a product in $\CPM_\lambda$ and in $\LPr_\lambda$.
  By condition (1), $\LPr(M, N)$ is locally $\lambda$-presentable and hence we may apply \cref{prop:left-induced} to equip it with a $\lambda$-combinatorial premodel category structure $\CPM(M, N)$ in which a morphism $t : F \to F'$ is a cofibration or an anodyne cofibration if and only if its image under $E$ is one.

  Now $Z$ and $E$ have been constructed so that for any morphism $t : F \to F'$, $Et$ is
  \begin{itemize}
  \item a cofibration if for every $f : A \to B$ in $I$, the induced map $F'A \amalg_{FA} FB \to FB'$ is a cofibration, and for every $f : A \to B$ in $J$, the induced map $F'A \amalg_{FA} FB \to FB'$ is an anodyne cofibration;
  \item an anodyne cofibration if for every $f : A \to B$ in either $I$ or $J$, the induced map $F'A \amalg_{FA} FB \to FB'$ is an anodyne cofibration.
  \end{itemize}
  By \cref{lemma:cpm-cof-gen} these conditions are equivalent to $t$ being a cofibration or an anodyne cofibration as defined by \cref{def:cpm-cof}.
\end{proof}

The premodel category structure on $\CPM$ is chosen so that the expected tensor--Hom adjunction holds.

\begin{proposition}\label{prop:cpm-tensor-hom}
  Let $M_1$, $M_2$ and $N$ be combinatorial premodel categories.
  The tensor--Hom adjunction
  \[
    \THom_\LPr(M_1 \otimes M_2, N) \eqv \THom_\LPr(M_1, \LPr(M_2, N))
  \]
  restricts to an adjunction
  \[
    \THom_\CPM(M_1 \otimes M_2, N) \eqv \THom_\CPM(M_1, \CPM(M_2, N)).
  \]
\end{proposition}

\begin{proof}
  We need to show that an adjunction of two variables $F_0 : M_1 \times M_2 \to N$ is a Quillen bifunctor if and only if the corresponding left adjoint $F : M_1 \to \CPM(M_2, N)$ is a left Quillen functor, where $F$ is given by the formula $F(A_1)(A_2) = F_0(A_1, A_2)$.
  By definition, $F$ is a left Quillen functor if and only if:
  \begin{enumerate}
  \item For every cofibration $f_1 : A_1 \to B_1$ of $M_1$, $F(f_1) : F(A_1) \to F(B_1)$ is a cofibration of $\CPM(M_2, N)$.
    That is, for every cofibration $f_2 : A_2 \to B_2$ of $M_2$, $F(B_1)(A_2) \amalg_{F(A_1)(A_2)} F(A_1)(B_2) \to F(B_1)(B_2)$ is a cofibration which is anodyne if $f_2$ is.
  \item For every anodyne cofibration $f_1 : A_1 \to B_1$ of $M_1$, $F(f_1) : F(A_1) \to F(B_1)$ is an anodyne cofibration of $\CPM(M_2, N)$.
    That is, for every cofibration $f_2 : A_2 \to B_2$ of $M_2$, $F(B_1)(A_2) \amalg_{F(A_1)(A_2)} F(A_1)(B_2) \to F(B_1)(B_2)$ is an anodyne cofibration.
  \end{enumerate}
  These are precisely the conditions that are required for the original adjunction of two variables $F_0$ to be a Quillen bifunctor.
\end{proof}

In fact this adjunction is the restriction to cofibrant objects of an ``enriched tensor--Hom adjunction''
\[
  \CPM(M_1 \otimes M_2, N) \eqv \CPM(M_1, \CPM(M_2, N)).
\]
This can be verified formally using the associativity of $\otimes$, or directly by an analysis of a cube induced by three morphisms $t : F \to F'$, $f_1 : A_1 \to B_1$, $f_2 : A_2 \to B_2$.

Let $M$ be a combinatorial premodel category and $A$ a small category.
Using \cref{prop:cpm-tensor-hom} we can identify the tensor $A \otimes M = M^{A^\op}_\proj$ with $\Set^{A^\op}_\proj \otimes M$, since
\begin{align*}
  \THom_\CPM(A \otimes M, N)
  & \eqv \TCat(A, \THom_\CPM(M, N)) \\
  & \eqv \TCat(A, \THom_\CPM(\Set, \CPM(M, N))) \\
  & \eqv \THom_\CPM(A \otimes \Set, \CPM(M, N)) \\
  & \eqv \THom_\CPM(\Set^{A^\op}_\proj, \CPM(M, N)) \\
  & \eqv \THom_\CPM(\Set^{A^\op}_\proj \otimes M, N).
\end{align*}
Thus tensors by categories are a special case of tensor products in $\CPM$.

Unlike the other results in this chapter, \cref{prop:cpm-internal-hom} does not provide information on the rank of combinatoriality of $\CPM(M, N)$ or about the relationship between $\CPM(M, N)$ and the sub-2-categories $\CPM_\lambda$.
We will need such information in order to carry out our variant of the small object argument in which we attach only solutions to lifting problems of some bounded ``size''.
The issue is with conditions (1) and (2) on $\lambda$ which appeared in the proof of \cref{prop:cpm-internal-hom}.
In order to control these conditions we need additional information about $M$, which we study in the next section.

\section{The size of a combinatorial premodel category}\label{sec:algebra-size}

Conditions (1) and (2) in the proof of \cref{prop:cpm-internal-hom} depend only on the underlying locally presentable categories involved.
The key property which controls these conditions is a bound on the size of a presentation for $M$, analogous to condition (3) in the same proof.

\begin{definition}
  Let $\mu \ge \lambda$ be regular cardinals.
  A category is \emph{$\mu$-small $\lambda$-presentable} if it is equivalent to a category of the form $\OO(A, E)$ with $A$ a $\mu$-compact category and $E$ a $\mu$-small set of morphisms between $\lambda$-compact objects of $\PP(A)$.
\end{definition}

If $M$ is a locally $\lambda$-presentable category, then $M$ admits a $\lambda$-presentation $\OO(A, E)$.
We call $\OO(A, E)$ a \emph{$\mu$-small $\lambda$-presentation} when $A$ and $E$ satisfy the conditions of the above definition.
Every $\lambda$-presentation is $\mu$-small for some $\mu$, so every locally $\lambda$-presentable category is $\mu$-small $\lambda$-presentable for some $\mu$.
A $\mu$-small $\lambda$-presentable category is also $\mu'$-small $\lambda'$-presentable for any $\lambda' \ge \lambda$ and $\mu' \ge \mu$.

\begin{remark}
  When $\lambda = \mu = \aleph_0$, the condition that $A$ is $\mu$-compact means that $A$ is finitely presented, not necessarily finite.
  It turns out to make no difference to the class of $\aleph_0$-small $\aleph_0$-presentable categories if we instead require $A$ to be finite, but we will not need this fact.
  In all other cases, $\mu$ is uncountable and then $A$ being $\mu$-compact is equivalent to $A$ being $\mu$-small, that is, having fewer than $\mu$ objects and morphisms.
  We choose to require $A$ to be $\mu$-compact rather than $\mu$-small in order to make the following fact obvious.
\end{remark}

\begin{proposition}\label{small-presentable-closure}
  The class of $\mu$-small $\lambda$-presentable categories is closed under $\mu$-small colimits, tensors by $\mu$-compact categories and tensor products.
\end{proposition}

\begin{proof}
  Closure under $\mu$-small colimits follows from \cref{prop:colimit-of-o-form} and closure under tensor products follows from the formula for $\OO(A_1, E_1) \otimes \OO(A_2, E_2)$.
  Closure under tensors by $\mu$-compact categories follows from closure under tensor products because $A \otimes M = \PP(A) \otimes M$ and $\PP(A) = \OO(A, \emptyset)$ is $\mu$-small $\lambda$-presentable when $A$ is $\mu$-compact.
\end{proof}

\begin{remark}
  The cardinal $\mu$ will play a role in the next chapter when we consider $\mu$-directed colimits in $\CPM_\lambda$.
  In order for these to have the expected properties we will need $\mu \ge \lambda$, so it is convenient to build in this assumption from the start.
  In this rest of this section we will only need to consider the case $\mu = \lambda$.
\end{remark}

\begin{proposition}\label{prop:lpr-rank}
  Let $\lambda$ be an uncountable regular cardinal.
  Suppose $M$ is $\lambda$-small $\lambda$-presentable and $N$ is locally $\lambda$-presentable.
  Then the category $\LPr(M, N)$ of left adjoints from $M$ to $N$ is locally $\lambda$-presentable.
  Moreover, an object $F$ of $\LPr(M, N)$ is $\lambda$-compact if and only if it is strongly $\lambda$-accessible as a functor $F : M \to N$.
\end{proposition}

\begin{proof}
  Choose a $\lambda$-small $\lambda$-presentation $M = \OO(A, E)$.
  By \cref{prop:out-of-perp}, $\LPr(\OO(A, E), N)$ is equivalent (via precomposition with $L : \PP(A) \to \OO(A, E)$) to the full subcategory of $\LPr(\PP(A), N)$ on the functors which invert the morphisms of $E$.
  In turn $\LPr(\PP(A), N)$ is equivalent to the diagram category $N^A$, which is locally $\lambda$-presentable.
  For each object $X$ of $\PP(A)$, let $H_X : \LPr(\PP(A), N) \to N$ denote the functor given by evaluation on $X$.
  For $X$ a representable presheaf $\yo a$, $H_X$ corresponds to the projection $\pi_a : N^A \to N$ which is strongly $\lambda$-accessible by \cref{prop:compact-diagram} because $A$ is $\lambda$-small.
  (We only assumed that $A$ is $\lambda$-compact, but $A$ is also $\lambda$-small because $\lambda$ is uncountable.)
  The functor $H_X$ commutes with colimits in the variable $X$ and the closure of the representables under $\lambda$-small colimits consists of all $\lambda$-compact objects of $\PP(A)$; hence $H_X$ is strongly $\lambda$-accessible for any $\lambda$-compact $X$.

  For each morphism $e$ in $E$, let $H_e : \LPr(\PP(A), N) \to N^{[1]}$ be the functor given by evaluation on $e$.
  These functors are strongly $\lambda$-accessible because $E$ consists of morphisms between $\lambda$-compact objects.
  Because $E$ has cardinality less than $\lambda$, the functors $H_e$ assemble to a strongly $\lambda$-accessible functor $H : \LPr(\PP(A), N) \to \prod_{e \in E} N^{[1]}$.
  We now construct the solid part of the diagram below in which the vertical functor sends each copy of $N$ to $N^{[1]}$ by sending an object $X$ to the constant diagram $X \to X$.
  \begin{equation*}
    \begin{tikzcd}
      \LPr(\OO(A, E), N) \ar[r, dotted] \ar[d, dotted, "L^*"'] & \prod_{e \in E} N \ar[d] \\
      \LPr(\PP(A), N) \ar[r, "H"'] & \prod_{e \in E} N^{[1]}
    \end{tikzcd}
    \tag{$*$}
  \end{equation*}
  The pseudopullback of this diagram is the category whose objects consist of a left adjoint from $\PP(A)$ to $N$ together with, for each $e \in E$, an isomorphism from the image of $e$ to a constant diagram.
  Up to equivalence this is the same as the full subcategory of $\LPr(\PP(A), N)$ on the functors which send each $e \in E$ to an isomorphism and therefore may be identified with $\LPr(\OO(A, E), N)$ as shown by the dotted arrows.

  Now the solid part of $(*)$ is a finite diagram of locally $\lambda$-presentable categories and strongly $\lambda$-accessible functors and so by \cref{prop:lpr-limits} its pullback $\LPr(\OO(A, E), N) = \LPr(M, N)$ is also a pullback in $\LPr_\lambda$.
  In particular $\LPr(M, N)$ is locally $\lambda$-presentable.

  It remains to identify the $\lambda$-compact objects of $\LPr(\OO(A, E), N)$.
  A $\lambda$-compact object corresponds to a morphism in $\LPr_\lambda$ from $\Set$, so since the completed diagram $(*)$ is a pullback, an object $F : \OO(A, E) \to N$ is $\lambda$-compact if and only if its images $L^*(F) = F \circ L$ in $\LPr(\PP(A), N)$ and in $\prod_{e \in E} N$ are both $\lambda$-compact.
  The second condition is redundant because $H : \LPr(\PP(A), N) \to \prod_{e \in E} N^{[1]}$ is strongly $\lambda$-accessible.
  From the equivalence $\LPr(\PP(A), N) \eqv N^A$ and \cref{prop:compact-diagram}, an object $F' : \PP(A) \to N$ of $\LPr(\PP(A), N)$ is $\lambda$-compact if and only if it sends every representable presheaf to a $\lambda$-compact object; in turn this is equivalent to $F'$ being strongly $\lambda$-accessible.
  Finally, by the last part of \cref{prop:out-of-perp}, $F' = F \circ L$ is strongly $\lambda$-accessible if and only if $F$ is.
\end{proof}

\begin{remark}
  The identification of $\LPr(\OO(A, E), N)$ as a pullback shows that $\OO(A, E)$ may also be described as a pushout in $\LPr$ as shown below.
  \[
    \begin{tikzcd}
      \coprod_{e \in E} \PP([1]) \ar[r] \ar[d] & \PP(A) \ar[d] \\
      \coprod_{e \in E} \PP(\{*\}) \ar[r] & \OO(A, E)
    \end{tikzcd}
  \]
  The top functor sends the generating morphism of each copy of $\PP([1])$ to the corresponding morphism $e$ of $\PP(A)$, so it belongs to $\LPr_\lambda$ if and only if $E$ consists of morphisms between $\lambda$-compact objects.
  When $\OO(A, E)$ is a $\mu$-small $\lambda$-presentation, then, $\OO(A, E)$ can be built up from objects of the form $\PP(A)$ for $A$ a $\mu$-compact category by the formation of $\mu$-small colimits in $\LPr_\lambda$.
  Therefore the $\mu$-small $\lambda$-presentable categories are precisely the closure of $\{\,\PP(A) \mid \mbox{$A$ a $\mu$-small category}\,\}$ in $\LPr_\lambda$ under $\mu$-small colimits.
  (In turn, $\PP(A) = A \otimes \Set$ and so we may also characterize the $\mu$-small $\lambda$-presentable categories as the closure of $\{\Set\}$ in $\LPr_\lambda$ under $\mu$-small \emph{weighted} colimits.)
\end{remark}

\begin{proposition}\label{prop:lpr-lam-tensor-hom}
  Let $M_1$, $M_2$ and $N$ be locally $\lambda$-presentable categories with $M_2$ $\lambda$-small $\lambda$-presentable.
  Then the adjunction
  \[
    \THom_\LPr(M_1 \otimes M_2, N) \eqv \THom_\LPr(M_1, \LPr(M_2, N))
  \]
  restricts to an adjunction
  \[
    \THom_{\LPr_\lambda}(M_1 \otimes M_2, N) \eqv \THom_{\LPr_\lambda}(M_1, \LPr(M_2, N)).
  \]
\end{proposition}

\begin{proof}
  According to \cref{prop:strongly-accessible-bifunctor}, we must show that an adjunction of two variables $F : M_1 \times M_2 \to N$ is strongly $\lambda$-accessible if and only if the corresponding functor $F' : M_1 \to \LPr(M_2, N)$ is strongly $\lambda$-accessible.
  This follows from the fact that the $\lambda$-compact objects of $\LPr(M_2, N)$ are the strongly $\lambda$-accessible functors.
\end{proof}

\Cref{prop:lpr-rank} can be recovered from this adjunction by taking $M_1 = \Set$.

We now apply these results to combinatorial premodel categories.

\begin{definition}
  Let $\mu \ge \lambda$ be regular cardinals.
  A $\lambda$-presentation $(A, E, I, J)$ is \emph{$\mu$-small} if $A$ is $\mu$-compact and $E$, $I$ and $J$ are $\mu$-small sets.
  A premodel category is \emph{$\mu$-small $\lambda$-combinatorial} if it admits a $\mu$-small $\lambda$-presentation.
  Explicitly, this means that its underlying category is $\mu$-small $\lambda$-presentable and it admits generating cofibrations and generating anodyne cofibrations which are $\mu$-small sets of morphisms between $\lambda$-compact objects.
\end{definition}

Every $\lambda$-combinatorial premodel category is $\mu$-small $\lambda$-combinatorial for some $\mu$.
A $\mu$-small $\lambda$-combinatorial premodel category is also $\mu'$-small $\lambda'$-combinatorial for any $\lambda' \ge \lambda$ and $\mu' \ge \mu$.

\begin{proposition}
  The class of $\mu$-small $\lambda$-combinatorial premodel categories is closed under $\mu$-small colimits, tensors by $\mu$-compact categories and tensor products.
\end{proposition}

\begin{proof}
  We checked the corresponding closure properties for the underlying locally presentable categories in \cref{small-presentable-closure}, so it remains to check the condition on the cardinalities of the generating (anodyne) cofibrations.
  For $\mu$-small colimits this follows from \cref{remark:colimit-gen}.
  For tensor products this follows from the formula for the generating (anodyne) cofibrations of the tensor products.
  For tensors by $\mu$-compact categories this follows from the identification of the tensor $A \otimes M$ with the tensor product $\Set^{A^\op}_\proj \otimes M$, together with the fact that $\Set^{A^\op}_\proj$ has fewer than $\mu$ generating cofibrations (one for each object of $A$) and no generating anodyne cofibrations.
\end{proof}

\begin{proposition}
  Let $\lambda$ be an uncountable regular cardinal and suppose $M$ is $\lambda$-small $\lambda$-combinatorial and $N$ is $\lambda$-combinatorial.
  Then $\CPM(M, N)$ is $\lambda$-combinatorial.
  Moreover, an object $F$ of $\CPM(M, N)$ is $\lambda$-compact if and only if it is strongly $\lambda$-accessible.
\end{proposition}

\begin{proof}
  We show that the conditions on $\lambda$ in the proof of \cref{prop:cpm-internal-hom} are met.
  Conditions (1) and (2) follow from \cref{prop:lpr-rank}, using the facts that $M$ is $\lambda$-small $\lambda$-presentable and $N$ is locally $\lambda$-presentable.
  Conditions (3) and (4) hold by assumption.
  Therefore $\CPM(M, N)$ is $\lambda$-combinatorial.
  The last claim also follows from \cref{prop:lpr-rank} because it is just a statement about the underlying category $\CPM(M, N) = \LPr(M, N)$.
\end{proof}

\begin{proposition}
  Let $M_1$, $M_2$ and $N$ be $\lambda$-combinatorial premodel categories and suppose that $M_2$ is $\lambda$-small $\lambda$-combinatorial.
  Then the adjunction
  \[
    \THom_\CPM(M_1 \otimes M_2, N) \eqv \THom_\CPM(M_1, \CPM(M_2, N))
  \]
  restricts to an adjunction
  \[
    \THom_{\CPM_\lambda}(M_1 \otimes M_2, N) \eqv \THom_{\CPM_\lambda}(M_1, \CPM(M_2, N)).
  \]
\end{proposition}

\begin{proof}
  Same as \cref{prop:lpr-lam-tensor-hom}.
\end{proof}

Hence $\CPM(M, N)$ is also an ``internal Hom'' for $\CPM_\lambda$ when $M$ is $\lambda$-small $\lambda$-combinatorial.

\section{$\CPM$ as a 2-multicategory}\label{sec:cpm-multicat}

We showed that the tensor product $\otimes : \CPM \times \CPM \to \CPM$ is unital and associative up to equivalence.
It is also obviously symmetric up to equivalence, since giving a Quillen bifunctor $F : M_1 \times M_2 \to N$ is the same as giving a Quillen bifunctor $F' : M_2 \times M_1 \to N$.
We would like to assert that $\otimes$ makes $\CPM$ into a symmetric monoidal 2-category.
However, the explicit definition of a symmetric monoidal 2-category (as presented for example in \cite{Stay}) is impractical for our purposes.

Instead, we will sketch the construction of a (strict) \emph{symmetric 2-multicategory} (or $\Cat$-valued operad) structure on $\CPM$.
For the definition of a symmetric 2-multicategory, see \cite[Definition~1.1]{Hormann}.
In short, we must give, for every $n \ge 0$, $n$-tuple of ``input'' objects $(M_1, \ldots, M_n)$ and ``output'' object $N$, a category $\THom_\CPM((M_1, \ldots, M_n), N)$ of $n$-ary morphisms from $(M_1, \ldots, M_n)$ to $N$.
These are related by composition functors, which are associative and unital.
The symmetric groups act on $\THom_\CPM((M_1, \ldots, M_n), N)$ by permuting the $M_i$ in a way compatible with composition.
This structure is \emph{strict}: all the expected equations between functors built from composition, units and the symmetry hold on the nose, not up to isomorphism.
For example, $\TCat$ is a symmetric 2-multicategory in which $\THom_\TCat((M_1, \ldots, M_n), N) = \TCat(M_1 \times \cdots \times M_n, N)$.

\begin{definition}
  Let $M_1$, \ldots, $M_n$ and $N$ be combinatorial premodel categories.
  We call a functor $F : M_1 \times \cdots \times M_n \to N$ a \emph{Quillen multifunctor} if:
  \begin{enumerate}
  \item
    $F$ preserves colimits in each argument.
  \item
    For any cofibrations $f_1 : A_1 \to B_1$ of $M_1$, \ldots, $f_n : A_n \to B_n$ of $M_n$, let $Q : [1]^n \to N$ be the cubical diagram built out of $F$ and the $f_i$.
    Write $\top$ for the terminal object of $Q$.
    Then the morphism $\colim_{q \ne \top} Q(q) \to Q(\top)$ is a cofibration, which is anodyne if any of $f_i$ are.
  \end{enumerate}
\end{definition}

For $n = 2$ this is the definition of a Quillen bifunctor, except that we have used the adjoint functor theorem in condition (1).
Similarly, for $n = 1$ this is the definition of a Quillen functor up to the adjoint functor theorem.
For $n = 0$ the diagram $Q$ consists just one object $F()$, and condition (2) requires that this object be cofibrant.

We define $\THom_\CPM((M_1, \ldots, M_n), N)$ to be the full subcategory of $\TCat(M_1 \times \cdots \times M_n, N)$ consisting the Quillen multifunctors.
The composition and symmetric action is inherited from the symmetric 2-multicategory structure of $\TCat$ (one must check that the composition of Quillen multifunctors is again a Quillen multifunctor).
By an argument similar to the binary case (\cref{prop:cpm-tensor}), any $n$-tuple $(M_1, \ldots, M_n)$ admits an $n$-fold tensor product $M_1 \otimes \cdots \otimes M_n$ equipped with a universal (up to equivalence) Quillen multifunctor from $M_1 \times \cdots \times M_n$.
One can check (for example using the formula for the tensor product) that these universal Quillen multifunctors induce equivalences
\[
  \begin{tikzcd}
    \THom_\CPM((M'_1, \ldots, M'_k, M_1 \otimes \cdots \otimes M_n, M''_1, \ldots, M''_l), N) \ar[d, "\eqv"] \\
    \THom_\CPM((M'_1, \ldots, M'_k, M_1, \ldots, M_n, M''_1, \ldots, M''_l), N)
  \end{tikzcd}
\]
for any $M'_1$, \ldots, $M'_k$, $M''_1$, \ldots, $M''_l$, $N$.

\chapter{Directed colimits in $\CPM_\lambda$}
\label{chap:size}

We would like to construct factorizations in $\CPM$ (and later $V\CPM$) using the small object argument.
However, there are two obstructions to doing so directly.
\begin{enumerate}
\item
  The Hom categories of $\CPM$ are not essentially small.
  This means that at each successor step in the small object argument, there is a proper class of lifting problems to solve and so we cannot simply adjoin a solution to each of these for set-theoretic reasons.
\item
  Less obviously, most objects of $\CPM$ are not $\mu$-small (in the sense required for the small object argument) for any $\mu$.

  For example, take the object $\Set_\init$.
  The functor $\THom_\CPM(\Set_\init, -) : \CPM \to \TCat$ sends a combinatorial premodel category $M$ to its underlying category.
  Let $\mu$ be a regular cardinal and define a diagram in $\CPM$ by $M_\alpha = \coprod_{i < \alpha} \Set$ for each $\alpha \le \mu$, with morphisms $M_\alpha \to M_\beta$ induced by the inclusion of index sets.
  Then in $\CPM$, $M_\mu = \colim_{\alpha < \mu} M_\alpha$.
  Concretely, the underlying category of $M_\alpha$ is the product $\prod_{i < \alpha} \Set$ and each functor $M_\alpha \to M_\beta$ is given by ``extension by $\emptyset$''.
  In particular any object $X$ of $M_\mu$ in the image of $M_\alpha \to M_\mu$ has $X_i = \emptyset$ for all $\alpha \le i < \mu$.
  Hence the object $Y$ with $Y_i = *$ for all $i$ is not in the image of any $M_\alpha \to M_\mu$ and so $M_\mu$ is not the colimit of the $M_\alpha$ in $\TCat$.

  Informally, the problem is that $\CPM$ has ``operations'' $\colim_I$ of arbitrary ``arity'' $|I|$, which can be used to build the object $Y$ of $M_\mu$ out of objects coming from $M_\alpha$ for \emph{every} $\alpha < \mu$, regardless of the size of $\mu$.
  Then $Y$ does not exist at any stage $M_\alpha$ for $\alpha < \mu$.
\end{enumerate}
These difficulties do not apply to the subcategory $\CPM_\lambda$, however.
\begin{enumerate}
\item
  The Hom categories of $\CPM_\lambda$ are essentially small.
  (A strongly $\lambda$-accessible functor $F : M \to N$ between locally $\lambda$-presentable categories is determined up to isomorphism by its restriction $F^\lambda : M^\lambda \to N^\lambda$, and $M^\lambda$ and $N^\lambda$ are essentially small.)
\item
  For every object $M$ of $\CPM_\lambda$, $\THom_{\CPM_\lambda}(M, -) : \CPM_\lambda \to \TCat$ preserves $\mu$-directed colimits for sufficiently large $\mu$.
  For $M = \Set$, for example, we can take $\mu = \lambda$.
  The counterexample above does not work because $\THom_{\CPM_\lambda}(\Set_\init, -)$ sends a $\lambda$-combinatorial premodel category to its category of $\lambda$-compact objects and the object $Y$ defined earlier is not $\lambda$-compact.
\end{enumerate}
Proving statement (2) is the goal of this chapter.
More specifically we will show that for any $\mu$-small $\lambda$-combinatorial premodel category $M$, $\THom_{\CPM_\lambda}(M, -) : \CPM_\lambda \to \TCat$ preserves $\mu$-directed colimits.

In \cref{chap:small} we will explain how to reduce the construction of factorizations in $\CPM$ to the construction of factorizations in $\CPM_\lambda$ under favorable conditions on the generating morphisms.

\begin{remark}
  By ``$\mu$-directed colimits'', we mean conical colimits indexed on an ordinary $\mu$-directed poset.
  Our arguments would also apply to $\mu$-filtered colimits but at the cost of more complicated notation.
  In actuality the only case we will really need is that of colimits indexed on an ordinal, for use in the small object argument.
\end{remark}

\section{Short cell complexes}

Let $M$ be a locally $\lambda$-presentable category equipped with a weak factorization system $(\sL, \sR)$ which is generated by a set $I$ of morphisms between $\lambda$-compact objects.
Suppose $i : A \to B$ is a morphism belonging to $\sL$.
By the small object argument, we may factor $i$ as an $I$-cell complex $i' : A \to B'$ followed by a morphism of $\sR$.
Then the retract argument shows that $i$ is a retract of $i'$ in the slice category $M_{A/}$.
\begin{equation*}
  \begin{tikzcd}
    A \ar[r, equals] \ar[d, "i"'] & A \ar[r, equals] \ar[d, "i'"] & A \ar[d, "i"] \\
    B \ar[r] & B' \ar[r] & B
  \end{tikzcd}
  \tag{$*$}
\end{equation*}

The cell complexes produced by the small object argument are generally not very efficient for this purpose.
When the objects $A$ and $B$ are $\lambda$-compact, we can hope to choose the $I$-cell complex $i'$ to have fewer than $\lambda$ cells.
This is possible using the \emph{fat small object argument} of \cite{FSOA}, a version of the ``good trees'' of \cite[section~A.1.5]{HTT}.

\begin{proposition}\label{prop:short-retract}
  With notation as above, suppose that $i : A \to B$ is a morphism between $\lambda$-compact objects.
  Then we can write $i$ as a retract of $i' : A \to B'$ as shown in $(*)$ for an $I$-cell complex $i'$ which has fewer than $\lambda$ cells.
  (That is, $i'$ is a transfinite composition of length less than $\lambda$ of pushouts of morphisms of $I$.)
\end{proposition}

\begin{proof}
  We apply \cite[Lemma~4.19]{FSOA} with $\mathcal{X} = I$.
  The morphism $i$ is a retract of an $I$-cell complex so it belongs to $\cof_\lambda\,\mathcal{X}$.
  Then \cite[Lemma~4.19]{FSOA} states that $i$ can be written as a retract of a transfinite composition of length less than $\lambda$ of pushouts of morphisms of $I$.
  By inspection of the proof this retract can be taken to be a retract in $M_{A/}$ as shown in $(*)$.
\end{proof}

The proof of \cite[Lemma~4.19]{FSOA} relies on the following fact which we will also use separately.

\begin{proposition}\label{prop:cell-directed-colimit}
  Let $i : A \to B$ be an $I$-cell complex (with no assumptions on the objects $A$ and $B$).
  Then $i$ can be expressed as a $\lambda$-directed colimit in $M_{A/}$ of morphisms $i_\alpha : A \to B_\alpha$ each of which is an $I$-cell complex with fewer than $\lambda$ cells.
\end{proposition}

\begin{proof}
  By \cite[Theorem~4.11]{FSOA} $i$ is the composite of a $\lambda$-good $\lambda$-directed diagram $D : P \to M$ with links that are pushouts of morphisms of $I$.
  For each $\alpha$ in $P$ let $i_\alpha : A \to B_\alpha$ be the morphism $D(\bot) \to D(\alpha)$.
  Since $P$ is $\lambda$-good, each $i_\alpha$ is the composition of a $\lambda$-good sub-diagram of $P$ of cardinality less than $\lambda$.
  Then by \cite[Remark~4.6]{FSOA}, $i_\alpha$ can also be written as a transfinite composition of length less than $\lambda$ of pushouts of morphisms of $I$.
\end{proof}

\section{Recognizing directed colimits of categories}

Let $D$ be a directed poset and let $(C_i)_{i \in D} \to C$ be a cocone in $\TCat$ indexed on $D$.
We write $F^i_j : C_i \to C_j$ (for each $i \le j$) for the functors that make up the original diagram $(C_i)_{i \in D}$ and $F^j : C_j \to C$ for the cocone functors.

\begin{notation}
  Throughout this chapter we will work as though all the cocones of categories we consider are \emph{strict}.
  For example, here we will assume $F^i_i = \id_{C_i}$, $F^i_k = F^j_k \circ F^i_j$ and $F^i = F^j \circ F^i_j$.
  We can justify this using standard strictification results for pseudofunctors, since all the notions we consider are equivalence-invariant.
  However, this strictification is really just a notational device to avoid writing an excess of coherence isomorphisms.
\end{notation}

Because the theory of categories is finitary, we can detect directed colimits of categories as described below.

\begin{proposition}\label{prop:dcolimit-cat}
  The cocone $(C_i)_{i \in D} \to C$ is a colimit in $\TCat$ if and only if:
  \begin{enumerate}
  \item[(a)]
    Every object of $C$ is isomorphic to the image of some object of some $C_i$ under the functor $F^i : C_i \to C$.
  \item[(b)]
    For every $i \in D$ and every pair of objects $A_i$ and $B_i$ of $C_i$, the cocone
    \[
      (\Hom_{C_j}(F^i_j A_i, F^i_j B_i))_{j \ge i} \to \Hom_C(F^i A_i, F^i B_i)
    \]
    is a colimit of sets.
    
    In turn, this means that
    \begin{enumerate}
    \item[{\rm(b1)}] every morphism $f : F^i A_i \to F^i B_i$ is the image under $F^j$ of some $f_j : F^i_j A_i \to F^i_j B_j$, and
    \item[{\rm(b2)}] every equality $F^i f = F^i g$ for $f$, $g : A_i \to B_i$ holds already at some stage, that is, there exists $j \ge i$ such that $F^i_j f = F^i_j g$.
    \end{enumerate}
  \end{enumerate}
\end{proposition}

In practice we will mainly check a variant of this condition.

\begin{definition}
  A category $C$ is \emph{idempotent complete} if for any object $X$ of $C$ and any idempotent $e : X \to X$ (so $e^2 = e$), there exists a retract $X \xrightarrow{r} Y \xrightarrow{i} X$ (so $ri = \id_Y$) with $ir = e$.
  We say the object $Y$ (as a retract of $X$) \emph{splits} the idempotent $e$.
\end{definition}

\begin{remark}\label{remark:splitting-unique}
  In this situation, the diagram
  \[
    \begin{tikzcd}
      X \ar[r, shift left, "e"] \ar[r, shift right, "\id_X"'] & X \ar[r, "r"] & Y
    \end{tikzcd}
  \]
  is a coequalizer: given a map $f : X \to Z$ such that $f = fe$, the map $fi : Y \to Z$ is the unique map whose composition with $r$ is $f$.
  Conversely, any coequalizer $r : X \to Y$ of $e$ and $\id_X$ yields a retract of $X$ splitting $e$; the map $i : Y \to X$ is the one obtained by applying the universal property of the coequalizer to $e : X \to X$.
  It follows that the data $(Y, r, i)$ is determined by $X$ and $e$ up to unique isomorphism.

  In particular, suppose $F : C \to C'$ is any functor with $C$ idempotent complete, and suppose that the object $Y'$ of $C'$ is a retract of an object $X'$ with associated idempotent $e' : X' \to X'$.
  Suppose that (up to isomorphism) $X'$ and $e' : X' \to X'$ are the image under $F$ of some object $X$ of $C$ and some idempotent $e : X \to X$.
  Let $Y$ be the retract of $X$ which splits $e$.
  Then $FY$ is also a retract of $FX$ splitting $Fe$, and so $FY$ is isomorphic to the original object $Y'$.
\end{remark}

\begin{remark}
  Let $M$ and $N$ be objects of $\CPM_\lambda$.
  Then $\THom_{\CPM_\lambda}(M, N)$ is idempotent complete.
  Indeed, by \cref{prop:cpm-internal-hom}, there is a premodel category structure $\CPM(M, N)$ on the category of all left adjoints from $M$ to $N$ for which $\THom_{\CPM_\lambda}(M, N)$ can be identified with the full subcategory on those objects which both are cofibrant (i.e., are left Quillen functors) and preserve $\lambda$-compact objects when viewed as functors from $M$ to $N$.
  The category $\CPM(M, N)$ is cocomplete, hence in particular idempotent complete, and both of the above conditions are preserved under passage to retracts.
\end{remark}

\begin{proposition}\label{prop:dcolimit-cat-retract}
  In \cref{prop:dcolimit-cat}, suppose that each of the categories $C_i$ is idempotent complete.
  Then (a) can be replaced by
  \begin{enumerate}
  \item[(a$'$)] every object of $C$ is a retract of the image of some object of some $C_i$ under the functor $F^i : C_i \to C$.
  \end{enumerate}
\end{proposition}

\begin{proof}
  We must check that conditions (a$'$) and (b) imply the original condition (a).
  Let $X$ be any object of $C$.
  By (a$'$), $X$ is a retract of $F^i X_i$ for some $i \in D$ and some object $X_i$ of $C_i$.
  Let $e : X \to X$ be the associated idempotent.
  By (b1), we can find $j \ge i$ and a morphism $e_j : F^i_j X_i \to F^i_j X_i$ whose image under $F^j : C_j \to C$ equals $e$.
  The morphism $e_j$ might not satisfy $e_j^2 = e_j$, but we know that $F^j(e_j^2) = (F^j e_j)^2 = e^2 = e = F^j e_j$ and so, by (b2), we can find $k \ge j$ such that $F^j_k (e_j^2) = F^j_k e_j$.
  Set $e_k = F^j_k e_j : F^i_k X_i \to F^i_k X_i$; then $e_k$ is an idempotent with $F^k e_k = F^k F^j_k e_j = F^j e_j = e$.
  Since $C_k$ is idempotent complete, there is a retract $Y_k$ of $F^i_k X_i$ splitting $e_k$.
  Then $F^k : C_k \to C$ sends $Y_k$ to a retract of $F^i X_i$ splitting $e$, hence (by \cref{remark:splitting-unique}) to an object isomorphic to the original object $X$.
  Thus, we have verified condition (a) and hence can apply \cref{prop:dcolimit-cat}.
\end{proof}

We will make use of the following arguments several times.
They give conditions under which a diagram or a cell complex in $C$ factors through $C_j$ for some $j$.

\begin{lemma}\label{lemma:small-diagram-factors}
  Suppose $(C_i)_{i \in D} \to C$ is a cocone in $\TCat$ indexed on a $\lambda$-directed poset $D$ which satisfies condition (b) of \cref{prop:dcolimit-cat}.
  Let $X : K \to C$ be a diagram in $C$ indexed on a $\lambda$-compact category $K$ such that, for each $k \in K$, the object $X(k)$ of $C$ is in the essential image of $F^i : C_i \to C$ for some $i \in D$.
  Then the entire diagram $X : K \to C$ factors (up to isomorphism) through $F^j : C_j \to C$ for some $j \in D$.
\end{lemma}

Recall that a category is $\lambda$-compact if and only if it admits a $\lambda$-small presentation, that is, one by fewer than $\lambda$ objects, fewer than $\lambda$ morphisms and fewer than $\lambda$ relations.
In particular, any $\lambda$-small category is $\lambda$-compact.

\begin{proof}
  First, we may replace each object of the diagram $X$ by one which is in the actual image of $F^i : C_i \to C$ for some $i$.
  We will show the resulting functor $X : K \to C$ factors on the nose through $F^j : C_j \to C$ for some $j$.
  Pick a $\lambda$-small presentation for $K$.
  Since $K$ has fewer than $\lambda$ objects and $D$ is $\lambda$-directed, we can choose a single $i_1 \in D$ such that every object $X(k)$ belongs to the image of $F^{i_1} : C_{i_1} \to C$.
  Choose some $X_{i_1} : \Ob K \to C_{i_1}$ so that $F^{i_1} \circ X_{i_1}$ is the object part of $X$.
  Then by condition (b1), for each generating morphism $g : k \to l$ of $K$, we may find a $j \ge i_1$ such that the morphism $X(g) : X(k) \to X(l)$ of $C$ is the image of some $f_j : F^{i_1}_j X_{i_1}(k) \to F^{i_1}_j X_{i_1}(l)$ under $F^j : C_j \to C$.
  Since $K$ has fewer than $\lambda$ generating morphisms, we can choose a single $i_2 \in D$ such that this condition holds for $j = i_2$ for every generating morphism $g$.
  Then, for each generating morphism $g : k \to l$, choose a morphism $F^{i_1}_{i_2} X_{i_1}(k) \to F^{i_1}_{i_2} X_{i_1}(l)$ sent by $F^{i_2}$ to $X(g) : X(k) \to X(l)$.
  These morphisms collectively define a functor $X_{i_2}$ from the free category $\hat K$ on the generating morphisms of $K$ to $C_{i_2}$ whose image under $F^{i_2} : C_{i_2} \to C$ agrees with $\hat K \to K \xrightarrow{X} C$; they may not yet satisfy the required relations to define a functor from $K$ itself.
  However, each of these relations holds for $X : K \to C$ and so, by condition (b2), also holds for $F^{i_2}_j X_{i_2} : \hat K \to C_j$ for some $j \ge i_2$.
  Since the presentation for $K$ has fewer than $\lambda$ relations, we may find a single $i_3 \in D$ such that all of the relations hold for $F^{i_2}_{i_3} X_{i_2} : \hat K \to C_{i_3}$.
  Set $j = i_3$ and let $X_{i_3} : K \to C_{i_3}$ be the functor induced by $F^{i_2}_{i_3} X_{i_2} : \hat K \to C_{i_3}$; then $F^{i_3} \circ X_{i_3} = X$ as desired.
\end{proof}

\begin{lemma}\label{lemma:cell-complex-factors}
  Suppose $(C_i)_{i \in D} \to C$ is a cocone of cocomplete categories and colimit-preserving functors indexed on a $\lambda$-directed poset $D$ and, for each $i$, $I_i$ is a set of morphisms between $\lambda$-compact objects of $C_i$ such that $F^i_j I_i \subset I_j$ for each $i \le j$.
  Set $I = \bigcup_{i \in D} F^i I_i$.
  Assume that the cocone $(C_i)_{i \in D} \to C$ satisfies the following version of condition (b):
  \begin{enumerate}
  \item[{\rm(b$'$)}] For every $i \in D$ and every pair of objects $A_i$ and $B_i$ of $C_i$ with $A_i$ $\lambda$-compact, the cocone
    \[
      (\Hom_{C_j}(F^i_j A_i, F^i_j B_i))_{j \ge i} \to \Hom_C(F^i A_i, F^i B_i)
    \]
    is a colimit of sets.
  \end{enumerate}
  Suppose $f : F^i A_i \to B_i$ is an $I$-cell complex containing fewer than $\lambda$ cells.
  Then there exists $j \ge i$ and an $I_j$-cell complex $f_j : F^i_j A_i \to B_j$ of the same length as $f$ with $F^j f_j = f$ (up to isomorphism).
\end{lemma}

\begin{proof}
  Write $f$ as the transfinite composition $X : \gamma \to C$ of morphisms of $I$, with $\gamma < \lambda$.
  By transfinite induction on $\beta \le \gamma$ we construct $j_\beta \in D$, $X_\beta \in \Ob C_{j_\beta}$ and morphisms $\colim_{\alpha < \beta} F^{j_\alpha}_{j_\beta} X_\alpha \to X_\beta$ such that the diagram $(F^{j_\alpha}_{j_\beta} X_\alpha)_{\alpha \le \beta}$ is an $I_{j_\beta}$-cell complex whose image under $F^{j_\beta}$ is isomorphic to the restriction of $X$ to $\beta$.
  \begin{itemize}
  \item
    For $\beta = 0$ we set $j_0 = i$ and $X_0 = A_i$.
  \item
    At a limit step, we set $j_\beta$ to be some upper bound in $D$ for $\{\,j_\alpha \mid \alpha < \beta\,\}$.
    Such an upper bound exists because $D$ is $\lambda$-directed and $\beta \le \gamma < \lambda$.
    Then we define $X_\beta = \colim_{\alpha < \beta} F^{j_\alpha}_{j_\beta} X_\beta$.
  \item
    At a successor step $\beta = \alpha+1$, the morphism $X(\alpha) \to X(\beta)$ is a pushout of some morphism $A \to B$ of $I$, which in turn is the image of some morphism $A_k \to B_k$ of $I_k$ under $F^k$.
    Increasing $k$ if necessary we may assume that $k \ge j_\alpha$.
    \[
      \begin{tikzcd}
        F^k A_k \ar[r] \ar[d, "\varphi"'] & F^k B_k \ar[d] \\
        X(\alpha) \ar[r] & X(\beta)
      \end{tikzcd}
    \]
    Up to isomorphism, $X(\alpha) = F^{j_\alpha} X_\alpha = F^k (F^{j_\alpha}_k X_\alpha)$.
    The object $A_k$ is $\lambda$-compact so we may apply condition (b$'$) to obtain $j_\beta$ and a map $\varphi_{j_\beta} : F^k_{j_\beta} A_k \to F^{j_\alpha}_{j_\beta} X_\alpha$ whose image under $F^{j_\beta}$ is the attaching map $\varphi : F^k A_k \to X(\alpha)$.
    Then we define $X_\beta$ by forming the pushout
    \[
      \begin{tikzcd}
        F^k_{j_\beta} A_k \ar[r] \ar[d, "\varphi_{j_\beta}"'] & F^k_{j_\beta} B_k \ar[d] \\
        F^{j_\alpha}_{j_\beta} X_\alpha \ar[r] & X_\beta
      \end{tikzcd}
    \]
    in $C_{j_\beta}$.
    The resulting morphism $F^{j_\alpha}_{j_\beta} X_\alpha \to X_\beta$ is a pushout of a morphism of $I_{j_\beta}$, and it is sent by $F^{j_\beta}$ to $X(\alpha) \to X(\beta)$ because $F$ preserves pushouts.
  \end{itemize}
  Taking $\beta = \gamma$ and $j = j_\gamma$ produces the desired $I_j$-cell complex $f_j$.
\end{proof}

\section{\texorpdfstring{$\lambda$}{\textlambda}-directed colimits in \texorpdfstring{$\LPr_\lambda$}{LPr\_\textlambda}}

Our goal in this section is to prove that the ``$\lambda$-compact objects'' functor $(-)_\lambda : \LPr_\lambda \to \TCat$ preserves $\lambda$-directed colimits.
We will build up to this result in a series of steps.

\subsection{$\lambda$-directed colimits of presheaf categories}

Let $D$ be a $\lambda$-directed poset and let $(A_i)_{i \in D} \to A$ be a colimit cocone in $\TCat$ indexed on $D$.
We assume as before that the diagram consisting of $A$ and the $A_i$ is a strict functor to $\TCat$ and write $J^i_j : A_i \to A_j$ and $J^i : A_i \to A$ for the functors which make up this diagram.

The pseudofunctor $\PP : \TCat \to \LPr$ takes this cocone to a cocone $(\PP(A_i))_{i \in D} \to \PP(A)$ whose functors $F^i_j : \PP(A_i) \to \PP(A_j)$ and $F^i : \PP(A_i) \to \PP(A)$ are given by left Kan extension along $J^i_j$ and $J^i$ respectively.
Since restriction along the Yoneda embedding $\yo : A \to \PP(A)$ induces an equivalence $\LPr(\PP(A), N) \eqv \TCat(A, N)$ for any $A$, this resulting cocone is again a colimit.
By replacing the categories $\PP(A_i)$ and $\PP(A)$ by equivalent ones, we will again assume that the $F^i_j$ and $F^i$ make up a strict functor to $\LPr$.

The categories $\PP(A_i)$ and $\PP(A)$ and the functors $F^i_j$ and $F^i$ also belong to $\LPr_\lambda$ for any regular cardinal $\lambda$, and $(\PP(A_i))_{i \in D} \to \PP(A)$ is also a colimit cocone in $\LPr_\lambda$.
In particular, we can apply $(-)_\lambda : \LPr_\lambda \to \TCat$ to obtain another (strict) cocone in $\TCat$.
We will show that it is a colimit by applying \cref{prop:dcolimit-cat}.
We first establish a slightly stronger version of the condition (b).

\begin{proposition}\label{prop:presheaf-b}
  Let $X_i$ and $Y_i$ be objects of $\PP(A_i)$ such that $X_i$ is $\lambda$-compact.
  Then
  \[
    (\Hom_{\PP(A_j)}(F^i_j X_i, F^i_j Y_i))_{j \ge i} \to \Hom_{\PP(A)}(F^i X_i, F^i Y_i)\]
  is a colimit cocone.
\end{proposition}

\begin{proof}
  Let $\mathcal{S}$ denote the class of objects $X_i$ of $C_i$ for which the conclusion holds for every object $Y_i$ of $\PP(A_i)$.
  We first show that $\mathcal{S}$ is closed under $\lambda$-small colimits.
  Indeed, each of $\Hom_{\PP(A_j)}(F^i_j -, F^i_j Y_i)$ and $\Hom_{\PP(A)}(F^i -, F^i_j Y_i)$ sends colimits to limits, because $F^i_j$ and $F^i$ are left adjoints.
  Moreover, colimits over the $\lambda$-directed poset $\{\,j \in D \mid j \ge i\,\}$ preserve $\lambda$-small limits, and so the conclusion is preserved under $\lambda$-small colimits.
  Thus, it suffices to show that the representable objects belong to $\mathcal{S}$.
  Hence we may assume without loss of generality that $X_i$ is representable.

  Next, the functors $F^i_j$ and $F^i$ are defined as left Kan extensions (of $J^i_j$ and $J^i$ respectively) and so preserve representables.
  It follows that the functors $\colim_{j \ge i} \Hom_{\PP(A_j)}(F^i_j X_i, F^i_j -)$ and $\Hom_{\PP(A)}(F^i X_i, F^i -)$ preserve \emph{all} colimits.
  Since every object $Y_i$ of $\PP(A_i)$ is a colimit of representables, we may also assume without loss of generality that $Y_i$ is representable.

  Thus, suppose $X_i = \yo x_i$ and $Y_i = \yo y_i$ for objects $x_i$ and $y_i$ of $A_i$.
  Then we can compute $F^i_j X_i = F^i_j (\yo x_i) = \yo (J^i_j x_i)$ and $F^i_j Y_i = \yo (J^i_j y_i)$, so by the Yoneda lemma
  \[
    \Hom_{\PP(A_j)}(F^i_j X_i, F^i_j Y_i) =
    \Hom_{\PP(A_j)}(\yo (J^i_j x_i), \yo (J^i_j y_i)) =
    \Hom_{A_j}(J^i_j x_i, J^i_j y_i).
  \]
  By a similar argument $\Hom_{\PP(A)}(F^i X_i, F^i Y_i) = \Hom_{A}(J^i x_i, J^i y_i)$.
  Thus it remains only to check that $(\Hom_{A_j}(J^i_j x_i, J^i_j y_i))_{j \ge i} \to \Hom_A(J^i x_i, J^i y_i)$ is a colimit cocone; and this follows from (the ``only if'' direction of) \cref{prop:dcolimit-cat}.
\end{proof}

\begin{proposition}\label{prop:presheaf-colimit}
  The cocone $(\PP(A_i)_\lambda)_{i \in D} \to \PP(A)_\lambda$ is a colimit in $\TCat$.
\end{proposition}

\begin{proof}
  We apply \cref{prop:dcolimit-cat}.
  We checked condition (b) above, so it remains to check condition (a).
  Let $\mathcal{S}$ denote the full subcategory of $\PP(A)_\lambda$ on the objects which are in the essential image of $F^i : \PP(A_i)_\lambda \to \PP(A)_\lambda$ for some $i \in D$.
  Then $\mathcal{S}$ contains the representables, by applying \cref{prop:dcolimit-cat} to the original colimit cone $(A_i)_{i \in D} \to A$.
  It remains, then, to verify that $\mathcal{S}$ is closed under $\lambda$-small colimits.

  Suppose $X : K \to \mathcal{S} \subset \PP(A)_\lambda$ is a $\lambda$-small diagram.
  We need to show that $\colim_K X$ again belongs to $\mathcal{S}$.
  We can apply \cref{lemma:small-diagram-factors} since the condition on objects holds by the definition of $\mathcal{S}$, and we already checked condition (b).
  Therefore $X$ can be expressed (up to isomorphism) as the image of a diagram $X_j : K \to \PP(A_j)_\lambda$ under $F^j : \PP(A_j)_\lambda \to \PP(A)_\lambda$ for some $j$.
  The category $\PP(A_j)_\lambda$ has $\lambda$-small colimits and these are preserved by $F^j$, so $\colim_K X$ is isomorphic to $F^j(\colim_K X_j)$ and hence belongs to $\mathcal{S}$.
\end{proof}

\subsection{General $\lambda$-directed colimits in $\LPr_\lambda$}

Now suppose that, in addition to the $\lambda$-directed colimit cocone $(A_i)_{i \in D} \to A$ in $\TCat$, we are given, for each $i \in D$, a set of morphisms $E_i$ between $\lambda$-compact objects of $\PP(A_i)$.
Then, for each $i$, the orthogonal subcategory $\OO(A_i, E_i) \subset \PP(A_i)$ is a locally $\lambda$-presentable category and the reflector $L_i : \PP(A_i) \to \OO(A_i, E_i)$ is a left adjoint which preserves $\lambda$-compact objects, i.e., a morphism of $\LPr_\lambda$.
Suppose furthermore that for each $i \le j$ in $D$, the left Kan extension $F^i_j : \PP(A_i) \to \PP(A_j)$ of each $J^i_j : A_i \to A_j$ fits into a square as shown below, for morphisms $G^i_j$ of $\LPr_\lambda$.
\[
  \begin{tikzcd}
    \PP(A_i) \ar[r, "F^i_j"] \ar[d, "L_i"'] & \PP(A_j) \ar[d, "L_j"] \\
    \OO(A_i, E_i) \ar[r, "G^i_j"'] & \OO(A_j, E_j)
  \end{tikzcd}
\]
By \cref{prop:out-of-perp}, this means that each map of $E_i$ is sent by $L_j \circ F^i_j$ to an isomorphism and so the subcategory $\OO(A_j, E_j)$ is not affected by enlarging $E_j$ to include the images of the maps $E_i$ under $F^i_j$.
Thus, by replacing $E_j$ by $\bigcup_{i \le j} F^i_j E_i$ for each $j$, we leave the subcategories $\OO(A_j, E_j)$ unmodified while ensuring that $F^i_j E_i \subset E_j$ for each $i \le j$ in $D$.
We will henceforth assume that this last condition is satisfied.

Define $E = \bigcup_i F^i E_i$; then $E$ is a set of morphisms between $\lambda$-compact objects of $\PP(A)$.
Write $L : \PP(A) \to \OO(A, E)$ for the reflector; $L$ is a morphism of $\LPr_\lambda$.
For each $i$, the composition $\PP(A_i) \xrightarrow{F^i} \PP(A) \xrightarrow{L} \OO(A, E)$ factors through $L_i : \PP(A_i) \to \OO(A_i, E_i)$, inducing a morphism $G^i : \OO(A_i, E_i) \to \OO(A, E)$.
As usual, we will assume for notational convenience that the functors $G^i_j$ and $G^i$ form a strict diagram $(\OO(A_i, E_i))_{i \in D} \to \OO(A, E)$.
By \cref{prop:colimit-of-o-form}, this cocone is a colimit in $\LPr_\lambda$ (and in $\LPr$).

Our goal is to show that $(-)_\lambda : \LPr_\lambda \to \TCat$ preserves this colimit.
As in the previous subsection, we will first verify a stronger version of condition (b).
We write $U_i : \OO(A_i, E_i) \to \PP(A_i)$ and $U : \OO(A, E) \to \PP(A)$ for the inclusions, which are right adjoint to $L_i$ and $L$.

\begin{proposition}\label{prop:orthogonal-colimit-b}
  Let $X_i$ and $Y_i$ be objects of $\OO(A_i, E_i)$ such that $X_i$ is $\lambda$-compact.
  Then
  \[
    (\Hom_{\OO(A_j, E_j)}(G^i_j X_i, G^i_j Y_i))_{j \ge i} \to
    \Hom_{\OO(A, E)}(G^i X_i, G^i Y_i)
  \]
  is a colimit cocone.
\end{proposition}

\begin{proof}
  Let $\mathcal{S}$ denote the class of objects $X_i$ of $\OO(A_i, E_i)$ for which the conclusion holds for every object $Y_i$ of $\OO(A_i, E_i)$.
  As in the proof of \cref{prop:presheaf-b}, the class $\mathcal{S}$ is closed under $\lambda$-small colimits, hence in particular under retracts.
  Let $X_i$ be any $\lambda$-compact object of $\OO(A_i, E_i)$.
  Since $\OO(A_i, E_i)$ is a reflective subcategory of $\PP(A_i)$, $X_i$ is isomorphic to $L_i(U_iX_i)$.
  Now $U_iX_i$ can be written as a $\lambda$-filtered colimit of $\lambda$-compact objects of $\PP(A_i)$; this colimit is preserved by $L_i$ and $X_i$ is $\lambda$-compact, so $X_i$ is a retract of $L_iX'_i$ for some $\lambda$-compact object $X'_i$ of $\PP(A_i)$.
  It suffices, then, to check that objects of the form $L_iX'_i$ for $\lambda$-compact $X'_i$ belong to $\mathcal{S}$, so set $X_i = L_iX'_i$.
  We may also write $Y_i = L_iY'_i$, where $Y'_i = U_iY_i$.
  The cocone in question then has the form
  \[
    (\Hom_{\OO(A_j, E_j)}(L_j F^i_j X'_i, L_j F^i_j Y'_i))_{j \ge i} \to
    \Hom_{\OO(A, E)}(L F^i X'_i, L F^i Y'_i)
  \]
  or, using the adjunction between the reflectors and inclusions and writing $R_j = U_j L_j$ and $R = U L$,
  \begin{equation}
    (\Hom_{\PP(A_j)}(F^i_j X'_i, R_j F^i_j Y'_i))_{j \ge i} \to
    \Hom_{\PP(A)}(F^i X'_i, R F^i Y'_i).
    \tag{$*$}
  \end{equation}

  Suppose first that, for some $j \ge i$, the morphisms $f_j : F^i_j X'_i \to R_j F^i_j Y'_i$ and $g_j : F^i_j X'_i \to R_j F^i_j Y'_i$ have the same image in $\Hom_{\PP(A)}(F^i X'_i, R F^i Y'_i)$.
  This means that $F^j f_j$ and $F^j g_j$ become equal after being composed with $\eta : F^j R_j F^i_j Y'_i \to R F^j R_j F^i_j Y'_i$.
  Now by \cref{prop:cell-directed-colimit}, $\eta$ can be written as a $\lambda$-filtered colimit of morphisms with domain $F^j R_j F^i_j Y'_i$, each of which is an $E^+$-cell complex of length less than $\lambda$.
  Since $F^j f_j$ and $F^j g_j$ have $\lambda$-compact domain $F^j F^i_j X'_i$, they agree after composing with one of these morphisms $h$.
  By \cref{lemma:cell-complex-factors}, we may find, for some $k \ge j$, an $E_k^+$-cell complex $h_k$ in $\PP(A_k)$ with domain $F^j_k R_j F^i_j Y'_j$ which is carried by $F^k$ onto the $E^+$-cell complex $h$ (up to isomorphism).
  Thus, the morphisms $h_k \circ F^j_k f_j$ and $h_k \circ F^j_k g_j$ of $\PP(A_k)$ are sent by $F^k$ to equal morphisms of $\PP(A)$.
  By \cref{prop:presheaf-b} again, there exists $l \ge k$ such that $F^k_l$ already sends them to equal morphisms, that is, so that $F^k_l h_k \circ F^j_l f_j = F^k_l h_k \circ F^j_l g_j$.
  Now $F^k_l h_k$ is an $E_l^+$-cell complex, so $R_l$ sends it to an isomorphism, and therefore $F^j_l f_j$ and $F^j_l g_j$ induce the same morphism of $\OO(A_l, E_l)$.

  Next, let $f : F^i X'_i \to R F^i Y'_i$ be any morphism.
  We must show that $f$ is the image of a morphism $F^i_j X'_i \to R_j F^i_j Y'_i$ for some $j \ge i$.
  Again, $\eta : F^i Y'_i \to R F^i Y'_i$ may be written as a $\lambda$-filtered colimit of morphisms with domain $F^i Y'_i$, each of which is an $E^+$-cell complex of length less than $\lambda$.
  Since $F^i X'_i$ is $\lambda$-compact, $f$ factors through one of these morphisms $h$ as shown below.
  \[
    \begin{tikzcd}
      F^i X'_i \ar[rd, "g"'] \ar[rrd, "f"] \\
      & Z \ar[r] & R F^i Y'_i \\
      F^i Y'_i \ar[ru, "h"] \ar[rru, "\eta"']
    \end{tikzcd}
  \]
  As before, we may realize $h$ as the image (up to isomorphism), for some $k \ge i$, of an $E_k^+$-cell complex $h_k : F^i_k Y'_i \to Z_k$ in $\PP(A_k)$.
  Replace $Z$, if necessary, by the isomorphic object $F^k Z_k$.
  Both $\eta$ and $h$ are $E^+$-cell complexes so $R$ sends the map $Z \to RF^iY'_i$ to an isomorphism.
  Thus it suffices to show that $g : F^i X'_i \to Z = F^k Z_k$ is the image of a map $g_l : F^i_l X'_i \to F^k_l Z_k$ for some $l \ge k$.
  This follows by applying \cref{prop:presheaf-b} once more.

  Thus, we have shown that $(*)$ and therefore also the original cocone are colimits.
\end{proof}

\begin{proposition}\label{prop:orthogonal-colimit}
  The cocone $(\OO(A_i, E_i)_\lambda)_{i \in D} \to \OO(A, E)_\lambda$ is a colimit in $\TCat$.
\end{proposition}

\begin{proof}
  As in the proof of \cref{prop:presheaf-colimit}, we apply \cref{prop:dcolimit-cat-retract}; it remains to check that every $\lambda$-compact object of $\OO(A, E)$ is a retract of the image under $G^i$ of a $\lambda$-compact object of $\OO(A_i, E_i)$ for some $i$.
  But, as we showed at the start of the previous proof, any $\lambda$-compact object of $\OO(A, E)$ is a retract of the image under $L$ of some $\lambda$-compact object $X'$ of $\PP(A)$.
  By \cref{prop:presheaf-colimit}, $X'$ is (isomorphic to) the image of some object $X'_i$ under $F^i$ for some $i$, and then $X$ is a retract of $G^i X'_i$.
\end{proof}

\begin{proposition}\label{prop:lambda-compact-directed}
  $(-)_\lambda : \LPr_\lambda \to \TCat$ preserves $\lambda$-directed colimits.
\end{proposition}

\begin{proof}
  By \cref{prop:o-form}, any diagram in $\LPr_\lambda$ is equivalent to one of the form $(\OO(A_i, E_i))_{i \in D}$ considered in \cref{prop:orthogonal-colimit}.
\end{proof}

Furthermore, the strengthened form of condition (b) proved in \cref{prop:orthogonal-colimit-b} also holds for any colimit in $\LPr_\lambda$.

\begin{remark}
  There is an alternative, arguably more elementary proof of this result.
  By \cite{Bird}, the colimit of a diagram in $\LPr_\lambda$ (or $\LPr$) may be computed by taking the right adjoints of all the functors involved and computing the limit of the resulting diagram in $\TCat$.
  The resulting limit category is again locally $\lambda$-presentable and the cone morphisms out of the limit category have left adjoints which assemble to a colimit cocone on the original diagram.
  When the indexing category is $\lambda$-directed, it is possible to write down explicit formulas for these left adjoints which easily yield \cref{prop:orthogonal-colimit-b}.
  Furthermore, each object of the colimit category has a canonical expression as a colimit of images of objects in each of the original categories, which verifies condition (a$'$).

  This approach does not seem to extend to the problem of understanding $\lambda$-directed colimits of $\lambda$-combinatorial premodel categories, however, and so we have opted in this section to give an argument using the same methods which will appear again later in this chapter.
\end{remark}

\section{\texorpdfstring{$\mu$}{\textmu}-small \texorpdfstring{$\lambda$}{\textlambda}-presentable categories}

Recall that a \emph{$\mu$-small $\lambda$-presentable category} is one equivalent to a category of the form $\OO(A, E)$ for a $\mu$-compact category $A$ and $\mu$-small set $E$ of morphisms between $\lambda$-compact objects of $\PP(A)$.
As always, we assume that $\mu \ge \lambda$.
In this section, we show that if $M$ is a $\mu$-small $\lambda$-presentable category, then $\THom_{\LPr_\lambda}(M, -) : \LPr_\lambda \to \TCat$ preserves $\mu$-directed colimits.

We have already verified the simplest nontrivial case $M = \Set$.

\begin{proposition}
  $\THom_{\LPr_\lambda}(\Set, -) : \LPr_\lambda \to \TCat$ preserves $\lambda$-directed colimits.
\end{proposition}

\begin{proof}
  For any $N$, giving a left adjoint $F : \Set \to N$ is equivalent to giving an object $X$ of $N$, and $F$ preserves $\lambda$-compact objects if and only if $X$ is $\lambda$-compact.
  Thus $\THom_{\LPr_\lambda}(\Set, -) = (-)_\lambda$ and so this is merely a restatement of \cref{prop:lambda-compact-directed}.
\end{proof}

The case of $M$ a presheaf category is not much more difficult.

\begin{proposition}\label{prop:presheaf-compact}
  Let $A$ be a $\mu$-compact category.
  Then $\THom_{\LPr_\lambda}(\PP(A), -) : \LPr_\lambda \to \TCat$ preserves $\mu$-directed colimits.
\end{proposition}

\begin{proof}
  For any $N$, giving a left adjoint $F : \PP(A) \to N$ is equivalent to giving an ordinary functor $X : A \to N$, and $F$ preserves $\lambda$-compact objects if and only $X(a)$ is $\lambda$-compact for each object $a$ of $A$.
  Let $(N_i)_{i \in D} \to N$ be a colimit cocone in $\LPr_\lambda$.
  We already know that $((N_i)_\lambda)_{i \in D} \to N_\lambda$ is a colimit cocone in $\TCat$, and we will use \cref{prop:dcolimit-cat} to deduce that $(((N_i)_\lambda)^A)_{i \in D} \to (N_\lambda)^A$ is one as well.\footnote{
    When $\mu > \lambda$, the category $A$ can have $\lambda$ or more objects and then $(N_\lambda)^A$ is generally not equal to $(N^A)_\lambda$.
    Thus, we avoid the ambiguous notation $N^A_\lambda$.
  }

  For condition (b2), we must check that any two parallel morphisms (natural transformations) in $((N_i)_\lambda)^A$ which become equal in $(N_\lambda)^A$ already become equal in $((N_j)_\lambda)^A$ for some $j \ge i$.
  The category $A$ has fewer than $\mu$ objects, and by condition (b2) for $((N_i)_\lambda)_{i \in D} \to N_\lambda$, for each object $a$ of $A$, the values on $a$ of the two natural transformations in question become equal in $((N_j)_\lambda)^A$ for some $j \ge i$.
  Since $D$ is $\mu$-directed, we can choose $j \ge i$ large enough that this condition holds for every $a$ simultaneously; then the two natural transformations become equal as morphisms of $((N_j)_\lambda)^A$, verifying condition (b2).

  For condition (b1), we must check that for any two objects (diagrams) $X$ and $Y$ of $((N_i)_\lambda)^A$, any natural transformation $t$ between their images in $(N_\lambda)^A$ arises from a natural transformation between their images in $((N_j)_\lambda)^A$ for some $j \ge i$.
  The category $A$ has fewer than $\mu$ objects, so by condition (b1) for $((N_i)_\lambda)_{i \in D} \to N_\lambda$, we can choose $i_1 \ge i$ such that for each object $a$ of $A$, the value of $t$ on $a$ arises from a morphism between the images of $X(a)$ and $Y(a)$ in $((N_{i_1})_\lambda)^A$; choose such a morphism for each object $a$ of $A$.
  These form an ``unnatural transformation'' between the images of $X$ and $Y$ in $((N_{i_1})_\lambda)^A$, which need not yet be a natural transformation.
  However, we can choose a set of fewer than $\mu$ morphisms of $A$ which generate it as a category, and for each such morphism $f : a \to b$, the naturality square for $f$ commutes in $(N_\lambda)^A$ and thus (by condition (b2) for $((N_i)_\lambda)_{i \in D} \to N_\lambda$) also in $((N_j)_\lambda)^A$ for some $j \ge i_1$.
  Choosing $j$ large enough so that all these squares commute in $((N_j)_\lambda)^A$, the image of our ``unnatural transformation'' in $((N_j)_\lambda)^A$ is a natural transformation between the images of $X$ and $Y$ whose image is the original natural transformation $t$.

  Finally, for condition (a), suppose $X$ is any object (diagram) of $(N_\lambda)^A$.
  By condition (a) for $((N_i)_\lambda)_{i \in D} \to N_\lambda$, each object $X(a)$ is in the essential image of $(N_i)_\lambda$ for some $i \in D$.
  Then we can apply \cref{lemma:small-diagram-factors}.
\end{proof}

In order to handle the general case $C = \OO(A, E)$, we need a preliminary lemma.

\begin{lemma}\label{lemma:iso-in-colimit}
  Suppose $(C_i)_{i \in D} \to C$ is a directed colimit cocone in $\TCat$.
  Let $i \in D$ and let $f_i : X_i \to Y_i$ be a morphism of $C_i$.
  If the image of $f_i$ in $C$ is an isomorphism, then the image of $f_i$ in $C_j$ is already an isomorphism for some $j \ge i$.
\end{lemma}

\begin{proof}
  Write $f : X \to Y$ for the image of $f_i$ in $C$, and let $g : Y \to X$ be its inverse.
  Use condition (b1) of \cref{prop:dcolimit-cat} to find $i_1 \ge i$ for which $g$ is the image of a morphism $g_{i_1}$ between the images of $Y$ and $X$ in $C_{i_1}$, and then use condition (b2) twice to find $j \ge i_1$ for which the equations $gf = \id_X$ and $fg = \id_Y$ already hold for the images of $f_i$ and $g_{i_1}$ in $C_j$.
  Then the image of $f_i$ in $C_j$ is an isomorphism.
\end{proof}

Of course the functors making up the diagram $(C_i)_{i \in D}$ preserve isomorphisms, so if the conclusion of the above lemma holds for some $j$, then it also holds for any $j' \ge j$.

\begin{proposition}\label{prop:lpr-small-compact}
  If $M$ is a $\mu$-small $\lambda$-presentable category, then $\THom_{\LPr_\lambda}(M, -) : \LPr_\lambda \to \TCat$ preserves $\mu$-directed colimits.
\end{proposition}

\begin{proof}
  We may assume $M = \OO(A, E)$, where $A$ is a $\mu$-compact category and $E$ is a $\mu$-small set of morphisms between $\lambda$-compact objects of $\PP(A)$.
  Suppose $(N_i)_{i \in D} \to N$ is a $\mu$-directed colimit cocone in $\LPr_\lambda$.
  By \cref{prop:presheaf-compact}, we already know that
  \[
    (\THom_{\LPr_\lambda}(\PP(A), N_i))_{i \in D} \to \THom_{\LPr_\lambda}(\PP(A), N)
  \]
  is a colimit in $\TCat$.
  We must show that it remains a colimit when $\THom_{\LPr_\lambda}(\PP(A), -)$ is replaced by $\THom_{\LPr_\lambda}(\OO(A, E), -)$, its full subcategory on the functors which take each morphism of $E$ to an isomorphism.
  For condition (b) of \cref{prop:dcolimit-cat}, there is nothing new to prove since $\THom_{\LPr_\lambda}(\OO(A, E), -)$ is (equivalent to) a full subcategory of $\THom_{\LPr_\lambda}(\PP(A), -)$.
  Thus it suffices to check condition (a).

  Write $G^i_j : N_i \to N_j$ and $G^i : N_i \to N$ for the functors which make up the colimit cocone $(N_i)_{i \in D} \to N$.
  As usual, we will assume the $G^i_j$ and $G^i$ form a strict diagram.
  Let $F$ be an object of $\THom_{\LPr_\lambda}(\OO(A, E), N)$, that is, $F : \OO(A, E) \to N$ is a left adjoint which preserves $\lambda$-compact objects.
  We must show that $F$ is isomorphic to $G^i \circ F_i$ for some object $F_i$ of $\THom_{\LPr_\lambda}(\OO(A, E), N_i)$.
  The composition $F' = F \circ L : \PP(A) \to \OO(A, E) \to N$ is an object of $\THom_{\LPr_\lambda}(\PP(A), N)$, and so (by \cref{prop:presheaf-compact}) is isomorphic to the composition $G^i \circ F'_i : \PP(A) \to N_i \to N$ for some $i \in D$ and some object $F'_i$ of $\THom_{\LPr_\lambda}(\PP(A), N_i)$.
  If $F'_i$ happens to send each morphism of $E$ to an isomorphism of $N_i$, then $F'_i$ factors (essentially uniquely) through a functor $F_i : \OO(A, E) \to N_i$ and then $G^i \circ F_i$ is isomorphic to $F$.
  In general, this will not be the case, but we will show that $G^i_j \circ F'_i : \PP(A) \to N_i \to N_j$ does send the morphisms of $E$ to isomorphisms for some $j \ge i$, and then we can apply the same argument to $j$ and $F'_j = G^i_j \circ F'_i$ instead.

  Let $e : X \to Y$ be a morphism of $E$.
  Then $X$ and $Y$ are $\lambda$-compact objects of $\PP(A)$, and so $F'_i$ sends $e$ to a morphism $f_i$ between $\lambda$-compact objects of $N_i$.
  The diagram $(N_i)_{i \in D} \to N$ is a $\mu$-directed colimit in $\LPr_\lambda$ and so by \cref{prop:lambda-compact-directed} (because $\mu \ge \lambda$) the diagram $((N_i)_\lambda)_{i \in D} \to N_\lambda$ is also a colimit in $\TCat$.
  The image $G^i f_i$ of $f_i$ in $N$ is an isomorphism, because $G^i f_i = G^i F'_i e \iso F' e = F (L e)$ and $L : \PP(A) \to \OO(A, E)$ sends elements of $E$ to isomorphisms.
  Then by \cref{lemma:iso-in-colimit}, we can find some $j \ge i$ such that $G^i_j f_i = (G^i_j \circ F'_i) e$ is already an isomorphism in $N_j$.
  Since $E$ contains fewer than $\mu$ morphisms and $D$ is $\mu$-directed, we can choose a single $j \ge i$ for which this holds simultaneously for all $e \in E$.
  Then $G^i_j \circ F'_i$ sends the morphisms of $E$ to isomorphisms, as desired.
\end{proof}

\section{\texorpdfstring{$\mu$}{\textmu}-small \texorpdfstring{$\lambda$}{\textlambda}-combinatorial premodel categories}

Recall that a \emph{$\mu$-small $\lambda$-combinatorial premodel category} $M$ is one for which
\begin{enumerate}
\item the underlying category of $M$ is $\mu$-small $\lambda$-presentable;
\item $M$ admits generating cofibrations and generating anodyne cofibrations which are each sets of fewer than $\mu$ morphisms between $\lambda$-compact objects.
\end{enumerate}
Our aim is to show that if $M$ is $\mu$-small $\lambda$-combinatorial, then $\THom_{\CPM_\lambda}(M, -) : \CPM_\lambda \to \TCat$ preserves $\mu$-directed colimits.

Let $(N_i)_{i \in D}$ be a (strict) diagram in $\CPM_\lambda$, made up of left Quillen functors $G^i_j : N_i \to N_j$.
For each $i$, choose generating cofibrations $I_i$ and generating anodyne cofibrations $J_i$ for $N_i$, each of which is a set of morphisms between $\lambda$-compact objects of $N_i$.
Since $G^i_j : N_i \to N_j$ is a left Quillen functor and preserves $\lambda$-compact objects, we may arrange that $G^i_j I_i \subset I_j$ and $G^i_j J_i \subset J_j$ for each $i \le j$ in $D$.
Write $N$ for the colimit of the $(N_i)_{i \in D}$ in $\CPM_\lambda$, with cocone morphisms $G^i : N_i \to N$, which we assume also form a strict diagram together with the $G^i_j$.
The underlying locally $\lambda$-presentable category of $N$ is the colimit of the $N_i$ in $\LPr_\lambda$, and $N$ admits generating cofibrations and anodyne cofibrations given by $I = \bigcup_{i \in D} G^i I_i$ and $J = \bigcup_{i \in D} G^i J_i$.

The following preliminary lemma will play a role analogous to that of \cref{lemma:iso-in-colimit}.

\begin{lemma}\label{lemma:cof-in-colimit}
  Suppose that $(N_i)_{i \in D} \to N$ is a $\lambda$-directed colimit cocone in $\CPM_\lambda$.
  Let $i \in D$ and let $f_i : X_i \to Y_i$ be a morphism between $\lambda$-compact objects of $N_i$.
  If the image $G^i f_i : G^i X_i \to G^i Y_i$ of $f_i$ in $N$ is a cofibration (respectively, anodyne cofibration), then the image $G^i_j f_i$ of $f_i$ in $N_j$ is already a cofibration (respectively, anodyne cofibration) for some $j \ge i$.
\end{lemma}

\begin{proof}
  The argument is the same for cofibrations and for anodyne cofibrations, so we only treat the case of cofibrations.
  Using \cref{prop:short-retract}, $G^i f_i$ can be written as a retract of an $I$-cell complex $h$ with fewer than $\lambda$ cells as shown below.
  \begin{equation}
    \begin{tikzcd}
      G^i X_i \ar[r, equals] \ar[d, "G^i f_i"'] &
      G^i X_i \ar[r, equals] \ar[d, "h"] &
      G^i X_i \ar[d, "G^i f_i"] \\
      G^i Y_i \ar[r, "s"'] & Z \ar[r, "r"'] & G^i Y_i
    \end{tikzcd}
    \tag{$*$}
  \end{equation}
  Using \cref{lemma:cell-complex-factors}, we can construct, for some $j \ge i$, an $I_j$-cell complex $h_j : G^i_j X_i \to Z_j$ in $N_j$, of the same length as $h$, which is sent by $G^j$ to $h$ (up to isomorphism).
  In $N_j$, then, we have the following partial diagram of $\lambda$-compact objects, which is sent by $G^j$ onto the corresponding part of $(*)$.
  Furthermore, $h_j$ is an $I_j$-cell complex and therefore a cofibration.
  \[
    \begin{tikzcd}
      G^i_j X_i \ar[r, equals] \ar[d, "G^i_j f_i"'] &
      G^i_j X_i \ar[r, equals] \ar[d, "h_j"] &
      G^i_j X_i \ar[d, "G^i_j f_i"] \\
      G^i_j Y_i \ar[r, dotted] & Z_j \ar[r, dotted] & G^i_j Y_i
    \end{tikzcd}
  \]
  As $(N_j) \to N$ is a $\lambda$-directed colimit in $\LPr_\lambda$, $((N_j)_\lambda) \to N_\lambda$ is a colimit in $\TCat$ by \cref{prop:lambda-compact-directed}.
  Thus, we can,
  \begin{itemize}
  \item first, find $k \ge j$ and maps $s_k : G^i_k Y_i \to G^j_k Z_j$ and $r_k : G^j_k Z_j \to G^i_k Y_i$ sent by $G^k$ to $s$ and $r$;
  \item second, find $l \ge k$ such that after applying $G^k_l$ to the (not necessarily commuting) diagram formed by filling the dotted arrows with $s_k$ and $r_k$, both squares commute in $N_l$ and the bottom composition is the identity.
  \end{itemize}
  Then $G^i_l f_i$ is a retract of the cofibration $G^j_l h_j$ of $N_l$ and therefore itself a cofibration.
\end{proof}

Again, if the conclusion of the preceding lemma holds for a particular $j$, then it also holds for any $j' \ge j$, because $G^j_{j'}$ is a left Quillen functor.

\begin{proposition}
  If $M$ is a $\mu$-small $\lambda$-combinatorial premodel category, then $\THom_{\CPM_\lambda}(M, -) : \CPM_\lambda \to \TCat$ preserves $\mu$-directed colimits.
\end{proposition}

\begin{proof}
  Let $I$ and $J$ be generating cofibrations and generating anodyne cofibrations for $M$ which are sets of fewer than $\mu$ objects between $\lambda$-compact objects.
  Suppose $(N_i)_{i \in D} \to N$ is a $\mu$-directed colimit cocone in $\CPM_\lambda$; we must show that
  \[
    (\THom_{\CPM_\lambda}(M, N_i))_{i \in D} \to \THom_{\CPM_\lambda}(M, N)
  \]
  is a colimit cocone in $\TCat$.
  Since $(N_i)_{i \in D} \to N$ is also a ($\mu$-directed) colimit cocone in $\LPr_\lambda$, we know from \cref{prop:lpr-small-compact} that
  \[
    (\THom_{\LPr_\lambda}(M, N_i))_{i \in D} \to \THom_{\LPr_\lambda}(M, N)
  \]
  is a colimit cocone in $\TCat$.
  As in the proof of \cref{prop:lpr-small-compact}, $\THom_{\CPM_\lambda}(M, N)$ is a full subcategory of $\THom_{\LPr_\lambda}(M, N)$ for any $N$, namely the full subcategory on the functors which send the morphisms of $I$ to cofibrations of $N$ and the morphisms of $J$ to anodyne cofibrations of $N$.
  Thus, condition (b) of \cref{prop:dcolimit-cat} has already been checked and it remains to verify condition (a).
  The argument for this is the same as the end of the proof of \cref{prop:lpr-small-compact}, using the fact that $I$ and $J$ each contain fewer than $\mu$ morphisms between $\lambda$-compact objects of $M$ and using \cref{lemma:cof-in-colimit} in place of \cref{lemma:iso-in-colimit}.
\end{proof}

\begin{corollary}\label{cor:cpm-lam-small}
  For any object $M$ of $\CPM_\lambda$, $\THom_{\CPM_\lambda}(M, -) : \CPM_\lambda \to \TCat$ preserves $\mu$-directed colimits for all sufficiently large $\mu$.
\end{corollary}

\chapter{Combinatorial $V$-premodel categories}
\label{chap:modules}

In the past two chapters, we have shown that the 2-category $\CPM$ of combinatorial premodel categories enjoys the following algebraic properties.

\begin{enumerate}
\item
  $\CPM$ is a complete and cocomplete 2-category.
\item
  $\CPM$ is tensored and cotensored over small categories, with the tensor given by the projective premodel structure $I \otimes M = M^{I^\op}_\proj$ and the cotensor given by the injective premodel structure $M^I = M^I_\inj$.
\item
  In particular, for any small category $I$, the combinatorial premodel category $\Set^{I^\op}_\proj$ satisfies the universal property
  \[
    \Hom_\CPM(\Set^{I^\op}_\proj, N) \eqv \Hom_\TCat(I, N^\cof)
  \]
  with the equivalence given by precomposition with $\yo : I \to \Set^{I^\op}_\proj$ (which sends each object of $I$ to a cofibrant object of $\Set^{I^\op}_\proj$).
\item
  By adjoining additional generating cofibrations and anodyne cofibrations to $\Set^{I^\op}_\proj$, we can design combinatorial premodel categories for which $\CPM(M, N)$ is equivalent to the full subcategory of $\TCat(I, N^\cof)$ on the diagrams satisfying any desired requirements of the form that a morphism built out of the diagram by colimits is a cofibration or an anodyne cofibration.
  As an arbitrary example, there is some object $M$ of $\CPM$ for which $\CPM(M, N)$ consists of cofibrations $A \to B$ between cofibrant objects of $N$ such that the induced map $A \amalg B \to B \amalg B$ is an anodyne cofibration.
\item
  The 2-category $\CPM$ has a filtration by sub-2-categories $\CPM_\lambda$, $\lambda$ ranging over all regular cardinals, such that every small diagram belongs to $\CPM_\lambda$ for sufficiently large $\lambda$.
  The subcategories $\CPM_\lambda$ are ``locally full'' in the sense that, for any objects $M$ and $N$ of $\CPM_\lambda$, the Hom category $\CPM_\lambda(M, N)$ is a full subcategory of $\CPM(M, N)$ (namely, the full subcategory consisting of the functors which preserve $\lambda$-compact objects).
\item
  The colimit in $\CPM$ of any diagram in $\CPM_\lambda$ belongs to $\CPM_\lambda$ and is also a colimit in $\CPM_\lambda$.
  The corresponding statement for limits holds for $\lambda$-small diagrams (as long as $\lambda$ is uncountable).
\item
  The category $\LPr(M, N)$ of \emph{all} left adjoints from $M$ to $N$ has a combinatorial premodel category structure, denoted by $\CPM(M, N)$, for which the cofibrant objects are precisely the left Quillen functors.
  That is, $\THom_\CPM(M, N) = \CPM(M, N)^\cof$.
\item
  Let $M$ be any object of $\CPM$.
  Then for all sufficiently large (in particular, uncountable) regular cardinals $\lambda$, $M$ is \emph{$\lambda$-small $\lambda$-combinatorial}, which has the following consequences.
  \begin{itemize}
  \item
    For any $\lambda$-combinatorial premodel category $N$, the category $\CPM(M, N)$ is also $\lambda$-combinatorial.
    Moreover, its $\lambda$-compact objects are precisely the left adjoint functors from $M$ to $N$ which preserve $\lambda$-compact objects.
  \item
    The functor $\THom_\CPM(M, -) : \CPM \to \TCat$ preserves $\lambda$-directed (conical) colimits.
  \end{itemize}
\item
  For any objects $M_1$ and $M_2$ of $\CPM$, there is a universal Quillen bifunctor $M_1 \times M_2 \to M_1 \otimes M_2$.
  The tensor product $M_1 \otimes M_2$ satisfies a tensor--Hom adjunction $\CPM(M_1 \otimes M_2, N) \eqv \CPM(M_1, \CPM(M_2, N))$.
  The subcategory $\CPM_\lambda$ is closed under $\otimes$.
\end{enumerate}

However, $\CPM$ has the drawback that its objects are not all relaxed premodel categories, and so we don't know how to do homotopy theory with them.
This drawback can be corrected by working with $V$-premodel categories for, say, a monoidal model category $V$.
Accordingly, we desire a 2-category $V\CPM$ of ``combinatorial $V$-premodel categories'' satisfying analogous properties to the ones listed above.

There are two equivalent ways we could proceed.

\begin{enumerate}
\item
  We could define a $V$-premodel category to be a \emph{$V$-enriched} category $M$ which is tensored and cotensored over $V$, and for which the tensor $\otimes : V \times M \to M$ is an (ordinary) Quillen bifunctor.
  This approach has the advantage that it encodes the additional data related to $V$ into the underlying $V$-enriched category $M$, which is not much more complicated than an ordinary category.
  Then we could generalize the material of the previous two chapters to $V$-enriched premodel categories.

  This approach turns out to be inconvenient in practice because a premodel category structure is really a structure on an ordinary category.
  Furthermore it is unclear whether some of the background theory we rely on (such as the fat small object argument) has a meaningful enriched equivalent.
  One ends up having to relate notions of enriched category theory to their unenriched equivalents to an extent that eliminates the apparent simplicity of working with enriched categories.

\item
  Alternatively, we can use the 2-multicategory structure of $\CPM$ to define $V\CPM$ as the 2-category of modules over the monoid object $V$.
  In this approach the passage from combinatorial premodel categories to combinatorial $V$-premodel categories is a formal procedure which depends only on the algebraic structure of $V\CPM$.

  On a technical level, there are two ways to formalize this.
  One option is to use the theory of symmetric monoidal $\infty$-categories developed in \cite{HA}.
  The other is to just write out all the structure explicitly, since this is not that difficult for 2-categories.
  We will freely mix these two approaches.
\end{enumerate}

Because this chapter has little to do with premodel categories in particular, we will offer only an summary of the theory of $V\CPM$ and omit most details.

In \cref{sec:cpm-multicat} we explained how to construct a symmetric 2-multicategory structure on $\CPM$.
If we apply the functor sending a category to the nerve of its maximal subgroupoid to each Hom object of this multicategory structure we obtain a simplicial colored operad in the sense of \cite[Variant~2.1.1.3]{HA}.
The existence of universal Quillen multifunctors satisfying the property described at the end of \cref{sec:cpm-multicat} implies that the operadic nerve of this simplicial colored operad is a symmetric monoidal $\infty$-category $\CPM^\otimes$.

The underlying $\infty$-category of $\CPM^\otimes$ is only a $(2, 1)$-category, not a 2-category; in the process of constructing it we discarded the noninvertible 2-morphisms.
However, since $\CPM$ is tensored and cotensored over $\TCat$, the lost information is easily recovered: we can compute the space of strings of $n$ morphisms in $\THom_\CPM(M, N)$ as the space of objects of $\THom_\CPM([n] \otimes M, N)$.

The sub-2-category $\CPM_\lambda$ of $\CPM$ is closed under tensor products, as is the further sub-2-category of $\lambda$-small $\lambda$-combinatorial premodel categories.
Let $V$ be a combinatorial monoidal premodel category.
We will always assume that $\lambda$ is large enough that the following conditions on $V$ are satisfied.
\begin{enumerate}
\item $V$ is a $\lambda$-small $\lambda$-combinatorial premodel category.
\item The multiplication $\mu : V \otimes V \to V$ and the unit morphism $\Set \to V$ (corresponding to the object $1_V$) are strongly $\lambda$-accessible.
\end{enumerate}
Then $V$ is a monoid object in $\CPM_\lambda$.
Moreover the assumption that $V$ is $\lambda$-small $\lambda$-combinatorial means that for any $\mu \ge \lambda$ and $\mu$-small $\lambda$-combinatorial premodel category $M$, $V \otimes M$ is again $\mu$-small $\lambda$-combinatorial.

We define $V\CPM_\lambda$ to be the 2-category of modules over $V$ in $\CPM_\lambda$, and $V\CPM$ to be the 2-category of modules over $V$ in $\CPM$.
An object $M$ of $V\CPM$ belongs to $V\CPM_\lambda$ if and only if the underlying object of $\CPM$ belongs to $\CPM_\lambda$ (i.e., is $\lambda$-combinatorial) and the action $\otimes : V \times M \to M$ is strongly $\lambda$-accessible.
There is a free $V$-module functor $V \otimes_\Set - : \CPM \to V\CPM$ which is left adjoint to the forgetful functor $V\CPM \to \CPM$; the underlying object of $V \otimes_\Set M$ is $V \otimes M$.
We now discuss the list of items which opened this chapter in relation to $V\CPM_\lambda$ and $V\CPM$.

\begin{enumerate}
\item
  $V\CPM$ has limits and colimits.
  A diagram in $V\CPM$ is a limit or colimit if and only if its image under the forgetful functor to $\CPM$ is one.
\item
  The same holds for tensors and cotensors, which are given by projective and injective premodel category structures respectively.
  $V$ acts on these diagram categories carry componentwise.
\item
  The role of $\yo$ is now played by the $V$-valued Yoneda embedding $\yo : A \to V^{A^\op}$ which sends $a$ to the object $\yo a$ for which $(\yo a)_b$ is the coproduct of $\Hom_V(b, a)$ copies of $1_V$.
\item
  We can analogously adjoin additional generating cofibrations $I'$ and generating anodyne cofibrations $J'$ ``over $V$'' to a combinatorial $V$-premodel category $M$ to obtain a new combinatorial $V$-premodel category $M'$.
  This means that a morphism $F : M \to N$ of $V\CPM$ extends (essentially uniquely) to $M'$ if and only if $F$ sends the morphisms of $I'$ to cofibrations and the morphisms of $J'$ to anodyne cofibrations.
  Concretely, $M'$ has the same underlying category and action by $V$ as $M$.
  To describe its premodel category structure, choose generating cofibrations $I_V$ and anodyne cofibrations $J_V$ for $V$.
  Then $M'$ is formed from $M$ by adjoining additional generating cofibrations $I_V \bp I'$ and additional generating anodyne cofibrations $J_V \bp I' \cup I_V \bp J'$.
  Because $V$ is monoidal, adjoining these produces the smallest premodel category structure $M'$ which is a module over $V$ and contains $I'$ as cofibrations and $J'$ as anodyne cofibrations.
\item
  $V\CPM$ has a filtration by sub-2-categories $V\CPM_\lambda$ with the proviso that we only allow $\lambda$ sufficiently large compared to $V$, as described earlier.
\item
  $V\CPM_\lambda$ is closed in $V\CPM$ under colimits and $\lambda$-small limits (for $\lambda$ uncountable and also sufficiently large compared to $V$).
\item
  The category $\THom_{V\CPM}(M, N)$ can be computed as a limit of a truncated simplicial diagram
  \[
    \begin{tikzcd}
      \THom_{\CPM}(M, N) \ar[r, shift left=2] \ar[r, shift right=2] &
      \THom_{\CPM}((V, M), N) \ar[l] \ar[r, shift left=2] \ar[r] \ar[r, shift right=2] &
      \THom_{\CPM}((V, V, M), N)
    \end{tikzcd}
  \]
  and so we define $\CPM^V(M, N)$ to be given by the limit in $\CPM$ of the same diagram with $\THom_\CPM(-, -)$ replaced by $\CPM(-, -)$.
\item
  Let $M$ be an object of $V\CPM$ whose underlying premodel category is $\lambda$-small $\lambda$-combinatorial.
  Then $V \otimes M$ and $V \otimes V \otimes M$ are also $\lambda$-small $\lambda$-combinatorial.
  Therefore:
  \begin{itemize}
  \item
    If $N$ is an object of $V\CPM_\lambda$, then $\CPM^V(M, N)$ is a finite limit of $\lambda$-combinatorial premodel categories, hence $\lambda$-combinatorial.
    An object of $\CPM^V(M, N)$ is $\lambda$-compact when the corresponding functor $M \to N$ is strongly $\lambda$-accessible (since the induced functors $V \otimes M \to N$ and $V \otimes V \otimes M \to N$ are then automatically strongly $\lambda$-accessible as well).
  \item
    $\THom_{V\CPM}(M, -)$ is a finite limit of functors which preserve $\lambda$-directed colimits and therefore itself preserves $\lambda$-directed colimits.
  \end{itemize}
\item
  For $M_1 \in \CPM$, $M_2 \in V\CPM$, we define $M_1 \otimes M_2$ to be the tensor product in $\CPM$ equipped with the action
  \[
    V \otimes (M_1 \otimes M_2) \eqv M_1 \otimes (V \otimes M_2) \to M_1 \otimes M_2.
  \]

  Then there is an adjunction of two variables for $M_1 \in \CPM$, $M_2 \in V\CPM$, $N \in V\CPM$:
  \[
    \THom_{V\CPM}(M_1 \otimes M_2, N) \eqv \THom_\CPM(M_1, \CPM^V(M_2, N))
    \eqv \THom_{V\CPM}(M_2, N^{M_1}).
  \]
  Here $N^{M_1}$ is $\CPM(M_1, N)$ equipped with the action of $V$ given by $(K \otimes F)(A_1) = K \otimes (FA_1)$.
  The tensor product $M_1 \otimes M_2$ belongs to $V\CPM_\lambda$ if $M_1$ belongs to $\CPM_\lambda$ and $M_2$ belongs to $V\CPM_\lambda$.
\end{enumerate}

\chapter{Model 2-categories}

The 2-category $V\CPM$ cannot be a model category in the ordinary sense, because its underlying 1-category lacks limits and colimits.
Furthermore, we prefer not to distinguish between equivalent combinatorial $V$-premodel categories and so we would like the cofibrations and fibrations of $V\CPM$ to be equivalence-invariant.
This means that we cannot expect the lifting axioms of a model category to hold strictly in $V\CPM$.
Instead, we ought to relax the axioms of a model category by replacing equalities between morphisms by invertible 2-morphisms in an appropriate way.

A suitable framework for our purposes is provided by \cite{MG} which extends the notion of model category to the setting of $(\infty, 1)$-categories.
Our $V\CPM$ is a 2-category and not a $(2, 1)$-category, but the noninvertible 2-morphisms will not play any role related to the model category structure.
Thus, we define a model 2-category to be a complete and cocomplete 2-category whose underlying $(2, 1)$-category is equipped with a model $(\infty, 1)$-category structure in the sense of \cite{MG}.
For the sake of concreteness, we give an explicit definition specialized to $(2, 1)$-categories.

\section{Weak factorization systems on $(2, 1)$-categories}

In this section, we fix an ambient $(2, 1)$-category $\C$.
We assume that $\C$ is strict, since all of our examples arise from strict 2-categories, though only minor modifications are needed to handle bicategories or other models for $(2, 1)$-categories.

\begin{convention}
  We take as our starting point a $(2, 1)$-category in order to emphasize that the notions we introduce (notably weak factorization systems and model 2-category structures) involve only \emph{invertible} 2-morphisms.
  In practice, $\C$ will be the $(2, 1)$-category obtained from a 2-category by discarding the noninvertible 2-morphisms.
  We implicitly extend all of these notions to general 2-categories by this process.
  For example, a 1-morphism $P$ has the right lifting property with respect to a 1-morphism $I$ in a 2-category $\C$ if $P$ has the right lifting property with respect to $I$ (as defined below) in the $(2, 1)$-category obtained from $\C$ by discarding the noninvertible 2-morphisms.
\end{convention}

\begin{definition}\label{def:2-lp}
  Let $I : A \to B$ and $P : X \to Y$ be 1-morphisms of $\C$.
  We say that \emph{$P$ has the right lifting property with respect to $I$}, or that \emph{$I$ has the left lifting property with respect to $P$}, if for any square (commuting up to an invertible 2-morphism) of the form
  \[
    \begin{tikzcd}[row sep=large, column sep=large]
      A \ar[r, "H"] \ar[d, "I"'] &
      X \ar[d, "P"] \ar[ld, Rightarrow, shorten >=14pt, shorten <=14pt, "\alpha"] \\
      B \ar[r, "K"'] & Y
    \end{tikzcd}
  \]
  there exists a 1-morphism $L : B \to X$ and (invertible) 2-morphisms $\beta : H \to LI$, $\gamma : PL \to K$ such that the composition
  \[
    \begin{tikzcd}[row sep=large, column sep=large]
      A \ar[r, "H", ""'{name=H}] \ar[d, "I"'] & X \ar[d, "P"] \\
      B \ar[r, "K"', ""{name=K}] \ar[ru, "L"] & Y
      \ar[Rightarrow, from=H, to=2-1, shorten >=14pt, shorten <=14pt, shift right, "\beta"']
      \ar[Rightarrow, from=1-2, to=K, shorten >=14pt, shorten <=14pt, shift left, "\gamma"]
    \end{tikzcd}
  \]
  equals $\alpha$ (that is, $\alpha = \gamma I \circ P \beta$).
\end{definition}

When $\C$ is an ordinary category (viewed as a $(2, 1)$-category with only identity 2-morphisms), these definitions reduce to the ordinary notions because $\alpha$, $\beta$ and $\gamma$ must be equalities.

\begin{warning}
  The $(2, 1)$-category $\C$ has an associated ``homotopy category'' $\tau_{\le 1} \C$, the ordinary category with the same objects as $\C$ in which $\Hom_{\tau_{\le 1} \C}(X, Y)$ consists of the isomorphism classes of objects of $\THom_\C(X, Y)$.
  The statement that $P$ has the right lifting property with respect to $I$ in $\C$ is, in general, stronger than the corresponding statement about the images of $P$ and $I$ in $\tau_{\le 1} \C$; the latter statement does not include the condition that the two triangles in the completed diagram compose to the square in the original lifting problem.
\end{warning}


\begin{convention}
  The preceding warning notwithstanding, in order to simplify the exposition, we will never name or notate the 2-morphisms that appear in lifting problems.
  Instead, we adopt the convention that, in the context of a 2-category, a diagram like
  \[
    \begin{tikzcd}
      A \ar[r] \ar[d, "I"'] & X \ar[d, "P"] \\
      B \ar[r] \ar[ru, dotted, "L"] & Y
    \end{tikzcd}
  \]
  represents a square in $\C$ which commutes up to some specified (invertible) 2-morphism, and solving the lifting problem means factoring this square into two triangles which are glued along a 1-morphism $L : B \to X$.
\end{convention}

\begin{example}
  Take $\C = \TCat$ and let $I : \emptyset \to \{\star\}$ be the unique functor from the empty category to the terminal category.
  A lifting problem
  \[
    \begin{tikzcd}
      \emptyset \ar[r, "H"] \ar[d, "I"'] & X \ar[d, "P"] \\
      \{\star\} \ar[r, "K"'] & Y
    \end{tikzcd}
  \]
  amounts simply to giving an object $y$ of $Y$, namely the object $K(\star)$.
  The 2-morphism filling the square contains no data, because $\emptyset$ is an initial object.
  A lift
  \[
    \begin{tikzcd}
      \emptyset \ar[r, "H"] \ar[d, "I"'] & X \ar[d, "P"] \\
      \{\star\} \ar[r, "K"'] \ar[ru, "L"] & Y
    \end{tikzcd}
  \]
  amounts to an object $x$ of $X$ (arising from $L$) together with an isomorphism $Px \cong y$ (arising from $\gamma$).
  Again, $\beta$ contains no data; and the condition that the triangles $\beta$ and $\gamma$ compose to the original $\alpha$ is vacuously satisfied.
  Thus, $P$ has the right lifting property with respect to $I : \emptyset \to \{\star\}$ if and only if for every object $y$ of $Y$, there exists an object $x$ of $X$ and an isomorphism $Px \cong y$, that is, if and only if $P$ is essentially surjective.

  In this case, because the domain of $I$ is an initial object, the lifting property actually does reduce to a lifting property in $\tau_{\le 1} \TCat$ (though not to a lifting property in $\Cat$).
\end{example}

\begin{example}\label{ex:r-rlp-cat}
  Take $\C = \TCat$ and let $I : \{\star\} \to \{\star \to \star'\}$ be the functor sending $\star$ to $\star$.
  A lifting problem
  \[
    \begin{tikzcd}
      \{\star\} \ar[r, "H"] \ar[d, "I"'] & X \ar[d, "P"] \\
      \{\star \to \star'\} \ar[r, "K"'] & Y
    \end{tikzcd}
  \]
  amounts to
  \begin{enumerate}
  \item an object $x$ of $X$ (the object $H(\star)$),
  \item a morphism $g : y \to y'$ of $Y$ (the image under $K$ of the morphism $\star \to \star'$),
  \item and an isomorphism $\psi : Px \cong y$ (the component of $\alpha$ on $\star$).
  \end{enumerate}
  A solution of this lifting problem
  \[
    \begin{tikzcd}
      \{\star\} \ar[r, "H"] \ar[d, "I"'] & X \ar[d, "P"] \\
      \{\star \to \star'\} \ar[r, "K"'] \ar[ru, "L"] & Y
    \end{tikzcd}
  \]
  amounts to
  \begin{enumerate}
  \item a morphism $f : x_1 \to x'$ of $X$ (the image under $L$ of the morphism $\star \to \star'$),
  \item an isomorphism $\varphi : x \to x_1$ (the component of $\beta$ on $\star$),
  \item and vertical isomorphisms forming a commutative square
    \[
      \begin{tikzcd}
        Px_1 \ar[r, "Pf"] \ar[d, "\psi_1"'] & Px' \ar[d, "\psi'"] \\
        y \ar[r, "g"'] & y'
      \end{tikzcd}
    \]
    (the components and naturality square of $\gamma$),
  \item such that $\psi = \psi_1 \circ P \varphi$ (the condition relating $\alpha$, $\beta$ and $\gamma$).
  \end{enumerate}
  We want to determine which $P$ satisfy this lifting property.
  Clearly, we might as well choose $x_1$ to be $x$ and $\varphi : x \to x_1$ to be the identity, replacing $f : x_1 \to x'$ by $f \circ \varphi : x \to x_1$.
  Furthermore, we may assume that $\psi$ is the identity, by replacing $g : y \to y'$ by $g \circ \psi : Px \to y'$.
  Thus, we can reformulate the lifting property as follows:
  $P$ has the right lifting property with respect to $I$ if and only if for any object $x$ of $X$ and any morphism $g : Px \to y'$ of $Y$, there exists a morphism $f : x \to x'$ and an isomorphism $\psi' : x' \to y'$ making the square
  \[
    \begin{tikzcd}
      Px \ar[r, "Pf"] \ar[d, equals] & Px' \ar[d, "\psi'"] \\
      Px \ar[r, "g"'] & y'
    \end{tikzcd}
  \]
  commute.
\end{example}

\begin{convention}
  Conditions like this one will play an important role in the model category structure on $V\CPM$.
  We will refer to the morphism $f : x \to x'$ as a \emph{lift} of $g : Px \to y'$, sometimes inserting the phrase ``up to isomorphism'' as a reminder that $Pf$ is not required to be equal to $g$ but only isomorphic to it in the sense described above.

  One way to justify this terminology is as follows.
  Recall that a functor $P : X \to Y$ is said to be an \emph{isofibration} if for any object $x$ of $X$ and any isomorphism $P : Px \to y'$ in $Y$, there exists an isomorphism $f : x \to x'$ whose image under $P$ is \emph{strictly} equal to $g$.
  When $P$ is an isofibration, we can turn an ``up to isomorphism'' lift $f : x \to x'$ of $g : Px \to y'$ into a strict lift by composing $f$ with a strict lift of the isomorphism $\psi' : Px' \to y'$.
  Thus, the ``up to isomorphism'' lifting property is equivalent to the strict one when $P$ is an isofibration.
  Now, any functor $P : X \to Y$ can be factored as an equivalence $X \to X'$ followed by an isofibration $X' \to Y$.
  All the structure we consider is invariant under equivalence, so we may replace a functor by an isofibration wherever convenient in order to turn ``up to isomorphism'' lifting properties into strict ones.
\end{convention}

\begin{proposition}\label{prop:lifting-replete}
  Let $I : A \to B$ be a 1-morphism of $\C$ and suppose the 1-morphisms $P : X \to Y$ and $P' : X \to Y$ are isomorphic.
  Then $P$ has the right lifting property with respect to $I$ if and only if $P'$ does.
  A dual statement holds for the left lifting property.
\end{proposition}

\begin{proof}
  Let $\lambda : P \to P'$ be an isomorphism and suppose $P$ has the right lifting property with respect to $I$.
  Given a lifting problem for $P'$, we may attach $\lambda$ as shown below to construct a lifting problem for $P$.
  \[
    \begin{tikzcd}
      A \ar[r, "H"] \ar[d, "I"'] & X \ar[d, "P'"'] \ar[d, bend left=50, "P"] \\
      B \ar[r, "K"'] \ar[ru, dotted] & Y
    \end{tikzcd}
  \]
  A solution for this lifting problem for $P$ can be turned into a solution for the original lifting problem by composing the lower triangle with the inverse of $\lambda$.
\end{proof}

\begin{notation}
  For classes of 1-morphisms $\sL$, $\sR$ of $\C$, we write $\llp(\sR)$ for the class of 1-morphisms $I$ of $C$ with the left lifting property with respect to every $P \in \sR$, and $\rlp(\sL)$ for the class of 1-morphisms $P$ of $C$ with the right lifting property with respect to every $I \in \sL$.
\end{notation}

By \cref{prop:lifting-replete}, the classes $\llp(\sR)$ and $\rlp(\sL)$ are closed under replacing a 1-morphism by a parallel isomorphic one.
As usual, $\llp$ and $\rlp$ are each inclusion-reversing, and for any class $\sL$, $\sL \subset \llp(\rlp(\sL))$, so that $\rlp(\llp(\rlp(\sL))) = \rlp(\sL)$.

We next give an alternate description of the lifting condition.

\begin{definition}
  Let $I : A \to B$ and $P : X \to Y$ be two 1-morphisms of $\C$.
  We define the \emph{category of lifting problems} or \emph{category of squares} $\Sq(I, P)$ to be the pseudopullback
  \[
    \begin{tikzcd}
      \Sq(I, P) \ar[r] \ar[d] & \THom_\C(A, X) \ar[d, "P_*"] \\
      \THom_\C(B, Y) \ar[r, "I^*"'] & \THom_\C(B, X)
    \end{tikzcd}
  \]
  in $\TCat$.
  Concretely, an object of $\Sq(I, P)$ consists of
  \begin{enumerate}
  \item a 1-morphism $H : A \to X$,
  \item a 1-morphism $K : B \to Y$,
  \item and a 2-morphism $\alpha : P H \to K I$,
  \end{enumerate}
  that is, precisely the data defining a lifting problem as in \cref{def:2-lp}; while a morphism from $(H, K, \alpha)$ to $(H', K', \alpha')$ consists of 2-morphisms $\varphi : H \to H'$ and $\psi : K \to K'$ such that the square
  \[
    \begin{tikzcd}
      PH \ar[r, "\alpha"] \ar[d, "P\varphi"'] & KI \ar[d, "\psi I"] \\
      PH' \ar[r, "\alpha'"] & K'I
    \end{tikzcd}
  \]
  commutes.
  The commutative square
  \[
    \begin{tikzcd}
      \THom_\C(B, X) \ar[r, "I^*"] \ar[d, "P_*"'] & \THom_\C(A, X) \ar[d, "P_*"] \\
      \THom_\C(B, Y) \ar[r, "I^*"'] & \THom_\C(B, X)
    \end{tikzcd}
  \]
  induces a functor $\LL(I, P) : \THom_\C(B, X) \to \Sq(I, P)$, given explicitly by the formula
  \[
    \LL(I, P)(L) = (LI, PL, \id_{PLI}).
  \]
\end{definition}

\begin{remark}
  When $\C$ is a $(2, 1)$-category, $\Sq(I, P)$ is not merely a category but a groupoid.
  When $\C$ is a 2-category, the above construction produces a category $\Sq(I, P)$.
  The maximal subgroupoid of this category is the same as the groupoid we get by first discarding the noninvertible morphisms of $\C$ and then performing the above construction.
  We only really care about the invertible morphisms of $\Sq(I, P)$, so this slight ambiguity in its definition will not bother us.
\end{remark}

\begin{proposition}\label{prop:rlp-iff-ess-surj}
  $P$ has the right lifting property with respect to $I$ if and only if the functor $\LL(I, P)$ is essentially surjective.
\end{proposition}

\begin{proof}
  An isomorphism $(H, K, \alpha) \to \LL(I, P)(L)$ consists of an isomorphism $\beta : H \to LI$ and an isomorphism $\gamma^{-1} : K \to PL$ such that $\gamma^{-1}I \circ \alpha = \id_{PLI} \circ P \beta$, or equivalently $\alpha = \gamma I \circ P \beta$.
  In other words, for any object $(H, K, \alpha)$ of $\Sq(I, P)$, an object $L$ together with an isomorphism $(H, K, \alpha) \to \LL(I, P)(L)$ is precisely the same as a solution to the lifting problem $(H, K, \alpha)$.
\end{proof}

\begin{proposition}
  Equivalences have the left and right lifting property with respect to any 1-morphism.
\end{proposition}

\begin{proof}
  If $P : X \to Y$ is an equivalence, then so are the morphisms marked $\eqv$ in the diagram below.
  \[
    \begin{tikzcd}
      \THom_\C(B, X) \ar[rd, "{\LL(I, P)}"]
      \ar[rrd, bend left=10, "I^*"] \ar[rdd, bend right=10, "P_*"', "\eqv"] \\
      & \Sq(I, P) \ar[r] \ar[d, "\eqv"'] & \THom_C(A, X) \ar[d, "P_*", "\eqv"'] \\
      & \THom_\C(B, Y) \ar[r, "I^*"'] & \THom_C(B, X)
    \end{tikzcd}
  \]
  Then $\LL(I, P)$ is also an equivalence and in particular essentially surjective.
\end{proof}

\begin{example}
  Take $\C = \LPr$ and let $\PP(I) : \PP(\{\star\}) \to \PP(\{\star \to \star'\})$ be the left adjoint induced by the functor $I : \{\star\} \to \{\star \to \star'\}$ sending $\star$ to $\star$.
  For any small category $A$ and locally presentable category $M$ there is an equivalence $\THom_\LPr(\PP(A), M) \eqv \TCat(A, M)$.
  Hence if $P : X \to Y$ is any morphism of $\LPr$, then the functor $\LL(\PP(I), P)$ is equivalent to the functor $\LL(I, P)$ as computed in $\TCat$, and so $P$ has the right lifting property with respect to $\PP(I)$ if and only if the underlying functor of $P$ has the right lifting property with respect to $I$.
  We described such functors in \cref{ex:r-rlp-cat}.
\end{example}

The left and right lifting properties enjoy the same stability properties as in the ordinary 1-categorical case.
We will treat the case of retracts in detail and leave the remaining properties to the reader.

\begin{definition}
  A 1-morphism $F : A \to B$ is a \emph{retract} of a 1-morphism $F' : A' \to B'$ if there exists a diagram (in which each bounded region is filled by an invertible 2-morphism)
  \[
    \begin{tikzcd}
      A \ar[r] \ar[rr, bend left, "\id_A"] \ar[d, "F"'] & A' \ar[r] \ar[d, "F'"] & A \ar[d, "F"] \\
      B \ar[r] \ar[rr, bend right, "\id_B"'] & B' \ar[r] & B
    \end{tikzcd}
  \]
  such that the composite 2-morphism is the identity on $F : A \to B$.
\end{definition}

\begin{proposition}
  For any $\sL$ and $\sR$, the classes $\llp(\sR)$ and $\rlp(\sL)$ are closed under retracts.
\end{proposition}

\begin{proof}
  Suppose that $I' : A' \to B'$ belongs to $\llp(\sR)$ and $I : A \to B$ is a retract of $I'$.
  Given a lifting problem for $I$ with respect to some $P : X \to Y$ in $\sR$
  \[
    \begin{tikzcd}
      A \ar[r, "H"] \ar[d, "I"'] & X \ar[d, "P"] \\
      B \ar[r, "K"'] & Y
    \end{tikzcd}
  \]
  we can attach the diagram exhibiting $I$ as a retract $I'$ to form a diagram
  \[
    \begin{tikzcd}
      A \ar[r] \ar[rr, bend left, "\id_A"] \ar[d, "I"'] & A' \ar[r] \ar[d, "I'"] &
      A \ar[r, "H"] \ar[d, "I"'] & X \ar[d, "P"] \\
      B \ar[r] \ar[rr, bend right, "\id_B"'] & B' \ar[r] &
      B \ar[r, "K"'] & Y
    \end{tikzcd}
  \]
  in which the composition of all the 2-morphisms agrees with the one filling the original lifting problem.
  Using the lifting property of $I'$, we obtain a diagram
  \[
    \begin{tikzcd}
      A \ar[r] \ar[rr, bend left, "\id_A"] \ar[d, "I"'] & A' \ar[r] \ar[d, "I'"] &
      A \ar[r, "H"] & X \ar[d, "P"] \\
      B \ar[r] \ar[rr, bend right, "\id_B"'] & B' \ar[r] \ar[rru, "L"] &
      B \ar[r, "K"'] & Y
    \end{tikzcd}
  \]
  in which the two triangles compose to the composition of the rightmost two squares of the previous diagram.
  The composition $B \to B' \xrightarrow{L} X$, together with the composition of the 2-morphisms above it and the composition of the 2-morphisms below it, then provide a solution to the original lifting problem for $I$.
  The case of $\rlp(\sL)$ is dual.
\end{proof}

\begin{proposition}
  For any class $\sR$, the class $\llp(\sR)$ is closed under coproducts, pushouts, and transfinite compositions, and dually for $\rlp(\sL)$ for any class $\sL$.
\end{proposition}

\begin{proof}
  The proofs are analogous to those for the 1-categorical case, with equalities between morphisms replaced by invertible 2-morphisms.
  We need only check that the resulting 2-morphisms compose to the one in the original lifting problem.
  We leave the details to the reader.
\end{proof}

\begin{definition}
  If $\sL$ is a class of 1-morphisms, then an \emph{$\sL$-cell morphism} is a transfinite composition of pushouts of coproducts of morphisms belonging to $\sL$.
\end{definition}

By the preceding proposition, if $\sL \subset \rlp(\sR)$, then any $\sL$-cell morphism also has the left lifting property with respect to $\sR$.

\begin{definition}
  A \emph{factorization} of a 1-morphism $F : X \to Y$ consists of an object $Z$, a 1-morphism $L : X \to Z$, a 1-morphism $R : Z \to Y$ and an isomorphism $\alpha : F \iso RL$.
  If $L \in \sL$ and $R \in \sR$ then we call the factorization a ``factorization of a morphism in $\sL$ followed by a morphism in $\sR$''.
\end{definition}

\begin{remark}
  Unlike lifting properties, the existence of a factorization of $F$ as a morphism in $\sL$ followed by a morphism in $\sR$ \emph{can} be expressed in terms of the homotopy category $\tau_{\le 1} \C$.
\end{remark}

\begin{definition}
  A \emph{weak factorization system} on $\C$ is a pair $(\sL, \sR)$ of classes of 1-morphisms of $\C$ such that
  \begin{enumerate}
  \item $\sL = \llp(\sR)$ and $\sR = \rlp(\sL)$.
  \item Each 1-morphism of $\C$ admits a factorization as a morphism in $\sL$ followed by a morphism in $\sR$.
  \end{enumerate}
\end{definition}

\begin{example}
  Suppose $\C$ is an ordinary category, viewed as a $(2, 1)$-category with only identity 2-morphisms.
  Then a weak factorization system on $\C$ is the same as a weak factorization system on the original 1-category in the usual sense.
\end{example}

\begin{example}
  Let $\C$ be any $(2, 1)$-category.
  We claim that $(\Eqv, \All)$ and $(\All, \Eqv)$ are weak factorization systems on $\C$, where $\Eqv$ is the class of equivalences in $\C$.
  We already checked that equivalences have the left and right lifting property with respect to all 1-morphisms, so to verify the first condition on a weak factorization system, it remains to show that any 1-morphism with the left (or right) lifting property with respect to all 1-morphisms is an equivalence.
  Indeed, suppose $I : A \to B$ has the left (or right) lifting property with respect to itself.
  Then the square (filled with the identity 2-morphism on $I$)
  \[
    \begin{tikzcd}
      A \ar[r, "\id_A"] \ar[d, "I"'] & A \ar[d, "I"] \\
      B \ar[r, "\id_B"'] & B
    \end{tikzcd}
  \]
  admits a lift
  \[
    \begin{tikzcd}
      A \ar[r, "\id_A"] \ar[d, "I"'] & A \ar[d, "I"] \\
      B \ar[r, "\id_B"'] \ar[ru, "L"] & B
    \end{tikzcd}
  \]
  and in particular $LI$ and $IL$ are isomorphic to identity morphisms.
  Hence $I$ is an equivalence.
  (More specifically, the data of a solution to the above lifting problem is precisely what is needed to make $I$ into a ``half-adjoint equivalence''.)
  The factorization condition is obvious, so both $(\Eqv, \All)$ and $(\All, \Eqv)$ are weak factorization systems.
\end{example}

\begin{example}\label{ex:vcpm-wfs}
  Let $\II$ be a set (or more generally a class) of morphisms of $\C$.
  Set $\sR = \rlp(\II)$ and $\sL = \llp(\sR)$.
  Then $\sR = \rlp(\sL)$, so if every morphism admits a factorization as a morphism in $\sL$ followed by a morphism in $\sR$, then $(\sL, \sR)$ is a weak factorization system on $\C$.
  In this case we call it the weak factorization system generated by $\II$.

  As in the 1-categorical case, under certain hypotheses on $\C$ and $\II$, such factorizations may be constructed using the small object argument.
  However, we have already noted in \cref{chap:size} that the small object argument cannot be applied directly in the case of most interest to us, the 2-category $V\CPM$.
  The next chapter will describe how to adapt the small object argument to $V\CPM$ under a mild additional assumption on the set $\II$.
\end{example}

\section{Premodel and model category structures}

Now, fix a complete and cocomplete 2-category $\C$.

\begin{remark}
  If $\C$ is complete and cocomplete as a 2-category, then the $(2, 1)$-category obtained by discarding the noninvertible 2-morphisms of $\C$ is again complete and cocomplete, with the same limits and colimits.
  Moreover, a model 2-category structure on $\C$ will be the same thing as a model 2-category structure on its underlying $(2, 1)$-category.
  The only additional requirement related to the full 2-categorical structure of $\C$ is that $\C$ itself admits all limits and colimits.
\end{remark}

\begin{definition}
  A \emph{premodel 2-category structure} on $\C$ is a pair of weak factorization systems $(\sC, \sAF)$ and $(\sAC, \sF)$ on $\C$ such that $\sAC \subset \sC$.
\end{definition}

As in the 1-categorical case, we call the morphisms of $\sC$ cofibrations, the morphisms of $\sAC$ anodyne cofibrations, the morphisms $\sF$ fibrations, and the morphisms of $\sAF$ anodyne fibrations.
We will use premodel 2-category structures for technical purposes in the course of the next chapter, but we are primarily interested in the notion of a model 2-category structure.

\begin{definition}
  Let $\sW$ be a class of 1-morphisms of $\C$.
  We say that $\sW$ satisfies the \emph{two-out-of-three property} if in any diagram
  \[
    \begin{tikzcd}[column sep=tiny]
      & Y \ar[rd, "G"] \\
      X \ar[ru, "F"] \ar[rr, "H"'{name=H}] & & Z
      \ar[from=1-2, to=H, shorten <=6pt, shorten >=6pt, Rightarrow, "\alpha", "\sim"']
    \end{tikzcd}
  \]
  if two of $F$, $G$, and $H$ belong to $\sW$, then so does the third.
\end{definition}

\begin{remark}
  Assume that $\sW$ also contains all identity morphisms (as will follow from the other requirements on the weak equivalences of a model 2-category structure).
  Then if $F : X \to Y$ and $F' : X \to Y$ are isomorphic parallel 1-morphisms, then $F$ belongs to $\sW$ if and only if $F'$ does.
  Hence $\sW$ consists of a union of isomorphism classes in $\THom_\C(X, Y)$ for each pair of objects $X$ and $Y$ of $\C$.
  Moreover, the two-out-of-three condition on such a class $\sW$ is equivalent to the ordinary two-out-of-three condition for the image of $\sW$ in the homotopy category $\tau_{\le 1} \C$.
\end{remark}

\begin{definition}\label{def:model-2-cat}
  A \emph{model 2-category structure} on $\C$ consists of classes of 1-morphisms $\sW$, $\sC$, and $\sF$ such that
  \begin{enumerate}
  \item $\sW$ satisfies the two-out-of-three property and is closed under retracts;
  \item $(\sC, \sF \cap \sW)$ and $(\sC \cap \sW, \sF)$ are weak factorization systems on $\C$.
  \end{enumerate}
\end{definition}

\begin{remark}
  Since $\sC \cap \sW$ is the left class of a weak factorization system, $\sW$ automatically must contain all the equivalences of $\C$ and is therefore also closed under replacing a 1-morphism by a parallel isomorphic one, as explained above.
\end{remark}

\begin{remark}
  \cite[Definition~1.1.1]{MG} defines a model $(\infty, 1)$-category by the traditional list of axioms.
  It is easy to see that in the case of a $(2, 1)$-category, \cref{def:model-2-cat} is a repackaged form of the same structure.
  Hence the theory of model $(\infty, 1)$-categories developed in \cite{MG} also applies to model 2-categories.

  In particular, we will make use of the notion of a Quillen bifunctor and the analogue of \cref{prop:quillen-bifunctor} for premodel 2-categories.
  This is treated in \cite[section~5.4]{MG}.

  In the opposite direction, the axioms of a model $(\infty, 1)$-category structure only involve the homotopy $(2, 1)$-category.
  (The retract and lifting axioms involve the existence of 3-morphisms between certain 2-morphisms; the two-out-of-three and factorization axioms only involve the existence of 2-morphisms between 1-morphisms.)
\end{remark}

\chapter{The large small object argument}
\label{chap:small}

Let $V$ be a combinatorial monoidal premodel category.
In this chapter we explain how to adapt the small object argument to the 2-category $V\CPM$.
We will use this ``large small object argument'' in the next chapter to construct the factorizations required for our model 2-category structure on $V\CPM$.

Let $\II$ be a set of morphisms of $V\CPM$.
We would like to use the small object argument to factor any morphism $F : M \to N$ of $V\CPM$ as an $\II$-cell complex followed by a morphism with the right lifting property with respect to $\II$.
(We use a bold letter for $\II$ to distinguish it from a set of generating cofibrations $I$ for some particular \emph{object} of $V\CPM$.)
As we explained in \cref{chap:size}, we cannot perform the small object argument in $V\CPM$ directly, both because the Hom categories of $V\CPM$ are large and because the objects of $V\CPM$ are not small in the required sense.
However, the sub-2-category $V\CPM_\lambda$ has neither of these deficiencies.
The usual small object argument applied to $V\CPM_\lambda$ then produces a factorization of $F : M \to N$ in which the first morphism is an $\II$-cell complex.
However, the second morphism of the factorization is a functor $F' : M' \to N$ which, in general, has a restricted right lifting property with respect to $\II$: we only know that we can solve lifting problems in which the horizontal morphisms belong to $V\CPM_\lambda$.

We would like to show that the functor $F'$ constructed in this way has the full right lifting property with respect to $\II$, that is, for all lifting problems in $V\CPM$.
In general this need not be the case.
As we will describe below, one can exhibit $\II$ and $F : M \to N$ for which there is \emph{no} factorization of $F$ of the expected form.
However, we will give a simple condition on $\II$ under which the restricted right lifting property implies the full one.
For such $\II$, then, the factorization constructed by applying the small object argument inside $V\CPM_\lambda$ does fulfill the original conditions.

\begin{nonexample}
  Take $V = \Set$ and $\II = \{0 \to \Set\}$.
  Consider the problem of factoring a morphism $0 \to N$ as an $\II$-cell complex $0 \to M$ followed by a functor $F : M \to N$ with the right lifting property with respect to $\II$.
  The latter condition means that $F^\cof : M^\cof \to N^\cof$ is essentially surjective.
  Clearly $M$ must have the form $\Set^{\oplus S}$ for some set $S$ and so the essential image of the functor $F : M \to N$ consists of all objects which can be written as a coproduct of copies of some set of fixed objects $(A_i)_{i \in S}$.
  For $N = \Kan$, however, there can be no set of (cofibrant) objects which generates all (cofibrant) objects under coproducts, as there exist connected simplicial sets of arbitrarily large cardinality.

  This example shows that some condition on $\II$ is needed to construct factorizations in $V\CPM$.
\end{nonexample}

\section{Reedy presentations}

In this chapter and the next, Reedy (or injective) premodel category structures on particular finite inverse categories with additional anodyne cofibrations adjoined turn out to play a major role.
We introduce a special notation for these and review their key properties.

\begin{notation}
  Let $D$ be a direct category.
  We regard $D$ as a Reedy category in which every nonidentity morphism increases degree (so $D^+ = D$).
  For any $V$-premodel category $M$, we write $M[D]_\R$ for the diagram category $M^{D^\op}_\Reedy$, where as usual $D^\op$ is equipped with the opposite Reedy category structure.
\end{notation}

\begin{remark}
  Since $D^\op$ is an inverse category, $M^{D^\op}_\Reedy$ is also the injective premodel category structure $M^{D^\op}_\inj$ by \cref{prop:reedy-direct-proj}.
  However, we will not use this fact often and therefore we mainly refer to $M[D]_\R$ as having the Reedy premodel category structure.
\end{remark}

The object $M[D]_\R$ of $V\CPM$ is also $\Set[D]_\R \otimes M$ and therefore obeys the adjointness relation
\[
  \THom_{V\CPM}(M[D]_R, N) \eqv \THom_{V\CPM}(M, N^D_\Reedy).
\]
Most often, $M$ will be $V$ itself and then we have
\[
  \THom_{V\CPM}(V[D]_R, N) \eqv \THom_{V\CPM}(V, N^D_\Reedy) \eqv (N^D_\Reedy)^\cof,
\]
the category of Reedy cofibrant diagrams of shape $D$ in $N$.
Here $D$ is direct, so the latching object $L_d X$ for $X$ in $N^D$ is the colimit over all nonidentity morphisms $d' \to d$ of $X_{d'}$.

\begin{notation}
  For $M$ an object of $V\CPM$, we write $M\langle A_1 \ancto B_1, A_2 \ancto B_2, \ldots \rangle$ for the object $M'$ of $V\CPM$ obtained by adjoining the morphisms of $S = \{A_1 \to B_1, A_2 \to B_2, \ldots\}$ as new generating anodyne cofibrations over $V$.
  This means that a morphism $F : M \to N$ of $V\CPM$ factors (essentially uniquely) through $M \to M'$ if and only if $F$ sends each morphism of $S$ to an anodyne cofibration of $N$.
  Concretely, writing $I_V$ for a set of generating cofibrations for $V$, this means that the underlying combinatorial premodel category of $M'$ is obtained from that of $M$ by adjoining $I_V \bp S$ as new generating anodyne cofibrations.
\end{notation}

In practice we will only use this notation when $M = V[D]_\R$ and the morphisms $A_i \to B_i$ are the images under the Yoneda embedding $\yo : D \to V[D]_\R$ of certain morphisms of $D$.
Then the category $\THom_{V\CPM}(V[D]_R\langle S \rangle, N)$ is equivalent to the full subcategory of $N^D_\Reedy$ on those diagrams in which the morphisms corresponding to elements of $S$ are not just cofibrations but anodyne cofibrations.

\section{The \texorpdfstring{$\lambda$}{\textlambda}-small small object argument}

Fix a regular cardinal $\lambda$ (as always, we assume that $\lambda$ is large enough that $V$ is $\lambda$-small $\lambda$-combinatorial monoidal) and a set $\II$ of morphisms of $V\CPM_\lambda$.
In this section, we will work entirely within $V\CPM_\lambda$.
To simplify the notation, we will write $\THom_\lambda(M, N)$ for the category $\THom_{V\CPM_\lambda}(M, N)$.

Recall that for morphisms $I : A \to B$ and $P : M \to N$ of $V\CPM_\lambda$, there is a category defined as the (pseudo)pullback of the square
\[
  \begin{tikzcd}
    \cdot \ar[r] \ar[d] & \THom_\lambda(A, N) \ar[d, "F_*"] \\
    \THom_\lambda(B, M) \ar[r, "I^*"'] & \THom_\lambda(A, M)
  \end{tikzcd}
\]
which we will denote by $\Sq_\lambda(I, P)$ and call the \emph{category of ($\lambda$-small) squares} from $I$ to $P$.
The commutative square with $\cdot$ above replaced by $\THom_\lambda(B, N)$ induces a functor $\LL_\lambda(I, P) : \THom_\lambda(B, N) \to \Sq_\lambda(I, P)$.
By \cref{prop:rlp-iff-ess-surj}, this functor is essentially surjective exactly when $P$ satisfies the right lifting property with respect to $I$ in $V\CPM_\lambda$.
In later sections, we will also describe this situation by saying that $P$ satisfies the ``$\lambda$-small right lifting property'' with respect to $I$.

\begin{proposition}\label{prop:small-small-object-argument}
  Let $P : M \to N$ be a morphism of $V\CPM_\lambda$.
  Then there exists a factorization (up to isomorphism) $M \xrightarrow{J} M' \xrightarrow{P'} N$ of $P$ in which
  \begin{enumerate}
  \item
    $J$ is an $\II$-cell map;
  \item
    $P'$ has the right lifting property (in $V\CPM_\lambda$) with respect to $\II$.
  \end{enumerate}
\end{proposition}

\begin{proof}
  The proof is just the usual small object argument.
  We may choose a regular cardinal $\gamma$ large enough so that the domains of the members of $\II$ are $\gamma$-small in $V\CPM_\lambda$.
  We construct a sequence of factorizations $M = M_0 \to M_1 \to \cdots \to M_\gamma \to N$ through morphisms $P_\alpha : M_\alpha \to N$ as follows.
  \begin{itemize}
  \item
    For $\alpha = 0$ we take $M_0 = M$ and $P_0 = P : M \to N$.
  \item
    At a limit stage $\beta$, we set $M_\beta = \colim_{\alpha < \beta} M_\alpha$ and let $P_\beta : M_\beta \to N$ be the morphism induced by the $P_\alpha$ for $\alpha < \beta$.
  \item
    At a successor stage, we define $M_\alpha \to M_{\alpha+1}$ as follows.
    Let $\mathbf{S}_\alpha$ denote the category of all lifting problems of the form
    \[
      \begin{tikzcd}
        A \ar[r] \ar[d] & M_\alpha \ar[d, "P_\alpha"] \\
        B \ar[r] & N
      \end{tikzcd}
    \]
    in which the morphism $A \to B$ belongs to $\II$.
    More precisely, $\mathbf{S}_\alpha = \coprod_{I \in \II} \Sq_\lambda(I, P_\alpha)$.
    The category $\mathbf{S}_\alpha$ is essentially small; choose a set $S_\alpha$ of representatives of its isomorphism classes.
    We then form a pushout
    \[
      \begin{tikzcd}
        \coprod_{s \in S_\alpha} A_s \ar[r] \ar[d] & M_\alpha \ar[d] \ar[rdd, "P_\alpha"] \\
        \coprod_{s \in S_\alpha} B_s \ar[r] \ar[rrd] & M_{\alpha+1} \ar[rd, dotted] \\
        & & N
      \end{tikzcd}
    \]
    and define the dotted morphism $P_{\alpha+1} : M_{\alpha+1} \to N$ using the morphisms $B_s \to N$ of each square $s \in S$.
  \end{itemize}
  We terminate the construction at $M' = M_\gamma$, yielding a factorization $M = M_0 \xrightarrow{J} M_\gamma = M' \xrightarrow{P'} N$.
  By definition, $J$ is an $\II$-cell map.
  Moreover, given any $\alpha < \gamma$ and any lifting problem
  \[
    \begin{tikzcd}
      A \ar[r] \ar[d, "I"'] & M_\alpha \ar[d, "P_\alpha"] \\
      B \ar[r] & N
    \end{tikzcd}
  \]
  with $I \in \II$, we may replace it by an isomorphic lifting problem which belongs to the set $S_\alpha$.
  Then by construction there exists a diagram
  \[
    \begin{tikzcd}
      A \ar[r] \ar[d, "I"'] &
      \coprod_{s \in S} A_s \ar[r] \ar[d] &
      M_\alpha \ar[d] \ar[rd, "P_\alpha"] \\
      B \ar[r] &
      \coprod_{s \in S} B_s \ar[r] &
      M_{\alpha+1} \ar[r, "P_{\alpha+1}"'] & N
    \end{tikzcd}
  \]
  composing to the original square.

  It remains to verify that $P' : M' \to N$ has the right lifting property with respect to $\II$ in $V\CPM_\lambda$.
  Consider any lifting problem
  \[
    \begin{tikzcd}
      A \ar[r, "H"] \ar[d, "I"'] & M' \ar[d, "P'"] \\
      B \ar[r, "K"'] & N
    \end{tikzcd}
  \]
  with $I \in \II$.
  Since $A$ is $\gamma$-small in $V\CPM_\lambda$, the morphism $H$ factors (up to isomorphism) through $M_\alpha \to M' = \colim_{\alpha < \gamma} M_\alpha$ for some $\alpha < \gamma$, and so we can factor the lifting problem as
  \[
    \begin{tikzcd}
      A \ar[r] \ar[rr, bend left=40, "H"] \ar[d, "I"'] &
      M_\alpha \ar[r] \ar[rd, "P_\alpha"] & M' \ar[d, "P'"] \\
      B \ar[rr, "K"'] & & N
    \end{tikzcd}
  \]
  and then as
  \[
    \begin{tikzcd}
      A \ar[r] \ar[rr, bend left=40, "H"] \ar[d, "I"'] &
      M_\alpha \ar[r] \ar[d] \ar[rd, "P_\alpha"] & M' \ar[d, "P'"] \\
      B \ar[r] \ar[rr, bend right=40, "K"'] &
      M_{\alpha+1} \ar[r, "P_{\alpha+1}"'] & N
    \end{tikzcd}
  \]
  by the construction of $M_{\alpha+1}$.
  Moreover, the two compositions of 2-morphisms
  \[
    \begin{tikzcd}
      M_\alpha \ar[r] \ar[d] \ar[rd, "F_\alpha"] & M' \ar[d, "P'"] \\
      M_{\alpha+1} \ar[r, "P_{\alpha+1}"'] & N
    \end{tikzcd}
    \quad\mbox{and}\quad
    \begin{tikzcd}
      M_\alpha \ar[r] \ar[d] & M' \ar[d, "P'"] \\
      M_{\alpha+1} \ar[r, "P_{\alpha+1}"'] \ar[ru] & N
    \end{tikzcd}
  \]
  agree by the coherence conditions on the data which define $P' = P_\gamma : M' \to N$ as induced by the $P_\alpha : M_\alpha \to N$.
  The composition $L : B \to M_{\alpha+1} \to M'$ (together with the obvious pastings of invertible 2-morphisms) then provides a solution to the original lifting problem.
\end{proof}

\begin{remark}
  This is just the small object argument for $(\infty, 1)$-categories \cite[Proposition~1.3.6]{MG} specialized to $V\CPM_\lambda$.
\end{remark}

\begin{remark}
  The inclusion $V\CPM_\lambda \to V\CPM$ preserves all colimits, so the functor $J$ is also an $\II$-cell map in $V\CPM$.
  However, $P'$ need not satisfy the right lifting property with respect to $\II$ in the full 2-category $V\CPM$.
\end{remark}



\section{The extensible right lifting property}

Let $\II$ be a set of morphisms of $V\CPM$.
In this section, we will describe a condition on $\II$ under which (for sufficiently large $\lambda$) any morphism $P : M \to N$ in $V\CPM_\lambda$ with the $\lambda$-small right lifting property with respect to $\II$ automatically has the full right lifting property with respect to $\II$, so that the factorization obtained by applying the small object argument in $V\CPM_\lambda$ also serves as a factorization of the desired form in $V\CPM$.
The condition on $\II$ is easy to check in practice and we will demonstrate it with examples in the next section.
The proof that our condition on $\II$ is in fact adequate will be postponed to the end of the chapter.

We begin by revisiting the reformulation of the right lifting property in terms of the $\Sq$ construction.
Let $I : A \to B$ and $P : M \to N$ be morphisms of $V\CPM$.
Recall that there is an associated functor $\THom(B, M) \to \Sq(I, P)$ (where for brevity we write $\THom$ for $\THom_{V\CPM}$) and that $P$ has the right lifting property with respect to $I$ if and only if this associated functor is essentially surjective.

The category $\THom(B, M)$ is also the full subcategory of cofibrant objects of a combinatorial premodel category $\CPM^V(B, M)$ whose underlying category is $\THom_{V\LPr}(B, M)$, the category of all left adjoint $V$-module functors from $B$ to $M$.
That is,
\[
  \THom(B, M) = \CPM^V(B, M)^\cof = \THom(V, \CPM^V(B, M)).
\]
In the same way, we can construct a combinatorial premodel category $\Sq\CPM^V(I, P)$ by forming the pullback in the square below.
We obtain an induced left Quillen functor $\LL(I, P)$ shown by the dotted arrow.
\[
  \begin{tikzcd}
    \CPM^V(B, M) \ar[rd, dotted] \ar[rrd, "I^*"] \ar[rdd, "P_*"'] \\
    & \Sq\CPM^V(I, P) \ar[r] \ar[d] & \CPM^V(A, M) \ar[d, "P_*"] \\
    & \CPM^V(B, N) \ar[r, "I^*"'] & \CPM^V(A, N)
  \end{tikzcd}
\]
Applying the (limit-preserving) functor $(-)^\cof = \THom(V, -)$ to this diagram recovers the diagram defining $\Sq(I, P)$ and the functor $\THom(B, M) \to \Sq(I, P)$.
Thus, we can identify the latter functor with $\LL(I, P)^\cof :\CPM^V(B, M)^\cof \to \Sq\CPM^V(I, P)^\cof$.
In particular, $P$ has the right lifting property with respect to $I$ if and only if $\LL(I, P)^\cof$ is essentially surjective.

We now define a modified lifting property which imposes a stronger condition on the left Quillen functor $\LL(I, P)$.

\begin{definition}\label{def:extensible}
  Let $F : M \to N$ be a left Quillen functor between premodel categories.
  We call $F$ \emph{extensible} if for every cofibrant object $m$ of $M$ and every cofibration $g : Fm \to n'$ of $N$, there exists a cofibration $f : m \to m'$ of $M$ lifting $g$ (up to isomorphism).
\end{definition}

\begin{proposition}\label{prop:ess-surj-of-extensible}
  If $F : M \to N$ is extensible, then $F^\cof : M^\cof \to N^\cof$ is essentially surjective.
\end{proposition}

\begin{proof}
  Take $m$ to be the initial object of $M$ in the definition of an extensible functor.
  Then it says that every cofibrant object $n'$ of $N$ lifts to a cofibrant object $m'$ of $M$.
\end{proof}

\begin{definition}
  We say that $P : M \to N$ satisfies the \emph{extensible right lifting property} with respect to $I : A \to B$ if $\LL(I, P)$ is extensible.
\end{definition}

The condition that a left Quillen functor is extensible can itself be expressed as a right lifting property.
(When every morphism of $M$ and $N$ is a cofibration, this condition is the one we considered in \cref{ex:r-rlp-cat}.)

\begin{definition}
  We define $E : V[\star]_\R \to V[\star \to \star']_\R$ to be the left Quillen $V$-functor induced by the functor sending $\star$ to $\star$.
\end{definition}

Let $F : M \to N$ be a left Quillen $V$-functor.
By the universal properties of the $V$-premodel categories $V[\star]_\R$ and $V[\star \to \star']_\R$, a lifting problem
\[
  \begin{tikzcd}
    V[\star]_\R \ar[r] \ar[d, "E"'] & M \ar[d, "F"] \\
    V[\star \to \star']_\R \ar[r] & N
  \end{tikzcd}
\]
amounts to specifying
\begin{enumerate}
\item a cofibrant object $m$ of $M$ (corresponding to the morphism $V[\star]_\R \to M$),
\item a cofibration $g : n \to n'$ between cofibrant objects of $N$ (corresponding to the morphism $V[\star \to \star']_\R \to N$), and
\item an isomorphism $Fm \iso n$, which we may as well take to be an equality.
\end{enumerate}
The existence of a lift $V[\star \to \star']_\R \to M$ corresponds to a cofibration $f : m \to m'$ which is sent by $F$ to $f$ (up to isomorphism).
Thus, $F$ has the right lifting property with respect to $E$ if and only if the underlying left Quillen functor of $F$ is extensible.

Taking $V = \Set$, we conclude that an ordinary left Quillen functor $F : M \to N$ is extensible if and only if it has the right lifting property with respect to the functor $E^\Set : \Set[\star]_\R \to \Set[\star \to \star']_\R$.

The extensible right lifting property can also be expressed as an ordinary right lifting property, as we explain next.

\begin{definition}\label{def:morphism-extension}
  Let $I : A \to B$ be a morphism of $V\CPM$.
  We define the \emph{extension} of $I$ to be the morphism $\E I$ of $V\CPM$ constructed as the induced map in the diagram below.
  \[
    \begin{tikzcd}
      A \ar[r] \ar[d, "I"'] &
      A[\star \to \star']_\R \ar[d] \ar[rdd, bend left=10, "{I[\star \to \star']_\R}"] \\
      B \ar[r] \ar[rrd, bend right=10] &
      B \amalg_A A[\star \to \star']_\R \ar[rd, "\E I"] \\
      & & B[\star \to \star']_\R
    \end{tikzcd}
  \]
  Here the morphism $A = A[\star]_\R \to A[\star \to \star']_\R$ is induced by the functor sending $\star$ to $\star$, and similarly for the morphism $B \to B[\star \to \star']_\R$.
\end{definition}

\begin{remark}
  If $I : A \to B$ belongs to $V\CPM_\lambda$, then so does $\E I$, since $V\CPM_\lambda$ is closed in $V\CPM$ under pushouts and the formation of Reedy premodel categories with finite index categories.
\end{remark}

\begin{proposition}\label{prop:erlp-iff-rlp-e}
  A morphism $P : M \to N$ has the extensible right lifting property with respect to $I : A \to B$ if and only if $P$ has the (ordinary) right lifting property with respect to $\E I$.
\end{proposition}

\begin{proof}
  
  We claim that lifting problems of the form
  \[
    \begin{tikzcd}
      \Set[\star]_\R \ar[r] \ar[d, "E^\Set"'] & \CPM^V(B, M) \ar[d, "{\LL(I, P)}"] \\
      \Set[\star \to \star']_\R \ar[r] & \CPM^V(A, M) \times_{\CPM^V(A, N)} \CPM^V(B, N)
    \end{tikzcd}
  \]
  and their solutions correspond precisely to lifting problems of the form
  \[
    \begin{tikzcd}
      B \amalg_A A[\star \to \star']_\R \ar[r] \ar[d, "\E I"'] & M \ar[d, "P"] \\
      B[\star \to \star']_\R \ar[r] & N
    \end{tikzcd}
  \]
  and their solutions.
  We repeatedly use the adjunction relation
  \begin{align*}
    \THom(V[D]_\R, \CPM^V(M, N))
    & = (\CPM^V(M, N)^D_\Reedy)^\cof \\
    & = \CPM^V(M[D]_\R, N)^\cof \\
    & = \THom(M[D]_\R, N),
  \end{align*}
  for various $D$, $M$ and $N$.
  Starting with the top square:
  \begin{itemize}
  \item
    A morphism $\Set = \Set[\star]_\R \to \CPM^V(B, M)$ corresponds to a morphism $B \to M$.
  \item
    A morphism $\Set[\star \to \star']_\R \to \CPM^V(A, M) \times_{\CPM^V(A, N)} \CPM^V(B, N)$ consists of compatible morphisms $\Set[\star \to \star']_\R \to \CPM^V(A, M)$ and $\Set[\star \to \star']_\R \to \CPM^V(B, M)$, which correspond to morphisms $A[\star \to \star']_\R \to M$ and $B[\star \to \star']_\R \to N$.
    The compatibility corresponds to an isomorphism filling the part of the bottom square shown below.
    \[
      \begin{tikzcd}
        A[\star \to \star']_\R \ar[r] \ar[d, "\E I"'] & M \ar[d, "P"] \\
        B[\star \to \star']_\R \ar[r] & N
      \end{tikzcd}
    \]
  \item The isomorphism filling the top square corresponds to two isomorphisms, between
    \begin{itemize}
    \item the compositions $A \to B \to M$ and $A \to A[\star \to \star']_\R \to M$, and
    \item the compositions $B \to M \to N$ and $B \to B[\star \to \star']_\R \to N$.
    \end{itemize}
    The first isomorphism allows to construct the top morphism $B \amalg_A A[\star \to \star']_\R \to M$ of the bottom square out of the morphisms $B \to M$ and $A[\star \to \star']_\R \to M$.
    The second isomorphism fills the remaining part of the bottom square.
    \[
      \begin{tikzcd}
        B \ar[r] \ar[d, "\E I"'] & M \ar[d, "P"] \\
        B[\star \to \star']_\R \ar[r] & N
      \end{tikzcd}
    \]
  \item
    The condition that the two isomorphisms
    \[
      \begin{tikzcd}
        \Set[\star]_\R \ar[r] \ar[d, "E^\Set"'] & \CPM^V(B, M) \ar[d, "{\LL(I, P)}"] \\
        \Set[\star \to \star']_\R \ar[r] & \CPM^V(A, M)
      \end{tikzcd}
      \quad\mbox{and}\quad
      \begin{tikzcd}
        \Set[\star]_\R \ar[r] \ar[d, "E^\Set"'] & \CPM^V(B, M) \ar[d, "{\LL(I, P)}"] \\
        \Set[\star \to \star']_\R \ar[r] & \CPM^V(B, N)
      \end{tikzcd}
    \]
    agree after composing with the morphisms to $\CPM^V(A, N)$ is equivalent to the condition that the two isomorphisms
    \[
      \begin{tikzcd}
        B \ar[r] \ar[d, "\E I"'] & M \ar[d, "P"] \\
        B[\star \to \star']_\R \ar[r] & N
      \end{tikzcd}
      \quad\mbox{and}\quad
      \begin{tikzcd}
        A[\star \to \star']_\R \ar[r] \ar[d, "\E I"'] & M \ar[d, "P"] \\
        B[\star \to \star']_\R \ar[r] & N
      \end{tikzcd}
    \]
    agree after precomposing with the morphisms from $A$.
  \end{itemize}
  Thus, the two lifting problems encode the same data.
  Similarly, a lift $\Set[\star \to \star']_\R \to \CPM^V(B, M)$ corresponds to a lift $B[\star \to \star']_\R \to M$.
  (We omit the detailed verification that all of the compatibility conditions on the lifts correspond precisely, which is similar to the above.)
\end{proof}

\begin{remark}\label{remark:extension-and-po-product}
  The extensible right lifting property is a kind of ``enriched lifting property'' and \cref{prop:erlp-iff-rlp-e} is another instance of the multiplicative structure of weak factorization systems which underlies \cref{prop:quillen-bifunctor}.
  $\E I$ is the pushout product $E^\Set \bp I$ induced by the tensor product $\otimes : \CPM \times V\CPM \to V\CPM$.
\end{remark}

Using the $\E$ construction, we can now formulate the condition we will impose on $\II$ to guarantee the existence of factorizations.

\begin{definition}\label{def:self-extensible}
  A set $\II$ of morphisms of $V\CPM$ is \emph{self-extensible} if, for each $I \in \II$, the morphism $\E I$ is an $\II$-cell morphism.
\end{definition}

\begin{remark}
  A union of self-extensible sets is evidently again self-extensible.
\end{remark}

We will prove later in this chapter that when $\II$ is self-extensible, we can construct factorizations in $V\CPM$ by applying the small object argument in $V\CPM$ for sufficiently large $\lambda$.
Roughly speaking, the reason is that (unlike essential surjectivity) the extensibility of a morphism $F : M \to N$ of $\CPM_\lambda$ can be checked on the $\lambda$-compact objects.
By applying this fact to the functor $\LL(I, P) : \CPM^V(B, M) \to \Sq\CPM^V(I, P)$, we can relate extensivity of $\LL(I, P)$ (hence in particular essential surjectivity of $\LL(I, P)^\cof$) to the $\lambda$-small right lifting property which is the output of the $\lambda$-small small object argument.

For this purpose, we will actually need the following technical variant of \cref{prop:erlp-iff-rlp-e}.

\begin{definition}
  A left Quillen functor $F : M \to N$ between premodel categories is \emph{$\lambda$-extensible} if for every \emph{$\lambda$-compact} cofibrant object $m$ of $M$ and every cofibration $g : Fm \to n'$ of $N$ to a \emph{$\lambda$-compact} object $n'$, there exists a cofibration $f : m \to m'$ to a \emph{$\lambda$-compact} object $m'$ of $M$ lifting $g$ (up to isomorphism).

  A morphism $P : M \to N$ of $V\CPM_\lambda$ has the \emph{$\lambda$-small extensible right lifting property} with respect to $I$ if $\LL(I, P)$ is $\lambda$-extensible.
\end{definition}

\begin{proposition}\label{prop:erlp-iff-rlp-e-lambda}
  Let $\lambda$ be an uncountable regular cardinal and suppose that $I : A \to B$ is a morphism of $V\CPM_\lambda$ with $A$ and $B$ $\lambda$-small $\lambda$-combinatorial $V$-premodel categories.
  Then, for any morphism $P : M \to N$ of $V\CPM_\lambda$:
  \begin{enumerate}
  \item The left Quillen functor $\LL(I, P)$ belongs to $\CPM_\lambda$.
  \item $P$ has the $\lambda$-small extensible right lifting property if and only if $P$ has the (ordinary) $\lambda$-small right lifting property with respect to $\E I$.
  \end{enumerate}
\end{proposition}

\begin{proof}
  For $D$ a $\lambda$-small direct category, a left Quillen functor from $\Set[D]_\R$ to a $\lambda$-combinatorial premodel category $N$ is strongly $\lambda$-accessible if and only if the corresponding (Reedy cofibrant) diagram $D \to N$ consists of $\lambda$-compact objects of $N$.
  It follows that a morphism $F : M \to N$ of $\CPM_\lambda$ is $\lambda$-extensible if and only if it has the right lifting property with respect to $E^\Set$ in $\CPM_\lambda$.
  
  Furthermore, when $M$ is a $\lambda$-small $\lambda$-combinatorial $V$-premodel category and $N$ is any object of $V\CPM_\lambda$, the premodel category $\CPM^V(M, N)$ is again $\lambda$-combinatorial and its $\lambda$-compact objects of $\CPM^V(M, N)$ are precisely those left adjoint $V$-functors which are strongly $\lambda$-accessible, i.e., preserve $\lambda$-compact objects.
  In particular, when $D$ is a $\lambda$-small direct category, so that $V[D]_\R$ and $M[D]_\R$ are $\lambda$-small $\lambda$-combinatorial, the adjunction
  \begin{align*}
    \THom(V[D]_\R, \CPM^V(M, N))
    & = (\CPM^V(M, N)^D_\Reedy)^\cof \\
    & = \CPM^V(M[D]_\R, N)^\cof \\
    & = \THom(M[D]_\R, N),
  \end{align*}
  which played the central role in the proof of \cref{prop:erlp-iff-rlp-e}, restricts to an adjunction
  \begin{align*}
    \THom_\lambda(V[D]_\R, \CPM^V(M, N))
    & = (\CPM^V(M, N)^D_\Reedy)^\cof_\lambda \\
    & = \CPM^V(M[D]_\R, N)^\cof_\lambda \\
    & = \THom_\lambda(M[D]_\R, N),
  \end{align*}
  where as earlier we write $\THom_\lambda$ for $\THom_{V\CPM_\lambda}$.
  Now because $A$ and $B$ are assumed to be $\lambda$-small $\lambda$-combinatorial, the pullback square defining $\Sq\CPM^V(I, P)$ belongs to $V\CPM_\lambda$ and therefore is also a pullback square in $V\CPM_\lambda$ since $V\CPM_\lambda$ is closed in $V\CPM$ under finite limits.
  Then the same argument as in the proof of \cref{prop:erlp-iff-rlp-e} shows that $\LL(I, P)$ has the right lifting property with respect to $E^\Set$ in $\CPM_\lambda$ if and only if $P$ has the right lifting property with respect to $\E I$ in $V\CPM_\lambda$.
\end{proof}

\begin{remark}
  There are alternative, \emph{a priori} slightly weaker conditions on $\II$ that are also sufficient for constructing factorizations.
  For example, we can require that any morphism with the right lifting property with respect to $\II$ also has the extensible right lifting property with respect to $\II$.
  However, in our application, it will turn out to be most convenient to verify the condition of \cref{def:self-extensible} anyways.
\end{remark}

\section{Checking self-extensibility}\label{sec:check-se}

We next demonstrate how to verify the condition of self-extensibility in two examples.
These examples are not chosen arbitrarily---they will be turn out to be the generating cofibrations of the model 2-category structure on $V\CPM$ constructed in the next chapter.

\subsection{Self-extensibility of $\{E\}$}

\newcommand{\Dsquare}{
  \begin{tikzpicture}[baseline={([yshift=-.5ex]current bounding box.center)}, inner sep=2pt, scale=0.75]
    \node at (0, 0) (0){$0$};
    \node at (1, 0) (1){$1$};
    \node at (0, -1) (0'){$0'$};
    \node at (1, -1) (1'){$1'$};
    \draw[->](0)--(1);
    \draw[->](0')--(1');
    \draw[->](0)--(0');
    \draw[->](1)--(1');
  \end{tikzpicture}
}
\newcommand{\Dcorner}{
  \begin{tikzpicture}[baseline={([yshift=-.5ex]current bounding box.center)}, inner sep=2pt, scale=0.75]
    \node at (0, 0) (0){$0$};
    \node at (1, 0) (1){$1$};
    \node at (0, -1) (0'){$0'$};
    \draw[->](0)--(1);
    \draw[->](0)--(0');
  \end{tikzpicture}
}  

We begin with the functor $E$ itself.
In order to avoid notational confusion caused by the two different roles played by $E$, we will instead verify the self-extensibility of $\{E'\}$ for the following morphism $E'$ which is obviously identical to $E$.

\begin{notation}
  Write $E' : V[0]_\R \to V[0 \to 1]_\R$ for the morphism of $V\CPM$ induced by the functor sending $0$ to $0$.
\end{notation}

We first compute the morphism $\E E'$ using \cref{def:morphism-extension}.
An iterated Reedy premodel category is equivalent to the Reedy premodel category on the product of the index categories, so the codomain of $\E E'$ is
\[
  (V[0 \to 1]_\R)[\star \to \star']_\R = V\Big[\Dsquare\Big]_\R.
\]
More generally, we have $(V[D]_\R)[\star \to \star']_\R = V[D \times \{\star \to \star'\}]_\R$.
Above we have adopted the following notational convention.

\begin{notation}\label{not:times-e}
  Suppose that $D$ is some explicit category whose objects we have named, such as $D = \{0 \to 1\}$ above.
  Then we name the objects of the product category $D \times \{\star \to \star'\}$ as follows.
  The object $(d, \star)$ is assigned the same name as $d$, while the object $(d, \star')$ is assigned the name of $d$ with a prime.
  This convention reflects the fact that we regard $D$ as contained in $D \times \{\star \to \star'\}$ via the functor $d \mapsto (d, \star)$.
\end{notation}


Turning now to the domain of $\E E'$, we claim that the square below is a pushout.
\[
  \begin{tikzcd}
    V[0]_\R \ar[r] \ar[d] & V[0 \to 0']_\R \ar[d] \\
    V[0 \to 1]_\R \ar[r] & V\Big[\Dcorner\Big]_\R
  \end{tikzcd}
\]
Indeed, we can verify the required universal property using the universal property of the objects $V[-]_\R$.
A morphism from the pushout of the square to an object $N$ of $V\CPM$ is supposed to correspond to
\begin{enumerate}
\item a cofibration $n \to n'$ between cofibrant objects of $N$,
\item a cofibration $n_0 \to n_1$ between cofibrant objects of $N$, and
\item an isomorphism $n \iso n_0$.
\end{enumerate}
The category of such morphisms is equivalent to the category of diagrams
\[
  \begin{tikzcd}
    n_0 \ar[r] \ar[d] & n_1 \\
    n'_0
  \end{tikzcd}
\]
of cofibrant objects and cofibrations of $N$, and this in turn is equivalent to $\THom(V\Big[\Dcorner\Big]_\R, N)$.

\begin{proposition}
  $\E E'$ is the morphism
  \[
    V\Big[\Dcorner\Big]_\R \to V\Big[\Dsquare\Big]_\R
  \]
  induced by the obvious inclusion of indexing categories.
\end{proposition}

\begin{proof}
  We have already identified the domain and codomain of $\E E'$.
  Identifying $\E E'$ itself is also best done in terms of universal properties.
  Specifically, we must show that precomposing a morphism $F : V\Big[\Dsquare\Big]_\R \to N$ with $\E E'$ corresponds to taking a Reedy cofibrant diagram
  \[
    \begin{tikzcd}
      n_0 \ar[r] \ar[d] & n_1 \ar[d] \\
      n'_0 \ar[r] & n'_1
    \end{tikzcd}
  \]
  in $N$ and discarding the part of the diagram involving $n'_1$.
  This follows easily from unfolding the definition of $\E E'$ and using the universal properties of all the constructions involved.
\end{proof}

Next we will show that there exists a pushout square
\begin{equation}
  \begin{tikzcd}
    V[0]_\R \ar[r] \ar[d, "E'"'] & V\Big[\Dcorner\Big]_\R \ar[d, "\E E'"] \\
    V[0 \to 1]_R \ar[r] & V\Big[\Dsquare\Big]_\R
  \end{tikzcd}
  \tag{$*$}
\end{equation}
which exhibits $\E E'$ as an $E'$-cell morphism.
To construct the top morphism of this pushout square, first recall that the equivalence between $\THom(V[D]_\R, N)$ and Reedy cofibrant $D$-indexed diagrams in $N$ is induced by a functor $\yo : D \to V[D]_\R$, the ($V$-valued) Yoneda embedding.
In particular, there is a ``tautological'' Reedy cofibrant diagram
\[
  \begin{tikzcd}
    \yo 0 \ar[r] \ar[d] & \yo 1 \\
    \yo 0'
  \end{tikzcd}
\]
in $V\Big[\Dcorner\Big]_\R$.
Form the pushout $p = \yo 1 \amalg_{\yo 0} \yo 0'$ of this diagram; then $p$ is cofibrant.
Moreover, for a morphism $F : V\Big[\Dcorner\Big]_\R \to N$ corresponding to a Reedy cofibrant diagram
\[
  \begin{tikzcd}
    n_0 \ar[r] \ar[d] & n_1 \\
    n'_0
  \end{tikzcd}
\]
in $N$, $F$ sends the object $p$ to the pushout $n_1 \amalg_{n_0} n'_0$.

Because $p$ is cofibrant, we obtain a morphism $V = V[0]_\R \to V\Big[\Dcorner\Big]_\R$ corresponding to the object $p$.
The functor $\E E'$ sends $p$ to the object $p' = \yo 1 \amalg_{\yo 0} \yo 0'$ of $V\Big[\Dsquare\Big]_\R$.
Since the tautological diagram of the latter category
\[
  \begin{tikzcd}
    \yo 0 \ar[r] \ar[d] & \yo 1 \ar[d] \\
    \yo 0' \ar[r] & \yo 1'
  \end{tikzcd}
\]
is Reedy cofibrant, the induced map $f : p' \to \yo 1'$ is a cofibration between cofibrant objects.
We can therefore define a morphism $V[0 \to 1]_\R \to V\Big[\Dsquare\Big]_\R$ sending the generating morphism $0 \to 1$ to $f : p' \to \yo 1'$, and this defines a square of the form $(*)$ above.

\begin{proposition}
  The square $(*)$ constructed in this way is a pushout in $V\CPM$.
\end{proposition}

\begin{proof}
  Again this is easily verified using universal properties.
  A morphism $V\Big[\Dsquare\Big]_\R \to N$ is the same as a Reedy cofibrant diagram
  \[
    \begin{tikzcd}
      n_0 \ar[r] \ar[d] & n_1 \ar[d] \\
      n'_0 \ar[r] & n'_1
    \end{tikzcd}
  \]
  in $N$.
  By the inductive description of Reedy cofibrant diagrams, giving such a diagram is the same as giving a Reedy cofibrant diagram
  \[
    \begin{tikzcd}
      n_0 \ar[r] \ar[d] & n_1 \\
      n'_0
    \end{tikzcd}
  \]
  together with an object $n'_1$ and a cofibration (latching map)
  \[
    n_1 \amalg_{n_0} n'_0 \to n'_1
  \]
  and this is equivalent to the data encoded in a square
  \[
    \begin{tikzcd}
      V[0]_\R \ar[r, "\yo 0 \mapsto p"] \ar[d, "E'"'] & V\Big[\Dcorner\Big]_\R \ar[d] \\
      V[0 \to 1]_\R \ar[r] & N
    \end{tikzcd}
  \]
  in $V\CPM$.
\end{proof}

Switching back to the original notation $E$, we record the result of this subsection below.

\begin{proposition}
  $\E E$ is an $E$-cell morphism, and so $\{E\}$ is self-extensible.
\end{proposition}

\begin{proof}
  This follows from the above construction of $\E E'$ as a pushout of $E'$.
\end{proof}

\begin{remark}\label{remark:z-e}
  The ``unit morphism'' $Z : 0 \to V$ of $V\CPM$ is also an $E$-cell morphism.
  In fact, there is a pushout square
  \[
    \begin{tikzcd}
      V[0]_\R \ar[r] \ar[d, "E"'] & 0 \ar[d, "Z"] \\
      V[0 \to 0']_\R \ar[r] & V
    \end{tikzcd}
  \]
  which already implicitly appeared in the proof of \cref{prop:ess-surj-of-extensible}.
  In the other direction, it is easy to see that $E = \E Z$.
  Therefore $\{E\}$ is a kind of ``self-extensible closure'' of $\{Z\}$, in the sense that $\{E\}$ is self-extensible and the extensible right lifting properties with respect to $\{Z\}$ and $\{E\}$ are equivalent.

  More generally, if $I : A \to B$ is any morphism of $V\CPM$, then one can show by a similar argument that $\E (\E I)$ is an $\E I$-cell morphism.
  (Informally, because $\E E^\Set = E^\Set \bp E^\Set$ is a pushout of $E^\Set$, it follows that $\E (\E I) = E^\Set \bp (E^\Set \bp I) = (E^\Set \bp E^\Set) \bp I$ is also a pushout of $\E I = E^\Set \bp I$, where $\bp$ is as described in \cref{remark:extension-and-po-product}.)
  Thus for any set $\II$, $\E \II = \{\,\E I \mid I \in \II\,\}$ is self-extensible.
  Similarly, $I$ is also an $\E I$-cell morphism and so the extensible right lifting property for $\E \II$ is equivalent to that for $\II$.
  In particular, $\E \II$ is self-extensible.
  Thus (using the large small object argument, which we will finish later in this chapter) \emph{any} set $\II$ ``generates'' a weak factorization system on $V\CPM$ in which the morphisms of the right class have the extensible right lifting property with respect to $\II$.
  However, the morphisms of the left class of this weak factorization system are not necessarily (retracts of) $\II$-cell morphisms, but only $\E \II$-cell morphisms.

  In our application we will arrange that the generating cofibrations $\II$ and the generating acyclic cofibrations $\JJ$ for $V\CPM$ are already self-extensible, so we will not need to use the argument outlined in the preceding paragraph.
\end{remark}

\subsection{Self-extensibility of $\{\Lambda\}$}\label{subsec:lambda}

To verify self-extensibility of $\{E\}$ we only had to use facts about Reedy premodel category structures.
Our second example will also involve anodyne cofibrations.

\begin{notation}
  We write $\Lambda$ for the morphism of $V\CPM$
  \[
    \Lambda : V[0 \to 1]_\R \to V[0 \to 1]_\R\langle \yo 0 \ancto \yo 1 \rangle
  \]
  induced by the identity functor of $\{0 \to 1\}$.
\end{notation}

Recall that $V[0 \to 1]_\R\langle \yo 0 \ancto \yo 1 \rangle$ is the $V$-premodel category formed by adding the morphism $\yo 0 \to \yo 1$ to the generating anodyne cofibrations (over $V$) to $V[0 \to 1]_\R$.
Thus, whereas a morphism $V[0 \to 1]_\R \to N$ corresponds to a cofibration between cofibrant objects of $N$, a morphism $V[0 \to 1]_\R\langle \yo 0 \ancto \yo 1 \rangle \to N$ corresponds to an \emph{anodyne} cofibration between cofibrant objects of $N$.
A morphism $P : M \to N$ has the right lifting property with respect to $\Lambda$ if every cofibration $f : m_0 \to m_1$ of $M^\cof$ which is sent by $P$ to an anodyne cofibration of $N$ is already an anodyne cofibration in $M$.


As before, our first task is to calculate a presentation for the codomain
\[
  (V[0 \to 1]_\R\langle \yo 0 \ancto \yo 1 \rangle)[\star \to \star']_\R
\]
of $\E \Lambda$.
For any object $N$ of $V\CPM$, giving a morphism from this object to $N$ is equivalent to giving a morphism from $V[0 \to 1]_\R\langle \yo 0 \ancto \yo 1 \rangle$ to $N^{\{\star \to \star'\}}_\Reedy$.
In turn, this amounts to giving an anodyne cofibration between cofibrant objects of this premodel category $N^{\{\star \to \star'\}}_\Reedy$.
By definition of the Reedy premodel category structure, this amounts to a square
\[
  \begin{tikzcd}
    n_0 \ar[r] \ar[d] & n_1 \ar[d] \\
    n'_0 \ar[r] & n'_1
  \end{tikzcd}
\]
in $N$ which is Reedy cofibrant and in which the morphism $n_0 \to n_1$ and the induced morphism $n'_0 \amalg_{n_0} n_1 \to n'_1$ are anodyne cofibrations.
Thus, we can describe the codomain $(V[0 \to 1]_\R\langle \yo 0 \ancto \yo 1 \rangle)[\star \to \star']_\R$ of $\E \Lambda$ by the presentation $V\Big[\Dsquare\Big]_\R\langle \yo 0 \ancto \yo 1, \yo 0' \amalg_{\yo 0} \yo 1 \ancto \yo 1' \rangle$.

\begin{remark}
  More generally, given an object of $V\CPM$ of the form
  \[
    A = V[D]_\R\langle x_i \ancto y_i \mid i \in I \rangle,
  \]
  an analogous argument shows that
  \[
    A[\star \to \star']_\R =
    V[D \times \{\star \to \star'\}]_\R
    \langle x_i \ancto y_i, x'_i \amalg_{x_i} y_i \ancto y'_i \mid i \in I \rangle.
  \]
  Here $x'_i$ and $y'_i$ denote the images of $x_i$ and $y_i$ respectively under the ``prime'' functor induced by $d \mapsto (d, \star')$.
\end{remark}

The square
\[
  \begin{tikzcd}
    V[0 \to 1]_\R \ar[r] \ar[d, "\Lambda"'] & V\Big[\Dsquare\Big]_\R \ar[d] \\
    V[0 \to 1]_\R\langle \yo 0 \ancto \yo 1 \rangle \ar[r] &
    V\Big[\Dsquare\Big]_\R \langle \yo 0 \ancto \yo 1 \rangle
  \end{tikzcd}
\]
is a pushout, and so we can identify the morphism $\E \Lambda$ as the morphism
\[
  V\Big[\Dsquare\Big]_\R \langle \yo 0 \ancto \yo 1 \rangle \to
  V\Big[\Dsquare\Big]_\R \langle \yo 0 \ancto \yo 1, \yo 0' \amalg_{\yo 0} \yo 1 \ancto \yo 1' \rangle
\]
induced by the identity functor of $\Dsquare$.
(Hence, $\E \Lambda$ is the identity on the underlying $V$-module categories.)

Now, as in the previous subsection, let $f$ be the morphism $\yo 0' \amalg_{\yo 0} \yo 1 \to \yo 1'$ of the category $V\Big[\Dsquare\Big]_\R \langle \yo 0 \ancto \yo 1 \rangle$ induced by the tautological square.
Again, $f$ is a cofibration, and now $\E \Lambda$ takes $f$ to an anodyne cofibration.
Therefore, there is a square
\[
  \begin{tikzcd}
    V[0 \to 1]_\R \ar[r, "f"] \ar[d, "\Lambda"'] &
    V\Big[\Dsquare\Big]_\R \langle \yo 0 \ancto \yo 1 \rangle \ar[d, "\E \Lambda"] \\
    V[0 \to 1]_\R \langle \yo 0 \ancto \yo 1 \rangle \ar[r, "(\E \Lambda) f"'] &
    V\Big[\Dsquare\Big]_\R \langle \yo 0 \ancto \yo 1, \yo 0' \amalg_{\yo 0} \yo 1 \ancto \yo 1' \rangle
  \end{tikzcd}
\]
This square is a pushout, as can be seen immediately from the universal properties.
We have thus proved the following.

\begin{proposition}
  $\E \Lambda$ is a $\Lambda$-cell morphism, and so $\{\Lambda\}$ is self-extensible.
\end{proposition}

\begin{proof}
  The above square exhibits $\E \Lambda$ as a pushout of $\Lambda$.
\end{proof}

\section{The large small object argument}

Our final task in this chapter is to explain how to factor a given morphism $F : M \to N$ of $V\CPM$ as an $\II$-cell morphism followed by a morphism with the right lifting property with respect to $\II$ when $\II$ is self-extensible.
As mentioned earlier, the key technical advantage of the extensible right lifting property is that the property of being extensible can be checked on $\lambda$-compact objects.

\begin{proposition}\label{prop:extensible-of-lambda-extensible}
  Let $F : M \to N$ be a strongly $\lambda$-accessible left Quillen functor between $\lambda$-combinatorial premodel categories, with $\lambda$ an uncountable regular cardinal.
  If $F$ is $\lambda$-extensible, then $F$ is extensible.
\end{proposition}

\begin{proof}
  Suppose that $m$ is a cofibrant object of $M$ and $g : Fm \to n'$ is a cofibration of $N$.
  Applying \cite[A.1.5.12]{HTT} with $\overline S$ the class of all cofibrations of $N$, so that $S$ is the set of cofibrations between $\lambda$-compact objects of $N$, we can express $g : Fm \to n'$ as a transfinite composition
  \[
    Fm = n_0 \to n_1 \to \cdots \to n_\gamma = n'
  \]
  such that for each $\alpha < \gamma$, the morphism $n_\alpha \to n_{\alpha+1}$ is a pushout of a member of $S$, that is, a pushout of some cofibration between $\lambda$-compact objects.
  We will lift this composition to a transfinite composition of cofibrations
  \[
    m = m_0 \to m_1 \to \cdots \to m_\gamma = m'
  \]
  by transfinite induction, starting by setting $m_0 = m$.
  At a limit stage $\beta \le \gamma$, having already lifted the portion of the sequence preceding $n_\beta$, we have $n_\beta = \colim_{\alpha < \beta} n_\alpha$, so (as $F$ preserves colimits) we may take $m_\beta = \colim_{\alpha < \beta} m_\alpha$.
  This process terminates with the desired lift $m \to m'$ of $Fm \to n'$.
  Thus, it suffices to describe how to carry out each successor stage, that is, given $m_\alpha$ (which is cofibrant) and the cofibration $g_\alpha : Fm_\alpha = n_\alpha \to n_{\alpha+1}$, how to lift $g_\alpha$ to a cofibration $f_\alpha : m_\alpha \to m_{\alpha+1}$.
  This problem has the same form as the original one, except that we also know that the map $g_\alpha$ is the pushout of some cofibration $k : a \to a'$ between $\lambda$-compact objects of $N$.
  We thus return to the original notation, writing $m$ for $m_\alpha$, $n'$ for $n_{\alpha+1}$ and $g$ for $g_\alpha$.

  By \cite[Corollary~5.1]{FSOA} the cofibrant object $m$ of $M$ may be written as a $\lambda$-directed colimit $m = \colim_{i \in I} m_i$ of $\lambda$-compact cofibrant objects of $M$.
  Write $\sigma_j : m_j \to \colim_{i \in I} m_i = m$ for the cocone maps of this colimit.
  Now in the pushout square
  \[
    \begin{tikzcd}
      a \ar[r, "k"] \ar[d] & a' \ar[d] \\
      Fm \ar[r, "g"'] & n'
    \end{tikzcd}
  \]
  we have $Fm = F(\colim_{i \in I} m_i) = \colim_{i \in I} Fm_i$.
  Since $I$ is $\lambda$-directed and $a$ is $\lambda$-compact, the left morphism $a \to Fm$ factors through $F\sigma_j : Fm_j \to Fm$ for some $j \in I$.
  Set $n'_j = Fm_j \amalg_a a'$, factoring the above pushout square as a composition of two pushout squares as shown below.
  \[
    \begin{tikzcd}
      a \ar[r, "k"] \ar[d] & a' \ar[d] \\
      Fm_j \ar[r, "g_j"] \ar[d, "F\sigma_j"'] & n'_j \ar[d] \\
      Fm \ar[r, "g"'] & n'
    \end{tikzcd}
  \]
  The map $g_j : Fm_j \to n'_j$ is a pushout of $k$ and therefore a cofibration.
  Moreover, $n'_j$ is $\lambda$-compact because $m_j$ (and hence $Fm_j$), $a$, and $a'$ are.
  Therefore, applying the hypothesis on $F$ to $m_j$ and the cofibration $g_j$, we obtain a cofibration $f_j : m_j \to m'_j$ which lifts $g_j$, producing the diagram consisting of two pushout squares below.
  \[
    \begin{tikzcd}
      a \ar[r, "k"] \ar[d] & a' \ar[d] \\
      Fm_j \ar[r, "Ff_j"] \ar[d, "F\sigma_j"'] & Fm'_j \ar[d] \\
      Fm \ar[r, "g"'] & n'
    \end{tikzcd}
  \]
  Now we can define $m'$ and a cofibration $f : m \to m'$ by forming the pushout
  \[
    \begin{tikzcd}
      m_j \ar[r, "f_j"] \ar[d, "\sigma_j"'] & m'_j \ar[d] \\
      m \ar[r, "f"'] & m'
    \end{tikzcd}
  \]
  and because $F$ preserves pushouts, $f$ is a lift of $g$ (up to isomorphism).
\end{proof}

\begin{remark}
  The condition in the definition of a $\lambda$-extensible left Quillen functor that the lift $m'$ of $n'$ is a $\lambda$-compact object is not used in the above proof.
  We only included it in the definition of $\lambda$-extensible in order to relate $\lambda$-extensibility to a $\lambda$-small right lifting property.
\end{remark}

This fact allows us to relate the $\lambda$-small right lifting property with respect to a self-extensible $\II$ to the full right lifting property with respect to $\II$, provided $\lambda$ is sufficiently large.

\begin{lemma}\label{lemma:erlp-of-srlp}
  Let $\II$ be a self-extensible set of morphisms of $V\CPM$.
  Then for all sufficiently large regular cardinals $\lambda$, any morphism $P : M \to N$ in $V\CPM_\lambda$ with the $\lambda$-small right lifting property with respect to $\II$ has the extensible right lifting property with respect to $\II$.
  (In particular, such a morphism $P$ also has the ordinary right lifting property with respect to $\II$.)
\end{lemma}

\begin{proof}
  We will assume that $\lambda$ is an uncountable regular cardinal so large that all of the following conditions are satisfied.
  \begin{enumerate}
  \item
    Every morphism $I$ of $\II$ belongs to $V\CPM_\lambda$.
  \item
    For each morphism $I : A \to B$ of $\II$, both $A$ and $B$ are $\lambda$-small $\lambda$-combinatorial $V$-premodel categories.
  \item
    For each morphism $I$ of $\II$, the morphism $\E I$ is an $\II$-cell morphism inside $V\CPM_\lambda$.
    ($\E I$ is an $\II$-cell morphism in $V\CPM$ by the assumption that $\II$ is self-extensible.
    Once $\lambda$ is large enough that the diagram exhibiting $\E I$ as an $\II$-cell complex belongs to $V\CPM_\lambda$, the same diagram also exhibits $\E I$ as an $\II$-cell complex in $V\CPM_\lambda$ because colimits in $V\CPM_\lambda$ are the same as those in $V\CPM$.)
  \end{enumerate}

  Now, suppose $P : M \to N$ is a morphism of $V\CPM_\lambda$ with the $\lambda$-small right lifting property with respect to $\II$.
  By definition, this simply means that $P$ has the right lifting property with respect to $\II$ within $V\CPM_\lambda$.
  Since (by (3)) each morphism $\E I$ for $I \in \II$ is an $\II$-cell morphism within  $V\CPM_\lambda$, $P$ also has the right lifting property with respect to each $\E I$ within $V\CPM_\lambda$.
  By \cref{prop:erlp-iff-rlp-e-lambda}, $\LL(I, P)$ is then $\lambda$-extensible and belongs to $\CPM_\lambda$.
  Then by \cref{prop:extensible-of-lambda-extensible}, $\LL(I, P)$ is extensible and so $P$ has the extensible right lifting property with respect to $I$.
\end{proof}

We can now complete the large small object argument.

\begin{proposition}
  Let $\II$ be a self-extensible set of morphisms of $V\CPM$.
  Then any morphism $F : M \to N$ of $V\CPM$ admits a factorization as an $\II$-cell morphism followed by a morphism with the extensible right lifting property (hence also the ordinary right lifting property) with respect to $\II$.
\end{proposition}

\begin{proof}
  Choose a regular cardinal $\lambda$ large enough to apply \cref{lemma:erlp-of-srlp} to $\II$ and also large enough so that $F$ belongs to $V\CPM_\lambda$.
  The $\lambda$-small small object argument (\cref{prop:small-small-object-argument}) produces a factorization of $F$ as an $\II$-cell morphism followed by a morphism $F'$ with the $\lambda$-small right lifting property with respect to $\II$.
  By \cref{lemma:erlp-of-srlp}, $F'$ then also has the extensible right lifting property with respect to $\II$.
\end{proof}

\begin{corollary}
  Let $\II$ be a self-extensible set of morphisms of $V\CPM$.
  Then $(\llp(\rlp(\II)), \rlp(\II))$ is a weak factorization system on $V\CPM$.
\end{corollary}

\begin{proof}
  As we noted in \cref{ex:vcpm-wfs}, the only nontrivial part is the existence of factorizations and these have been constructed above.
\end{proof}

\chapter{The model 2-category $V\CPM$}
\label{chap:vcpmmodel}

Let $V$ be a tractable\footnote{
  A combinatorial model category is \emph{tractable} if it admits a set of generating cofibrations with cofibrant domains \cite[Corollary~2.7]{Barw}.
}
symmetric monoidal model category.
In this chapter, we construct a model 2-category structure on $V\CPM$ whose weak equivalences are the Quillen equivalences.
This structure is defined by explicit generating cofibrations and acyclic cofibrations.
An object $M$ of $V\CPM$ is fibrant if and only if every trivial cofibration of $M^\cof$ is an anodyne cofibration.
In particular, each $V$-\emph{model} category is a fibrant object of $V\CPM$.

This chapter is organized as follows.
The model 2-category structure we define on $V\CPM$ will be induced from a premodel 2-category structure on $\CPM$ along the right adjoint $V\CPM \to \CPM$ which forgets the $V$-module structure.
In other words, a morphism $F : X \to Y$ of $V\CPM$ is a fibration or an acyclic fibration if and only if the underlying left Quillen functor is a fibration or an anodyne fibration in $\CPM$.
We begin by introducing the sets $\II$ and $\JJ$ which generate the premodel 2-category structure on $\CPM$ and verifying that they satisfy the self-extensibility condition required for the large small object argument.
It is then a formal matter to transfer this premodel 2-category structure to $V\CPM$ for any combinatorial monoidal premodel category $V$.
It remains to show that, under the additional assumptions that $V$ is a tractable symmetric monoidal model category, this premodel 2-category structure on $V\CPM$ is actually a model 2-category structure whose weak equivalences are precisely the Quillen equivalences.
The two-out-of-three and retract axioms for the weak equivalences are obvious by definition, so this amounts to verifying that the equalities $\sAC = \sC \cap \sW$ and $\sAF = \sF \cap \sW$ hold in $V\CPM$.

The fibrations and anodyne fibrations are characterized by lifting properties which allow the condition $\sAF = \sF \cap \sW$ to be checked directly.
The generating cofibrations $\II$ and anodyne cofibrations $\JJ$ have been chosen so that the (anodyne) fibrations of $V\CPM$ resemble the (acyclic) fibrations of the fibration category of cofibration categories constructed by Szumi\l o in \cite{Sz}, allowing us to imitate the proof of the corresponding fact for cofibration categories.
We must however adjust the conditions defining the fibrations to ensure that they can be stated as right lifting properties with respect to morphisms of $V\CPM$.
(The fibrations of a fibration category do not have an analogous requirement.)

For the condition $\sAC = \sC \cap \sW$ we argue indirectly.
By a standard argument, it suffices to check that an anodyne cofibration is a weak equivalence.
First, we show that any morphism $M \to X$ of $V\CPM$ to a fibrant object factors through a Quillen equivalence $M \to \hat M$ to a fibrant object.
For this we apply the small object argument using just one of the generating anodyne cofibrations, whose cell complexes are easy to understand.
Consequently, every morphism $M \to N$ can be approximated by a morphism $\hat M \to \hat N$ between fibrant objects.
Second, we construct an explicit ``mapping path category'' factorization of any morphism between fibrant objects as a Quillen equivalence followed by a fibration.
When the original morphism $M \to N$ is an anodyne cofibration it has the left lifting property with respect to the resulting fibration and this allows us to conclude that the original morphism is a weak equivalence by the two-out-of-six property.

As mentioned above, the generating cofibrations $\II$ and anodyne cofibrations $\JJ$ are chosen so as to determine lifting properties similar to those of \cite[Proposition~1.11]{Sz} and \cite[Definition~1.9]{Sz} respectively.
In the remainder of this introductory section, we present some \emph{a priori} considerations by which one might arrive at this particular set $\II$---and hence also, if not the specific set $\JJ$, at least the particular model 2-category structure we will consider, which is determined by $\II$ and the weak equivalences.

\begin{itemize}
\item
  The object $V$ is surely an obvious candidate for a cofibrant object of $V\CPM$.
  In other words, the morphism $0 \to V$ (which we named $Z$ in \cref{remark:z-e}) ought to be a cofibration.

\item
  We intend to construct factorizations using the large small object argument of \cref{chap:small}.
  We cannot apply the large small object argument to $\{Z\}$ because it is not self-extensible.
  Instead, we should take the morphism $\E Z = E$ as a generating cofibration.
  (Since $Z$ is a pushout of $E$, we do not need to also include $Z$ among the generating cofibrations.)

\item
  Suppose that $F : M \to N$ is an acyclic fibration in $V\CPM$ between two \emph{model} $V$-categories.
  We don't yet know what the acyclic fibrations of $V\CPM$ should be, but we do know that they have to be Quillen equivalences.
  Consider a lifting problem of the form below, involving $F$ and the morphism $\Lambda$ from \cref{subsec:lambda}.
  \[
    \begin{tikzcd}
      V[0 \to 1]_\R \ar[r, "f"] \ar[d, "\Lambda"'] & M \ar[d, "F"] \\
      V[0 \to 1]_\R\langle \yo 0 \ancto \yo 1 \rangle \ar[r] & N
    \end{tikzcd}
  \]
  A square of this form corresponds to a cofibration $f$ between cofibrant objects of $M$ which is sent by $F$ to an anodyne cofibration.
  Since $N$ is a model category, an anodyne cofibration of $N$ is just an acyclic cofibration.
  Because $F$ is a Quillen equivalence, it reflects weak equivalences between cofibrant objects.
  Therefore, $f$ is also an acyclic cofibration in $M$ or equivalently an anodyne cofibration.
  Hence, $f$ also defines a functor $V[0 \to 1]_\R\langle \yo 0 \ancto \yo 1 \rangle \to M$, which provides a lift in the original square.

  At this point, we can only make this argument under the assumption that $M$ and $N$ are model categories.
  Still, it provides evidence that the morphism $\Lambda$ ought to be a cofibration in $V\CPM$.
\end{itemize}

These observations suggest that $\II = \{E, \Lambda\}$ is the smallest reasonable choice for the generating cofibrations of a model 2-category structure on $V\CPM$, and these will indeed be the generating cofibrations of the structure we construct in this chapter.

\begin{remark}
  The weak equivalences of $V\CPM$ will be the Quillen equivalences.
  In order to define what it means for a left Quillen functor $F : M \to N$ to be a Quillen equivalence we need $M$ and $N$ to be relaxed premodel categories, so that they have well-behaved homotopy theories.
  Since $V$ is a monoidal model category, every $V$-premodel category is relaxed.
  The $V$-module structures play no other role in defining the weak equivalences of $V\CPM$.
  Specifically, whether a morphism $F : M \to N$ in $V\CPM$ is a Quillen equivalence depends only on the underlying left Quillen functor between relaxed premodel categories, and not on the $V$-module structures of $M$, $N$, or $F$.
  Similarly, the condition for a morphism $F : M \to N$ of $V\CPM$ to be a fibration or anodyne fibration will not depend on the $V$-module structures, only on the underlying left Quillen functor.
  For these reasons, we will sometimes blur the distinction between morphisms of $V\CPM$ and ordinary left Quillen functors when discussing weak equivalences and fibrations.
\end{remark}

\begin{remark}
  We will use the hypothesis that $V$ is a model category (or at least a relaxed premodel category) in order to define the weak equivalences of $V\CPM$, and the remaining hypotheses on $V$ to show that anodyne cofibrations are Quillen equivalences.
  However, it will be useful to carry out the construction of the \emph{premodel} 2-category structure on $V\CPM$ in greater generality.
  In particular, we will initially want to include the case $V = \Set$.
  Thus, until further notice, we assume only that $V$ is a \emph{combinatorial monoidal premodel category}, not necessarily a model category, symmetric monoidal, or tractable.
\end{remark}

\section{The premodel 2-category structure on $\CPM$}

We begin in the case $V = \Set$.
In this section, we'll construct a particular premodel 2-category structure on $\CPM$.
In the next section, we will describe how to transfer this premodel 2-category structure to $V\CPM$ for a general monoidal combinatorial premodel category $V$ along the forgetful functor $V\CPM \to \CPM$.
Our main theorem is that when $V$ is a tractable symmetric monoidal model category, this transferred structure is a model 2-category structure on $V\CPM$ whose weak equivalences are the Quillen equivalences.

In order to define a premodel 2-category structure on $\CPM$ we simply have to give a nested pair of weak factorization systems.
Naturally we intend to apply the large small object argument to construct the required factorizations, so our plan is as follows.
\begin{itemize}
\item
  Write down two particular sets $\II$ and $\JJ$ of morphisms of $\CPM$.
  The choices of these sets will be justified in the rest of this chapter.
\item
  Define a morphism to be an anodyne fibration (respectively, fibration) if it has the right lifting property with respect to $\II$ (respectively, $\JJ$), and a cofibration (respectively, anodyne cofibration) if it has the left lifting property with respect to all anodyne fibrations (respectively, fibrations).
\item
  Verify that each member of $\JJ$ is an $\II$-cell morphism, so that every anodyne fibration is a fibration.
\item
  Verify that each of the sets $\II$ and $\JJ$ is self-extensible.
\item
  Apply the large small object argument of \cref{chap:small} to factor an arbitrary morphism of $\CPM$ into an $\II$-cell morphism (hence a cofibration) followed by an acyclic fibration, or a $\JJ$-cell morphism (hence an anodyne cofibration) followed by a fibration.
\end{itemize}
Furthermore, we will verify that this premodel 2-category structure on $\CPM$ is ``monoidal'' in the sense that $\II \bp \II$ consists of cofibrations and $\II \bp \JJ$ (and $\JJ \bp \II$) consists of anodyne cofibrations.
Here $\bp$ is taken with respect to the tensor product $\otimes$ on $\CPM$.
This will be required for the analysis of the path category construction on $V\CPM$ near the end of the chapter.

\begin{definition}
  We write $\II = \{E, \Lambda\}$ (with notation as in \cref{sec:check-se}).
  A morphism of $\CPM$ is an \emph{anodyne fibration} if it has the right lifting property with respect to $\II$ and a \emph{cofibration} if it has the left lifting property with respect to all anodyne fibrations.
\end{definition}

We have already verified in \cref{sec:check-se} that $\II$ is self-extensible, so it generates a weak factorization system on $\CPM$.

We now turn to the generating anodyne cofibrations.

\begin{definition}
  We write
  \[
    \Sigma_L :
    \Set[0 \to 1 \to 2]_\R\langle \yo 0 \ancto \yo 2, \yo 1 \ancto \yo 2 \rangle \to
    \Set[0 \to 1 \to 2]_\R\langle \yo 0 \ancto \yo 1, \yo 1 \ancto \yo 2 \rangle
  \]
  and
  \[
    \Sigma_R :
    \Set[0 \to 1 \to 2]_\R\langle \yo 0 \ancto \yo 2, \yo 0 \ancto \yo 1 \rangle \to
    \Set[0 \to 1 \to 2]_\R\langle \yo 0 \ancto \yo 1, \yo 1 \ancto \yo 2 \rangle
  \]
  for the morphisms of $\CPM$ induced by the identity functor of $\{0 \to 1 \to 2\}$.
\end{definition}
Each of $\Sigma_L$ and $\Sigma_R$ adjoins a single generating anodyne cofibration: $\yo 0 \ancto \yo 1$ in the case of $\Sigma_L$ and $\yo 1 \ancto \yo 2$ in the case of $\Sigma_R$.
(We omit $\yo 0 \ancto \yo 2$ from the generating anodyne cofibrations of the codomains of $\Sigma_L$ and $\Sigma_R$ because it is redundant, being the composition of $\yo 0 \ancto \yo 1$ and $\yo 1 \ancto \yo 2$.)

\newcommand{\Dtwosquare}{
  \begin{tikzpicture}[baseline={([yshift=-.5ex]current bounding box.center)}, inner sep=2pt, scale=0.75]
    \node at (0, 0) (0){$0$};
    \node at (1, 0) (1){$1$};
    \node at (2, 0) (2){$2$};
    \node at (0, -1) (0'){$0'$};
    \node at (1, -1) (1'){$1'$};
    \node at (2, -1) (2'){$2'$};
    \draw[->](0)--(1);
    \draw[->](1)--(2);
    \draw[->](0')--(1');
    \draw[->](1')--(2');
    \draw[->](0)--(0');
    \draw[->](1)--(1');
    \draw[->](2)--(2');
  \end{tikzpicture}
}

\begin{proposition}\label{prop:sigma-l-self-extensible}
  $\E \Sigma_L$ is a pushout of $\Sigma_L$, so $\{\Sigma_L\}$ is self-extensible.
\end{proposition}

\begin{proof}
  Using the techniques of \cref{sec:check-se}, we compute that $\E \Sigma_L$ is the morphism
  \[
    \E \Sigma_L : \Set\Big[\Dtwosquare\Big]_\R\langle S \rangle \to
    \Set\Big[\Dtwosquare\Big]_\R\langle S, \yo 0' \amalg_{\yo 0} \yo 1 \ancto \yo 1' \rangle
  \]
  induced by the identity, where
  \[
    S = \{\yo 0 \ancto \yo 1, \yo 1 \ancto \yo 2,
    \yo 0' \amalg_{\yo 0} \yo 2 \ancto \yo 2',
    \yo 1' \amalg_{\yo 1} \yo 2 \ancto \yo 2'\}.
  \]
  Thus, let $f : \yo 0' \amalg_{\yo 0} \yo 1 \to \yo 1'$ be the morphism in the domain of $\E \Sigma_L$ which becomes an anodyne cofibration in the codomain of $\E \Sigma_L$.
  This $f$ is a cofibration between cofibrant objects.
  To express $\E \Sigma_L$ as a pushout of $\Sigma_L$, it suffices to exhibit an anodyne cofibration out of $\yo 1'$ whose composition with $f$ is also an anodyne cofibration.
  For this, we can take $g : \yo 1' \to \yo 2'$.
  It is the composition of the anodyne cofibrations $\yo 1' \ancto \yo 1' \amalg_{\yo 1} \yo 2$ (a pushout of $\yo 1 \ancto \yo 2$) and $\yo 1' \amalg_{\yo 1} \yo 2 \ancto \yo 2'$, and the composition $gf : \yo 0' \amalg_{\yo 0} \yo 1 \to \yo 2'$ is also the composition of the anodyne cofibrations $\yo 0' \amalg_{\yo 0} \yo 1 \ancto \yo 0' \amalg_{\yo 0} \yo 2$ (a pushout of $\yo 1 \ancto \yo 2$) and $\yo 0' \amalg_{\yo 0} \yo 2 \ancto \yo 2'$.
\end{proof}

\begin{proposition}
  $\E \Sigma_R$ is a pushout of $\Sigma_R$, so $\{\Sigma_R\}$ is self-extensible.
\end{proposition}

\begin{proof}
  In a similar manner, we compute $\E \Sigma_R$ as the morphism
  \[
    \E \Sigma_R : \Set\Big[\Dtwosquare\Big]_\R \langle S \rangle \to
    \Set\Big[\Dtwosquare\Big]_\R \langle S, \yo 1' \amalg_{\yo 1} \yo 2 \ancto \yo 2' \rangle
  \]
  induced by the identity, for
  \[
    S = \{\yo 0 \ancto \yo 1, \yo 1 \ancto \yo 2,
    \yo 0' \amalg_{\yo 0} \yo 1 \ancto \yo 1',
    \yo 0' \amalg_{\yo 0} \yo 2 \ancto \yo 2'\}.
  \]
  Setting $g : \yo 1' \amalg_{\yo 1} \yo 2 \to \yo 2'$, it suffices to find an anodyne cofibration $f$ between cofibrant objects of the domain of $\E \Sigma_R$ such that $gf$ is again an anodyne cofibration.
  Form the diagram below, in which every square is a pushout.
  \[
    \begin{tikzcd}
      \yo 0 \ar[r] \ar[d] & \yo 1 \ar[r] \ar[d] & \yo 2 \ar[d] \\
      \yo 0' \ar[r] \ar[rd] & \cdot \ar[r] \ar[d] & \cdot \ar[d, "f"] \\
      & \yo 1' \ar[r] \ar[rd] & \cdot \ar[d, "g"] \\
      & & \yo 2'
    \end{tikzcd}
  \]
  The map labeled $f$ is an anodyne cofibration because it is a pushout of $\yo 0' \amalg_{\yo 0} \yo 1 \ancto \yo 1'$, and its composition with $g$ is $\yo 0' \amalg_{\yo 0} \yo 2 \ancto \yo 2'$.
\end{proof}

\newcommand{\Dpcorner}{
  \begin{tikzpicture}[baseline={([yshift=-.5ex]current bounding box.center)}, inner sep=2pt, scale=0.75]
    \node at (0, 0) (00){$00$};
    \node at (1, 0) (01){$01$};
    \node at (0, -1) (10){$10$};
    \draw[->](00)--(01);
    \draw[->](00)--(10);
  \end{tikzpicture}
}

\newcommand{\Dpsquare}{
  \begin{tikzpicture}[baseline={([yshift=-.5ex]current bounding box.center)}, inner sep=2pt, scale=0.75]
    \node at (0, 0) (00){$00$};
    \node at (1, 0) (01){$01$};
    \node at (0, -1) (10){$10$};
    \node at (1, -1) (11){$11$};
    \draw[->](00)--(01);
    \draw[->](10)--(11);
    \draw[->](00)--(10);
    \draw[->](01)--(11);
  \end{tikzpicture}
}

\newcommand{\Dccorner}{
  \begin{tikzpicture}[baseline={([yshift=-.5ex]current bounding box.center)}, inner sep=2pt, scale=1.3]
    \node at (0, 0) (00){$00$};
    \node at (1, 0) (01){$01$};
    \node at (0, -1) (10){$10$};
    \node at (1, -1) (11){$11$};
    \node at (0.5, 0.5) (00'){$00'$};
    \node at (1.5, 0.5) (01'){$01'$};
    \draw[->](00)--(01);
    \draw[->](10)--(11);
    \draw[->](00)--(10);
    \draw[->](01)--(11);
    \draw[->](00')--(01');
    \draw[->](00)--(00');
    \draw[->](01)--(01');
  \end{tikzpicture}
}

\newcommand{\Dcube}{
  \begin{tikzpicture}[baseline={([yshift=-.5ex]current bounding box.center)}, inner sep=2pt, scale=1.3]
    \node at (0, 0) (00){$00$};
    \node at (1, 0) (01){$01$};
    \node at (0, -1) (10){$10$};
    \node at (1, -1) (11){$11$};
    \node at (0.5, 0.5) (00'){$00'$};
    \node at (1.5, 0.5) (01'){$01'$};
    \node at (0.5, -0.5) (10'){$10'$};
    \node at (1.5, -0.5) (11'){$11'$};
    \draw[->](00)--(01);
    \draw[->](10)--(11);
    \draw[->](00)--(10);
    \draw[->](01)--(11);
    \draw[->](00')--(01');
    \draw[->](10')--(11');
    \draw[->](00')--(10');
    \draw[->](01')--(11');
    \draw[->](00)--(00');
    \draw[->](01)--(01');
    \draw[->](10)--(10');
    \draw[->](11)--(11');
  \end{tikzpicture}
}

\begin{definition}
  We write
  \[
    \Psi : \Set[00 \to 01]_R \to
    \Set\Big[\Dpsquare\Big]_\R \langle \yo 10 \ancto \yo 11, \yo 01 \ancto \yo 11 \rangle.
  \]
\end{definition}

\begin{proposition}\label{prop:psi-self-extensible}
  $\{\Sigma_R, \Psi\}$ is self-extensible.
\end{proposition}

\begin{proof}
  We will show that $\E \Psi$ is the composition of a pushout of $\Psi$ and a pushout of $\Sigma_R$.
  This is sufficient since we already showed that $\Sigma_R$ is self-extensible above.

  $\E \Psi$ is the morphism
  \[
    \E \Psi : \Set\Bigg[\Dccorner\Bigg]_\R \langle S_1 \rangle \to
    \Set\Bigg[\Dcube\Bigg]_\R \langle S_2 \rangle
  \]
  for
  \begin{align*}
    S_1 & = \{\yo 10 \ancto \yo 11, \yo 01 \ancto \yo 11\}, \\
    S_2 & = \{\yo 10 \ancto \yo 11, \yo 01 \ancto \yo 11,
    \yo 11 \amalg_{\yo 10} \yo 10' \ancto \yo 11',
    \yo 11 \amalg_{\yo 01} \yo 01' \ancto \yo 11'\}.
  \end{align*}
  We will first attach a copy of $\Psi$ to the domain of $\E \Psi$ so that the pushout is nearly the codomain of $\E \Psi$, but with slightly different generating anodyne cofibrations.
  In the domain of $\E \Psi$, form the pushouts $x_0 = \yo 10 \amalg_{\yo 00} \yo 00'$ and $x_1 = \yo 11 \amalg_{\yo 01} \yo 01'$.
  The objects $x_0$ and $x_1$ are cofibrant and by the gluing lemma for cofibrations \cite[Lemma~7.2.15]{Hi}, the induced map $h : x_0 \to x_1$ is a cofibration.
  We can then map the domain $\Set[00 \to 01]_\R$ of $\Psi$ to the domain of $\E \Psi$ along this map $h$.
  We claim that the pushout of $\Psi$ along this morphism has the form
  \[
    \begin{tikzcd}
      \Set[00 \to 01]_R \ar[r] \ar[d, "\Psi"'] &
      \Set\Bigg[\Dccorner\Bigg]_\R \langle S_1 \rangle \ar[d] \\
      \Set\Big[\Dpsquare\Big]_\R \langle \yo 10 \ancto \yo 11, \yo 01 \ancto \yo 11 \rangle \ar[r] &
      \Set\Bigg[\Dcube\Bigg]_\R \langle S_3 \rangle
    \end{tikzcd}
  \]
  for
  \[
    S_3 = \{\yo 10 \ancto \yo 11, \yo 01 \ancto \yo 11,
    \yo 10' \ancto \yo 11', \yo 11 \amalg_{\yo 01} \yo 01' \ancto \yo 11'\}.
  \]
  Ignoring for a moment the generating anodyne cofibrations, the claim amounts to the fact that giving a Reedy cofibrant cube in a premodel category $N$
  \[
    \begin{tikzcd}[row sep=1em, column sep=1em]
      & n'_{00} \ar[rr] \ar[dd] & & n'_{01} \ar[dd] \\
      n_{00} \ar[rr] \ar[dd] \ar[ru] & & n_{01} \ar[dd] \ar[ru] \\
      & n'_{10} \ar[rr] & & n'_{11} \\
      n_{10} \ar[rr] \ar[ru] & & n_{11} \ar[ru]
    \end{tikzcd}
  \]
  is equivalent to giving a Reedy cofibrant diagram
  \[
    \begin{tikzcd}[row sep=1em, column sep=1em]
      & n'_{00} \ar[rr] & & n'_{01} \\
      n_{00} \ar[rr] \ar[dd] \ar[ru] & & n_{01} \ar[dd] \ar[ru] \\
      \\
      n_{10} \ar[rr] & & n_{11}
    \end{tikzcd}
  \]
  together with a Reedy cofibrant square
  \[
    \begin{tikzcd}
      n_{10} \amalg_{n_{00}} n'_{00} \ar[r] \ar[d] &
      n_{11} \amalg_{n_{01}} n'_{01} \ar[d] \\
      n'_{10} \ar[r] & n'_{11}
    \end{tikzcd}
  \]
  which can be verified by assigning the cube a degree function in which the objects of highest degree are $10'$ and then $11'$.
  The claim about the generating anodyne cofibrations $S_3$ of the pushout follows, taking into account the fact that the bottom morphism of the above pushout square sends $\yo 10$ to $\yo 10'$, $\yo 01$ to $x_1 = \yo 11 \amalg_{\yo 01} \yo 01'$, and $\yo 11$ to $\yo 11'$.
  For convenience, let us call the combinatorial premodel category that we have just constructed $M$.

  Now $M$ does not have quite the same anodyne cofibrations as the codomain of $\E \Psi$; we need to adjoin the additional cofibration $g : \yo 11 \amalg_{\yo 10} \yo 10' \to \yo 11'$ as an anodyne cofibration.
  However, in $M$, the map $f : \yo 10' \to \yo 11 \amalg_{\yo 10} \yo 10'$ is the pushout of $\yo 10 \ancto \yo 11$, hence an anodyne cofibration, and its composition $gf : \yo 10' \to \yo 11'$ is an anodyne cofibration.
  Therefore, we can attach a copy of $\Sigma_R$ along $f$ and $g$ to adjoin $g$ as an anodyne cofibration, realizing $\E \Psi$ as the composition of a pushout of $\Psi$ and a pushout of $\Sigma_R$.
\end{proof}

\begin{definition}
  We write $\JJ = \{\Sigma_L, \Sigma_R, \Psi\}$.
  A morphism of $\CPM$ is a \emph{fibration} if it has the right lifting property with respect to $\JJ$ and an \emph{anodyne cofibration} if it has the left lifting property with respect to all fibrations.
\end{definition}

\begin{proposition}
  $\JJ$ is self-extensible.
\end{proposition}

\begin{proof}
  This follows from \cref{prop:sigma-l-self-extensible,prop:psi-self-extensible}.
\end{proof}

Thus $\JJ$ also generates a weak factorization system on $\CPM$.

\begin{proposition}
  Each member of $\JJ$ is an $\II$-cell morphism.
\end{proposition}

\begin{proof}
  The morphisms $\Sigma_L$ and $\Sigma_R$ each adjoin an anodyne cofibration which was already a cofibration between cofibrant objects, so they are each pushouts of $\Lambda \in \II$.
  For $\Psi$, we must first build the correct underlying locally presentable category by attaching two copies of $E$.
  First form a pushout
  \[
    \begin{tikzcd}
      \Set[0]_\R \ar[r] \ar[d, "E"'] & \Set[00 \to 01]_\R \ar[d] \\
      \Set[0 \to 0']_\R \ar[r] & \Set\Big[\Dpcorner\Big]_\R
    \end{tikzcd}
  \]
  in which the top morphism sends $\yo 0$ to $\yo 00$.
  Next, let $x$ be the cofibrant object $\yo 10 \amalg_{\yo 00} \yo 01$ of $\Set\Big[\Dpcorner\Big]_\R$ and form a pushout
  \[
    \begin{tikzcd}
      \Set[0]_\R \ar[r] \ar[d, "E"'] & \Set\Big[\Dpcorner\Big]_\R \ar[d] \\
      \Set[0 \to 0']_\R \ar[r] & \Set\Big[\Dpsquare\Big]_\R
    \end{tikzcd}
  \]
  in which the top morphism sends $\yo 0$ to $x$.
  We verified that both these squares are pushouts in \cref{sec:check-se}.
  Finally, attach two copies of $\Lambda$ along the cofibrations $\yo 10 \to \yo 11$ and $\yo 01 \to \yo 11$ between cofibrant objects, to impose the correct generating anodyne cofibrations.
  The composition of these morphisms is $\Psi$.
\end{proof}

In particular, every anodyne fibration is a fibration and we have thus constructed a premodel 2-category structure on $\CPM$.
To finish this section, we show that this premodel 2-category structure is monoidal.

\begin{proposition}\label{prop:cpm-monoidal}
  Each morphism of $\II \bp \II$ is an $\II$-cell morphism, and each morphism of $\II \bp \JJ$ is a $\JJ$-cell morphism.
\end{proposition}

\begin{proof}
  Recall that $\II = \{E, \Lambda\}$.
  We have already checked all the conditions involving $E \bp -$, because they amount to self-extensibility of $\II$ and $\JJ$.
  Thus, it remains to check the conditions involving $\Lambda \bp -$.
  For notational clarity, we will rename the index category appearing in the definition of $\Lambda$ so that
  \[
    \Lambda : \Set[\star \to \star']_\R \to
    \Set[\star \to \star']_\R \langle \yo \star \ancto \yo \star' \rangle,
  \]
  and we again adopt the convention of \cref{not:times-e} for naming the objects of a product category $D \times \{\star \to \star'\}$.

  By the formula for the tensor product of combinatorial premodel categories, for any combinatorial premodel category $M$, the tensor product $\Set[\star \to \star']_\R \langle \yo \star \ancto \yo \star' \rangle \otimes M$  differs from $\Set[\star \to \star'] \otimes M = M[\star \to \star']_\R$ in that the former has an additional generating anodyne cofibration $A' \amalg_A B \ancto B'$ for each generating anodyne cofibration $A \ancto B$ of $M$.
  Consequently, for any morphism $F : M \to N$ of $\CPM$, the morphism $\Lambda \bp F$ has the form $N[\star \to \star']_\R \langle S_1 \rangle \to N[\star \to \star']_\R \langle S_2 \rangle$ where
  \begin{align*}
    S_1 & = \{\, A' \amalg_A B \ancto B' \mid \mbox{$A \to B$ is the image under $F$ of a generating cofibration of $M$} \,\}, \\
    S_2 & = \{\, A' \amalg_A B \ancto B' \mid \mbox{$A \to B$ is a generating cofibration of $N$} \,\}.
  \end{align*}
  Now the morphisms $\Lambda$, $\Sigma_L$, and $\Sigma_R$ each adjoin only anodyne cofibrations; their domains and codomains have the same underlying categories and generating cofibrations.
  Therefore each of $\Lambda \bp \Lambda$, $\Lambda \bp \Sigma_L$, and $\Lambda \bp \Sigma_R$ is actually an equivalence in $\CPM$.

  It remains to consider $\Lambda \bp \Psi$.
  By the above formula, this morphism has the form
  \[
    \Set\Bigg[\Dcube\Bigg]_\R\langle S_1 \rangle \to
    \Set\Bigg[\Dcube\Bigg]_\R\langle S_2 \rangle
  \]
  where
  \begin{align*}
    S_1 & = \{\yo 10 \ancto \yo 11, \yo 01 \ancto \yo 11,
          \yo 11 \amalg_{\yo 10} \yo 10' \ancto \yo 11',
          \yo 11 \amalg_{\yo 01} \yo 01' \ancto \yo 11', \\
    & \qquad \qquad
      \yo 00 \ancto \yo 00',
      \yo 01 \amalg_{\yo 00} \yo 00' \ancto \yo 01\}, \\
    S_2 & = S_1 \cup \{\yo 10 \amalg_{\yo 00} \yo 00' \ancto \yo 10',
          y \ancto \yo 11'\}
  \end{align*}
  where $y$ denotes the colimit of the image under $\yo$ of the entire cube minus its terminal object $11'$.
  As in the proof of \cref{prop:psi-self-extensible}, write $x_0 = \yo 10 \amalg_{\yo 00} \yo 00'$ and $x_1 = \yo 11 \amalg_{\yo 01} \yo 01'$.
  Let $M$ denote the domain of $\Lambda \bp \Psi$.
  We will show that the morphism $\Lambda \bp \Psi$ can be obtained by attaching one copy of $\Sigma_L$ and one copy of $\Sigma_R$ to $M$.

  In the diagram
  \[
    \begin{tikzcd}
      \yo 10 \ar[d] & \yo 00 \ar[l] \ar[r] \ar[d] & \yo 00' \ar[d] \\
      \yo 11 & \yo 01 \ar[l] \ar[r] & \yo 01'
    \end{tikzcd}
  \]
  in $M$, the left morphism $\yo 10 \to \yo 11$ and the induced map $\yo 01 \amalg_{\yo 00} \yo 00' \to \yo 01'$ are each (generating) anodyne cofibrations, so by the gluing lemma \cite[Lemma~7.2.15]{Hi} the induced map between pushouts $x_0 \to x_1$ is also an anodyne cofibration.
  Now in the diagram
  \[
    \begin{tikzcd}
      x_0 \ar[r] \ar[d] & x_1 \ar[d] \\
      \yo 10' \ar[r] & \yo 11'
    \end{tikzcd}
  \]
  the right morphism $x_1 \to \yo 11'$ is a generating anodyne cofibration, and the bottom morphism $\yo 10' \to \yo 11'$ is also an anodyne cofibration because $\yo 10 \to \yo 11$ and $\yo 11 \amalg_{\yo 10} \yo 10' \to \yo 11'$ are generating anodyne cofibrations.
  Form the pushout in the above square to produce
  \[
    \begin{tikzcd}
      x_0 \ar[r] \ar[d, "*"'] & x_1 \ar[d] \ar[rdd, bend left=10] \\
      \yo 10' \ar[r] \ar[rrd, bend right=10] & y \ar[rd, "*"] \\
      & & \yo 11'
    \end{tikzcd}
  \]
  in which the two maps belonging to $S_2 \setminus S_1$ are marked by $*$s.
  We can make the map $x_0 \to \yo 10'$ into an anodyne cofibration by attaching a copy of $\Sigma_L$ using the additional map $\yo 10' \to \yo 11'$, because the other way around the square $x_0 \to x_1 \to \yo 11'$ is a composition of two anodyne cofibrations; and we can make the map $y \to \yo 11'$ into an anodyne cofibration by attaching a copy of $\Sigma_R$ using the additional map $\yo 10' \to y$, which is an anodyne cofibration because it is a pushout of the anodyne cofibration $x_0 \to x_1$.
  Hence $\Lambda \bp \Psi$ is a $\JJ$-cell morphism.
\end{proof}

\section{Transferring the structure to $V\CPM$}

Now let $V$ be a monoidal combinatorial premodel category.
Recall from \cref{chap:modules} that there is a left adjoint $V \otimes_{\Set} - : \CPM \to V\CPM$ to the functor $V\CPM \to \CPM$ which forgets $V$-module structures.
On the objects of the form $\Set[D]_\R\langle S \rangle$ we have been considering, it is simply given by the formula $V \otimes_{\Set} \Set[D]_\R\langle S \rangle = V[D]_\R \langle S \rangle$.

\begin{definition}
  Renaming the previously defined $\II$ and $\JJ$ to $\II^\Set$ and $\JJ^\Set$, we define $\II$ and $\JJ$ to be the images of $\II^\Set$ and $\JJ^\Set$ respectively under $V \otimes_{\Set} -$.
  Similarly we rename the old $\Sigma_L$ to $\Sigma_L^\Set$ and write $\Sigma_L = V \otimes_{\Set} \Sigma_L^\Set$, and so on.
\end{definition}

If $V = \Set$ then $V \otimes_{\Set} -$ is the identity, so this notation extends the previous one.
In the general case, the effect is simply to replace all occurrences of ``$\Set$'' with ``$V$'' in the definitions of $E$, $\Lambda$, $\Sigma_L$, $\Sigma_R$, and $\Psi$.

\begin{proposition}
  With this new notation, the sets $\II$ and $\JJ$ are self-extensible in $V\CPM$.
\end{proposition}

\begin{proof}
  The left adjoint $V \otimes_{\Set} -$ preserves colimits and commutes with the formation of Reedy premodel categories, hence it also commutes with the $\E$ construction.
  Thus, it sends all the cell complexes which witness the self-extensibility of $\II^\Set$ and $\JJ^\Set$ to ones for $\II$ and $\JJ$.
  (Alternatively, we can just repeat the same proofs with $\Set$ replaced by $V$ as needed.)
\end{proof}

\begin{definition}
  A morphism of $V\CPM$ is an anodyne fibration (respectively, fibration) if it has the right lifting property with respect to $\II$ (respectively, $\JJ$), and a cofibration (respectively, anodyne cofibration) if it has the left lifting property with respect to all anodyne fibrations (respectively, fibrations).
\end{definition}

\begin{proposition}
  Every morphism of $V\CPM$ admits a factorization as a cofibration followed by an anodyne fibration, and also a factorization as an anodyne cofibration followed by a fibration.
\end{proposition}

\begin{proof}
  Because $\II$ and $\JJ$ are self-extensible, we can apply the large small object argument in $V\CPM$.
\end{proof}

\begin{proposition}
  A morphism of $V\CPM$ is a fibration (respectively, anodyne fibration) if and only if its underlying left Quillen functor is a fibration (respectively, anodyne fibration) in $\CPM$.
\end{proposition}

\begin{proof}
  This follows from the adjunction between $V \otimes_{\Set} -$ and the forgetful functor $V\CPM \to \CPM$ along with the definitions of $\II$ and $\JJ$.
\end{proof}

In particular any anodyne fibration of $V\CPM$ is a fibration, and so we have constructed a premodel 2-category structure on $V\CPM$.

The 2-category $\CPM$ acts on $V\CPM$ via the tensor product $\otimes : \CPM \times V\CPM \to V\CPM$, which is part of an adjunction of two variables also involving the exponential $N^M$ (for $M$ in $\CPM$ and $N$ in $V\CPM$) and the $\CPM$-valued Hom $\CPM^V(M, N)$ (for $M$ and $N$ in $V\CPM$).

\begin{proposition}
  The tensor product $\otimes : \CPM \times V\CPM \to V\CPM$ is a Quillen bifunctor.
\end{proposition}

\begin{proof}
  Using the formula $M \otimes (V \otimes_{\Set} N) = V \otimes_{\Set} (M \otimes N)$ for $M$ and $N$ in $\CPM$, this reduces to \cref{prop:cpm-monoidal}.
\end{proof}

We record the following specific consequence for later use in the path category construction.

\begin{proposition}\label{prop:vcpm-exp-fibration}
  Let $I : A \to B$ be a cofibration of $\CPM$ and $M$ a fibrant object of $V\CPM$.
  Then the induced morphism $I^* : M^B \to M^A$ is a fibration.
  If $I$ is an anodyne cofibration, then $I^*$ is an anodyne fibration.
\end{proposition}

\begin{proof}
  This follows from the fact that $\otimes : \CPM \times V\CPM \to V\CPM$ is a Quillen bifunctor by the usual sort of adjunction argument.
\end{proof}

\section{The weak equivalences and the (anodyne) fibrations}

Henceforth we shall assume that $V$ is a model category (or at least a relaxed premodel category), so that every object of $V\CPM$ is a relaxed premodel category.
In this section we will verify that the Quillen equivalences and the (anodyne) fibrations of $V\CPM$ satisfy the expected relationship.
This section is an adaptation of \cite[Proposition~1.11]{Sz}.

\begin{proposition}
  A morphism $F : M \to N$ of $V\CPM$ is an anodyne fibration if and only if it satisfies the following two conditions.
  \begin{enumerate}
  \item[(AF1)] If $f$ is a cofibration between cofibrant objects of $M$ for which $Ff$ is an anodyne cofibration, then $f$ is already an anodyne cofibration of $M$.
  \item[(AF2)] $F$ is extensible (\cref{def:extensible}).
  \end{enumerate}
\end{proposition}

\begin{proof}
  The two conditions amount to the right lifting properties with respect to $\Lambda$ and $E$, respectively.
\end{proof}

In order to study the fibrations of $V\CPM$, we make the following auxiliary definitions.

\begin{definition}
  A left Quillen functor $F : M \to N$ between relaxed premodel categories is called \emph{saturated} if whenever $f : A \to B$ is a trivial cofibration of $M^\cof$ which is sent by $F$ to an anodyne cofibration, $f$ is already an anodyne cofibration of $M$.
\end{definition}

Informally, $F : M \to N$ is saturated if $M$ has as anodyne cofibrations the largest possible class consistent with its homotopy category and the condition that $F : M \to N$ be a left Quillen functor.

\begin{example}\label{ex:model-saturated}
  Suppose that $M$ is a \emph{model} category.
  Then any left Quillen functor $F : M \to N$ is saturated, as every trivial cofibration of $M^\cof$ is already an anodyne cofibration.
  In particular, the functor $M \to 0$ is saturated.
\end{example}

\begin{proposition}\label{prop:saturated-iff-lifting}
  Let $F : M \to N$ be a morphism of $V\CPM$.
  The following are equivalent:
  \begin{enumerate}
  \item
    $F$ is saturated.
  \item
    $F$ has the right lifting property with respect to $\Sigma_L$.
  \item
    $F$ has the right lifting property with respect to both $\Sigma_L$ and $\Sigma_R$.
  \end{enumerate}
\end{proposition}

\begin{proof}
  Obviously (3) implies (2).
  Suppose that $F$ is saturated, and consider a lifting problem
  \begin{equation*}
    \begin{tikzcd}
      V[0 \to 1 \to 2]_\R\langle \yo 0 \ancto \yo 2, \yo 1 \ancto \yo 2 \rangle \ar[r] \ar[d, "\Sigma_L"'] & M \ar[d, "F"] \\
      V[0 \to 1 \to 2]_\R\langle \yo 0 \ancto \yo 1, \yo 1 \ancto \yo 2 \rangle \ar[r] & N
    \end{tikzcd}
    \tag{$*$}
  \end{equation*}
  and write $f$ for the image of the map $\yo 0 \to \yo 1$ under the top morphism.
  To construct a lift, we must show that $f$ is an anodyne cofibration.
  Now $f$ is a cofibration between cofibrant objects of $M$, and it is a left weak equivalence by the two-out-of-three property.
  Moreover, $Ff$ is an anodyne cofibration and so $f$ is already an anodyne cofibration by the hypothesis on $F$.
  Hence $F$ has the right lifting property with respect to $\Sigma_L$, and the same argument applies to $\Sigma_R$.

  It remains to show that (2) implies (1), so suppose $F$ has the right lifting property with respect to $\Sigma_L$ and let $f : A_0 \to A_1$ be a cofibration in $M^\cof$ which is sent by $F$ to an anodyne cofibration of $M$.
  Using \cref{prop:tcof-iff-comp-acof}, we may choose an anodyne cofibration $g : A_1 \to A_2$ such that $gf : A_0 \to A_2$ is also an anodyne cofibration.
  Then $f$ and $g$ determine a lifting problem of the form $(*)$, and the lifting property implies that $f$ is an anodyne cofibration.
\end{proof}

\begin{definition}
  Let $N$ be a premodel category and let $f : B_{00} \to B_{01}$ be a cofibration between cofibrant objects of $N$.
  A \emph{pseudofactorization} of $f$ is an extension of $f$ to a Reedy cofibrant diagram of the form
  \[
    \begin{tikzcd}
      B_{00} \ar[r, "f"] \ar[d] & B_{01} \ar[d] \\
      B_{10} \ar[r] & B_{11}
    \end{tikzcd}
  \]
  in which the maps $B_{10} \to B_{11}$ and $B_{01} \to B_{11}$ are anodyne cofibrations.

  Let $F : M \to N$ be a left Quillen functor.
  Then $F$ has the \emph{pseudofactorization lifting property} if for each cofibration $f : A_{00} \to A_{01}$ between cofibrant objects of $M$, each pseudofactorization of $Ff$ in $N$ lifts (up to isomorphism) to a pseudofactorization of $f$ in $M$.
\end{definition}

Compared to \cite[Definition~1.9 (3)]{Sz}, we impose extra conditions on a pseudofactorization (the diagram must be Reedy cofibrant, and $B_{10} \to B_{11}$ and $B_{01} \to B_{11}$ are not merely trivial cofibrations but anodyne cofibrations) in order to arrange that the pseudofactorization lifting property is equivalent to the right lifting property with respect to the cofibration $\Psi$.

\begin{proposition}
  A morphism $F : M \to N$ of $V\CPM$ is a fibration if and only if it is saturated and has the pseudofactorization lifting property.
\end{proposition}

\begin{proof}
  We verified in \cref{prop:saturated-iff-lifting} that $F$ is saturated if and only if it has the right lifting property with respect to both $\Sigma_L$ and $\Sigma_R$, and by the definition of $\Psi$, a morphism $F : M \to N$ has the pseudofactorization lifting property if and only if it has the right lifting property with respect to $\Psi$.
\end{proof}

We now study the (anodyne) fibrations in $V\CPM$ when $V$ is relaxed, with the eventual aim of showing that the anodyne fibrations are exactly the fibrations which are also Quillen equivalences.

\begin{lemma}\label{lemma:terminal-plp}
  Let $M$ be a relaxed premodel category.
  Then the functor $M \to 0$ has the pseudofactorization lifting property.
\end{lemma}

\begin{proof}
  We just need to check that every cofibration $f : A \to B$ of $M^\cof$ admits some pseudofactorization.
  Because $M$ is left relaxed, we can choose an anodyne cylinder object $g : C \to D$ on the cofibrant object $f : A \to B$ of $M^{[1]}_\proj$.
  As in the proof of \cref{lemma:homotopy-rel-anodyne}, let $E$ be the ``relative anodyne cylinder object'' $E = D \amalg_C A$.
  We claim that
  \[
    \begin{tikzcd}
      A \ar[r, "f"] \ar[d, "f"'] & B \ar[d, "j'_1"] \\
      B \ar[r, "j'_0"'] & E
    \end{tikzcd}
  \]
  is a pseudofactorization of $f$, where $j'_0$ and $j'_1$ are the compositions of the induced map $j' : B \amalg_A B \to E$ with the two inclusions of $B$ in $B \amalg_A B$.
  We saw in the proof of \cref{lemma:homotopy-rel-anodyne} that $j'$ is a cofibration, so it remains to check that $j'_0$ and $j'_1$ are anodyne cofibrations.
  Indeed, these maps are pushouts of the anodyne cofibrations $j_0$ and $j_1$.
\end{proof}

\begin{proposition}\label{prop:fibrant-iff-saturated}
  An object $M$ of $V\CPM$ is fibrant (that is, the functor $M \to 0$ is a fibration) if and only if the functor $M \to 0$ is saturated.
  (Concretely this means that every trivial cofibration of $M^\cof$ is an anodyne cofibration.)
\end{proposition}

\begin{proof}
  Follows from \cref{lemma:terminal-plp}.
\end{proof}

\begin{example}
  Let $M$ be a $V$-model category.
  Then $M$ is fibrant in $V\CPM$, because the functor $M \to 0$ is saturated by \cref{ex:model-saturated}.
\end{example}

We next show that fibrations are ``extensible for anodyne cofibrations''.

\begin{lemma}\label{lemma:fib-lift-acof}
  Let $F : M \to N$ be a fibration in $V\CPM$.
  Then for any cofibrant object $A$ of $M$ and any \emph{anodyne} cofibration $g : FA \to B'$, there exists an anodyne cofibration $f : A \to B$ such that $Ff = g$.
\end{lemma}

\begin{proof}
  The diagram
  \[
    \begin{tikzcd}
      FA \ar[r, "\id_{FA}"] \ar[d, "g"'] & FA \ar[d, "g"] \\
      B \ar[r, "\id_B"'] & B
    \end{tikzcd}
  \]
  defines a pseudofactorization in $N$ of $\id_{FA} = F(\id_A)$.
  Use the pseudofactorization lifting property of $F$ to obtain a lift to a pseudofactorization in $M$ of $\id_A$.
  The right-hand vertical map of the resulting pseudofactorization is then a lift of $g : FA \to B$ which is an anodyne cofibration.
\end{proof}

In order to relate the anodyne fibrations to the Quillen equivalences, we recall the criterion due to Cisinski for determining when a cofibration functor $F : C \to D$ between cofibration categories induces an equivalence on homotopy categories which we mentioned in \cref{sec:cof-cats}.
For convenience, we recall the conditions here.
\begin{enumerate}
\item[(AP1)]
  For each morphism $f$ of $C$, $Ff$ is a weak equivalence if and only if $f$ is.
\item[(AP2)]
  Let $A$ be an object of $C$ and $g : FA \to B$ a morphism of $D$.
  Then there exists a morphism $f : A \to A'$ in $C$ and weak equivalences $B \to B'$ and $FA' \to B'$ making the diagram below commute.
  \[
    \begin{tikzcd}
      FA \ar[r, "Ff"] \ar[d, "g"'] & FA' \ar[d, "\sim"] \\
      B \ar[r, "\sim"'] & B'
    \end{tikzcd}
  \]
\end{enumerate}
Then by \cite[Th\'eor\`eme~3.19]{CisCD} $\Ho F : \Ho C \to \Ho D$ is an equivalence if and only if $F$ satisfies both (AP1) and (AP2).

\begin{lemma}
  Let $F : M \to N$ be a fibration in $V\CPM$.
  Then the following conditions are equivalent.
  \begin{enumerate}
  \item $F$ satisfies (AF1).
    That is, if $f$ is a cofibration in $M^\cof$ which $F$ sends to an anodyne cofibration in $N^\cof$, then $f$ is already an anodyne cofibration in $M^\cof$.
  \item If $f$ is a cofibration in $M^\cof$ which $F$ sends to a trivial cofibration in $N^\cof$, then $f$ is already a trivial cofibration in $M^\cof$.
  \item $F^\cof : M^\cof \to N^\cof$ satisfies (AP1).
    That is, if $f$ is any morphism in $M^\cof$ which $F$ sends to a left weak equivalence in $N^\cof$, then $f$ is already a left weak equivalence in $M^\cof$.
  \end{enumerate}
\end{lemma}

Note the converse direction of each of these conditions is automatic because any left Quillen functor $F$ preserves anodyne cofibrations, trivial cofibrations and left weak equivalences.

\begin{proof}
  We first prove that (1) and (2) are equivalent.
  Suppose $F$ satisfies (AF1) and $f : A \to B$ is a cofibration in $M^\cof$ which $F$ sends to a trivial cofibration in $N^\cof$.
  Applying \cref{prop:tcof-iff-comp-acof} to $Ff$, we may find an anodyne cofibration $g' : FB \to C'$ such that $g'(Ff) : FA \to C'$ is also an anodyne cofibration.
  Using \cref{lemma:fib-lift-acof}, choose a lift of $g'$ to an anodyne cofibration $f' : B \to B'$.
  Then $F$ maps the composition $f'f$ to an anodyne cofibration $g'(Ff)$.
  By assumption $F$ satisfies (AF1), so it follows that $f'f$ is also an anodyne cofibration.
  Therefore $f$ is a left weak equivalence by the two-out-of-three property.

  Conversely, suppose (2) holds and let $f$ be a cofibration which is sent by $F$ to an anodyne cofibration.
  By condition (2), $f$ is at least a trivial cofibration; but then since $F$ is saturated and $Ff$ is an anodyne cofibration, $f$ is an anodyne cofibration as well.

  We now show that (2) and (3) are equivalent.
  Clearly (2) implies (3).
  Suppose (2) holds and let $f : A \to B$ be any morphism of $M^\cof$ which is sent by $F$ to a left weak equivalence.
  Factor $f$ as a cofibration followed by a weak equivalence.
  The image of this factorization under $F$ is also a cofibration followed by a weak equivalence.
  Since $Ff$ is a left weak equivalence, this latter cofibration is actually a trivial cofibration by the two-out-of-three property.
  Therefore the original cofibration is also a trivial cofibration by condition (2), and so the original morphism $f$ is a left weak equivalence.
\end{proof}

\begin{lemma}\label{lemma:afib-iff-qeq-of-fib}
  Let $F : M \to N$ be a fibration in $V\CPM$.
  Then $F$ is an anodyne fibration if and only if it is a Quillen equivalence.
\end{lemma}

\begin{proof}
  Above, we showed that $F$ satisfies (AF1) if and only if it satisfies (AP1).
  Thus, it suffices to prove that $F$ satisfies (AF2) if and only if it satisfies (AP2).

  First suppose that $F$ satisfies (AF2) and let $A$ be an object of $M^\cof$ and $g : FA \to B$ a morphism of $N^\cof$.
  Factor $g$ into a cofibration $j : FA \to B'$ followed by a left weak equivalence $q : B' \to B$.
  Using (AF2), we may lift $j$ to a morphism $f : A \to A'$ of $M^\cof$, so that in particular $FA' = B'$.
  Then the square
  \[
    \begin{tikzcd}
      FA \ar[r, "Ff"] \ar[d, "g"'] & FA' \ar[d, equals] \\
      B \ar[r, "q"'] & B
    \end{tikzcd}
  \]
  verifies (AP2) for $g$.

  Conversely, suppose that $F$ satisfies (AP2).
  Let $A$ be a cofibrant object of $M$ and $g : FA \to B'$ a cofibration in $N$.
  We must lift $g$ to a cofibration $f : A \to A'$ in $M$.

  We begin by applying (AP2) to the map $g : FA \to B'$, obtaining a morphism $f_1 : A \to A_1$ of $M^\cof$ whose image in $N^\cof$ fits into a square of the form below.
  \[
    \begin{tikzcd}
      FA \ar[r, "Ff_1"] \ar[d, "g"'] & FA_1 \ar[d, "\sim"] \\
      B' \ar[r, "\sim"'] & B_1
    \end{tikzcd}
  \]
  Factor $f_1 : A \to A_1$ into a cofibration $f'_1 : A \to A'_1$ followed by a left weak equivalence $A'_1 \to A_1$.
  Then $F$ sends this left weak equivalence to a left weak equivalence of $N^\cof$, so $Ff'_1$ also fits into a square of the above form.
  Thus, replacing $f_1$ by $f'_1$, we may assume without loss of generality that $f_1$ is a cofibration.

  Next, form the pushout $FA_1 \amalg_{FA} B'$ and factor the induced morphism $FA_1 \amalg_{FA} B' \to B_1$ into a cofibration $FA_1 \amalg_{FA} B' \to B'_1$ followed by a left weak equivalence.
  By the two-out-of-three property, the resulting maps $FA_1 \to B'_1$ and $B' \to B'_1$ are also left weak equivalences.
  Thus, replacing $B_1$ by $B'_1$, we may assume without loss of generality that the induced map $FA_1 \amalg_{FA} B' \to B_1$ is also a cofibration.
  In particular, the maps $FA_1 \to B_1$ and $B' \to B_1$ are also cofibrations.

  Our square now nearly has the form of a pseudofactorization of $Ff_1$, except that the morphisms $B' \to B_1$ and $FA_1 \to B_1$ are only trivial cofibrations and not anodyne cofibrations.
  We can remedy this using \cref{prop:tcof-iff-comp-acof}.
  Choose anodyne cofibrations $B_1 \to B_2$ and $B_1 \to B_3$ such that the compositions $FA_1 \to B_1 \to B_2$ and $B' \to B_1 \to B_3$ are also anodyne cofibrations, and form the pushout of $B_1 \to B_2$ and $B_1 \to B_3$ as shown below.
  \[
    \begin{tikzcd}
      FA \ar[r, hookrightarrow, "Ff_1"] \ar[d, hookrightarrow, "g"'] &
      FA_1 \ar[d, hookrightarrow, "\sim"'] \ar[rd, hookrightarrow, "\an"] \\
      B' \ar[r, hookrightarrow, "\sim"] \ar[rd, hookrightarrow, "\an"'] &
      B_1 \ar[r, hookrightarrow, "\an"'] \ar[d, hookrightarrow, "\an"] &
      B_2 \ar[d, hookrightarrow, "\an"] \\
      & B_3 \ar[r, hookrightarrow, "\an"'] &
      B_4
    \end{tikzcd}
  \]
  Now the square formed by $FA$, $B'$, $FA_1$, and $B_4$ is a pseudofactorization of $Ff_1$, because $g$ is a cofibration, the induced map $FA_1 \amalg_{FA} B' \to B_4$ is the composition of the cofibration $FA_1 \amalg_{FA} B' \to B_1$ and the anodyne cofibration $B_1 \to B_4$, hence a cofibration, and the maps $FA_1 \to B_4$ and $B' \to B_4$ are compositions of anodyne cofibrations, hence anodyne cofibrations themselves.
  Applying the pseudofactorization lifting property of $F$, we obtain in particular a cofibration $f : A \to B$ which is a lift of $g$.
\end{proof}

\begin{proposition}\label{prop:afib-iff-fib-weq}
  A morphism of $V\CPM$ is an anodyne fibration if and only if it is a fibration and a Quillen equivalence.
\end{proposition}

\begin{proof}
  In view of the fact that any anodyne fibration is a fibration, this is simply a restatement of \cref{lemma:afib-iff-qeq-of-fib}.
\end{proof}

\section{The tame--saturated factorization}

To verify that the classes of morphisms we have defined form a model 2-category structure on $V\CPM$, the main remaining task is to show that anodyne cofibrations are Quillen equivalences.
Our first step in this direction is to show that $\Sigma_L$-cell morphisms are Quillen equivalences.
Consequently, applying the large small object argument to $\{\Sigma_L\}$ yields factorizations as Quillen equivalences followed by saturated morphisms (since, by \cref{prop:saturated-iff-lifting}, a morphism is saturated if and only if it has the right lifting property with respect to $\Sigma_L$).
In order to control the homotopy theory of a premodel category in terms of its generating anodyne cofibrations over $V$, we will \emph{assume that $V$ is tractable}.

\begin{definition}
  A morphism $F : M \to M'$ of $V\CPM$ is \emph{tame} if it is a pushout of a coproduct of copies of $\Sigma_L$.
\end{definition}

If $F : M \to M'$ is tame then $M'$ is obtained from $M$ by adjoining a set $K$ of new generating anodyne cofibrations over $V$, each of which was already a trivial cofibration in $M^\cof$.
In particular, $M'$ has the same underlying $V$-module category and the same cofibrations as $M$.
(More formally, $F$ is an equivalence of $V$-module categories.
We will assume that the underlying $V$-module category of $M'$ has been identified with that of $M$ via this equivalence, so that $F$ becomes the identity $V$-module functor.)

\begin{proposition}
  If $F : M \to M'$ is tame then the underlying (non-enriched) combinatorial premodel category $M'$ is also obtained from $M$ by adjoining a set of new generating anodyne cofibrations, each of which was already a trivial cofibration in $M^\cof$.
\end{proposition}

\begin{proof}
  The diagram expressing $F : M \to M'$ as a pushout of a coproduct of copies of $\Sigma_L$ is also a pushout in $\CPM$.
  Let $I_V$ denote a set of generating cofibrations for $V$ with cofibrant domains (such an $I_V$ exists because $V$ is tractable).
  The underlying left Quillen functor of $\Sigma_L$ is $V \otimes \Sigma_L^\Set$, so by the formula for the tensor product in $\CPM$ it adjoins generating anodyne cofibrations $I_V \bp\, (\yo 0 \to \yo 1)$.
  Then by the formula for the colimit in $\CPM$, $M'$ is obtained by adjoining as generating anodyne cofibrations the images of these maps in $M$, namely the maps $I_V \bp K$.
  By \cref{prop:bifunctor-trivial}, these maps are trivial cofibrations in $M$ because the maps of $K$ are.
\end{proof}

\begin{proposition}
  If $F : M \to M'$ is tame, then $M'$ has the same left weak equivalences as $M$.
\end{proposition}

\begin{proof}
  The functor $F$ certainly preserves left weak equivalences, so suppose $f : A \to B$ is a left weak equivalence in $M'$; we need to show $f$ is already a left weak equivalence in $M$.
  Consider any left Quillen functor $F : M \to N$ from $M$ to a model category $N$.
  By the preceding proposition $M'$ is obtained by adjoining some trivial cofibrations of $M$ as anodyne cofibrations.
  Now $F$ sends these trivial cofibrations of $M$ to acyclic (i.e., anodyne) cofibrations of $N$, so $F$ extends to $M'$ as a left Quillen functor.
  Since $f$ is a left weak equivalence in $M'$, $Ff$ is a weak equivalence in $N$.
  Therefore, $f$ is a left weak equivalence in $M$.
\end{proof}

\begin{proposition}
  If $F : M \to M'$ is tame, then $M'$ has the same trivial cofibrations as $M$.
\end{proposition}

\begin{proof}
  Immediate from the previous proposition, since $M'$ also has the same cofibrations as $M$.
\end{proof}

\begin{proposition}
  A transfinite composition of tame morphisms of $V\CPM$ is again tame.
\end{proposition}

\begin{proof}
  Let $M_0 \to M_1 \to \cdots \to M_\gamma$ be a transfinite composition of tame morphisms of $V\CPM$.
  We will show by transfinite induction that the morphism $M_0 \to M_\alpha$ is tame for each $\alpha \le \gamma$.
  That is, for each $\alpha$ there exists a set $K_\alpha$ of trivial cofibrations of $M_0^\cof$ such that $M_\alpha$ is obtained from $M_0$ by adjoining $K_\alpha$ as generating anodyne cofibrations over $V$.

  This is evident for $\alpha = 0$ (take $K_0 = \emptyset$).
  At a limit step, we have $M_\beta = \colim_{\alpha < \beta} M_\alpha$.
  All the $M_\alpha$ for $\alpha < \beta$ have the same underlying category and generating cofibrations, and so we may simply take $K_\beta = \bigcup_{\alpha < \beta} K_\alpha$.
  At a successor step, $M_{\alpha+1}$ is obtained from $M_\alpha$ by adjoining a set $K$ of trivial cofibrations of $M_\alpha^\cof$ as generating anodyne cofibrations over $V$.
  But $M_0 \to M_\alpha$ is tame, so by the previous proposition, $K$ also consists of trivial cofibrations of $M_0^\cof$.
  Hence we may take $K_{\alpha+1} = K_\alpha \cup K$.
\end{proof}

\begin{proposition}
  If $F : M \to M'$ is tame, then $F$ is a Quillen equivalence.
\end{proposition}

\begin{proof}
  The homotopy category of a relaxed premodel category is computed by taking the full subcategory of cofibrant objects and inverting the left weak equivalences.
  Since $M$ and $M'$ have the same underlying category, cofibrations and left weak equivalences, they have the same homotopy category and so by definition $F : M \to M'$ is a Quillen equivalence.
\end{proof}

\begin{proposition}
  Every morphism of $V\CPM$ admits a factorization as a Quillen equivalence followed by a saturated morphism.
\end{proposition}

\begin{proof}
  The set $\{\Sigma_L\}$ is self-extensible, so any morphism $F : M \to N$ can be factored as a $\Sigma_L$-cell morphism $M \to M'$ followed by a morphism with the right lifting property with respect to $\Sigma_L$.
  The first morphism is a Quillen equivalence by the preceding two propositions, and the second morphism is saturated by \cref{prop:saturated-iff-lifting}.
\end{proof}

In particular, for any object $M$ of $V\CPM$, the morphism $M \to 0$ has a factorization as a Quillen equivalence $M \to \hat M$ followed by a saturated morphism $\hat M \to 0$.
By \cref{prop:fibrant-iff-saturated}, $\hat M$ is actually fibrant and so $M \to \hat M$ is a fibrant replacement.

\section{The mapping path category construction}

The final ingredient needed for the model category structure on $V\CPM$ is an explicit construction of a factorization of a morphism between \emph{fibrant} objects as a weak equivalence followed by a fibration.
This construction is analogous to the classical mapping path space construction, and so we call it the mapping path category construction.
In order to carry out the construction, we will need to \emph{assume $V$ is symmetric monoidal}.

We begin by defining a kind of ``unit cylinder'' object.
\begin{definition}
  We write $C$ for the combinatorial premodel category
  \[
    C = \Set[0 \to 01 \leftarrow 1]_\R
    \langle \yo 0 \ancto \yo 01, \yo 1 \ancto \yo 01 \rangle.
  \]
\end{definition}
There is a left Quillen functor $I : \Set \oplus \Set = \Set[0, 1] \to C$ which is induced by the evident inclusion of the category $\{0, 1\}$ in $\{0 \to 01 \leftarrow 1\}$.
(More plainly, $I$ sends the generator of the first copy of $\Set$ to $\yo 0$ and the generator of the second copy to $\yo 1$.)

\begin{warning}
  $C$ is not really a cylinder object on $\Set$ in $\CPM$ because there is no left Quillen functor $C \to \Set$ making both compositions $\Set \to \Set \oplus \Set \xrightarrow{I} C \to \Set$ the identity.
  To see this, first note that a left Quillen functor from $C$ to any premodel category $M$ consists of a Reedy cofibrant diagram $A_0 \to A_{01} \leftarrow A_1$ in which each of $A_0 \to A_{01}$ and $A_1 \to A_{01}$ are anodyne cofibrations.
  In particular, the map $A_0 \amalg A_1 \to A_{01}$ is supposed to be a cofibration.
  Thus the diagram is a sort of cylinder object in $M$, albeit one which may have two different ``ends'' $A_0$ and $A_1$.
  Precomposition of the functor $C \to M$ with $I : \Set \oplus \Set \to C$ corresponds to recording the two ends $(A_0, A_1)$ of the cylinder.

  Now in order for $C$ to be a cylinder object for $\Set$ in $\CPM$, we would need to be able to find such a diagram in $\Set$ in which both of the ends $A_0$ and $A_1$ are one-element sets.
  But this is impossible as the only anodyne cofibrations in $\Set$ are the isomorphisms, so that $A_{01}$ must also be a one-element set, and then $A_0 \amalg A_1 \to A_{01}$ cannot be a cofibration (injection).

  Note, however, that after tensoring with a relaxed premodel category $V$, this problem goes away: we need $A_0$ and $A_1$ to each be the unit object $1_V$, and then we can take $A_{01}$ to be an anodyne cylinder object for $1_V$.
  We will not directly make use of this observation (because we have not defined the $V\CPM$-valued internal Hom in $V\CPM$) but we will perform an essentially equivalent construction.
\end{warning}

With this caveat, $C$ does otherwise resemble a cylinder object for $\Set$.

\begin{proposition}
  The left Quillen functor $I : \Set \oplus \Set \to C$ is a cofibration in $\CPM$, and each composition $\Set \to \Set \oplus \Set \xrightarrow{I} C$ is an anodyne cofibration.
\end{proposition}

\begin{proof}
  We can express $I$ as an $\II$-cell morphism by first forming the pushout
  \[
    \begin{tikzcd}
      \Set[\star]_\R \ar[r] \ar[d, "E"'] & \Set \oplus \Set \ar[d] \\
      \Set[\star \to \star']_R \ar[r] & \Set[0 \to 01 \leftarrow 1]_\R
    \end{tikzcd}
  \]
  in which the top morphism sends $\yo *$ to the object $\yo 0 \amalg \yo 1$, then attaching two copies of $\Lambda$ to make $\yo 0 \to \yo 01$ and $\yo 1 \to \yo 01$ into anodyne cofibrations.
  The two compositions $\Set \to \Set \oplus \Set \xrightarrow{I} C$ (which are equivalent), meanwhile, are pushouts of $\Psi$
  \[
    \begin{tikzcd}
      \Set[00 \to 01]_\R \ar[r] \ar[d, "\Psi"'] & \Set[1]_\R \ar[d] \\
      \Set\Big[\Dpsquare\Big]_\R \langle \yo 10 \ancto \yo 11, \yo 01 \ancto \yo 11 \rangle \ar[r] &
      \Set[0 \to 01 \leftarrow 1]_\R \langle \yo 0 \ancto \yo 01, \yo 1 \ancto \yo 01 \rangle
    \end{tikzcd}
  \]
  along the morphism which sends $\yo 00 \to \yo 01$ to the cofibration $\emptyset \to \yo 1$.
  (The bottom morphism sends $\yo 00$ to $\emptyset$, $\yo 01$ to $\yo 1$, $\yo 10$ to $\yo 0$ and $\yo 11$ to $\yo 01$.)
  Hence the compositions $\Set \to C$ are $\JJ$-cell morphisms.
\end{proof}

\begin{proposition}\label{prop:path-fibration}
  Let $M$ be a fibrant object of $V\CPM$.
  Then the induced map $I^* : M^C \to M \times M$ is a fibration, and the compositions of $I^*$ with the projections to each copy of $M$ are anodyne fibrations.
\end{proposition}

\begin{proof}
  This follows from the preceding proposition and \cref{prop:vcpm-exp-fibration}.
\end{proof}

As $C$ is described by a presentation, we can give an explicit description of $M^C$.
\begin{itemize}
\item The underlying category of $M^C$ consists of all diagrams $A_0 \to A_{01} \leftarrow A_1$ in $M$.
\item The action of $V$ is componentwise.
\item A morphism
  \[
    \begin{tikzcd}
      A_0 \ar[r] \ar[d] & A_{01} \ar[d] & A_1 \ar[l] \ar[d] \\
      A'_0 \ar[r] & A'_{01} & A'_1 \ar[l]
    \end{tikzcd}
  \]
  is a cofibration (respectively, anodyne cofibration) if
  \begin{itemize}
  \item $A_0 \to A'_0$, $A_1 \to A'_1$, and the corner map in the square
    \[
      \begin{tikzcd}
        A_0 \amalg A_1 \ar[r] \ar[d] & A_{01} \ar[d] \\
        A'_0 \amalg A'_1 \ar[r] & A'_{01}
      \end{tikzcd}
    \]
    are cofibrations (respectively, anodyne cofibrations), and
  \item $A'_0 \amalg_{A_0} A_{01} \to A'_{01}$ and $A'_1 \amalg_{A_1} A_{01} \to A'_{01}$ are anodyne cofibrations.
  \end{itemize}
\end{itemize}
The functor $I^* : M^C \to M \times M$ sends a diagram $A_0 \to A_{01} \leftarrow A_1$ to $(A_0, A_1)$.

\begin{proposition}
  Let $M$ be any object of $V\CPM$.
  Then there exists a morphism $J : M \to M^C$ of $V\CPM$ whose composition with $I^* : M^C \to M \times M$ is the diagonal $A \mapsto (A, A)$.
\end{proposition}

\begin{proof}
  Choose an anodyne cylinder object $B$ for the unit object of $V$, so that there is a cofibration $1_V \amalg 1_V \to B$ for which each composition $1_V \to 1_V \amalg 1_V \to B$ is an anodyne cofibration.
  Then the desired morphism is given by the formula $A \mapsto (A \to B \otimes A \leftarrow A)$, where the two maps $A \to B \otimes A$ are induced by tensoring the two inclusions $1_V \to B$ with $A$.
  This functor is evidently colimit-preserving, hence a left adjoint, and it has a canonical $V$-module functor structure coming from the symmetric monoidal structure of $V$.
  It is easily seen to be a left Quillen functor, using the fact that $\otimes : V \times M \to M$ is a Quillen bifunctor.
  Finally, its composition with $I^*$ is evidently the diagonal morphism.
\end{proof}

\begin{proposition}
  Let $M$ be a fibrant object of $V\CPM$.
  With $J$ as in the preceding proposition, $M \xrightarrow{J} M^C \xrightarrow{I^*} M \times M$ is a factorization of the diagonal $M \to M \times M$ as a Quillen equivalence followed by a fibration.
\end{proposition}

\begin{proof}
  By \cref{prop:path-fibration}, $I^* : M^C \to M \times M$ is a fibration.
  Moreover, its composition with either projection $M \times M \to M$ is an anodyne fibration, hence (by \cref{prop:afib-iff-fib-weq}) in particular a Quillen equivalence.
  The functor $J : M \to M^C$ is a one-sided inverse to either composition $M^C \xrightarrow{I^*} M \times M \to M$, hence also a Quillen equivalence.
\end{proof}

We call this factorization the \emph{path category factorization} of $M$.
Note that it depends on the choice of an anodyne cylinder for the unit object of $V$; we may assume that this choice has been fixed once and for all.

Now, by standard model category methods, we can construct a corresponding \emph{mapping path category factorization} of any morphism between fibrant objects of $V\CPM$.

\begin{proposition}
  Every morphism of $V\CPM$ between fibrant objects admits a factorization as a Quillen equivalence followed by a fibration.
\end{proposition}

\begin{proof}
  Let $M \to N$ be a left Quillen functor between fibrant $V$-premodel categories.
  Define the mapping path category $P_F$ of $F$ as the pullback of the top square in the diagram below.
  As the compositions $M \to M \times N \to N \times N$ and $M \to N \to N^C \to N \times N$ both equal $(F, F) : M \to N \times N$, there is an induced left Quillen functor $M \to P_F$ as shown.
  \[
    \begin{tikzcd}
      M \ar[r, "F"] \ar[rd] \ar[rdd, "{(\id, F)}"'] & N \ar[rd] \\
      & P_F \ar[r] \ar[d] & N^C \ar[d] \\
      & M \times N \ar[r, "F \times \id"'] \ar[d, "\pi_1"'] & N \times N \ar[d, "\pi_1"] \\
      & M \ar[r, "F"'] & N
    \end{tikzcd}
  \]
  As $N$ is fibrant and both squares in the diagram are pullbacks, we conclude that $P_F \to M \times N$ is a fibration and $P_F \to M$ is an anodyne fibration and in particular a Quillen equivalence.
  The composition $M \to P_F \to M$ is the identity, so $M \to P_F$ is also a Quillen equivalence by two-out-of-three.
  As $M$ is fibrant, the composition $P_F \to M
  \times N \to N$ is also a fibration, and then $M \to P_F \to N$
  provides the desired factorization of $F : M \to N$.
\end{proof}

\section{The model 2-category $V\CPM$}

We can now complete the proof that the classes of weak equivalences, cofibrations, and fibrations that we have defined make $V\CPM$ into a model 2-category.

\begin{lemma}\label{lemma:weq-of-acof}
  Any anodyne cofibration of $V\CPM$ is a Quillen equivalence.
\end{lemma}

\begin{proof}
  Let $F : M \to N$ be an anodyne cofibration.
  Using the tame--saturated factorization twice, first factor $N \to 0$ into a Quillen equivalence $N \to \hat N$ followed by a saturated morphism $\hat N \to 0$, and then factor the composition $M \to N \to \hat N$ into a Quillen equivalence $M \to \hat M$ followed by a saturated morphism $\hat M \to \hat N$.
  \[
    \begin{tikzcd}
      M \ar[r, "\sim"] \ar[d, "F"'] & \hat M \ar[d] \\
      N \ar[r, "\sim"'] & \hat N
    \end{tikzcd}
  \]
  Then $\hat M \to 0$ is also saturated and so by \cref{prop:fibrant-iff-saturated}, both $\hat M$ and $\hat N$ are fibrant.
  Thus, we can next use the mapping path category factorization to factor $\hat M \to \hat N$ into a Quillen equivalence followed by a fibration.
  The intermediate object is again fibrant, so we may simply replace $\hat M$ by it in the above diagram, thereby reducing to the case where $\hat M \to \hat N$ is a fibration.
  Now $F$ is an anodyne cofibration, so it has the left lifting property with respect to the fibration $\hat M \to \hat N$.
  The resulting lift $N \to \hat M$ implies that $F$ (as well as all the other morphisms in the diagram) is a Quillen equivalence, by two-out-of-six.
\end{proof}

\begin{theorem}
  Let $V$ be a symmetric monoidal model category.
  The weak equivalences, cofibrations, and fibrations defined above make $V\CPM$ into a model 2-category.
  A morphism is an anodyne cofibration if and only if it is both a cofibration and a weak equivalence.
\end{theorem}

\begin{proof}
  We first prove the last claim.
  An anodyne cofibration is a cofibration (because an anodyne fibration is a fibration) and also a weak equivalence by \cref{lemma:weq-of-acof}.
  Conversely, suppose $F : M \to N$ is a cofibration which is also a weak equivalence; we will use the retract argument to show that $F$ is an anodyne cofibration.
  Factor $F$ as an anodyne cofibration $M \to M'$ followed by a fibration $M' \to N$.
  By \cref{lemma:weq-of-acof}, the anodyne cofibration $M \to M'$ is a weak equivalence and so the fibration $M' \to N$ is a weak equivalence as well (by two-out-of-three), hence an anodyne fibration by \cref{prop:afib-iff-fib-weq}.
  Then $F$ has the left lifting property with respect to $M' \to N$.
  Construct a lift as shown below.
  \[
    \begin{tikzcd}
      M \ar[r] \ar[d, "F"'] & M' \ar[d] \\
      N \ar[r, equals] \ar[ru, dotted] & N
    \end{tikzcd}
  \]
  This diagram can rearranged to display $F$ as a retract of the anodyne cofibration $M \to M'$, hence itself an anodyne cofibration.
  \[
    \begin{tikzcd}
      M \ar[r, equals] \ar[d, "F"'] & M \ar[r, equals] \ar[d] & M \ar[d, "F"] \\
      N \ar[r, dotted] \ar[rr, bend right=40] & M' \ar[r] & N
    \end{tikzcd}
  \]

  Now we can verify that the Quillen equivalences $\sW$, the cofibrations $\sC$ and the fibrations $\sF$ form a model 2-category structure on $V\CPM$.
  We saw that $V\CPM$ is complete and cocomplete in \cref{chap:algebra}.
  The weak equivalences are closed under retracts and satisfy the two-out-of-three axiom because they are the morphisms sent by the functor $\Ho^L$ to equivalences of categories.
  Finally, $(\sC, \sF \cap \sW)$ and $(\sC \cap \sW, \sF)$ are weak factorization systems on $V\CPM$ because we showed that $\sF \cap \sW$ equals the class of anodyne fibrations and $\sC \cap \sW$ equals the class of anodyne cofibrations.
\end{proof}


\begin{singlespacing}
  \renewcommand{\bibname}{References}
  \setlength\bibitemsep{\baselineskip}
  \printbibliography
\end{singlespacing}

\end{document}